\documentclass[10pt,french,leqno]{article}

\newcount\Cor\Cor=0
\def\newcor{\global\advance\Cor by 1
\par\bigskip\noindent (\romannumeral\Cor) - }

\usepackage{lipsum} 
\usepackage{titlesec}
\titleformat{\subsection}[runin]
       {\normalfont\bfseries}
       {\thesubsection}
       {0.5em}
       {}
       []
       
\titleformat{\subsubsection}[runin]
       {\normalfont\bfseries}
       {\thesubsubsection}
       {0.5em}
       {}
       []

\usepackage[T1]{fontenc}
\usepackage[french]{babel}
\usepackage{enumitem}
\usepackage{bbm}
\usepackage{bm}
\usepackage{amsmath}
\usepackage{amssymb,url,xspace}
\usepackage{mathtools}
\usepackage{mathrsfs,euscript,color}
\usepackage[all]{xy}
\usepackage{xcolor}
\usepackage{stmaryrd}

\newtheorem{theorem}{\textsc{Th\'eor\`eme}}[subsection]
\newtheorem{proposition}[theorem]{\textsc{Proposition}}
\newtheorem{lemme}[theorem]{\textsc{Lemme}}
\newtheorem{corollary}[theorem]{\textsc{Corollaire}}
\newtheorem{remark}[theorem]{\textsc{Remarque}}

\newtheorem{definition}[theorem]{\textsc{D\'efinition}}

\newtheorem{hypothese}[theorem]{\textsc{Hypoth\`ese}}

\newtheorem{exemple}[theorem]{\textsc{Exemple}}

\newenvironment{demo}{\noindent \textit{DŽmonstration.}}

\catcode`\@=11
\catcode`\Ž=13
\catcode`\=13
\catcode`\ˆ=13
\catcode`\=13
\catcode`\=13
\catcode`\‰=13
\catcode`\™=13
\catcode`\"=13
\catcode`\•=13
\catcode`\ž=13
\catcode`\=13
\catcode`\'=13
\catcode`\Ë=13
\def Ž{\'e}
\def {\`e}
\def ˆ{\`a}
\def {\`u}
\def {\^e}
\def ‰{\^a}
\def ™{\^o}
\def "{\^\i}
\def •{\"\i}
\def ž{\^u}
\def {\c c}
\def '{\"e}
\def Ë{\`A}
\catcode`\:=13
\def :{~\string:}
\catcode`\;=13
\def ;{~\string;}
\catcode`\!=13
\def !{~\string!}
\def\cad{c'est-\`a-dire\ }
\def\ie{{\textit{i.e. }}}
\def\eg{{\textit{e.g. }}}

\def\bydef{\buildrel \mathrm{d\acute{e}f}\over{=}}

\def\wt#1{\widetilde{#1}}
\def\bs#1{\boldsymbol{#1}}
\def\wh{\widehat}

\def\mbb#1{\mathbb{#1}}
\def\bs#1{\boldsymbol{#1}}
\def\ES#1{\EuScript{#1}}

\def\ni{\noindent}

\def\ptf{\,.}
\def\vg{\,,}
\def\pvg{\,;}
\def\vgq{\,,\quad}

\def\guill#1{{\guillemotleft\,#1\,\guillemotright}}
\def\qhq#1{{\quad\hbox{#1}\quad}}

\makeatletter
\def\Ddots{\mathinner{\mkern1mu\raise\p@
\vbox{\kern7\p@\hbox{.}}\mkern2mu
\raise4\p@\hbox{.}\mkern2mu\raise7\p@\hbox{.}\mkern1mu}}
\makeatother

\def\K{\mathfrak{K}}
\def\F{{\bf F}}
\def\E{{\bf E}}
\def\C{\mathbb{C}}
\def\AF{\mathbb{A}_{\F}}
\def\AE{\mathbb{A}_{\E}}
\def\GLm{{\mathrm{GL}_m}}
\def\GLmd{{\mathrm{GL}_{md}}}
\def\GLn{{\mathrm{GL}_n}}
\def\U{\mathfrak{U}}
\def\G{{\bf G}}

\def\V{{\ES{V}}}
\def\Vinf{{\ES{V}_\infty}}
\def\Vfin{{\ES{V}_{\mathrm{fin}}}}
\def\AFfin{{\mathbb{A}_{\F,\mathrm{fin}}}}
\def\A{{\mathbb{A}}}
\def\Afin{{\mathbb{A}_{\mathrm{fin}}}}

\def\Hom{\mathrm{Hom}}
\def\Isom{\mathrm{Isom}}
\def\CG{C^\infty_{\mathrm{c}}(G)}
\def\CH{C^\infty_{\mathrm{c}}(H)}
\def\CHnat{C^\infty_{\mathrm{c}}(H^\natural)}
\def\CGreg{C^\infty_{\mathrm{c}}(G_{\mathrm{r\acute{e}g}})}

\def\CHGreg{C^\infty_{\mathrm{c}}(H\cap G_{\mathrm{r\acute{e}g}})}
\def\d{\mathrm{d}}
\def\Greg{G_{\mathrm{r\acute{e}g}}}
\def\Hreg{H_{\mathrm{r\acute{e}g}}}
\def\HGreg{H\cap G_{\mathrm{r\acute{e}g}}}
\def\Agen{A^{\mathrm{g\acute{e}n}}}

\makeatletter
\newcommand{\subjclass}[2][2020]{%
  \let\@oldtitle\@title%
  \gdef\@title{\@oldtitle\footnotetext{#1 \emph{Mathematics subject classification.} #2}}%
}
\newcommand{\keywords}[1]{%
  \let\@@oldtitle\@title%
  \gdef\@title{\@@oldtitle\footnotetext{\emph{Key words and phrases.} #1.}}%
  }
\newcommand{\motscles}[1]{%
  \let\@@oldtitle\@title%
  \gdef\@title{\@@oldtitle\footnotetext{\emph{Mots-clŽs.} #1.}}%
}
\makeatother

\author{Martin Fatou \& Bertrand Lemaire\footnote{Institut de MathŽmatiques de Marseille, CNRS I2M (UMR 7373), Aix-Marseille UniversitŽ, Campus de Luminy, Case postale 907, 13288 Marseille, France.}}

\title{Induction automorphe: reprŽsentations unitaires et spectre rŽsiduel}

\begin{document}




\maketitle

\begin{abstract}
Soit $E/F$ une extension finie cyclique de corps locaux de caractŽristique nulle, de degrŽ $d$, et soit $\kappa$ un caractre de $F^\times$ de noyau $\mathrm{N}_{E/F}(E^\times)$. Pour $m\in \mathbb{N}^*$, on prouve que toute reprŽsentation irrŽductible unitaire de $\GLm(E)$ admet un $\kappa$-relvement ˆ $\GLmd(F)$ --- donnŽ par une identitŽ de caractres comme dans Henniart-Herb \cite{HH}. Soit $\E/\F$ une extension finie cyclique de corps de nombres, de degrŽ $d$, et soit $\K$ un caractre de $\AF^\times$ de noyau $\F^\times \mathrm{N}_{\E/\F}(\AE^\times)$. On prouve que toute reprŽsentation automorphe discrte de $\GLm(\AE)$ admet un $\K$-relvement (fort) ˆ $\GLmd(\AF)$, \ie compatible aux applications de relvement locales. On dŽcrit l'image et les fibres de ces applications de relvement locale et globale. En local, on traite aussi les reprŽsentations elliptiques. 
\end{abstract}

\renewcommand{\abstractname}{Abstract}
\begin{abstract}
Let $E/F$ be a finite cyclic extension of local fields of characteristic zero, of degree $d$, and $\kappa$ be a character of $F^\times$ whose kernel is 
$\mathrm{N}_{E/F}(E^\times)$. For $m\in \mathbb{N}^*$, we prove that every irreducible unitary representation of $\GLm(E)$ has a $\kappa$-lift to $\GLmd(F)$ --- given by a character identity as in Henniart-Herb \cite{HH}. 
Let $\E/\F$ be a finite cyclic extension of number fields, of degree $d$, and $\K$ be a character of $\AF^\times$ whose kernel is $\F^\times \mathrm{N}_{\E/\F}(\AE^\times)$. We prove that every automorphic discrete representation of $\GLm(\AE)$ has a (strong) $\K$-lift to $\GLmd(\AF)$, \ie compatible with the local lifting maps. We describe the image and the fibres of these local and global lifting maps. Locally, we also treat the elliptic representations.
 \end{abstract}

\tableofcontents

\section{Introduction}

\subsection{Induction automorphe: Žtat des lieux.} \label{rŽsultats connus} L'induction automorphe pour le groupe linŽaire gŽnŽral est, 
au mme titre que le changement de base et la correspondance de Jacquet-Langlands, un modle du principe de fonctorialitŽ de Langlands. 
C'est aussi un cas particulier de transfert endoscopique qui gŽnŽralise l'article fondateur de Labesse-Langlands \cite{LL} sur la $L$-indiscernability pour $\mathrm{SL}_2$. Commenons par rappeler les rŽsultats Žtablis avant cet article. 

Soit $\E/\F$ une extension finie cyclique de corps de nombres, de degrŽ $d$. On note $\A=\AF$ l'anneau des adles de $\F$. 
Fixons un caractre $\K$ de $\A^\times$ de noyau 
$\F^\times \mathrm{N}_{\E/\F}(\AE^\times)$. 
Pour $m\in \mathbb{N}^*$, l'induction automorphe (globale) associe ˆ toute repr\'esentation automorphe cuspidale unitaire 
$\tau$ de $\GLm(\AE)$ une repr\'esentation automorphe $\pi$ de $\GLmd(\AF)$ induite de cuspidale unitaire, caractŽrisŽe par le fait qu'en presque toute place finie $v$ de $\F$, le facteur $L(\pi_v,s)$ est le produit 
des facteurs $L(\tau_w,s)$ o $w$ parcourt les places de $\E$ au-dessus de $v$. Elle est $\K$-stable, \ie $(\K\circ \det)\otimes \pi \simeq \pi$. 
Par la correspondance (conjecturale) de Langlands, cette opŽration doit correspondre ˆ l'induction 
de $\E$ ˆ $\F$ des reprŽsentations galoisiennes; d'o son nom. On dira que $\pi$ est un \textit{$\K$-relvement faible} de $\tau$. 

Soit $E/F$ une extension finie cyclique de corps locaux de caracatŽristique nulle, de degrŽ $d$; \ie $F$ est une extension finie de $\mathbb{Q}_p$ pour un nombre premier $p$, ou bien (si $d>1$) $E/F \simeq \mathbb{C}/ \mathbb{R}$. Soit $\kappa$ un caractre de $F^\times$ de noyau $\mathrm{N}_{E/F}(E^\times)$. Ë toute reprŽsentation irrŽductible gŽnŽrique unitaire $\tau$ de $\GLm(E)$, l'induction automorphe (locale) associe une reprŽsentation irrŽductible gŽnŽrique unitaire $\kappa$-stable $\pi$ de $\GLmd(F)$, 
caractŽrisŽe par une identitŽ entre la fonction-caractre $\Theta_\pi^A$ de $\pi$ tordue par un isomorphisme $A$ entre $\kappa\pi= (\kappa\circ \det)\otimes \pi$ et $\pi$, 
et la fonction-caractre $\Theta_\tau$ de $\tau$. On dira que $\pi$ est un \textit{$\kappa$-relvement} de $\tau$. Cette application de $\kappa$-relvement 
s'Žtend naturellement au cas d'une $F$-algbre cyclique de degrŽ fini $E$; \eg $E/F= \E_v/\F_v$ pour une place $v$ de $\F$.

L'existence d'un $\kappa$-relvement pour les reprŽsentations irrŽductibles gŽnŽriques unitaires de $\GLm(E)$ a ŽtŽ Žtablie par Henniart-Herb \cite{HH} pour les corps locaux non archimŽdiens, puis par Henniart \cite{H1} pour les corps locaux archimŽdiens. 
L'existence d'un $\K$-relvement faible pour les reprŽsentations automorphes cuspidales unitaires de $\GLm(\AE)$ a ŽtŽ Žtablie par Henniart \cite{H2}.  
Il dŽcrit aussi l'image et les fibres de cette application de $\K$-relvement faible et prouve qu'elle est compatible aux applications de relvement locales: tout $\K$-relvement faible $\pi$ de $\tau$ est un \textit{$\K$-relvement (fort)} au sens o pour toute place $v$ de $\F$, la composante locale $\pi_v$ de $\pi$ en $v$ est un 
$\K_v$-relvement de $\tau_v$. 

\subsection{Rappels sur le changement de base.} Continuons avec les notations de \ref{rŽsultats connus}. 
Soit $\sigma$ un gŽnŽrateur du groupe de Galois de $\E/\F$, resp. $E/F$.  

Les rŽsultats de \cite{H2} sont obtenus ˆ partir des travaux d'Arthur-Clozel \cite{AC} sur le changement de base. Pour $n\in \mathbb{N}^*$, le changement de base (global) est une application de relvement de certaines reprŽsentations automorphes de $\GLn(\A)$ ˆ 
$\GLn(\AE)$ qui, suivant l'heuristique de Langlands, correspond ˆ la restriction de $\F$ ˆ $\E$ des 
reprŽsentations galoisiennes; \cad ˆ l'opŽration \guill{inverse} de l'induction automorphe. On a aussi une application de relvement locale de $\GLn(F)$ ˆ $\GLn(E)$, 
dŽfinie par une identitŽ de caractres (relation de Shintani). 

Les rŽsultats de \cite{AC} sur le changement de base, complŽtŽs par Henniart dans \cite{H2}, 
concernent --- comme ceux sur l'induction automorphe dŽcrits en \ref{rŽsultats connus} --- les reprŽsentations automorphes cuspidales unitaires et, en local, les reprŽsentations irrŽductibles gŽnŽriques unitaires. Ils ont ŽtŽ gŽnŽralisŽs par Badulescu-Henniart \cite{BH} aux reprŽsentations automorphes discrtes et, en local, aux reprŽsentations irrŽductibles unitaires. PrŽcisŽment, toute reprŽsentation irrŽductible unitaire de $\GLn(F)$ admet un relvement ˆ la Shintani, appelŽ ici \textit{$\sigma$-relvement}, qui est une reprŽsentation irrŽductible unitaire $\sigma$-stable de $\GLn(E)$; et toute reprŽsentation automorphe discrte de $\GLn(\A)$ admet un $\sigma$-relvement (fort), 
\ie compatible aux applications de relvement locales, qui est une reprŽsentation automorphe $\sigma$-stable de $\GLn(\AE)$ 
induite de discrte. Dans \cite{BH} sont aussi dŽcrites l'image et les fibres de cette application de $\sigma$-relvement. 
La mŽthode utilisŽe dans \cite{BH} est basŽe sur le principe local-global et 
la comparaison --- Žtablie dans \cite{AC} --- des parties discrtes de deux formules des traces: l'une, ordinaire, 
pour $\GLn(\A)$ et l'autre, $\sigma$-tordue, pour $\GLn(\AE)$. Comme dans \cite{AC} et \cite{H2}, des considŽrations sur les p™les des 
fonctions $L$ de paires permettent de sŽparer les reprŽsentations automorphes apparaissant dans la formule de comparaison.

\subsection{Description des rŽsultats.}\label{description des rŽsultats} Cet article gŽnŽralise les rŽsultats rappelŽs en \ref{rŽsultats connus}. 
En local, on Žtend l'application de $\kappa$-relvement ˆ toutes les reprŽsentations irrŽductibles unitaires.   
En global, on Žtend l'application de $\K$-relvement (fort) ˆ tout le spectre discret de $\GLm(\AE)$, \ie on associe 
ˆ toute reprŽsentation automorphe discrte $\tau$ de $\GLm(\AE)$ une reprŽsentation automorphe 
$\pi$ de $\GLmd(\A)$ induite de discrte telle qu'en toute place $v$ de $\F$, $\pi_v$ soit un $\K_v$-relvement de $\tau_v$. 
Comme dans \cite{H2} on dŽcrit l'image et les fibres de cette application de $\K$-relvement. 

La dŽmonstration prŽsentŽe ici est calquŽe sur celle de \cite{BH}. Comme dans \textit{loc.~cit.}, l'ingrŽdient principal est la comparaison des parties discrtes 
de deux formules des traces: l'une, ordinaire, pour $\GLm(\AE)$ et l'autre, $\K$-tordue, pour $\GLn(\A)$. La formule de comparaison se dŽduit 
de la stabilisation de la formule des traces tordue Žtablie par M\oe glin et Waldspurger dans \cite{MW2}. En effet $G'=\mathrm{Res}_{\E/\F}(\mathrm{GL}_{m/\E})$ est ˆ isomorphisme prs 
l'unique groupe endoscopique 
elliptique de $(\mathrm{GL}_{md/ \F},\K\circ \det)$; et la formule en question n'est autre que le dŽveloppement endoscopique de la partie discrte de la formule 
des traces pour $(\mathrm{GL}_{md/ \F},\K\circ\det)$. Des considŽrations sur les fonctions $L$ de paires permettent lˆ encore de sŽparer 
les reprŽsentations automorphes apparaissant dans la formule de comparaison. 

Comme l'ont fait Badulescu et Henniart pour le changement de base \cite{BH}, on Žtend aussi l'application de $\kappa$-relvement ˆ toutes les reprŽsentations irrŽductibles elliptiques, \ie dont le caractre n'est pas identiquement nulle sur l'ouvert des ŽlŽments semi-simples rŽguliers elliptiques. Pour les corps locaux non archimŽdiens, le rŽsultat a ŽtŽ obtenu dans \cite{Fa}. On le complte ici (description de l'image et des fibres) en incluant le cas des corps locaux archimŽdiens. 

On dŽcrit aussi comme dans \cite{H2} le lien entre les opŽrations d'induction automorphe et de changement de base. 

Observons que pour les corps de fonctions, \ie les corps globaux ou locaux de caractŽristique $p>0$,  
on ne dispose pas de la formule de comparaison des parties discrtes des 
formules des traces: la formule des traces n'a pas encore ŽtŽ stabilisŽe, mme dans le cas trs simple qui nous occupe ici. 
En admettant cette formule de comparaison, on peut adapter la dŽmonstration pour obtenir les mmes rŽsultats qu'en caractŽristique nulle.  

\subsection{Sur la dŽmonstration.} La dŽmonstration des rŽsultats dŽcrits en \ref{description des rŽsultats} --- ˆ l'exception de celui sur les reprŽsentations elliptiques --- est basŽe sur le principe local-global: on prouve simultanŽment les rŽsultats locaux et globaux. 

Pour une reprŽsentation automorphe discrte $\Pi$ de $\GLm(\AE)$, on obtient facilement, gr‰ce ˆ la formule de comparaison, l'existence d'un $\K$-relvement faible $\pi$ de $\Pi$ de la forme $\pi=\delta \times \K \delta \times \cdots \times \K^{k-1}\delta$ (induite parabolique normalisŽe) o 
$k\geq 1$ est un entier divisant $n=dm$ et $\delta$ est une reprŽsentation automorphe discrte 
$\K^k$-stable de $\mathrm{GL}_{\frac{n}{k}}(\A)$. De plus $k$ est le cardinal, notŽ $x(\delta)$, de l'ensemble des classes d'isomorphisme de reprŽsentations 
$\K^i \delta$ ($i\in \mathbb{Z}$). En admettant par rŽcurrence sur $m$ que les rŽsultats globaux que l'on cherche ˆ dŽmontrer soient vrais pour tout entier $1\leq m'< m$, on prouve que $k$ est aussi le cardinal, notŽ $r(\Pi)$, du stabilisateur de la classe d'isomorphisme de $\Pi$ dans le groupe $\mathrm{Gal}(\E/\F)$; en particulier $\pi$ est discrte (\ie $k=1$) si et seulement $\Pi$ est $\mathrm{Gal}(\E/\F)$-rŽgulire. De plus la dŽmonstration assure que pour toute place $v$ dans un ensemble fini $S$ de places de $\F$ (contenant les places archimŽdiennes), la composante locale $\pi_v$ de $\pi$ est un $\K_v$-relvement de $\Pi_v$. 

Ensuite on rŽalise une reprŽsentation de Speh (locale) $u$ de $\GLm(E)$ comme une composante locale en une place $v_0$ de $\F$ inerte dans $\E$ d'une reprŽsentation automorphe discrte $\Pi$ de $\GLm(\AE)$; \cad que $\E_{v_0}/\F_{v_0}$ et $\Pi_{v_0}\simeq u$. Cela permet de prouver l'existence d'un $\kappa$-relvement pour toutes les reprŽsentations de Speh puis, par compatibilitŽ entre les opŽrations locales de $\kappa$-relvement et d'induction parabolique, 
l'existence d'un $\kappa$-relvement pour toutes les reprŽsentations irrŽductibles unitaires (d'aprs la classification de ces dernires ˆ partir des reprŽsentations de Speh). RŽciproquement, on prouve que toute reprŽsentation irrŽductible unitaire $\kappa$-stable de $\GLn(F)$ est un $\kappa$-relvement.

Puisque les composantes locales d'une reprŽsentation automorphe discrte $\Pi$ de $\GLm(\AE)$ sont irrŽductibles unitaires, on en dŽduit que le $\K$-relvement faible de $\Pi$ est un $\K$-relvement (fort). RŽciproquement, toute reprŽsentation automorphe de $\GLn(\A)$ de la forme $\pi= \delta \times \K\delta \times \cdots \times \K^{k-1}\delta$ pour un entier $k\geq 1$ divisant $n$ et une reprŽsentation automorphe discrte $\K^k$-stable $\delta$ de $\mathrm{GL}_{\frac{n}{k}}(\A)$ telle que $x(\delta)=k$, est un $\K$-relvement d'une reprŽsentation automorphe discrte $\Pi$ de $\GLm(\AE)$; et si $\Pi'$ est une reprŽsentation automorphe discrte de $\GLm(\AE)$ ayant $\pi$ pour $\K$-relvement, alors $\Pi' \simeq {^\gamma\Pi}$ pour un ŽlŽment $\gamma\in \mathrm{Gal}(\E/\F)$.
 
\subsection{Organisation de l'article.} Les sections \ref{le cas n-a}--\ref{le cas a} regroupent la thŽorie locale. Dans la section \ref{le cas n-a}, on rappelle, pour une extension finie $F/\mathbb{Q}_p$, les rŽsultats de la thŽorie des reprŽsentations de $\GLn(F)$ que nous utiliserons par la suite ainsi que, pour une extension finie cyclique $E/F$, les rŽsultats connus de $\kappa$-relvement pour les reprŽsentations de $\GLm(E)$. On traite ensuite 
le cas d'une $F$-algbre cyclique $E$ de degrŽ fini (section \ref{le cas d'une algbre cyclique}) puis le cas particulier des reprŽsentations sphŽriques quand $E/F$ est non ramifiŽe (section \ref{reprŽsentations sphŽriques}). Observons que, contrairement au cas du changement de base, on ne peut pas dŽmontrer d'emblŽe l'existence d'un $\kappa$-relvement pour les reprŽsentations irrŽductibles sphŽriques unitaires; c'est pourquoi nous introduisons la notion de $\kappa$-relvement faible (\ref{def kapparelfaible}). Dans la section \ref{le cas a}, on traite le cas des corps locaux archimŽdiens. 

Les principaux rŽsultats locaux et globaux sont ŽnoncŽs dans la section \ref{ŽnoncŽ des rŽsultats} et dŽmontrŽs dans la section \ref{demo}. Le cas des reprŽsentations elliptiques est traitŽ ˆ part (section \ref{kappa-rel elliptiques}). La formule de comparaison des parties discrtes des formules des traces, essentielle ˆ la dŽmonstration, est Žtablie dans la section \ref{comparaison part disc}. La section \ref{sur l'op d'entrel IK} fait le lien entre l'opŽrateur \guill{physique} $I_\K$  
qui ˆ une forme automorphe $\phi$ sur $G(\A)$ associe la forme automorphe $\K\phi= (\K\circ \det)\cdot \phi$, et l'opŽrateur global construit ˆ partir des opŽrateurs locaux normalisŽs. 

La compatibilitŽ entre les opŽrations globales de $\K$-relvement (induction automorphe) et de $\sigma$-relvement (changement de base) 
est Žtablie dans la section \ref{comp CB}. Enfin dans la section \ref{Žtat des lieux CF}, on dŽcrit brivement l'Žtat des lieux pour les corps de fonctions, \ie les corps globaux et locaux de caractŽristique $p>0$.

\subsection{Notations.} Dans cet article, $E/F$ dŽsigne une extension finie cyclique 
de corps locaux (commutatifs) de caractŽristique nulle ou, plus gŽnŽralement, une $F$-algbre cyclique $E$ de degrŽ fini; et $\E/\F$ dŽsigne 
une extension finie cyclique de corps de nombres. On note $d\geq 1$ le degrŽ de $E/F$, resp. $\E/\F$. On note $\Gamma(E/F)$, resp. $\Gamma(\E/\F)$, le groupe de Galois de $E/F$, resp. $\E/\F$. Lorsque nŽcessaire, on fixera un gŽnŽrateur $\sigma$ de ce groupe de Galois.

Si $F$ est non archimŽdien, \ie est une extension finie du corps $\mathbb{Q}_p$ des nombres $p$-adiques pour un nombre premier $p$, on note $\mathfrak{o}_F$ l'anneau des entiers de $F$, $\mathfrak{p}_F$ l'idŽal maximal de $\mathfrak{o}_F$, $\kappa_F$ le corps rŽsiduel $\mathfrak{o}_F/ \mathfrak{p}_F$ et $q_F$ le cardinal de $\kappa_F$. On note $\nu= \nu_F$ le caractre de $F^\times$ donnŽ par la valeur absolue $\vert \,\vert_F$ sur $F$ normalisŽe par $\vert \varpi\vert_F=q_F^{-1}$ 
pour une (\ie pour toute) uniformisante $\varpi$. Si $F$ est archimŽdien, \ie est isomorphe ˆ $\mathbb{R}$ ou $\mathbb{C}$, on note $\nu=\nu_F$ le caractre de $F^\times$ donnŽ par la valeur absolue usuelle $\vert\,\vert_F$ sur $F$. 

On note $\V=\V_{\F}$ l'ensemble des places de $\F$, $\Vinf=\V_{\F,\infty}$ l'ensemble des places archimŽdiennes et 
$\Vfin= \V_{\F,\mathrm{fin}}$ l'ensemble des places finies. 
On note $\A=\AF$ l'anneau des adles de $\F$. Si $v\in \V$, on note $\F_v$ le complŽtŽ de $\F$ en $v$ et on pose $\E_v= \F_v \otimes_{\F} \E=\prod_{w\vert v}\E_w$ o $w$ parcourt les places de $\E$ au-dessus de $v$. Ainsi $\E_v$ est une $\F_v$-algbre cyclique de groupe de Galois canoniquement identifiŽ ˆ $\Gamma(\E/\F)$. 
On pose $\F_\infty = \F\otimes_{\mathbb{Q}}\mathbb{R}= \prod_{v\in \EuScript{V}_\infty}\F_v$ et on note $\Afin=\AFfin$ l'anneau des adles finis de $\F$. 
Si $v\in \Vfin$, on pose $\mathfrak{o}_v = \mathfrak{o}_{\F_v}$, $\mathfrak{p}_v = \mathfrak{p}_{\F_v}$, etc. 

Les notations suivantes dŽpendent du corps de base (local ou global). On les rappellera au dŽbut de chaque section. 
Pour $F$ (local), on pose $G=\GLn(F)$ et $H=\GLm(E)$ o $n\geq 1$ et $m\geq 1$ sont des entiers fixŽs; avec le plus souvent $n=md$. Pour $\F$ (global), on pose ${\bf G}=\GLn(\F)$ et ${\bf H}=\GLm(\E)$; et on note $G$, resp. $G'$, le $\F$-groupe algŽbrique rŽductif connexe dŽfini par  
$G=\mathrm{GL}_{n/\F}$, resp. $G'=\mathrm{Res}_{\E/\F}(\mathrm{GL}_{m/ \E})$, o $\mathrm{Res}_{\E/\F}$ dŽsigne le foncteur restriction des scalaires. Ainsi $G'(\F)={\bf H}$ et $G'(\A)=\GLm(\AE)$. 

On note $\mathbb{U}$ le groupe des nombre complexes de module $1$. On appelle \textit{caractre} (d'un groupe topologique) un homomorphisme continu dans 
$\mathbb{C}^\times$, et \textit{caractre unitaire} un homomorphisme continu dans $\mathbb{U}$. 
Toutes les reprŽsentations sont supposŽes complexes, \ie ˆ valeurs dans le groupe des automorphismes d'une espace vectoriel sur $\mathbb{C}$. 
Toutes les induites paraboliques sont supposŽes normalisŽes. L'espace d'une reprŽsentation $\pi$ est notŽ $V_\pi$.

Sauf mention expresse du contraire, toutes les fonctions considŽrŽes sont ˆ valeurs complexes.

\section{La thŽorie locale (cas non archimŽdien)}\label{le cas n-a}

Dans cette section, $F$ est une extension finie de $\mathbb{Q}_p$. Toutes les reprŽsentations 
sont supposŽes lisses, \ie telles que le stabilisateur de tout vecteur dans l'espace de la reprŽsentation soit ouvert.

\subsection{Induction parabolique.}\label{IP}
Pour $k\in \mathbb{N}^*$, on note $G_k$ le groupe $\mathrm{GL}_k(F)$. Fixons un entier $n\geq 1$ et posons $G= G_n$. Soit $A_0= (F^\times)^n\subset G$ le tore maximal (dŽployŽ) formŽ des matrices diagonales et soit $P_0= A_0\ltimes U_0\subset G$ le sous-groupe de Borel formŽ des matrices triangulaires supŽrieures. 

Si $P\subset G$ est un sous-groupe parabolique, on note $U_P$ son radical unipotent; ainsi $U_0=U_{P_0}$. On appelle \textit{facteur de Levi} une composante de Levi d'un sous-groupe parabolique. 
Un sous-groupe parabolique, resp. facteur de Levi, de $G$ est dit \textit{semi-standard} s'il contient $A_0$. Si $P$ est un sous-groupe parabolique semi-standard 
de $G$, on note $M_P$ sa composante de Levi semi-standard et $A_P=A_{M_P}$ le centre de $M_P$; ainsi $M_{P_0}=A_{P_0}=A_0$. 
Un sous-groupe parabolique, resp. facteur de Levi, de $G$ est dit \textit{standard} s'il contient $P_0$, resp. si c'est la composante de Levi semi-standard d'un sous-groupe parabolique standard de $G$. Les facteurs de Levi standard de $G$ sont donc les sous-groupes de la forme $G_{n_1}\times \cdots \times G_{n_r}$ pour des entiers $n_i\geq 1$ tels que $\sum_{i=1}^r n_i =n$. L'application $P\mapsto M_P$ est une bijection entre l'ensemble (fini) des sous-groupes paraboliques standard de $G$ et l'ensemble des facteurs de Levi standard de $G$; on note $M \mapsto P_M= MU_0$ l'application rŽciproque. 

Si $k\geq 1$ est un entier divisant $n$, on note $P_k$ le sous-groupe parabolique standard de $G$ de composante de Levi standard $L_k= G_{\frac{n}{k}}\times \cdots \times G_{\frac{n}{k}}$ ($k$ fois). 

Si $P\subset G$ est un sous-groupe parabolique et $\tau$ est une reprŽsentation d'une composante de Levi $M$ de $P$, on note $i_P^G(\tau)$ l'induite parabolique (normalisŽe) de $\tau$ ˆ $G$ suivant $P$. Si de plus $P$ est standard, $M= M_P=G_{n_1}\times \cdots \times G_{n_r}$ et $\tau= \tau_1\otimes \cdots \otimes \tau_r$ o $\tau_i$ est une reprŽsentation $\tau_i$ de $G_{n_i}$ ($i=1,\ldots ,r$), on pose 
$$\tau_1\times \cdots \times \tau_r = i_P^G(\tau)\ptf$$ 

Si $\chi$ est un caractre de $F^\times$ et $\pi$ est une reprŽsentation de $G_l$ pour un entier 
$l\geq 1$, on pose 
$$\chi\pi=
(\chi\circ \det)\otimes \pi\ptf$$
Si $\pi$ est une reprŽsentation de $G_{\frac{n}{k}}$ pour un entier $k\geq 1$ divisant $n$, on pose 
$$R(\pi,k)= \nu^{\frac{k-1}{2}}\pi \times \nu^{\frac{k-1}{2}-1}\pi \times \cdots \times \nu^{-\frac{k-1}{2}}\pi$$ et 
$$\wt{R}(\pi,k)= \nu^{k-1}\pi\times \nu^{k-2}\pi \times \cdots \times \pi\ptf$$ On a donc 
$$\wt{R}(\pi,k)\simeq  \nu^{\frac{k-1}{2}}R(\pi,k)\ptf$$ 
Si $\pi$ est irrŽductible et si $R(\pi,k)$ a un unique quotient irrŽductible --- \eg si $\pi$ est tempŽrŽe ou, plus gŽnŽralement, gŽnŽrique unitaire (voir plus loin) ---, on note $u(\pi,k)$ ce quotient. 

\subsection{Classifications (rappels).}\label{classif (rappels)} On rappelle brivement les diverses classifications dont nous aurons besoin par la suite. 

\vskip1mm
{\bf $\bullet$ ReprŽsentations de carrŽ intŽgrable (cf. \cite{BZ,Z}).} Pour un entier $k\geq 1$ divisant $n$ et une reprŽsentation irrŽductible cuspidale unitaire 
$\rho$ de $G_{\frac{n}{k}}$, on note $\delta(\rho,k)$ l'unique sous-reprŽsentation irrŽductible de $R(\rho,k)$.  
C'est une une \textit{reprŽsentation de carrŽ intŽgrable (modulo le centre)} de $G=G_n$, \cad que tout coefficient de $\delta(\rho,k)$ est de carrŽ intŽgrable sur 
$Z\backslash G$; o $Z=Z_n$ est le centre de $G_n$. RŽciproquement, si $\delta$ est une reprŽsentation irrŽductible de carrŽ intŽgrable de 
$G$, il existe un entier $k\geq 1$ divisant $n$ et une reprŽsentation irrŽductible 
cuspidale unitaire $\rho$ de $G_{\frac{n}{k}}$ telle que $\delta \simeq \delta(\rho,k)$; 
l'entier $k$ et la classe d'isomorphisme de $\rho$ sont dŽterminŽs par la classe d'isomorphisme de $\delta$. 

Une reprŽsentation irrŽductible $\pi$ de $G$ est dite \textit{essentiellement de carrŽ intŽgrable (modulo le centre)} s'il existe un caractre 
$\chi$ de $F^\times$ tel que $\chi\pi$ soit de carrŽ intŽgrable. 
Si $\tilde{\delta}$\footnote{Ë ne pas confondre avec la contragrŽdiente de $\delta$, que nous noterons $\check{\delta}$.} est une reprŽsentation irrŽductible essentiellement de carrŽ intŽgrable de $G$, il existe un entier $k\geq 1$ et une reprŽsentation irrŽductible cuspidale $\rho$ de $G_{\frac{n}{k}}$ tels que $\tilde{\delta}$ soit isomorphe ˆ l'unique sous-reprŽsentation irrŽductible $\tilde{\delta}(\rho,k)$ 
de $\wt{R}(\rho,k)$; l'entier $k$ et la classe d'isomorphisme de $\rho$ sont dŽterminŽs par la classe d'isomorphisme de 
$\delta$. L'ensemble $[\rho,\nu^{k-1}\rho]=\{\rho, \nu\rho, \ldots, \nu^{k-1}\rho\}$ est le \textit{segment de Zelevinsky} de base $\rho$ et de longueur $k$. Observons que la reprŽsentation $\tilde{\delta}(\rho,k)$ est de carrŽ intŽgrable si et seulement si $\rho'=\nu^{\frac{k-1}{2}}\rho$ est unitaire, auquel cas 
$\tilde{\delta}(\rho,k)\simeq \delta(\rho',k)$. 

\vskip1mm
{\bf $\bullet$ Classification de Langlands.} Soit $M= G_{n_1}\times \cdots \times G_{n_r}$ un facteur de Levi standard de $G$. 
Si $\pi_M$ est une reprŽsentation irrŽductible unitaire de $M$, alors l'induite parabolique $i_{P_M}^G(\pi_M)$ est irrŽductible \cite{Be}.

Si pour $i=1,\ldots ,r$, $\delta_i$ est une reprŽsentation irrŽductible de carrŽ intŽgrable  de $G_{n_i}$, alors la reprŽsentation $\delta_1\times \cdots \times \delta_r$ de $G$ est irrŽductible et tempŽrŽe. 
Toute reprŽsentation irrŽductible tempŽrŽe $\pi$ de $G$ est isomorphe ˆ un tel produit $\delta_1\times \cdots \times \delta_r$; l'entier $r$ et les classes 
d'isomorphisme des $\delta_i$ sont, ˆ permutation prs, determinŽs par la classe d'isomorphisme de $\pi$. 

Soient $\alpha_1,\ldots , \alpha_r$ des nombres rŽels tels que $\alpha_1> \cdots > \alpha_r$, et soit $\pi_i$ une reprŽsentation irrŽductible tempŽrŽe 
de $G_{n_i}$ ($i=1,\ldots ,r$). La reprŽsentation $\nu^{\alpha_1}\pi_1\times \cdots \times \nu^{\alpha_r}\pi_r$, dite \textit{standard}, a un unique quotient irrŽductible, appelŽ   
\textit{quotient de Langlands}; on le note $L(\nu^{\alpha_1}\pi_1,\ldots ,\nu^{\alpha_r}\pi_r)$. Toute reprŽsentation irrŽductible $\pi$ de $G$ est isomorphe ˆ un tel 
quotient $L(\nu^{\alpha_1}\pi_1,\ldots ,\nu^{\alpha_r}\pi_r)$; l'entier $r$, les rŽels $\alpha_1,\ldots , \alpha_r$ et les classes d'isomorphisme des 
$\pi_i$ sont dŽterminŽs par la classe d'isomorphisme de $\pi$. 

\vskip1mm
{\bf $\bullet$ ReprŽsentations elliptiques (cf. \cite{Ba}).} Une reprŽsentation irrŽductible $\pi$ de $G$ est dite \textit{elliptique} si sa fonction-caractre $\Theta_\pi$ (voir \ref{def kappa-rel}) n'est pas identiquement nulle sur l'ouvert de $G$ formŽ des ŽlŽments semi-simples rŽguliers elliptiques. Toute reprŽsentation irrŽductible essentiellement de carrŽ intŽgrable de $G$ est elliptique; ˆ l'opposŽ une induite parabolique irrŽductible $i_P^G(\tau)$ pour un sous-groupe parabolique propre $P$ de $G$ et une reprŽsentation irrŽductible $\tau$ d'une composante de Levi de $P$, ne peut pas tre elliptique. 

D'aprs \cite[2.5]{Ba}, les reprŽsentations elliptiques de $G$ sont les sous-quotients irrŽductibles des reprŽsentations $\wt{R}(\rho,k)$ o $k\geq 1$ est un entier divisant $n$ et $\rho$ est une reprŽsentation irrŽductible cuspidale de $G_{\frac{n}{k}}$. On sait \cite{Z} que la reprŽsentation $\wt{R}(\rho,k)$ a exactement $2^{k-1}$ sous-quotients irrŽductibles, chacun apparaissant avec multiplicitŽ $1$. Ces sous-quotients sont paramŽtrŽs par l'ensemble $\mathcal{P}(\{1,\ldots ,k-1\})$ des parties de $\{1,\ldots ,k-1\}$ ou, ce qui revient au mme, par l'ensemble 
$\ES{L}_{\mathrm{st}}(L_k)$ des facteurs de Levi standard de $G$ contenant $L_k$. Pour $M\in \ES{L}_{\mathrm{st}}(L_k)$, on commence par former l'induite parabolique $$\wt{R}^M(\rho,k)=i_{P_k\cap M}^M(\nu^{k-1}\rho \otimes \nu^{k-2}\rho\otimes\cdots \otimes \rho)$$ suivant le sous-groupe parabolique standard $P_k\cap M$ de $M$. Cette reprŽsentation $\wt{R}^M(\rho,k)$ de $M$ a une unique sous-reprŽsentation irrŽductible, notŽe $\tilde{\delta}^M(\rho,k)$, qui est essentiellement de carrŽ intŽgrable. On forme ensuite l'induite parabolique $$\wt{R}(\rho,k;M)= i_P^G(\tilde{\delta}^M(\rho,k))$$ suivant le sous-groupe parabolique standard $P$ de $G$ de composante de Levi $M$. Cette reprŽsentation $\wt{R}(\rho,k;M)$ de $G$ est standard; elle a donc un unique quotient, que l'on note 
$\wt{u}(\rho,k; M)$. Par construction $\wt{R}(\rho,k;M)$ est une sous-reprŽsentation de $\wt{R}(\rho,k;L_k)= \wt{R}^G(\rho,k)= \wt{R}(\rho,k)$. Plus gŽnŽralement on a $$\wt{R}(\rho,k;M)\subset \wt{R}(\rho,k;M')\qhq{si} M'\subset M\ptf$$ L'application $M \mapsto \wt{u}(\rho,k;M)$ est une bijection de $\ES{L}_{\mathrm{st}}(L_k)$ sur l'ensemble des sous-quotients irrŽductibles de $\wt{R}(\rho,k)$. 

\vskip1mm
{\bf $\bullet$ ReprŽsentations unitaires (\cite{Be,T1}).} Pour un entier $q\geq 1$ divisant $n$ et une reprŽsentation irrŽductible tempŽrŽe $\pi$ de $G_{\frac{n}{q}}$, la reprŽsentation $R(\pi,q)$ de $G$ 
a un unique quotient irrŽductible, notŽ $u(\pi,q)$. Observons que $u(\pi,q)$ n'est autre que le quotient de Langlands 
$L(\nu^{\frac{q-1}{2}}\pi, \nu^{\frac{q-1}{2}-1}\pi, \ldots , \nu^{-(\frac{q-1}{2})}\pi)$. Si $\pi$ est de carrŽ intŽgrable, la reprŽsentation 
$u(\pi,q)$ est dite \textit{de Speh}. Dans ce cas, pour tout nombre rŽel $\alpha$ tel que $0<\alpha < \frac{1}{2}$, on pose 
$$u(\pi,q;\alpha)= \nu^{\alpha} u(\pi,q)\times \nu^{-\alpha} u(\pi,q)\pvg$$ 
c'est une reprŽsentation irrŽductible de $G_{2n}$. 

Soit $\U$ l'ensemble des classes d'isomorphisme de toutes les reprŽsentations du type $u(\delta,q)$ ou $u(\delta,q;\alpha)$; o 
$\delta$ est une reprŽsentation irrŽductible de carrŽ intŽgrable de $G_a$ pour un entier $a\geq 1$, $q\in \mathbb{N}^*$ et $\alpha \in\; ]0,\frac{1}{2}[$. 
D'aprs Tadi\'c \cite{T1}, on a: 
\begin{itemize}
\item tout ŽlŽment de $\U$ est unitaire (\ie unitarisable); 
\item tout produit d'ŽlŽments de $\U$ est irrŽductible; 
\item toute reprŽsentation irrŽductible unitaire de $G$ est isomorphe ˆ un produit d'ŽlŽments de $\U$, dŽterminŽ de manire unique ˆ permutation prs des facteurs du produit.
\end{itemize}

\vskip1mm
{\bf $\bullet$ ReprŽsentations gŽnŽriques unitaires.} D'aprs Zelevinsky \cite[9.7]{Z}, une reprŽsentation irrŽductible $\pi$ de $G$ est \textit{gŽnŽrique} (voir \ref{norm op entrelac}) si et seulement si elle est isomorphe ˆ une induite parabolique irrŽductible $\nu^{\alpha_1}\delta_1\times \cdots \times \nu^{\alpha_r}\pi_r$ pour des reprŽsentations irrŽductibles de carrŽ intŽgrable $\delta_1,\ldots , \delta_r$ et des nombres rŽels $\alpha_1,\ldots ,\alpha_r$. On en dŽduit qu'une reprŽsentation irrŽductible gŽnŽrique de $G$ est elliptique si et seulement si elle est essentiellement de carrŽ intŽgrable. Si l'on Žcrit $\pi$ comme un quotient de Langlands 
$L(\wt{\pi}_1,\ldots , \wt{\pi}_s)$ o $\wt{\pi}_i$ est une reprŽsentation irrŽductible essentiellement tempŽrŽe ($i=1,\ldots , s)$, il s'ensuit que $\pi$ est gŽnŽrique si et seulement si l'induite parabolique 
$\wt{\pi}_1\times \cdots \times \wt{\pi}_s$ est irrŽductible (auquel cas $\pi \simeq \wt{\pi}_1\times \cdots \times \wt{\pi}_s$). Notons 
$\U_{\mathrm{g\acute{e}n}}\subset \U$ le sous-ensemble formŽ des (classes d'isomorphisme de) reprŽsentations du type $u(\delta,q=1)=\delta$ ou $u(\delta,q=1; \alpha)= \nu^\alpha\delta \times \nu^{-\alpha}\delta$; o 
$\delta$ est une reprŽsentation irrŽductible de carrŽ intŽgrable de $G_a$ pour un entier $a\geq 1$ et $\alpha \in\; ]0,\frac{1}{2}[$. D'aprs la classification des reprŽsentations irrŽductibles unitaires, on a: 
\begin{itemize}
\item tout ŽlŽment de $\U_{\mathrm{g\acute{e}n}}$ est (unitaire et) gŽnŽrique; 
\item toute reprŽsentation irrŽductible unitaire gŽnŽrique de $G$ est isomorphe ˆ un produit d'ŽlŽments de $\U_{\mathrm{g\acute{e}n}}$, dŽterminŽ de manire unique ˆ permutation prs des facteurs du produit.
\end{itemize}

\vskip1mm
{\bf $\bullet$ ReprŽsentations sphŽriques unitaires.} Une reprŽsentation irr\'eductible $\pi$ de $G$ est dite 
\textit{sphŽrique} si elle a un vecteur non nul fixŽ par le sous-groupe ouvert compact maximal $K= \GLn(\mathfrak{o}_F)$ de $G$. Si l'on Žcrit $\pi$ 
comme un quotient de Langlands 
$L(\wt{\pi}_1,\ldots , \wt{\pi}_s)$ o $\wt{\pi}_i$ est une reprŽsentation irrŽductible essentiellement tempŽrŽe ($i=1,\ldots , s)$, alors 
$\pi$ est sphŽrique si et seulement si $s=n$ et $\wt{\pi}_i$ est un caractre non ramifiŽ de $G_1=F^\times$ ($i=1,\ldots , n$). 
En combinant cette classification avec celle des reprŽsentations irrŽductibles unitaires de Tadi\'c \cite{T1}, 
on obtient qu'une reprŽsentation irrŽductible sphŽrique unitaire $\pi$ de $G$ est isomorphe ˆ 
$$\xi_1\times \cdots \times \xi_k \times (\eta_1\times \eta'_1)\times \cdots \times (\eta_l \times \eta'_l) $$ 
o: 
\begin{itemize}
\item pour $i= 1, \ldots , k$, $\xi_i$ est un caractre non ramifiŽ unitaire de $G_{n_i}$ pour un entier $n_i\geq 1$, donc de la forme 
$\xi_i= \nu^{a_i} \circ \det$ pour un $a_i\in \mathbb{U}$; 
\item pour $j=1,\ldots ,l$, $\eta_j$ est un caractre non ramifiŽ de $G_{m_j}$ pour un entier $m_j\geq 1$, de la forme $\eta_j = \nu^{\alpha_j + i\beta_j}\circ \det$ 
avec $\alpha_j \in\; ]0,\frac{1}{2}[$ et $\beta_j\in \mathbb{R}$, et $\eta'_j$ est le caractre non ramifiŽ $\nu^{-\alpha_j + i\beta_j}\circ \det$ de $G_{m_j}$;
\item $\sum_{i=1}^k n_i + 2\sum_{j=1}^l m_j= n$.
\end{itemize}
\`A permutation prs, les caractres $\chi_i$, $\eta_j$ et $\eta'_j$ sont dŽterminŽs par $\pi$. De plus, la reprŽsentation (irrŽductible sphŽrique unitaire) $\pi$ est gŽnŽrique si et seulement si $n_1= \cdots = n_k = m_1= \cdots = m_l=1$. 

Observons que les caractres non ramifiŽs unitaires $\xi$ de $G$ --- par exemple la reprŽsentation triviale --- 
correspondent dans la description ci-dessus au cas $k=1$ et $l=0$. Un tel $\xi$ s'Žcrit 
$\xi = \nu^{a}\circ \det$ pour un $a \in \mathbb{U}$. On rŽalise aussi $\xi$ comme le quotient de Langlands 
$$u(\nu^a\!, n)= L(\nu^{\frac{n-1}{2}+a}\!, \nu^{\frac{n-1}{2}-1+a}\!,\ldots , \nu^{-\frac{n-1}{2}+a})\ptf$$

\subsection{Facteurs de transfert et fonctions concordantes.}\label{facteurs de transfert} Soit $E$ une extension finie cyclique de $F$, de degrŽ $d$. 
On note $\mathrm{N}_{E/F}: E^\times \rightarrow F^\times$ l'application norme et $X(E/F)$ le groupe des caractres de $F^\times/ \mathrm{N}_{E/F}(E^\times)$. 
On fixe un gŽnŽrateur $\kappa$ de $X(E/F)$, \cad un caractre de $F^\times$ de noyau $\mathrm{N}_{E/F}(E^\times)$. Observons que le caractre $\nu_E=\vert\,\vert_E$ de $E^\times$ est donnŽ par $\nu_E= \nu_F\circ \mathrm{N}_{E/F}$. Pour $k\in \mathbb{N}^*$, on note $H_k$ le groupe $\mathrm{GL}_k(E)$. 

Fixons un entier $m\geq 1$ et posons $H=H_m$. On note $A_{0,H}=(E^\times)^m\subset H$ le tore maximal formŽ des matrices diagonales et $P_{0,H}= A_{0,H}\ltimes U_{0,H}\subset H$ le sous-groupe de Borel formŽ des matrices triangulaires supŽrieures. On dŽfinit comme plus haut (en remplaant $(P_0,A_0)$ par $(P_{0,H},A_{0,H})$) les notions de sous-groupe parabolique, resp. facteur de Levi, \textit{semi-standard} de $H$, et sous-groupe parabolique, resp. facteur de Levi, \textit{standard} de $H$. Si $P'$ est un sous-groupe parabolique semi-standard de $H$, on note 
$M_{P'}$ sa composante de Levi semi-standard et $A_{P'}=A_{M_{P'}}$ le centre de $M_{P'}$. Les facteurs de Levi standard de $H$ sont donc les sous-groupes de la forme $H_{m_1}\times\cdots \times H_{m_r}$ pour des entiers $m_i\geq 1$ tels que $\sum_{i=1}^r m_i=m$. 

On suppose de plus que $n=md$, avec toujours $G=G_n$. 
On fixe une base de $E^m=E\times \cdots \times E$ vu comme $F$-espace vectoriel, ce qui donne un $F$-isomorphisme $\varphi$ de 
$E^m$ sur $F^n$ et un plongement de groupes, encore notŽ $\varphi$, de $H$ dans $G$. Changer de base revient ˆ remplacer $\varphi$ par $\mathrm{Int}_g\circ \varphi$ pour un ŽlŽment $g\in G$. On identifie $H$ au sous-groupe $\varphi(H)$ de $G$. 

Pour $x \in G$, on dŽfinit le facteur discriminant $D_G(x)\in F$ de la manire habituelle:  
$$\det((T-1)\mathrm{Id}_{\mathfrak{g}}+ \mathrm{Ad}_G(x)\,\vert \, \mathfrak{g})= D_G(x)T^n + \cdots \qhq{o} \mathfrak{g}=\mathrm{M}_n(F)\ptf$$ 
Soit $\Greg$ l'ensemble des $x\in G$ tels que $D_G(x)\neq 0$. C'est l'ensemble des ŽlŽments (absolument) semi-simples rŽguliers de $G$, \ie 
ceux qui ont $n$-valeurs propres distinctes dans une cl™ture algŽbrique de $F$. 
Pour $y\in H$, on dŽfinit de la mme manire (en remplaant $F$ par $E$) le facteur discriminant $D_H(y)\in E$ et l'ensemble $\Hreg$ des $y\in H$ tels que 
$D_H(y)\neq 0$. On a l'inclusion $$\HGreg \subset \Hreg\ptf$$ 

Pour $x\in \HGreg$ est dŽfini en \cite[3.2]{HH} un facteur de transfert\footnote{On peut aussi voir $H$ comme le groupe des points $F$-rationnel de $\mathrm{Res}_{E/F}(\mathrm{GL}_{m/E})$, qui est un groupe endoscopique elliptique de $(\mathrm{GL}_{n/F}, \kappa\circ \det)$, et utiliser les facteurs de transfert de Kottwitz-Shelstad normalisŽs comme en \cite[5.3]{KS}. \`A une constante prs, ces derniers co\"{\i}ncident avec ceux de \cite{HH} (cf. \cite[remark 1.1]{HI}).} 
$$\Delta(x)=\kappa(e\wt{\Delta}(x))\in \mathbb{C}^\times$$ o $\wt{\Delta}(x)$ est un ŽlŽment de $E^\times$ qui dŽpend du choix d'un gŽnŽrateur $\sigma$ de $\Gamma(E/F)$  
et $e$ est un ŽlŽment de $E^\times$ vŽrifiant $\sigma e= (-1)^{m(d-1)}e$. Comme (par construction) 
$\sigma\wt{\Delta}(x)= (-1)^{m(d-1)}\wt{\Delta}(x)$, l'ŽlŽment $e\wt{\Delta}(x)$ appartient ˆ $F^\times$ et la quantitŽ $\kappa(e\wt{\Delta}(x))$ est bien dŽfinie. 
Observons que l'ŽlŽment $\wt{\Delta}(x)\;(\in E)$ est dŽfini pour tout $x\in H$, et que l'on a $$H\cap \Greg= \{x\in H\,\vert\, \wt{\Delta}(x)\neq 0\}\ptf$$

\begin{remark}\label{facteurs de transfert, propriŽtŽ et choix}
\textup{\begin{enumerate}
\item[(i)]Le caractre $\chi_0=\kappa^{md(d-1)/2}$ de $F^\times$ est trivial si $m(d-1)$ est pair, et vaut $\kappa^{d/2}$ sinon. Pour $\gamma \in \HGreg$ et $z\in F^\times$, 
en identifiant naturellement $z$ ˆ un ŽlŽment du centre de $G$, on a $\wt{\Delta}(z \gamma)= z^{m^2d(d-1)/2}\wt{\Delta}(\gamma)$; 
d'o l'on dŽduit que $$\Delta(z\gamma)= \chi_0(z)\Delta(\gamma)\ptf$$ \item[(ii)]Remplacer $\sigma$ par un autre gŽnŽrateur de $\Gamma(E/F)$ revient ˆ multiplier la fonction $\wt{\Delta}$ par un signe. 
D'autre part la condition $\sigma e= (-1)^{m(d-1)}e$ ne dŽpend pas du choix de $\sigma$ et dŽtermine $e$ ˆ multiplication prs par un Žlement de $F^\times$. 
En particulier si l'extension $E/F$ est non ramifiŽe, on peut prendre pour $e$ une unitŽ de $E$; en ce cas le facteur de transfert $\Delta$ ne dŽpend ni du choix de $\sigma$, ni de celui de $e\in \mathfrak{o}_E^\times$.  \hfill $\blacksquare$
\end{enumerate}
}\end{remark}

On fixe une mesure de Haar $\d g$ sur le groupe localement compact $G$. Pour $\gamma\in \Greg$, le centralisateur de $\gamma$ dans $G$ est un tore maximal de $G$, notŽ $T_\gamma$, et on note $\d t_\gamma$ la mesure de Haar sur $T_\gamma$ qui donne le volume $1$ au sous-groupe compact maximal de $T_\gamma$. 
Observons que si $\gamma_1$ et $\gamma_2$ sont deux ŽlŽments de $\Greg$ qui sont conjuguŽs dans $G$, les mesures $\d t_{\gamma_1}$ sur $T_{\gamma_1}$ et $\d t_{\gamma_2}$ sur $T_{\gamma_2}$ se correspondent par n'importe quelle conjugaison dans $G$ transformant $\gamma_1$ en $\gamma_2$. On note $\CG$ l'espace vectoriel des fonctions localement constantes et ˆ support compact sur $G$. On note encore $\kappa$ le caractre $\kappa\circ\det$ de $G$. 
Pour $f\in \CG$ et $\gamma\in \Greg$, on pose 
$$I^G_\kappa(f,\gamma)= \vert D_G(\gamma)\vert_F^{\frac{1}{2}}\int_{T_\gamma\backslash G}f(g^{-1}\gamma g)\kappa (g)\frac{\d g}{\d t_\gamma} \qhq{si} 
T_\gamma \subset \ker(\kappa)$$ (l'intŽgrale est absolument convergente) et 
$I^G_\kappa(f,\gamma)= 0$ sinon. La distribution $f\mapsto I^G_\kappa(f,\gamma)$ sur $G$ est $\kappa$-invariante au sens o pour tout $x\in G$, on a 
$$I^G_\kappa(f\circ \mathrm{Int}_x, \gamma)= \kappa (x) I^G_\kappa(f,\gamma)\qhq{pour toute fonction}f\in \CG \ptf$$
De mme on fixe une mesure de Haar $\d h$ sur $H$. Pour $\delta\in \Hreg$, le centralisateur de $\delta$ dans $H$ est un tore maximal de $H$, notŽ $T'_\delta$, et on note $\d t'_\delta$ la mesure de Haar sur $T'_\delta$ qui donne le volume $1$ au sous-groupe compact maximal de $T'_\delta$. 
On note $\CH$ l'espace vectoriel des fonctions localement constantes et ˆ support compact sur $H$. Pour $\phi\in \CH$ et $\delta\in \Hreg$, on pose 
$$I^H(\phi, \delta)= \vert D_H(\delta)\vert_E^{\frac{1}{2}}\int_{T'_\delta\backslash H}\phi(h^{-1}\delta h) \frac{\d h}{d t'_\delta}$$ (l'intŽgrale est absolument convergente). 

\begin{remark}\label{sur le discriminant de Weyl}
\textup{On a dŽfini le facteur discriminant $D_H(\delta)$ en considŽrant $H$ comme le groupe des points $E$-rationnels d'un groupe algŽbrique dŽfini sur 
$E$ (en l'occurence $\mathrm{GL}_{m/E}$). On peut aussi considŽrer $H$ comme le groupe des points 
$F$-rationnels du $F$-groupe algŽbrique (rŽductif connexe) $G'=\mathrm{Res}_{E/F}(\mathrm{GL}_{m/E})$. Fixons une cl™ture 
algŽbrique $\overline{F}$ de $F$. Pour $\delta \in G'$ semi-simple rŽgulier, 
on a le facteur discriminant
$$D_{G'}(\delta)= \textstyle{\det_{\overline{F}}}(1-\mathrm{Ad}_\delta\,\vert \, \mathfrak{g}'(\overline{F})/\mathfrak{t}'_\delta(\overline{F}))\in \overline{F}$$ o $\mathfrak{g}'= \mathrm{Lie}(G')$ et $\mathfrak{t}'_\delta$ est le centralisateur de $\delta$ dans $\mathfrak{g}'$. Si de plus $\delta$ est $F$-rationnel, 
\ie si $\delta$ appartient ˆ $ \Hreg$, alors $D_{G'}(\delta)$ appartient ˆ $F$ et co\"{\i}ncide avec 
$\det_F(1-\mathrm{Ad}_\delta\,\vert \, \mathfrak{g}'(F)/\mathfrak{t}'_\delta(F))$; en ce cas on a l'ŽgalitŽ 
$$\vert D_H(\delta)\vert_E = \vert D_{G'}(\delta)\vert_F\ptf\eqno{\blacksquare}$$}
\end{remark}

\begin{definition}\label{fonctionsconcor}
\textup{
On dit que deux fonctions $f\in \CG$ et $\phi\in \CH$ sont \textit{concordantes} si pour tout $\gamma\in \HGreg$, on a l'ŽgalitŽ 
$$\Delta(\gamma) I^G_\kappa(f,\gamma)= I^H(\phi,\gamma)\ptf$$
}
\end{definition} D'aprs Waldspurger \cite[cor.~1.7]{W2}, on a le 

\begin{theorem}\label{existence transfert}
\begin{enumerate}
\item[(i)]Soit une fonction $f\in\CG$. Il existe une fonction $\phi\in \CH$ qui concorde avec $f$.
\item[(ii)]Soit une fonction $\phi\in \CH$ telle que pour tout $\gamma \in \HGreg$ et tout $x\in G$ tels que $\gamma'=x^{-1}\gamma x$ appartienne ˆ $ H$, on ait 
l'ŽgalitŽ  $$I^H(\phi,\gamma')= \kappa(x)^{-1} \frac{\Delta(\gamma')}{\Delta(\gamma)} I^H(\phi,\gamma)\ptf$$ 
Il existe une fonction $f\in \CG$ qui concorde avec $\phi$.
\end{enumerate}
\end{theorem}

\begin{remark}\label{propriŽtŽs fonct concor}
\textup{\begin{enumerate}
\item[(i)]Pour que $f\in \CG$ et $\phi\in \CH$ soient concordantes, il est nŽcessaire que pour tout $\gamma \in \HGreg$ et tout $x\in G$ tels que l'ŽlŽment 
$\gamma'=x^{-1}\gamma x$ appartienne ˆ $ H$, l'ŽgalitŽ du point (ii) soit vŽrifiŽe. Si de plus $x^{-1}Hx =H$, cette ŽgalitŽ se simplifie en \cite[cor.~4.3]{HH} 
 $$I^H(\phi,\gamma')= I^H(\phi,\gamma)\pvg$$ 
observons aussi que $\vert D_H(\gamma')\vert_E^{\frac{1}{2}}= \vert D_H(\gamma)\vert_E^{\frac{1}{2}}$. En particulier pour que $f$ et $\phi$ soient concordantes, il est nŽcessaire que la fonction $\gamma \mapsto I^H(\phi,\gamma)$ sur $\HGreg$ soit invariante sous l'action de $\Gamma(E/F)$ (\cite[4.4, remark 2]{HH}). 
\item[(ii)] Si l'on se limite aux fonctions ˆ support dans les ŽlŽments semi-simples $G$-rŽguliers, \cad ˆ $f\in \CGreg$ et 
$\phi\in \CHGreg$, le thŽorme est beaucoup plus facile ˆ dŽmontrer (cf. \cite[theorem~3.8]{HH}). \hfill $\blacksquare$
\end{enumerate}
 }
\end{remark}

\subsection{$\kappa$-relvement des reprŽsentations.}\label{def kappa-rel}
Pour une reprŽsentation $\pi$ de $G$ et une fonction $f\in \CG$, on note $\pi(f)=\pi(f,\d g)$ l'opŽrateur sur l'espace $V_\pi$ de $\pi$ dŽfini par 
$$\pi(f)(v)= \int_G f(g)\pi(v)\d g\qhq{pour tout}v\in V_\pi\ptf$$ En notant $\check{\pi}$ la reprŽsentation contragrŽdiente de $\pi$, on a donc
$$\langle \pi(f)(v),\check{v}\rangle =\int_G f(g) \langle v, \check{\pi}(\check{v}) \rangle \d g\qhq{pour tout}(v,\check{v})\in V_\pi \times V_{\check{\pi}}\ptf$$ 
On dŽfinit de la mme manire, pour toute reprŽsentation $\tau$ de $H$ et toute fonction $\phi\in \CH$, l'opŽrateur $\tau(\phi)= \tau(\phi, \d h)$. 

Une reprŽsentation $\pi$ de $G$ est dite \textit{$\kappa$-stable} si $\kappa\pi= (\kappa\circ \det)\otimes \pi$ est isomorphe ˆ $\pi$, \ie s'il existe un $\mathbb{C}$-automorphisme $A$ de $V_\pi$ tel que $$A\circ \kappa\pi(g)= \pi(g)\circ A \qhq{pour tout} g\in G\ptf$$ 
Si $\pi$ est irrŽductible et $\kappa$-stable, alors (lemme de Schur) l'espace $\Hom_G(\kappa\pi,\pi)$ des opŽrateurs d'entrelacement entre $\kappa\pi$ et $\pi$ est de dimension $1$. 

\begin{definition}\label{defkapparel}\textup{Soit $\tau$ une reprŽsentation irrŽductible de $H$. Soit $\pi$ une reprŽsentation irrŽductible $\kappa$-stable de $G$ et soit $A\in \Isom_G(\kappa\pi,\pi)$. On dit que $\pi$ est un \textit{$\kappa$-relvement} de $\tau$ s'il existe un constante $c(\tau,\pi,A)\in \mathbb{C}^\times$ telle que 
$$\mathrm{tr}(\pi(f)\circ A)= c(\tau,\pi,A)\mathrm{tr}(\tau(\phi))\leqno{(1)}$$ pour toute paire de fonctions concordantes 
$(f,\phi)\in \CGreg \times \CHGreg$.
}
\end{definition}

On sait d'aprs Harish-Chandra \cite{HC2} (cf. \cite[2.7]{HL2}) que si $\pi$ est une reprŽsentation irrŽductible $\kappa$-stable de $G$ et si $A\in \mathrm{Isom}_G(\kappa\pi,\pi)$, la distribution $f \mapsto \mathrm{tr}(\pi(f)\circ A_\pi)$ sur $G$ est donnŽe sur l'ouvert $\Greg$ par une fonction localement constante $\Theta_\pi^A $: pour toute fonction $f\in \CGreg$, on a $$\mathrm{tr}(\pi(f)\circ A)= \int_G \Theta_\pi^A(g)f(g) \d g$$  
(l'intŽgrale est absolument convergente). De la mme manire, pour toute reprŽsentation irrŽductible $\tau$ de $H$, la distribution $\phi\mapsto \mathrm{tr}(\tau(\phi))$ sur $H$ est donnŽe sur l'ouvert $\Hreg$ par une fonction localement constante 
$\Theta_\tau$.  
La notion de $\kappa$-relvement s'exprime aussi comme une ŽgalitŽ entre ces (fonctions-)caractres. Pour $\gamma\in \HGreg$, soit $X(\gamma)\subset \HGreg$ un ensemble (fini) de reprŽsentants des classes de $H$-conjugaison dans $H$ qui rencontrent la classe de $G$-conjugaison de $\gamma$. 
D'aprs la formule d'intŽgration de Weyl (cf. \cite[3.11]{HH}), l'ŽgalitŽ (1) est Žquivalente ˆ: pour tout $\gamma\in \HGreg$, on a 
$$\vert D_G(\gamma)\vert_F^{\frac{1}{2}}\Theta_\pi^A(\gamma)= c(\tau,\pi,A)\sum_{\gamma'\in X(\gamma)} \kappa(x_{\gamma'})^{-1}\Delta(\gamma')\vert D_H(\gamma')\vert_E^{\frac{1}{2}}\Theta_\tau(\gamma')\leqno{(2)}$$ o, pour chaque $\gamma'$, on a choisi un 
$x_{\gamma'}\in G$ tel que $\mathrm{Int}_{x_{\gamma'}}(\gamma')=\gamma$; la validitŽ de l'ŽgalitŽ ci-dessus ne dŽpend pas du choix des $x_{\gamma'}$. 

\begin{remark}\label{rem fonc-car kappa-rel}
\textup{
\begin{enumerate}
\item[(i)] Pour $\gamma \in \Greg$, on a $\Theta_\pi^A(\gamma)=0$ si $\gamma$ n'est pas $G$-conjuguŽ ˆ un ŽlŽment de $H$ (cf. \cite[rem.~2.7]{HL2}).
\item[(ii)] Pour $\gamma \in \HGreg$ elliptique (dans $H$ ou dans $G$, c'est la mme chose), l'ŽgalitŽ (2) se simplifie et devient (cf. 
\cite[3.11,\,rem.~2, p.~145]{HH})
$$\vert D_G(\gamma)\vert_F^{\frac{1}{2}}\Theta_\pi^A(\gamma)= c(\tau,\pi,A) \Delta(\gamma)\vert D_H(\gamma)\vert_E^{\frac{1}{2}}\sum_{\sigma'\in \Gamma(E/F)}\Theta_\tau(\sigma'\gamma)\ptf$$
\item[(iii)]
On sait aussi d'aprs Harish-Chandra \cite{HC1} (cf. \cite[2.7]{HL2}) que la distribution $f\mapsto \mathrm{tr}(\pi(f)\circ A)$ sur $G$ est localement intŽgrable: 
l'ŽgalitŽ 
$$\mathrm{tr}(\pi(f)\circ A)= \int_G \Theta_\pi^A(g)f(g) \d g$$ est vraie pour toute fonction $f\in \CG$. On en dŽduit 
que $\pi$ est un $\kappa$-relvement de $\tau$ si et seulement si l'ŽgalitŽ (1) est vŽrifiŽe pour toute paire de fonctions concordantes 
$(f,\phi)\in \CG \times \CH$.\hfill$\blacksquare$
\end{enumerate}
}
\end{remark}

La notion de $\kappa$-relvement ne dŽpend que des classes d'isomorphisme de $\tau$ et $\pi$. De plus  
un $\kappa$-relvement $\pi$ de $\tau$, s'il existe, est unique ˆ isomorphisme prs \cite[cor.~3.11]{HH}; c'est pourquoi on Žcrira aussi de manire un peu impropre  \guill{le} $\kappa$-relvement. Observons que la notion de $\kappa$-relvement ne dŽpend pas non plus du choix de $\kappa$ (dŽfinissant $E$): si $\pi$ est un $\kappa$-relvement de $\tau$, alors c'est un 
$\kappa'$-relvement de $\tau$ pour tout gŽnŽrateur $\kappa'$ de $X(E/F)$. 

\begin{lemme}[{\cite[4.5]{HH}}]\label{[HH] 4.5}
Soit $\tau$ une reprŽsentation irrŽductible de $H$, de caractre central $\omega_\tau: E^\times \rightarrow \mathbb{C}^\times$. Supposons que $\tau$ admette un $\kappa$-relvement $\pi$. 
\begin{enumerate}
\item[(i)] Le caractre central 
$\omega_\pi$ de $\pi$ est donnŽ par 
$\omega_\pi= \kappa^{m\frac{d(d-1)}{2}}\omega_\tau\vert_{F^\times}$. 
\item[(ii)] Pour tout caractre $\xi$ de $F^\times$, $\xi\pi$ est un $\kappa$-relvement de $(\xi\circ \mathrm{N}_{E/F})\tau$.
\item[(iii)] La contragrŽdiente $\check{\pi}$ de $\pi$ est un $\kappa$-relvement de $\check{\tau}$.
\item[(iv)] Pour tout $\sigma'\in \Gamma(E/F)$, $\pi$ est un $\kappa$-relvement de $\tau^{\sigma'}=\tau\circ \sigma'$.
\end{enumerate}
\end{lemme}

Pour la notion de reprŽsentation irrŽductible gŽnŽrique \textit{ˆ segments $\Gamma(E/F)$-rŽguliers}, on renvoie ˆ \cite[2.4]{HH}. Toute reprŽsentation irrŽductible gŽnŽrique unitaire de $H$ est ˆ segments $\Gamma(E/F)$-rŽguliers.

\begin{proposition}[{\cite[2.4]{HH}}]\label{proposition HH} 
\begin{enumerate}
\item[(i)] Toute reprŽsentation irrŽductible gŽnŽrique ˆ segments $\Gamma(E/F)$-rŽgulier $\tau$ de $H$ admet un $\kappa$-relvement $\pi$ qui est une reprŽsentation irrŽductible gŽnŽrique de $G$.
\item[(ii)] Toute reprŽsentation irrŽductible gŽnŽrique $\kappa$-stable $\pi$ de $G$ est le $\kappa$-relvement d'une reprŽsentation irrŽductible gŽnŽrique 
ˆ segments $\Gamma(E/F)$-rŽguliers $\tau$ de $H$.
\end{enumerate}
\end{proposition}

\begin{remark}\label{Ichino-Iraga}
\textup{Il est prouvŽ en \cite{HL2} que la constante de proportionnalitŽ $c(\tau,\pi,A)$ 
ne dŽpend pas des reprŽsentations irrŽductibles gŽnŽriques $\tau$ et $\pi$ telles que 
$\pi$ soit un $\kappa$-relvement de $\tau$, pourvu que l'opŽrateur $A\in \mathrm{Isom}_G(\kappa\pi,\pi)$ soit normalisŽ ˆ l'aide d'une fonctionnelle de Whittaker 
$\lambda: V_\pi\rightarrow \mathbb{C}$ dŽfinie via le choix d'un caractre additif non trivial $\psi$ de $F$. La constante $c(\tau,\pi,A)$ 
ne dŽpend alors que de $\psi$ et des choix auxilliaires effectuŽs au dŽpart (choix d'une $F$-base de $E^m$ dŽfinissant le plongement $\varphi : H \hookrightarrow G$, choix d'un gŽnŽrateur 
$\sigma$ de $\Gamma(E/F)$ et d'un ŽlŽment $e\in E^\times$ tel que $\sigma e = (-1)^{m(d-1)}e$). Il est aussi prouvŽ (\textit{loc.\,cit.}) comment elle dŽpend de $\psi$ et de ces choix; en particulier elle est de module $1$. 
Ichino et Iraga ont ensuite prouvŽ dans \cite{HI} que si l'on normalise les facteurs de transfert $\Delta(\gamma)$ comme dans Kottwitz-Shelstad \cite[5.3]{KS} --- c'est la normalisation dite \guill{de Whittaker}, dŽfinie ˆ l'aide du caractre $\psi$ ---, cette constante vaut $1$.\hfill $\blacksquare$}
\end{remark}

L'intŽrt de considŽrer les reprŽsentations irrŽductibles qui sont gŽnŽriques unitaires est fondamental pour les applications globales: on sait en effet que si $\F$ est un corps de nombres et si $\pi$ est une reprŽsentation automorphe cuspidale de $\GLn(\A)$, alors en toute place $v\in \Vfin$, la composante locale $\pi_v$ est gŽnŽrique (unitaire); on prŽvoit bien sžr qu'elle soit tempŽrŽe (c'est la conjecture de Ramanujan-Petterson). Le cas crucial dans la proposition \ref{proposition HH} est celui o la reprŽsentation $\tau$ de $H$ est de carrŽ intŽgrable. En ce cas le $\kappa$-relvement de $\tau$ est 
de la forme $ \pi_1\times \kappa \pi_1\times \cdots \times \kappa^{d_1-1}\pi_1$ pour un entier $d_1$ divisant $n$ et une reprŽsentation irrŽductible de carrŽ intŽgrable $\pi_1$ de $G_{n_1}$, $n_1=\frac{n}{d_1}$, telle que $\langle \kappa^{d_1} \rangle$ soit le stabilisateur de la classe d'isomorphisme de $\pi_1$ dans $X(E/F)$. Par induction parabolique et torsion par les caractres de $H$, cela entra"ne que le $\kappa$-relvement d'une reprŽsentation irrŽductible essentiellement tempŽrŽe de $H$ est une reprŽsentation (irrŽductible, $\kappa$-stable) essentiellement tempŽrŽe de $G$. 

Le $\kappa$-relvement des reprŽsentations irrŽductibles essentiellement tempŽrŽes de $H$ permet de dŽfinir, via la classification de Langlands, 
la notion suivante d'\textit{induite automorphe}: une reprŽsentation irrŽductible $\pi\simeq L(\wt{\pi}_1,\ldots , \wt{\pi}_r)$ de $G$ avec $\wt{\pi}_i$ essentiellement 
tempŽrŽe ($i=1,\ldots ,r$) est 
l'induite automorphe d'une reprŽsentation irrŽductible 
$\tau \simeq L(\wt{\tau}_1,\ldots , \wt{\tau}_s)$ de $H$ avec $\wt{\tau}_j$ essentiellement tempŽrŽe ($j=1,\ldots ,s$) si et seulement si $r=s$ et $\wt{\pi}_i$ est un $\kappa$-relvement de $\wt{\tau}_i$ ($i= 1,\ldots ,r$). L'unicitŽ du quotient de Langlands assure que $\pi$ est $\kappa$-stable; et comme pour la notion de $\kappa$-relvement, la notion d'induite automorphe ne dŽpend que des classes d'isomorphisme de $\tau$ et $\pi$. Par dŽfinition, toute reprŽsentation irrŽductible $\tau$ de $H$ admet une induite automorphe $\pi$; mais $\pi$ 
n'est en gŽnŽral pas un $\kappa$-relvement de $\tau$, comme le montre l'exemple suivant. 

\begin{exemple}
\textup{
Prenons le cas o $m=d=2$ et $n=4$. Notons $\bs{1}_E$ le caractre trivial de $\mathrm{GL}_1(E)=E^\times$. La sŽrie principale non ramifiŽe 
$\tau=\nu_E\times \bs{1}_E$ de $H= \mathrm{GL}_2(E)$ est rŽductible: elle s'insre dans la suite exacte courte 
$$0 \rightarrow \tilde{\delta}(\nu_E,2) \rightarrow \tau \rightarrow L(\nu_E, \bs{1}_E)\rightarrow 0\ptf$$  
La reprŽsentation $\tilde{\delta}(\nu_E,2)=\nu_E^{1/2} \delta(\bs{1}_E,2)$ est essentiellement de carrŽ intŽgrable, et le quotient de Langlands $L(\nu_E, \bs{1}_E)$ est le caractre non ramifiŽ $\nu_E^{1/2}\circ \det_E $ de $H$. Le caractre trivial $\bs{1}_E$ de $E^\times$ admet pour $\kappa$-relvement la sŽrie principale (irrŽductible) 
$\rho = \bs{1}_F \times \kappa$ de $\mathrm{GL}_2(F)$; et le caractre $\nu_E= \nu_F\circ \mathrm{N}_{E/F}$ de $E^\times$ admet pour $\kappa$-relvement la sŽrie principale 
$\nu_F \times \kappa\nu_F\simeq \nu_F\rho$ de $\mathrm{GL}_2(F)$. Par suite le quotient de Langlands $L(\nu_F\rho,\rho)$ est l'induite automorphe de $L(\nu_E,\bs{1}_E)$. Or $L(\nu_F\rho, \rho)$ est isomorphe ˆ l'induite parabolique 
$L(\nu_F,\bs{1}_F)\times \kappa L(\nu_F,\bs{1}_F)$ o $L(\nu_F,\bs{1}_F)$ est le caractre non ramifiŽ 
$\nu_F^{1/2}\circ \det_F$ de $\mathrm{GL}_2(F)$. Le caractre $\mathrm{tr}(L(\nu_F\rho,\rho))$ s'annule sur le sous-ensemble $G_{\mathrm{ell}}\subset G_{\mathrm{r\acute{e}g}}$ formŽ des ŽlŽments elliptiques. Comme $H\cap G_{\mathrm{ell}}$ est non vide et que le caractre $\nu_E^{1/2}\circ \det_E$ de $\mathrm{GL}_2(E)$ ne s'annule nulle part, la reprŽsentation $L(\nu_F\rho,\rho)$ ne peut pas tre un $\kappa$-relvement de $L(\nu_E,\bs{1}_E)$. \hfill$\blacksquare$}
\end{exemple}

\subsection{$\kappa$-relvement et induction parabolique.}\label{compatibilitŽ avec l'IP} Soit $M= G_{n_1}\times \cdots \times G_{n_r}$ un facteur de Levi standard de $G= \GLn(F)$. Si $\pi_M=\pi_1\otimes\cdots \otimes \pi_r$ est une reprŽsentation de $M$, on note $\kappa\pi_M$ la reprŽsentation $\kappa\pi_1\otimes \cdots \otimes \kappa\pi_r$  et on dit que $\pi_M$ est $\kappa$-stable si $\kappa\pi_M\simeq \pi_M$.

Pour $i=1,\ldots, r$, soit $\pi_i$ une reprŽsentation $\kappa$-stable (irrŽductible ou non) $\pi_i$ de $G_{n_i}$ et soit $A_i\in \mathrm{Isom}_G(\kappa\pi_i,\pi_i)$. 
La reprŽsentation $\pi_M=\pi_1\otimes \cdots \otimes \pi_r$ de $M$ est $\kappa$-stable et l'opŽrateur $B= A_1\otimes \cdots \otimes A_r$ sur l'espace $V_{\pi_M}= V_{\pi_1}\otimes \cdots \otimes V_{\pi_r}$ de $\pi_M$ appartient ˆ $\mathrm{Isom}_M(\kappa\pi_M , \pi_M)$. Soit $\pi= \pi_1\times \cdots \times \pi_r \;(=i_P^G(\pi_M))$ avec $P=P_M$; c'est une reprŽsentation $\kappa$-stable de $G$. On dŽfinit comme suit un opŽrateur 
$$A=A_1\times \cdots \times A_r \;(=i_P^G(B))\in  \mathrm{Isom}_G(\kappa\pi,\pi)\ptf$$ 
Rappelons que l'espace $V_\pi$ de $\pi$ est formŽ des fonctions $\varphi: G \rightarrow V_{\pi_M}$ tels que 
pour tout $(m,u,g)\in M\times U_P\times G$, on ait
$$\varphi(mug)= \bs{\delta}_P(m)^{\frac{1}{2}} \pi_M(m)(\varphi(g))\pvg$$ o 
$\bs{\delta}_P: P \rightarrow \mathbb{R}_+^*$ est le caractre-module de $P$ dŽfini par $\d (pup^{-1}) = \bs{\delta}_P(p)\d u$ pour une (\ie pour toute) mesure de Haar ˆ droite ou ˆ gauche $\d u$ sur $U_P$. La reprŽsentation $\pi$ agit sur $V_\pi$ par translations ˆ droite: pour $\varphi\in V_\pi$, on a 
$$(\pi(x)\varphi)(g)= \varphi(gx)\qhq{pour tous} x,\,g\in G\ptf$$ La reprŽsentation $\pi' =\kappa\pi$ agit sur le mme espace $V_\pi$: pour $\varphi\in V_\pi$ et 
$x,\, g\in G$, on a $(\pi'(x) \varphi)(g) =\kappa(x)\varphi(gx)$. Pour $\varphi\in V_\pi$, on pose
$$(A\varphi)(g)= \kappa(g)B(\varphi(g))\qhq{pour tout} g\in G\ptf$$ 
Pour $\varphi\in V_\pi$ et $x,\, g\in G$, on a bien $$(A\circ \pi'(x)\varphi)(g)= \kappa(g)B(\pi'(x)\varphi)(g)) = \kappa(gx)B(\varphi(gx))= (\pi(x)\circ A\varphi )(g)\vg$$ 
\ie $A\in \mathrm{Isom}_G(\pi'\!,\pi)$.

\begin{proposition}[{\cite[3.7]{HL2}}]\label{compatibilitŽ IP-IA}
Supposons que pour $i=1,\ldots , r$, la reprŽsentation $\pi_i$ de $G_{n_i}$ soit un $\kappa$-relvement d'une reprŽsentation irr\'eductible $\tau_i$ de $H_{m_i}$, avec $n_i=m_id$, et que les reprŽsentations 
$\tau= \tau_1\times\cdots\times \tau_r$ de $H$ et $\pi=\pi_1\times \cdots \times \pi_r$ soient irrŽductibles. Alors $\pi$ est un $\kappa$-relvement de $\tau$ et il existe une racine de l'unitŽ $\zeta$, qui ne dŽpend ni des $\tau_i$ ni des $\pi_i$, telle que 
$$c(\tau,\pi,A)= \zeta \prod_{i=1}^r c(\tau_i,\pi_i,A_i)\ptf$$  
\end{proposition}

\subsection{OpŽrateurs d'entrelacement et multiplicitŽ $1$.}\label{op entrelac et mult 1} Soit $\pi$ une reprŽsentation de $G$ de longueur finie et soit $\pi_0$ un sous-quotient irrŽductible de $\pi$ de multiplicitŽ $1$ (dans $\pi$), \cad que le coefficient de la classe d'isomorphisme de $\pi_0$ dans la semi-simplifiŽe de $\pi$ est $1$. Par dŽfinition du sous-quotient, il existe une paire $(W\subset X)$ de sous-espaces vectoriels $G$-stables de $V=V_\pi$ telle que la reprŽsentation $\pi_{X/W}$ de $G$ d'espace $X/W$ dŽduite de $\pi$ (par restriction et passage au quotient) soit isomorphe ˆ $\pi_0$. 
Si de plus $X$ est maximal pour cette propriŽtŽ, ce que l'on suppose dŽsormais, alors on a le 

\begin{lemme}[{\cite[7.1]{BH}}]\label{[BH,7.1]}
\begin{enumerate} 
\item[(i)] Pour toute paire $(W'\subset X')$ de sous-espaces $G$-stables de $V$ telle que $\pi_{X'/W'}\simeq \pi_0$, on a 
$X'\subset X$ et $W'\subset W$, et l'inclusion $X'\subset X$ induit un isomorphisme $\pi_{X'/W'}\simeq \pi_{X/W}$.
\item[(ii)] Pour toute paire $(W'\subset X')$ de sous-espaces $G$-stables de $V$ telle que $\pi_{X'/W'}\simeq \pi_0$ avec $X'$ maximal, on a $X'=X$ et $W'=W$. 
\item[(iii)] Pour tout $\phi \in \mathrm{Aut}_G(\pi)$, on a $\phi (X)= X$ et $\phi(W)=W$. 
\end{enumerate} 
\end{lemme}

La paire $(W\subset X)$ avec $X$ maximal et $\pi_{X/W}\simeq \pi_0$ est donc dŽterminŽe de manire unique. 
En particulier si $\phi\in \mathrm{Aut}_G(\pi)$, alors $\phi$ induit (par restriction et passage au quotient) un opŽrateur 
$\phi_{X/W}\in \mathrm{Aut}_G(\pi_{X/W})$. Choisissons un isomorphisme $\beta \in \mathrm{Isom}_G(\pi_0, \pi_{X/W})$. 
Pour $\phi \in \mathrm{Aut}_G(\pi)$, posons 
$$\phi_{\pi_0}= \beta^{-1}\circ \phi_{X/W} \circ \beta \in \mathrm{Aut}_G(\pi_0)\ptf$$ 
D'aprs le lemme de Schur (pour $\pi_0$), $\phi_{\pi_0}$ ne dŽpend pas du choix de $\beta$. 

\begin{remark}
\textup{Soit $\pi^\star$ une reprŽsentation de $G$. Si $u\in \mathrm{Hom}_G(\pi,\pi^\star)$, pour toute paire $(W'\subset X')$ de sous-espaces $G$-stables de $V$ telle que $\pi_{X'/W'}\simeq \pi_0$, 
$u$ induit un $G$-morphisme $u_{X'/W'}:X'/W' \rightarrow u(X')/u(W')$. En le composant avec 
un $\beta'\in \mathrm{Isom}_{G}(\pi_0,\pi_{X'/W'})$, on obtient un opŽrateur $$u_{X'/W'}\circ \beta' \in \mathrm{Hom}_G(\pi_0, \pi^\star_{u(X')/u(W')})$$ 
qui est un isomorphisme si et seulement si $u(V_\pi)$ contient un sous-quotient irrŽductible isomorphe ˆ $\pi_0$. 
D'autre part on a des $G$-morphismes canoniques 
$i:X'/W'\rightarrow X/W$ et $j:u(X')/u(W') \rightarrow u(X)/u(W)$ qui sont tous les deux des isomorphismes, et 
le $G$-morphisme $j\circ u_{X'/W'}\circ i^{-1}$ 
co\"{\i}ncide avec $u_{X/W}$.
\hfill $\blacksquare$ }
\end{remark}

La reprŽsentation $\kappa \pi$ est elle aussi de longueur finie. Tout sous-espace de $V\;(=V_\pi = V_{\kappa\pi})$ est $G$-stable (pour l'action de $\pi$) si et seulement 
s'il est $G$-stable pour l'action de $\kappa\pi$; 
et si $(W'\subset X')$ est une paire de sous-espaces $G$-stables de $V$, alors on a $\kappa \pi_{X'/W'}= (\kappa\pi)_{X'/W'}$. 
Soit $\pi_0$ un sous-quotient irrŽductible de $\pi$ (on n'a pas besoin ici de la multiplicitŽ $1$). Si $\beta\in \mathrm{Isom}_G(\pi_0, \pi_{X'/W'})$ pour une paire $W'\subset X'$ de sous-espaces 
$G$-sables de $V$, alors $\beta \in \mathrm{Isom}_G(\kappa\pi_0 , \kappa\pi_{X'/W'})$. Supposons de plus que $\pi$ soit $\kappa$-stable et soit $A\in \mathrm{Isom}_G(\kappa\pi,\pi)$. L'application 
$$(W'\subset X') \mapsto (A(W')\subset A(X'))$$ est une bijection de l'ensemble des paires $(W'\subset X')$ de sous-espaces $G$-stables de $V$ telles que 
$\pi_{X'/W'}\simeq  \kappa \pi_0$ sur l'ensemble des paires $(W''\subset X'')$ de sous-espaces $G$-stables de $V$ telles que $\pi_{X''/W''}\simeq \pi_0$. 

Soit $\pi_0$ un sous-quotient irrŽductible $\kappa$-stable de $\pi$ de multiplicitŽ $1$. Soit $(W\subset X)$ l'unique paire de sous-espaces $G$-stables de $V$ telle que $\pi_{X/W}\simeq \pi_0$ avec $X$ maximal. Choisissons un isomorphisme $\beta\in \mathrm{Isom}_G(\kappa \pi_0, \pi_{X/W})$. 
La paire $(A(W)\subset A(X))$ vŽrifiant $\pi_{A(X)/A(W)}\simeq \pi_0$, on a $A(X)=X$ et $A(W)=W$. 
L'opŽrateur $A$ induit (par restriction et passage au quotient) un isomorphisme $A_{X/W}\in \mathrm{Isom}_G(\kappa \pi_{X/W}, \pi_{X/W})$. 
L'opŽrateur $$A_{\pi_0}= \beta^{-1} \circ A_{U/W} \circ \beta \in \mathrm{Isom}_G(\kappa\pi_0,\pi_0)$$ 
ne dŽpend pas du choix de $\beta$ (lemme de Schur). On dira que $A_{\pi_0}$ est l'opŽrateur (sur l'espace de $\pi_0$) \textit{dŽduit de $A$ 
par multiplicitŽ $1$}. 

\vskip1mm
{\bf $\bullet$ Induction parabolique et multiplicitŽ $1$.} Soit $M= G_{n_1}\times\cdots \times G_{n_r}$ un facteur de Levi standard de $G$. 
Soit $\pi_M=\pi_1\otimes \cdots \otimes \pi_r$ une reprŽsentation de $M$. On suppose que $\pi_M$ est $\kappa$-stable, \cad que $\kappa\pi_i\simeq \pi_i$ 
($i=1,\ldots , r$). Soit $B= B_1\otimes \cdots \otimes B_r\in \mathrm{Isom}_G(\kappa\pi_M,\pi_M)$. Notons $\pi$ l'induite parabolique 
$\pi_1\times \cdots \times \pi_r$ et $A\in \mathrm{Isom}_G(\kappa\pi,\pi)$ l'opŽrateur $B_1\times\cdots \times B_r$. 
Supposons de plus que $\pi_M$ soit de longueur finie (ce qui entra"ne que $\pi$ l'est aussi) et que $\pi$ ait un sous-quotient irŽductible $\pi_0$ qui soit $\kappa$-stable et de multiplicitŽ $1$. Alors $A$ dŽfinit comme plus haut un opŽrateur $A_{\pi_0}\in \mathrm{Isom}_G(\kappa\pi_0,\pi_0)$, que l'on notera aussi $B_{\pi_0}$. On dire que $B_{\pi_0}$ est l'opŽrateur (sur l'espace de $\pi_0$) \textit{dŽduit de $B$ par induction parabolique et multiplicitŽ $1$}. 

Soient $M_1\subset M_2$ deux facteurs de Levi standard de $G$. Pour $i=1,\,2$, 
on pose $P_i = P_{M_i}$ et $U_i= U_{Pi}$. Alors $P_1\subset P_2$ et $P_{1,2}= P_1\cap M_2$ est un sous-groupe parabolique (standard) de 
$M_2$. On pose $U_{1,2}= U_{P_{1,2}}\;(=U_1\cap M_2)$. On note $i_{P_{1,2}}^{M_2}$ le foncteur induction parabolique normalisŽe dŽfini par rapport au sous-groupe parabolique $P_{1,2}$ de $M_2$. 
Soit $\pi_{M_1}$ une reprŽsentation de $M_1$. Posons $$\pi_{M_2}= i_{P_{1,2}}^{M_2}(\pi_{M_1})\vgq \pi'= i_{P_2}^G(\pi_{M_2})\qhq{et}\pi = i_{P_1}^G(\pi_{M_1})\ptf$$ On sait (par transitivitŽ des foncteurs induction parabolique) que $\pi'\simeq \pi$. PrŽcisons cet isomorphisme. L'espace $V_{\pi_{M_2}}$ de $\pi_{M_2}$ est formŽ des fonctions $\varphi_2: M_2 \rightarrow V_{\pi_{M_1}}$ telles que 
pour tout $(m_1,u_1,m_2)\in M_1\times (U_1\cap M_2)\times M_2$, on ait  
$$\varphi_2(m_1u_1m_2)= \bs{\delta}_{P_{1,2}}(m_1)^{\frac{1}{2}}\pi_1(m_1)(\varphi_2(m_2))\pvg$$ et $\pi_{M_2}$ agit sur cet espace par translations ˆ droite.  
L'espace $V_{\pi'} $ de $\pi'$ est formŽ des fonctions $\varphi': G \rightarrow V_{\pi_{M_2}}$ telles que pour tout 
$(m_2,u_2,g)\in M_2\times U_2\times G$, on ait 
$$\varphi'(m_2u_2g) = \bs{\delta}_{P_2}(m_2)^{\frac{1}{2}}\pi_2(m_2)(\varphi'(g))\pvg$$ et $G$ agit sur cet espace par translations ˆ droite. 
L'espace $V_{\pi}$ de $\pi$ et l'action de $G$ sur $V_\pi$ sont dŽfinis de la mme manire (en remplaant l'indice $2$ par $1$). 
Soit $\bs{h}: V_{\pi'} \rightarrow V_\pi$ l'application dŽfinie par 
$$(\bs{h}\varphi')(g)= \varphi'(g)(1)\qhq{pour tout} g\in G\ptf$$ Puisque pour $p_1\in P_1$, en Žcrivant $p_1=p_{1,2}u_2$ avec $p_{1,2}\in P_{1,2}$ et $u_2\in U_2$, on a $\bs{\delta}_{P_1}(p_1)= \bs{\delta}_{P_{1,2}}(p_{1,2})\bs{\delta}_{P_2}(p_1)$, la fonction 
$\bs{h}\varphi'$ est bien dans l'espace $V_\pi$. De plus (vŽrification facile) $\bs{h}$ entrelace $\pi'$ et $\pi$, \ie $\bs{h}\in \mathrm{Hom}_G(\pi'\!,\pi)$. 
Soit $\bs{h}': V_\pi \rightarrow V_{\pi'}$ l'application dŽfinie par 
$$(\bs{h}'\varphi)(g)(m_2)= \bs{\delta}_{P_2}(m_2)^{-\frac{1}{2}}\varphi(m_2g) \qhq{pour tout} (g,m_2)\in G\times M_2\ptf$$ 
On a (vŽrifications faciles) $\bs{h'}\in \mathrm{Hom}_G(\pi,\pi')$, $\bs{h}'\circ \bs{h}= \mathrm{Id}_{V_{\pi'}}$ et $ \bs{h}\circ \bs{h}'= \mathrm{Id}_{V_{\pi}}$. Par consŽquent $\bs{h}\in \mathrm{Isom}_G(\pi'\!,\pi)$ et $\bs{h}'=\bs{h}^{-1}$. 
Supposons de plus que $\pi_{M_1}$ soit 
$\kappa$-stable et soit $B_1\in \mathrm{Isom}_G(\kappa\pi_1,\pi_1)$. Posons 
$$B_2= i_{P_{1,2}}^{M_2}(B_1)\vgq A'= i_{P_2}^G(B_2)\qhq{et}A=i_{P_1}^G(B_1)\ptf$$

\begin{lemme}
\begin{enumerate}
\item[(i)]$\bs{h} \circ A'= A \circ \bs{h}\;(\in \mathrm{Isom}_G(\kappa \pi',\pi))$. 
\item[(ii)]Supposons que $\pi_{M_1}$ soit de longueur finie et que $\pi\;(=i_{P_1}^G(\pi_{M_1}))$ ait un sous-quotient irrŽductible $\pi_0$ qui soit $\kappa$-stable et de multiplicitŽ $1$. Soit $\pi_2$ le sous-quotient irrŽductible de $\pi_{M_2}$ tel que $\pi_0$ soit isomorphe ˆ un sous-quotient irrŽductible de 
$i_{P_2}^G(\pi_{2})$. Si $\pi_2$ est $\kappa$-stable, on a l'ŽgalitŽ 
$$((B_1)_{\pi_2}))_{\pi_0}= (B_1)_{\pi_0}\ptf$$   
\end{enumerate} 
\end{lemme}

\begin{demo} 
Prouvons (i). Pour $\varphi'\in V_{\pi'}$ et $g\in G$, on a 
$$((A\circ \bs{h})\varphi')(g)=\kappa(g) B_1(\bs{h}\varphi'(g))=\kappa(g) B_1(\varphi'(g)(1))$$ 
et $$((\bs{h}\circ A')\varphi')(g) = A'\varphi'(g)(1)= \kappa(g) B_2(\varphi'(g))(1)= \kappa(g)\kappa(1)B_1(\varphi'(g)(1))\ptf$$ 

Prouvons (ii). Observons que (par dŽfinition) $\pi_2$ est de multiplicitŽ $1$ dans $\pi_{M_2}$. Soit $(W_2\subset X_2)$ la paire de sous-espaces $G$-stables de $V_{\pi_{M_2}}$ avec $X_2$ maximal telle que $\pi_2\simeq \pi_{X_2/W_2}$. On a donc 
une suite exacte courte de $M_2$-modules 
$$0 \rightarrow W_2 \rightarrow X_2 \rightarrow V_{\pi_2} \rightarrow 0\ptf$$ 
Supposons de plus que $\pi_2$ soit $\kappa$-stable. 
Alors on a (cf. plus haut) $B_2(X_2)=X_2$ et $B_2(W_2)=W_2$. Posons $\wt{\pi}=i_{P_2}^G(\pi_2)$; c'est un sous-quotient de $\pi$ et $\pi_0$ est un sous-quotient irrŽductible de $\wt{\pi}$ de multplicitŽ $1$.  
Par exactitude du foncteur induction parabolique $i_{P_2}^G$, on a la suite exacte courte de $G$-modules 
$$0 \rightarrow i_{P_2}^G(W_2) \rightarrow i_{P_2}^{G}(X_2) \buildrel q \over{\longrightarrow} V_{\wt{\pi}}\rightarrow 0\ptf $$  
Soit $(W\subset X)$ la paire de sous-espaces vectoriels $G$-stables 
de $V_{\wt{\pi}}$ avec $X$ maximal telle que $\pi_0 \simeq \wt{\pi}_{X/W}$. On a la cha"ne d'inclusions 
$$i_{P_2}^G(W_2) \subset q^{-1}(W) \subset q^{-1}(X) \subset i_{P_2}^G(X_2)$$ et le $G$-isomorphisme 
$\overline{q} : i_{P_2}^G(X_2)/ i_{P_2}^G(W_2)\rightarrow V_{\wt{\pi}}$ dŽduit de $q$ par passage au quotient 
induit par restriction des $G$-isomorphismes 
$$\overline{q}_X: q^{-1}(X)/i_{P_2}^G(W_2) \rightarrow X\qhq{et} \overline{q}_W: q^{-1}(W)/i_{P_2}^G(W_2)\rightarrow W\ptf$$ 
Posons $V'= i_{P_2}^G(X_2)\subset V_{\pi'}$ et $V= \bs{h}V'\subset V_{\pi}$. Puisque $B_2(X_2)=X_2$, on a $A'(V')=V'$; et d'aprs le point (i), on a 
$$A(V)=V\qhq{et}\bs{h}\circ A' \vert_{V'}= A\circ \bs{h}\vert_{V'}\ptf\leqno{(*)}$$ 
D'aprs \ref{[BH,7.1]}\,(i), $q^{-1}(X)$ est le sous-espace maximal de 
$V'$ admettant $\pi_0$ comme quotient. Sur $Y=q^{-1}(X)/q^{-1}(W)$, on a l'opŽrateur $(i_{P_2}^G(B_1)_{\pi_2})_{Y}$ dŽduit 
de $i_{P_2}^G((B_1)_{\pi_2})$, ou ce qui revient au mme de $A'=i_{P_2}^G(B_2)$, par restriction et passage au quotient. 
Il dŽfinit l'opŽrateur $((B_1)_{\pi_2})_{\pi_0}$ via n'importe quel $G$-isomorphisme $Y\simeq V_{\pi_0}$. 
D'autre part comme $A'(X)=X$ et $A'(W)=W$, on a aussi $A(\bs{h}X)=\bs{h}X$ et $A(\bs{h}W)= \bs{h}W$. Sur $\bs{h}X/\bs{h}W$, on a l'opŽrateur $A_{\bs{h}X/\bs{h}W}$ dŽduit de $A$ par restriction et passage au quotient. Il dŽfinit l'opŽrateur $(B_1)_{\pi_0}$ via n'importe quel 
isomorphisme $\bs{h}X/\bs{h}W\simeq \pi_0$. On conclut gr‰ce ˆ l'ŽgalitŽ 
$A'=\bs{h}^{-1}\circ A \circ \bs{h}$ (point (i)).\hfill$\square$
\end{demo}
 
\subsection{Normalisation des opŽrateurs d'entrelacement.}\label{norm op entrelac} On commence par dŽfinir, pour une reprŽsentation irrŽductible 
$\kappa$-stable gŽnŽrique $\pi$ de $G$, un opŽrateur normalisŽ $\Agen_\pi\in \mathrm{Isom}_G(\kappa\pi,\pi)$. Gr‰ce ˆ la dŽcomposition de Langlands, on dŽfinit 
ensuite, pour toute reprŽsentation irrŽductible $\kappa$-stable $\pi$ de $G$, un opŽrateur normalisŽ $A_\pi\in \mathrm{Isom}_G(\kappa\pi,\pi)$. Enfin on vŽrifie que ces deux normalisations sont cohŽrentes: si 
$\pi$ est gŽnŽrique, alors $A_\pi=\Agen_\pi$. 

\vskip1mm
{\bf $\bullet$ DŽfinition de $\Agen_\pi= \Agen_{\pi,\psi}$ ($\pi$ gŽnŽrique).} On fixe un caractre additif non trivial $\psi$ de $F$. Il dŽfinit un caractre $\theta_\psi$ de $U_0$: 
$$\theta_\psi(u) = \psi(\textstyle{\sum_{i=1}^{n-1}}u_{i,i+1})\qhq{pour} u= (u_{i,j})_{1\leq i,j\leq n}\in U_0\ptf$$ Si $\pi$ est une reprŽsentation de $G$, on note $\ES{D}(\pi)=\ES{D}(\pi,\psi)$ l'espace des fonctions $\lambda: V_\pi\rightarrow \mathbb{C}$ telles que 
$$\lambda(\pi(u)v)= \theta_\psi(u)\lambda(v)\qhq{pour tout} (u,v)\in U_0\times V_\pi\ptf$$ 
Un ŽlŽment de $\ES{D}(\pi)\smallsetminus \{0\}$ est appelŽ \textit{fonctionnelle de Whittaker} pour $\pi$ (relativement ˆ $\psi$). On sait que $\dim_{\mbb{C}}(\ES{D}(\pi))\leq 1$ si $\pi$ est irrŽductible. 
La reprŽsentation $\pi$ est dite \textit{gŽnŽrique} s'il existe une fonctionnelle de Whittaker pour $\pi$, \ie si $\dim_{\mathbb{C}}(\ES{D}(\pi))>0$; auquel cas 
$\dim_{\mathbb{C}}(\ES{D}(\pi))=1$ si $\pi$ est irrŽductible. 
La notion de gŽnŽricitŽ ne dŽpend pas de $\psi$. 

Si $\pi$ est une reprŽsentation gŽnŽrique de $G$, alors pour toute fonctionnelle de Whittaker 
$\lambda\in \ES{D}(\pi)$, on a $$\lambda(\kappa\pi(u)(v))= \kappa(u)\lambda(\pi(u)(v))= \theta_\psi(u)\lambda(v)\qhq{pour tout}(u,v)\in U_0\times V_\pi\ptf$$ 
Par consŽquent $\lambda$ est aussi une fonctionnelle de Whittaker pour $\kappa\pi$. Si de plus $\pi$ est irrŽductible et $\kappa$-stable, on peut normaliser l'opŽrateur $A\in \mathrm{Isom}_G(\kappa\pi,\pi)$ en imposant la condition $\lambda \circ A= \lambda$ pour une (\ie pour toute) fonctionnelle de Whittaker 
$\lambda$ pour $\pi$. On note $$\Agen_\pi= \Agen_{\pi,\psi}$$ l'opŽrateur $A$ ainsi normalisŽ. Cette \textit{normalisation Whittaker} est compatible aux isomorphismes: si $\pi'$ est une reprŽsentation irrŽductible gŽnŽrique $\kappa$-stable de $G$ isomorphe ˆ $\pi$, 
et si $\phi\in \mathrm{Isom}_G(\pi,\pi')$, l'application $\ES{D}(\pi) \rightarrow \ES{D}(\pi'),\, \lambda \mapsto \lambda \circ \phi^{-1}$ est un isomorphisme. On en dŽduit 
l'ŽgalitŽ $$\Agen_{\pi'}\circ \phi = \phi \circ \Agen_\pi\ptf$$

\begin{remark}
\textup{La normalisation $\Agen_\pi= \Agen_{\pi,\psi}$ dŽpend de $\psi$. PrŽcisŽment, tout autre caractre additif non trivial de $F$ est de la forme 
$a\mapsto \psi_a(x)= \psi(ax)$ pour un $a\in F^\times$. En notant $t_a$ l'ŽlŽment $(1,a, \ldots ,a^{n-1})$ de $A_0= (F^\times)^n$, on vŽrifie que 
l'application $\ES{D}(\pi,\psi)\rightarrow \ES{D}(\pi,\psi_a),\, \lambda \mapsto \lambda \circ \pi(t_a^{-1})$ est un isomorphisme. On en dŽduit l'ŽgalitŽ 
$$\Agen_{\pi,\psi_a}= \kappa(t_a)\Agen_{\pi,\psi}= \kappa(a)^{\frac{n(n-1)}{2}}\Agen_{\pi,\psi}\ptf$$ 
Comme $n= md$, $\kappa(t_a)=1$ si $m$ est pair ou si $d$ est impair, et $\kappa(t_a)= \pm 1$ dans tous les cas. \hfill $\blacksquare$}
\end{remark}

\vskip1mm
{\bf $\bullet$ DŽfinition de $A_\pi=A_{\pi,\psi}$.} Soit $\pi$ une reprŽsentation irrŽductible et $\kappa$-stable de $G$. On \'ecrit $\pi\simeq L(\wt{\pi}_1,\ldots , \wt{\pi}_r)$ pour des reprŽsentations irrŽductibles essentiellement tempŽrŽes $\wt{\pi}_i$ de $G_{n_i}$, avec $\sum_{i=1}^r n_i= n$. Pour $i=1,\ldots ,r$, la reprŽsentation $\wt{\pi}_i$ est gŽnŽrique et $\kappa$-stable (par unicitŽ du quotient de Langlands), donc munie d'un opŽrateur normalisŽ 
$\Agen_{\wt{\pi}_i}\in \mathrm{Isom}_{G_{n_i}}(\kappa\wt{\pi}_i, \wt{\pi}_i)$. La reprŽsentation $\pi_M= \wt{\pi}_1\otimes \cdots \otimes \wt{\pi}_r$ de 
$M=G_{n_1}\times \cdots \times G_{n_r}$ est irrŽductible et $\kappa$-stable, et l'opŽrateur 
$\Agen_{\pi_M}\bydef \Agen_{\wt{\pi}_1}\otimes \cdots \otimes \Agen_{\wt{\pi}_r}$ appartient ˆ $\mathrm{Isom}_M(\kappa\pi_M, \pi_M)$. Posons 
$$\rho=i_{P_M}^G(\pi_M)\qhq{et}A= i_{P_M}^G(\Agen_{\pi_M})\in \mathrm{Isom}_G(\kappa \rho,\rho)\ptf$$ Notons $\overline{\rho}= L(\wt{\pi}_1,\ldots , \wt{\pi}_r)$ l'unique quotient irrŽductible de $\rho$. Puisque $\rho$ est $\kappa$-stable, $\overline{\rho}$ l'est aussi et $A$ induit par passage au quotient un opŽrateur 
$\overline{A}\in \mathrm{Isom}_G(\kappa \overline{\rho}, \overline{\rho})$. 

\begin{remark}\textup{Le quotient de Langlands est un cas particulier d'induction parabolique et multiplicitŽ $1$: avec les notations de \ref{op entrelac et mult 1}, on a $\overline{A}= (\Agen_{\pi_M})_{\overline{\rho}}$.\hfill $\blacksquare$}
\end{remark}

Choisissons un isomorphisme 
$u\in \mathrm{Isom}_G(\pi,\overline{\rho})= \mathrm{Isom}_G(\kappa\pi,\kappa\overline{\rho})$. L'opŽrateur 
$u^{-1}\circ \overline{A}\circ u \in \mathrm{Isom}_G(\kappa\pi , \pi)$ est bien dŽfini, \cad qu'il ne dŽpend pas du choix de 
$u$ (d'aprs le lemme de Schur). On le note $$A_\pi = A_{\pi,\psi}\ptf$$ Comme pour la normalisation Whittaker, cette normalisation est compatible aux isomorphismes: 
si $\pi'$ est une reprŽsentation irrŽductible gŽnŽrique $\kappa$-stable de $G$ isomorphe ˆ $\pi$, 
et si $\phi\in \mathrm{Isom}_G(\pi,\pi')$, on a 
l'ŽgalitŽ $$A_{\pi'}\circ \phi = \phi \circ A_\pi\ptf$$

\vskip1mm
{\bf $\bullet$ L'ŽgalitŽ $A_\pi= \Agen_\pi$ ($\pi$ gŽnŽrique).} Soit $\pi$ une reprŽsentation irrŽductible $\kappa$-stable et g\'enŽrique de $G$. On veut prouver l'ŽgalitŽ $A_\pi=\Agen_\pi$. On Žcrit 
$\pi\simeq L(\wt{\pi},\ldots ,\wt{\pi}_r)$ o $\wt{\pi}_i$ est une reprŽsentation irrŽductible essentiellement tempŽrŽe de $G_{n_i}$, avec $\sum_{i=1}^r n_i=n$. Puisque $\pi$ est gŽnŽrique, l'induite 
$\wt{\pi}_1\times \cdots \times \wt{\pi}_r$ est irrŽductible, et $\pi\simeq \wt{\pi}_1\times \cdots \times \wt{\pi}_r$. Puisque les opŽrateurs $\Agen_\pi$ et $A_\pi$ sont compatibles aux isomorphismes, on peut supposer que 
$\pi= \wt{\pi}_1\times \cdots \times \wt{\pi}_r$. Pour prouver l'ŽgalitŽ $A_\pi=\Agen_\pi$, il suffit de prouver que $\lambda \circ A_\pi =\lambda$ 
pour une (\ie pour toute) fonctionnelle de Whittaker $\lambda$ pour $\pi$. Pour $i=1,\ldots , r$, choisissons une fonctionnelle de Whittaker 
$\lambda_i$ pour $\wt{\pi}_i$. Notons $\lambda_M$ la fonctionnelle linŽaire $\lambda_1\otimes\cdots\otimes \lambda_r$ sur l'espace 
de $\pi_M=\wt{\pi}_1\otimes\cdots \otimes \wt{\pi}_r$, avec $M=G_{n_1}\times \cdots \times G_{n_r}$. 
Posons $M'= G_{n_r}\times G_{n_{r-1}}\times\cdots \times G_{n_1}$ et $P'=P_{M'}$.  
Pour $f\in V_\pi$, posons 
$$\lambda(f)= \int_{U_{P'}}\lambda_M(f(wu))\overline{\theta_\psi(u)}\d u$$ o $\d u$ est une mesure de Haar sur $U_{P'}$ et $w$ est l'ŽlŽment anti-diagonal par blocs dŽfini par $$w= \left(\begin{array}{cccc} 0 & \ldots & 0 & 1_{n_1}\\
\vdots & & 1_{n_2} & 0 \\
0 & \Ddots && \\
1_{n_r} & 0 & \ldots & 0
\end{array}\right)\ptf$$ 
Observons que $wM' w^{-1}= M$ et que l'intŽgrale dŽfinit un ŽlŽment de $\ES{D}(\pi,\psi)$, pourvu qu'elle converge.  
D'aprs \cite{JS2}, l'intŽgrale est absolument convergente et dŽfinit une fonctionnelle de Whittaker pour $\pi$; \ie $\lambda\in \ES{D}(\pi,\psi)\smallsetminus \{0\}$. 
Pour $u\in U_{P'}$, on a 
$$A_\pi(f)(wu)=\kappa(wu) \Agen_{\pi_M}(f(wu))= \Agen_{\pi_M}(f(wu))$$ 
car $ \det(wu)=1$; et comme (par dŽfinition) $\lambda_M \circ \Agen_{\pi_M}= \lambda_M$, on obtient que 
$\lambda_M(A_\pi(f)(wu))=\lambda_M(f(wu))$. D'o l'ŽgalitŽ cherchŽe: $\lambda\circ A_\pi = \lambda$.

\begin{remark}\label{op unit}
\textup{Soit $\pi$ une reprŽsentation irrŽductible $\kappa$-stable de $G$. 
\begin{enumerate}
\item[(i)] Puisque $A_\pi^d\;(= A_\pi\circ \cdots \circ A_\pi)$ appartient ˆ $\mathrm{Aut}_G(\pi)$, d'aprs le lemme de Schur, on a $A_\pi^d= \mathrm{Id}_{V_\pi}$ 
(c'est Žvident si $\pi$ est gŽnŽrique; en gŽnŽral, on se ramne au cas gŽnŽrique via la dŽcomposition de Langlands).
\item[(ii)] Si $\pi$ est unitaire et si $(\cdot,\cdot)$ est un produit scalaire hermitien $G$-invariant sur $V_\pi$, alors $(\cdot,\cdot)$ est $A_\pi$-invariant 
(cf. \cite[2.2, remark]{BH}). \hfill $\blacksquare$
\end{enumerate}
}
\end{remark}

\section{Extension au cas d'une algbre cyclique}\label{le cas d'une algbre cyclique}

Pour les applications globales, on Žtend la thŽorie de la section \ref{le cas n-a} 
au cas d'une $F$-algbre cyclique $E$ de degrŽ fini $d$ (on suppose toujours que $F$ est un corps local non archimŽdien).  

\subsection{Notations.}\label{notations ext cyclique} On pose $\Gamma=\Gamma(E/F)\;(\simeq \mathbb{Z}/d\mathbb{Z})$. 
L'algbre $E$ s'Žcrit comme un produit de corps $E=E_1\times \cdots \times E_r$ o $E_i/F$ est une extension cyclique de degrŽ $s= \frac{d}{r}$. On suppose fixŽ un gŽnŽrateur $\sigma$ de $\Gamma$ tel que $\sigma E_i=E_{i+1}$ ($i=1,\ldots , r-1$) et $\sigma E_r= E_1$. 
Pour $i=1,\ldots , r$, on a $\sigma^{r} E_i=E_i$ et $\sigma^r$ est un gŽnŽrateur du groupe $\Gamma_i= \Gamma(E_i/F)$.  
Le corps $F$, plongŽ diagonalement dans $E$, co\"{\i}ncide avec l'ensemble des ŽlŽments de $E$ fixŽs par $\Gamma$. 

On dispose d'une application norme $\mathrm{N}_{E/F}: E^\times \rightarrow F^\times$ dŽfinie, comme dans le cas o $E=E_1$, par $\mathrm{N}_{E/F}(x)= \prod_{i=0}^{d-1}\sigma^i x$. Pour $x= (x_1,\ldots , x_r)$, $x_i \in E_i^\times$, on a $\mathrm{N}_{E/F}(x)= \mathrm{N}_{E_1/F}(\prod_{i=1}^r \sigma^{r-i+1} x_i)$. En particulier le sous-groupe des normes 
$\mathrm{N}_{E/F}(E^\times)\subset F^\times$ co\"{\i}ncide avec $\mathrm{N}_{E_1/F}(E_1^\times)$. Soit $\kappa$ un caractre de $F^\times$ de noyau $\mathrm{N}_{E/F}(E^\times)$, \ie 
un caractre de $F^\times$ correspondant ˆ l'extension $E_1/F$ par la thŽorie du corps de classes local. 

On fixe un entier $m\geq 1$. On pose $H=\GLm(E)$ et $G=\GLn(F)$, $n=md$. On a donc 
$H= \prod_{i=1}^r\GLm(E_i)$. On fixe une $F$-base de $E^m$, ce qui donne un plongement $\varphi$ du groupe $H$ dans le groupe $G$. 
Comme dans le cas o $E=E_1$, changer de base revient ˆ conjuguer $\varphi$ par un ŽlŽment de $G$. On identifie $H$ au sous-groupe $\varphi(H)$ de $G$. 

\subsection{Facteurs de transfert et fonctions concordantes.}\label{facteurs de transfert bis} On dŽfinit 
les fonctions discriminant $D_G: G \rightarrow F$ et $D_H: H \rightarrow E$ comme en \ref{facteurs de transfert}. On note $\Greg$ l'ensemble des 
$x\in G$ tels que $D_G(x)\neq 0$, et $\Hreg\;(\supset H \cap \Greg)$ l'ensemble des $y\in H$ tels que $D_H(y)\neq 0$. Pour $x\in H$, on dŽfinit 
l'ŽlŽment $\wt{\Delta}(x)\in E$ comme en \cite[2.2]{HL2}; il dŽpend du choix de $\sigma$ et vŽrifie 
$$\sigma \wt{\Delta}(x)= \wt{\Delta}(\sigma x) = (-1)^{m(d-1)}\wt{\Delta}(x)\ptf$$ On choisit un ŽlŽment $e=(e_1,\ldots ,e_r)\in E^\times$ tel que $\sigma e = (-1)^{m(d-1)}e$ et pour $x\in H\cap \Greg$, on pose 
$$\Delta(x)= \kappa (e \wt{\Delta}(x))\in \mathbb{C}^\times\ptf$$
La remarque \ref{facteurs de transfert, propriŽtŽ et choix} est encore vraie ici; en particulier si la $F$-algbre cyclique $E/F$ est non ramifiŽe, \ie si l'extension 
$E_1/F$ est non ramifiŽe, le facteur de transfert $\Delta$ ne dŽpend ni du choix de $\sigma$ ni du choix de $e\in \mathfrak{o}_E^\times = \prod_{i=1}^r \mathfrak{o}_{E_i}^\times$. 

\vskip1mm
On fixe des mesures de Haar $\d g$ sur $G$ et $\d h$ sur $H$. Pour $f\in \CG$ et $\gamma\in \Greg$, la $\kappa$-intŽgrale orbitale $I^G_\kappa(f,\gamma)$ a dŽjˆ ŽtŽ dŽfinie (cf. \ref{facteurs de transfert}). Pour $\phi\in \CH$ et $\delta\in \Hreg$, on dŽfinit comme en \ref{facteurs de transfert} l'intŽgrale orbitale  
$I^H(\phi, \delta)$, o la valeur absolue normalisŽe $\vert\, \vert_E$ sur $E$ est donnŽe par $\vert x\vert_E = \vert x_1\vert_ {E_1}\cdots \vert x_r\vert_{E_r}$ pour 
$x=(x_1,\ldots ,x_r)\in E_1\times \cdots \times E_r$. La notion de fonctions \textit{concordantes} est la mme qu'en \ref{facteurs de transfert}. 
Observons que comme dans le cas le cas o $E=E_1$, pour que $f\in \CG$ et $\phi\in \CH$ soient concordantes, il est nŽcessaire que pour tout $\gamma \in \HGreg$ et tout $x\in G$ tels que l'ŽlŽment 
$\gamma'=x^{-1}\gamma x$ appartienne ˆ $ H$, on ait l'ŽgalitŽ $I^H(\phi,\gamma')= \kappa(x)^{-1} \frac{\Delta(\gamma')}{\Delta(\gamma)} I^H(\phi,\gamma)$.

\vskip1mm
Posons $H^\natural= \GLm(E_1)^{r}\;(= \GLm(E_1)\times\cdots \times \GLm(E_1))$ et $H_1= \mathrm{GL}_{mr}(E_1)$. 
L'application $$x=(x_1,\ldots , x_r) \mapsto \iota(x)=(x_1,\sigma^{-1} x_2, \ldots , \sigma^{1-r}x_r) $$ 
est un isomorphisme de $E$ sur $E^\natural=(E_1)^{r}$. Il induit un isomorphisme $h\mapsto \iota(h)$ de $H$ sur $H^\natural$. Pour $\phi\in \CH$, on note $\iota(\phi)\in \CHnat$ la fonction dŽfinie par $$\iota(\phi)(y)= \phi(\iota^{-1}(y))\ptf$$ 
Le groupe $H^\natural$ est un facteur de Levi standard de $H_1$. La $F$-base de $E^m$ fixŽe pour dŽfinir le plongement 
$\varphi$ de $H$ dans $G$, identifiŽe via $\iota$ ˆ une $F$-base de $(E_1)^{rm}$, dŽfinit un plongement $\varphi_1$ de $H_1$ dans $G$ dont la restriction 
ˆ $H^\natural$ co\"{\i}ncide avec $\varphi \circ \iota^{-1}$. On identifie 
$H_1$ au sous-groupe $\varphi_1(H_1)$ de $G$. 

Pour $x\in H_1$, on dŽfinit comme en \cite[3.2]{HH}, ˆ l'aide du gŽnŽrateur $\sigma_1\bydef \sigma^r$ de $\Gamma(E_1/F)$, un ŽlŽment $\wt{\Delta}_1(x)\in E_1$ qui vŽrifie $\sigma_1 \wt{\Delta}_1(x)= (-1)^{rm(s-1)}\wt{\Delta}_1(x)$. Les ŽlŽments de $H_1$ tels que $\wt{\Delta}_1(x)\neq 0$ sont prŽcisŽment ceux de $H_1\cap \Greg$. Pour dŽfinir le facteur de transfert $\Delta_1$, il nous faut choisir un ŽlŽment $e'_1\in E_1^\times$ tel que $\sigma_1e'_1= (-1)^{rm(s-1)}e'_1$; on pose alors 
$$\Delta_1(x)= \kappa(e'_1\wt{\Delta}_1(x))\qhq{pour tout}x\in H_1\cap \Greg\ptf$$ 
La condition $\sigma e = (-1)^{m(d-1)}e$ est Žquivalente ˆ 
$$
\left\{
\begin{array}{ll}e_{i+1} = (-1)^{im(d-1)}\sigma^i e_1 & (i=1,\ldots , r-1)\\
\sigma_1e_1= (-1)^{rm(d-1)} e_1
\end{array}\right..$$
Or comme $d=sr$, on a $$rm(d-1)\equiv rm(s-1)\quad(\mathrm{mod}\; 2\mathbb{Z})\ptf$$ Par consŽquent $(-1)^{rm(d-1)}= (-1)^{rm(s-1)}$ et l'on peut prendre $e'_1=e_1$. 
D'aprs la preuve de \cite[7.8, lemma]{HH}, pour $x\in H\cap \Greg$, l'ŽlŽment $\wt{\Delta}(x)$ appartient ˆ $\wt{\Delta}_1(\iota(x))N_{E_1/F}(E_1^\times)$. 
On en dŽduit que $$\Delta(x) = \Delta_1(\iota(x))\qhq{pour tout} x\in H\cap \Greg \ptf\leqno{(1)}$$ 
Comme en \ref{facteurs de transfert}, deux fonctions $f\in \CG$ et $\phi_1\in C^\infty_{\mathrm{c}}(H_1)$ sont dites \textit{concordantes} (relativement ˆ $E_1/F$) si 
$$\Delta_1(\gamma_1) I^G_\kappa(f,\gamma_1)= I^{H_1}(\phi_1,\gamma_1)\qhq{pour tout}\gamma_1\in H_1\cap \Greg \pvg$$ o la distribution $I^{H_1}(\cdot , \gamma_1)$ sur $H_1$ est dŽfinie comme plus haut via le choix d'une mesure de Haar $\d h_1$ sur $H_1$. 

Posons $K_1=\mathrm{GL}_{mr}(\mathfrak{o}_{E_1})$ et notons $P_1=H^\natural\ltimes U_1$ le sous-groupe parabolique standard de $H_1$ de composante de Levi $H^\natural$. 
Soit $(\phi_1)^{P_1}\in C^\infty_{\mathrm{c}}(H^\natural)$ le \textit{terme constant} de $\phi_1$ suivant $(K_1,P_1)$, dŽfini par 
$$(\phi_1)^{P_1}(\delta^\natural)= \bs{\delta}_{P_1}(y)^{\frac{1}{2}}\iint_{K_1\times U_1}\phi(k_1^{-1}\delta^\natural u_1 k_1)\d u_1 \d k_1\qhq{pour tout} \delta^\natural\in H^\natural\pvg$$ o 
les mesures de Haar $\d k_1$ sur $K_1$ et $\d u_1$ sur $U_1$ sont normalisŽes par 
$$\mathrm{vol}(K_1,\d k_1)=1=\mathrm{vol}(U_1\cap K_1,\d u_1)\ptf$$
Soit $\d h^\natural$ la mesure de Haar sur $H^\natural$ dŽduite de celle sur $H$ via $\iota$. On normalise la mesure 
de Haar $\d h_1$ sur $H_1$ de telle manire que pour tout 
$\delta^\natural\in H^\natural\cap (H_1)_{\mathrm{r\acute{e}g}}$, on ait a formule de descente: 
$$I^{H_1}(\phi_1,\delta^\natural)= I^{H^\natural\!}((\phi_1)^{P_1}\!,\delta^\natural)\ptf\leqno{(2)}$$ 

\begin{lemme}\label{existence transfert bis}
\begin{enumerate}
\item[(i)] Soit $(f,\phi_1)\in \CG\times \mathrm{C}^\infty_{\mathrm{c}}(H_1)$ une paire de fonctions concordantes (relativement ˆ $E_1/F$). La fonction $\phi =\iota^{-1}((\phi_1)^{P_1})\in \CH$ concorde avec $f$ (relativement ˆ $E/F$).
\item[(ii)] Soit $(\phi,\phi_1)\in \CH\times \mathrm{C}^\infty_{\mathrm{c}}(H_1)$ une paire de fonctions vŽrifiant les deux conditions suivantes:
\begin{itemize} 
\item pour tout $\gamma\in H \cap \Greg$, on a l'ŽgalitŽ $I^{H}(\phi,\gamma)= I^H(\iota^{-1}((\phi_1)^{P_1}),\gamma)$;
\item pour tout $\gamma_1 \in H_1 \cap \Greg$ et tout $x\in G$ tels que 
$\gamma'_1=x^{-1}\gamma_1 x\in H_1$, on a l'ŽgalitŽ $I^{H_1}(\phi_1,\gamma'_1)= \kappa(x)^{-1} \frac{\Delta_1(\gamma'_1)}{\Delta_1(\gamma_1)} I^{H_1}(\phi_1,\gamma_1)$. 
\end{itemize}
Alors il existe une fonction $f\in \CG$ qui concorde avec $\phi$. 
\end{enumerate}
 \end{lemme} 

\begin{demo} Soit $f\in \CG$. Si une fonction $\phi_1\in C^\infty_{\mathrm{c}}(H_1)$ concorde avec $f\in \CG$ (relativement ˆ $E_1/F$), d'aprs (1) et (2), 
pour tout $\gamma\in H\cap \Greg$, on a 
$$\Delta(\gamma)I^G_\kappa(f,\gamma)= I^{H'\!}((\phi_1)^{P_1}\!,\iota(\gamma))= I^H(\iota^{-1}((\phi_1)^{P_1}),\gamma)\ptf$$ 
Cela prouve (i). Cela prouve aussi (ii) car la seconde condition de (ii) assure que $\phi_1$ concorde avec une fonction $f\in \CG$ (relativement ˆ $E_1/F$). \hfill $\square$ 
\end{demo}

\subsection{$\kappa$-relvement.}\label{kappa-rel bis} 
La notion de $\kappa$-relvement (relativement ˆ $E/F$) de $H$ ˆ $G$  est dŽfinie comme en \ref{defkapparel}. C'est la composition de deux opŽrations: une induction parabolique 
de $H\simeq H^\natural=\GLm(E_1)^r$ ˆ $H_1=\mathrm{GL}_{mr}(E_1)$ composŽe avec un $\kappa$-relvement de $H_1$ ˆ $G$ (relativement ˆ $E_1/F$). 
PrŽcisŽment, si $\tau= \tau_1\otimes \cdots \otimes \tau_r$ est une reprŽsentation irrŽductible de $H$, on note $\tau^\natural= \tau^\natural_1\otimes \cdots \otimes \tau^\natural_r$ la reprŽsentation 
$\tau \circ \iota^{-1}$ de $H^\natural$, avec 
$\tau^\natural_i=(\tau_{i})^{\sigma^{i-1}}$. 
Soit $$\rho=\tau_1^\natural\times \cdots \times \tau^\natural_r\;(=i_{P_1}^{H_1}(\tau^\natural))$$ l'induite parabolique de 
$\tau^\natural$ ˆ $H_1$. En normalisant les mesures de Haar sur $H_1$ et $H^\natural$ comme en 
\ref{facteurs de transfert bis}\,(2), pour toute fonction $\phi_1\in C^\infty_{\mathrm{c}}(H_1)$, 
on a la formule de descente: 
$$\mathrm{tr}(\rho(\phi_1))= \mathrm{tr}(\tau^\natural((\phi_1)^{P_1}))\ptf\leqno{(1)}$$  

\begin{lemme}[{\cite[prop.~3.8]{HL2}}]\label{[HL] prop 3.8}
Supposons que la reprŽsentation $\rho$ de $H_1$ soit irrŽductible. Soit $\pi$ une reprŽsentation irrŽductible $\kappa$-stable de $G$. Pour que $\pi$ soit un $\kappa$-relvement de $\tau$ (relativement ˆ $E/F$), il faut et il suffit que ce soit un 
$\kappa$-relvement de $\rho$ (relativement ˆ $E_1/F$); auquel cas, avec les donnŽes auxiliaires fixŽes comme en \ref{facteurs de transfert bis}, on a l'ŽgalitŽ  $$c(\tau,\pi,A)= c(\rho,\pi,A)\qhq{pour tout} A\in \mathrm{Isom}(\kappa\pi,\pi)\ptf$$
\end{lemme}

\begin{demo} La preuve de \textit{loc.\,cit.} utilise l'ŽgalitŽ entre fonctions-caractres \ref{def kappa-rel}\,(2). On donne ici une preuve analogue n'utilisant que l'ŽgalitŽ \ref{def kappa-rel}\,(1). 
Pour toute fonction $\phi_1\in C^\infty_{\mathrm{c}}(H_1)$, on a l'ŽgalitŽ:
$$\mathrm{tr}(\rho(\phi_1))= \mathrm{tr}(\tau^\natural((\phi_1)^{P_1}))= \mathrm{tr}(\tau(\phi))\qhq{avec} \phi=\iota^{-1}((\phi_1)^{P_1})\ptf$$  
Pour toute paire de fonctions concordantes $(f,\phi_1)\in \CG\times C^\infty_{\mathrm{c}}(H_1)$ 
(relativement ˆ $E_1/F$), puisque (d'aprs \ref{existence transfert bis}\,(i)) la fonction $\iota^{-1}((\phi_1)^{P_1})$ concorde avec $f$ (relativement ˆ $E/F$), si $\pi$ est un $\kappa$-relvement de $\tau$, 
on a bien l'ŽgalitŽ 
$$\mathrm{tr}(\pi(f)\circ A)= c(\tau, \pi,A)\mathrm{tr}(\rho(\phi_1))\qhq{pour tout} A\in \mathrm{Isom}_G(\kappa\pi,\pi)\ptf$$ 
RŽciproquement (d'aprs \ref{existence transfert bis}\,(ii)), pour toute paire de fonctions concordantes $(f,\phi)$ (relativement ˆ $E/F$), il existe une fonction $\phi_1\in C^\infty_{\mathrm{c}}(H_1)$ qui concorde avec $f$ (relativement ˆ $E_1/F$) et telle que $\phi$ et $\phi^\natural= \iota^{-1}((\phi_1)^{P_1})$ ait les mmes intŽgrales orbitales semi-simples $G$-rŽgulires. Par consŽquent si $\pi$ est un $\kappa$-relvement de $\rho$, on a l'ŽgalitŽ 
$$\mathrm{tr}(\pi(f)\circ A) = c(\rho,\pi,A)\mathrm{tr}(\tau(\phi^\natural))\qhq{pour tout} A\in \mathrm{Isom}_G(\kappa\pi,\pi)\ptf$$ 
Puisque $\HGreg$ est ouvert dense dans $H$, d'aprs la formule d'intŽgration de Weyl, on a $\mathrm{tr}(\tau(\phi^\natural))= \mathrm{tr}(\tau(\phi))$. Cela prouve le lemme. 
\hfill$\square$
\end{demo}

\begin{remark}\label{variante inverse}
\textup{On peut intervertir les opŽrations d'induction parabolique et de $\kappa$-relvement dans 
la construction ci-dessus. Le plongement $\varphi$ de $H$ dans $G$ est donnŽ via le choix d'une $F$-base de $E^m$. Quitte ˆ changer cette $F$-base, \ie ˆ remplacer $H\;(=\varphi(H))$ par $gHg^{-1}$ pour un $g\in G$, on peut supposer que $\varphi = \varphi_1\times \cdots \times \varphi_r$ o $\varphi_i$ est un plongement de $\GLm(E_i)$ dans $\mathrm{GL}_{ms}(F)$. Alors $H$ s'identifie ˆ un sous-groupe de $L= (\mathrm{GL}_{ms}(F))^r$. Soit $\pi_L$ une reprŽsentation irrŽductible $\kappa$-stable de $L$, \ie $\pi_L=\pi_1\otimes \cdots \otimes \pi_r$ o $\pi_i$ est une reprŽsentation irrŽductible $\kappa$-stable de $\mathrm{GL}_{ms}(F)$. Supposons que pour $i=1,\ldots ,r$, $\pi_i$ soit un $\kappa$-relvement de 
$\tau_i$ (relativement ˆ $E_i/F$), \ie de $\tau_i^\natural$ (relativement ˆ $E_1/F$). Supposons aussi que l'induite parabolique 
$\pi'= \pi_1 \times \cdots \times \pi_r \;(= i_{P_L}^G(\pi_L))$ soit irrŽductible. Alors $\pi'$ est un $\kappa$-relvement de $\tau$. Cela rŽsulte des formules de descente parabolique pour les traces tordues $\mathrm{tr}(\pi'(f)\circ A)$ et les $\kappa$-intŽgrales orbitales $I^G_\kappa(f,\gamma)$, Žtablies en \cite[7.3.9]{HL4} et \cite[7.2.2]{HL4}. Observons que si la reprŽsentation $\rho= \tau_1^\natural\times \cdots \times \tau_r^\natural$ de $\mathrm{GL}_{mr}(E_1)$ est irrŽductible, cela rŽsulte aussi  de \ref{compatibilitŽ IP-IA} et \ref{[HL] prop 3.8}.\hfill $\blacksquare$
}
\end{remark}

\section{ReprŽsentations sphŽriques}\label{reprŽsentations sphŽriques}

Dans cette section, $F$ une extension finie de $\mathbb{Q}_p$ et $E$ est 
une $F$-algbre cyclique \textit{non ramifiŽe} de degrŽ fini $d$.

\subsection{Notations.} On pose $\Gamma=\Gamma(E/F)$ et on 
Žcrit $E=E_1\times \cdots \times E_r$ o 
$E_i/F$ est une extension non ramifiŽe de degrŽ $s= \frac{d}{r}$. On reprend les hypothses et les notations de \ref{notations ext cyclique}. En particulier, on a 
$\sigma E_i = E_{i+1}$ ($i=1,\ldots , r-1$) et $\sigma E_r = E_1$ pour un gŽnŽrateur $\sigma$ de $\Gamma$. 
On fixe un caractre $\kappa$ de $F^\times$ de noyau $\mathrm{N}_{E/F}(E^\times)= \mathrm{N}_{E_1/F}(E_1^\times)$ et on pose $\zeta=\kappa(\varpi)$ pour une (\ie pour toute) uniformisante $\varpi$ de $F^\times$. 

On pose $H=\GLm(E)= \prod_{i=1}^r \GLm(E_i)$ et $G=\GLn(F)$, $n=md$. On pose aussi 
$K_H= \GLm(\mathfrak{o}_E)=\prod_{i=1}^m\GLm(\mathfrak{o}_{E_i})$ et $K= \GLn(\mathfrak{o}_F)$. 
On suppose de plus que le plongement $\varphi$ de $H$ dans $G$ est dŽfini via le choix d'une $\mathfrak{o}_F$-base de $\mathfrak{o}_E^m=\prod_{i=1}^r\mathfrak{o}_{E_i}$. Cela entra"ne que $\varphi$ se restreint en un plongement de $K_H$ dans $K$. D'autre part l'isomorphisme $\iota$ de $E$ sur $(E_1)^r$ induit 
un isomorphisme, encore notŽ $\iota$, de $H$ sur $H^\natural = \GLm(E_1)^r$ qui se restreint en un isomorphisme de $K_H$ sur $K_{H^\natural}= \GLm(\mathfrak{o}_{E_1})^r$. 

On pose $H_1= \mathrm{GL}_{mr}(E_1)$ et $K_1= \mathrm{GL}_{mr}(\mathfrak{o}_{E_1})$. La $\mathfrak{o}_F$-base de $\mathfrak{o}_E^m$ dŽfinissant le plongement $\varphi$ de $H$ dans $G$, identifiŽe via $\iota$ ˆ une $\mathfrak{o}_F$-base de $\mathfrak{o}_{E_1}^{mr}$, dŽfinit un plongement $\varphi_1$ de $H_1$ dans $G$ qui envoie $K_1$ dans $K$.  

On suppose que l'ŽlŽment $e=(e_1,\ldots , e_r)\in E^\times$ tel que $\sigma e= (-1)^{m(d-1)}e$, utilisŽ pour dŽfinir le facteur de transfert $\Delta$, appartient 
ˆ $\mathfrak{o}_E^\times = \prod_{i=1}^r\mathfrak{o}_{E_i}^\times$. On suppose que les mesures de Haar $\d g$ sur $G(F)$ et $\d h$ sur $H$, utilisŽes pour dŽfinir les intŽgrales orbitales et les distributions-traces, sont normalisŽes par $\mathrm{vol}(K,\d g)=1 = \mathrm{vol}(K_H, \d h)$. Enfin on suppose que le caractre additif $\psi$ de $F$, utilisŽ pour dŽfinir les opŽrateurs d'entrelacement normalisŽs $A_{\pi,\psi}$, est de \textit{niveau $0$}, \ie trivial sur $\mathfrak{o}_F$ mais pas sur $\mathfrak{p}_F^{-1}$.

\subsection{ParamŽtrisation de Satake.}\label{param satake} 
Une reprŽsentation irrŽductible $\pi$ de $G$ est dite \textit{sphŽrique} si elle a un vecteur non nul fixŽ par $K$, \ie si $V_\pi^K \neq \{0\}$. 
On appelle \textit{sŽrie principale non ramifiŽe} de $G$ une reprŽsentation isomorphe ˆ $i_{A_0}^G(\chi)$ pour un caractre non ramifiŽ 
$\chi$ de $A_0=(F^\times)^n$. On sait \cite{S} que toute reprŽsentation irrŽductible sphŽrique de $G$ est isomorphe 
ˆ un sous-quotient d'une sŽrie principale non ramifiŽe, et que toute sŽrie principale non ramifiŽe de $G$ ˆ un unique sous-quotient irrŽductible sphŽrique.  
Pour $y=(y_1,\ldots ,y_n)\in (\mathbb{C}^\times)^n$, on note: 
\begin{itemize}
\item $\chi_y$ le caractre non ramifiŽ de $A_0$ dŽfini par $\chi_y(\varpi^{a_1},\ldots , \varpi^{a_n})= \prod_{i=1}^n y_i^{a_i}$ pour toute uniformisante 
$\varpi$ de $F$ et tout $n$-uplet $(a_1,\ldots , a_n)\in \mathbb{Z}^n$;
\item $\pi_y^*$ la sŽrie principale non ramifiŽe $i_{A_0}^G(\chi_y)$; 
\item $\pi_y$ l'unique sous-quotient irrŽductible sphŽrique de $\pi_y^*$.
\end{itemize}
L'application $y \mapsto \pi_y$ induit une bijection entre $(\mathbb{C}^\times)^n/\mathfrak{S}_n$ et l'ensemble $\mathrm{Irr}^K(G)$ des classes d'isomorphisme de reprŽsentations 
irrŽductibles sphŽriques de $G$. 

\begin{remark}\label{caractres nr}
\textup{Parmi les reprŽsentations irrŽductibles sphŽriques de $G$, on a les \textit{caractres non ramifiŽs}, \ie les reprŽsentations de la forme $\xi\circ \det$ pour un caractre non ramifiŽ $\xi$ de $F^\times$. Une telle reprŽsentation est l'unique quotient irrŽductible (\ie le quotient de Langlands) $u(\xi,n)=L(\nu^{\frac{n-1}{2}}\xi, \cdots , \nu^{-\frac{n-1}{2}}\xi)$ de l'induite parabolique 
$R(\xi,n)= \nu^{\frac{n-1}{2}}\xi\times \cdots \times \nu^{-\frac{n-1}{2}}\xi$. On a donc $u(\xi,n)= \pi_y$ avec 
$y=(q^{-\frac{n-1}{2}}\xi(\varpi), \ldots ,q^{\frac{n-1}{2}}\xi(\varpi))\in (\mathbb{C}^\times)^n$, o $q=q_F$ est le cardinal du corps rŽsiduel de $F$. 
En particulier $u(\xi=1,n)$ est la reprŽsentation triviale.\hfill $\blacksquare$}
\end{remark}

Pour $y=(y_1,\ldots , y_n)\in (\mathbb{C}^\times)^n$, la reprŽsentation $\kappa\pi_y= (\kappa\circ\det)\otimes \pi_y$ est isomorphe ˆ $\pi_{\zeta \cdot y}$ avec 
$\zeta \cdot y = (\zeta y_1,\ldots ,\zeta y_n)$. On en dŽduit que $\pi_y$ est $\kappa$-stable si et seulement l'image de $\zeta\cdot y$ dans $(\mathbb{C}^\times)^n/\mathfrak{S}_n$ co\"{\i}ncide avec celle de $y$. On note $\mathrm{Irr}^K_\kappa(G)$ le sous-ensemble de $\mathrm{Irr}^K(G)$ formŽ des reprŽsentations $\kappa$-stables.

\vskip1mm
De la mme manire, une reprŽsentation irrŽductible $\tau= \tau_1\otimes \cdots \otimes \tau_r$ de $H$ est dite \textit{sphŽrique} si elle 
a un vecteur non nul fixŽ par $K_H$, \ie si pour $i=1,\ldots , r$, on a $V_{\tau_i}^{K_{H_i}}\neq \{0\}$ avec $K_{H_i}= \GLm(\mathfrak{o}_{E_i})$. 

Supposons pour commencer que $E$ soit une extension de $F$ (\ie $E=E_1$). 
Pour $y\in (\mathbb{C}^\times)^m$, on dŽfinit comme plus haut le caractre non ramifiŽ 
$\chi_{E,y}$ de $A_{0,H}= (E^\times)^m$, on note $\Pi^*_y$ la sŽrie principale non ramifiŽe 
$i_{A_{0,H}}^H(\chi_{E,y})$ de $H$ et $\Pi_y$ son sous-quotient irrŽductible sphŽrique. Observons que la reprŽsentation 
$\Pi_y$ est $\Gamma$-stable: on a ${^\gamma\Pi_y}\simeq \Pi_y$ pour tout $\gamma\in \Gamma\;(=\mathrm{Gal}(E/F))$. L'application 
$y \mapsto \Pi_y$ induit une bijection entre $(\mathbb{C}^\times)^m/\mathfrak{S}_m$ et l'ensemble $\mathrm{Irr}^{K_H}(H)$ 
des classes d'isomorphisme de reprŽsentations 
irrŽductibles sphŽriques de $H$. Pour $y= (y_1,\ldots , y_m)\in (\mathbb{C}^\times)^m$, on choisit un ŽlŽment $t=(t_1,\ldots , t_m)\in (\mathbb{C}^\times)^m$ tel que 
$t^d = y$, \ie tel que $t_i^d= y_i$ ($i=1,\ldots , m$), et on note $\delta(y)\in (\mathbb{C}^\times)^n$ l'ŽlŽment dŽfini par 
$$\delta(y)=(t,\zeta t, \ldots , \zeta^{d-1}t)\ptf$$ L'image de $\delta(y)$ dans $(\mathbb{C}^\times)^n/\mathfrak{S}_n$ ne dŽpend pas du choix de $t$ puisque 
tout autre choix est de la forme $(\zeta^{b_1}t_1,\ldots , \zeta^{b_m}t_m)$ pour des entiers $b_i\in \{0,\ldots , d-1\}$. Elle ne dŽpend pas non plus du choix de $\kappa$ puisque changer de $\kappa$ revient ˆ remplacer $\zeta$ par une autre racine primitive 
$d$-ime de l'unitŽ. D'ailleurs elle ne dŽpend que de l'image de $y$ dans $(\mathbb{C}^\times)^m/\mathfrak{S}_m$. L'application 
$$\delta:( \mathbb{C}^\times)^m/\mathfrak{S}_m \rightarrow (\mathbb{C}^\times)^n/\mathfrak{S}_n$$ est donc bien dŽfinie; et elle est injective. 
Elle induit une application bijective $$\mathrm{Irr}^{K_H}(H)\rightarrow \mathrm{Irr}^K_\kappa(G)\,,\; \Pi_y \mapsto  \pi_{\delta(y)}\quad (y\in (\mathbb{C}^\times)^m/\mathfrak{S}_m)\ptf$$ 

Traitons maintenant le cas gŽnŽral (\ie $E=E_1\times \cdots \times E_r$ avec $rs=d$). Pour $y=(y(1),\ldots , y(r))\in (\mathbb{C}^\times)^{mr}$ avec $y(i)\in (\mathbb{C}^\times)^m$, on dŽfinit comme plus haut la reprŽsentation irrŽductible sphŽrique $\Pi_{y(i)}$ de $\GLm(E_i)$ ($i=1,\ldots , r$), et on note $\Pi_y$, resp. $\Pi_y^\natural$, la reprŽsentation irrŽductible sphŽrique de $H$, resp. $H^\natural$, dŽfinie par  
$$\Pi_y= \Pi_{y(1)} \otimes\cdots \otimes \Pi_{y(r)}\qhq{et} \Pi_y^\natural= \Pi_{y(1)} \otimes (\Pi_{y(2)})^{\sigma} \otimes \cdots \otimes (\Pi_{y(r)})^{\sigma^{r-1}}\ptf$$ 
La classe d'isomorphisme de $\Pi_y$, resp. $\Pi_y^\natural$, ne dŽpend que de l'image de 
$y$ dans $((\mathbb{C}^\times)^m/\mathfrak{S}_m)^r$ et l'application $y\mapsto \Pi_y$ induit une bijection de $((\mathbb{C}^\times)^m/\mathfrak{S}_m)^r$ sur 
$\mathrm{Irr}^{K_H}(H)$. 
Pour $y\in (\mathbb{C}^\times)^{mr}$, on choisit un ŽlŽment $t\in (\mathbb{C}^\times)^{mr}$ tel que $t^s= y$ et on pose $\delta(y)= (t, \zeta t, \ldots , \zeta^{s-1}t)\in 
(\mathbb{C}^\times)^n$. 
L'image de $\delta(y)$ dans $(\mathbb{C}^\times)^n/\mathfrak{S}_n$ est bien dŽfinie et elle ne dŽpend que de l'image de $y$ 
dans $((\mathbb{C}^\times)^m/\mathfrak{S}_m)^r$. Observons que si $s>1$ (\ie si $E\neq E_1$), l'application 
$$\delta : ((\mathbb{C}^\times)^m/\mathfrak{S}_m)^r\rightarrow (\mathbb{C}^\times)^n/\mathfrak{S}_n$$ n'est pas injective; ses fibres sont les 
$\mathfrak{S}_r$-orbites dans $((\mathbb{C}^\times)^m/\mathfrak{S}_m)^r$. Elle induit une application surjective 
$$\mathrm{Irr}^{K_H}(H)\rightarrow \mathrm{Irr}^K_\kappa(G)\,,\; \Pi_y \mapsto  \pi_{\delta(y)}\quad (y\in ((\mathbb{C}^\times)^m/\mathfrak{S}_m)^r)$$ dont les fibres 
sont les $\Gamma$-orbites dans $\mathrm{Irr}^{K_H}(H)$. 

\begin{remark}\label{kappa-rel faible et Sigma}
\textup{Un ŽlŽment $y\in (\mathbb{C}^\times)^{mr}$ dŽfinit aussi une reprŽsentation irrŽductible sphŽrique $\Sigma_y$ de $H_1= \mathrm{GL}_{mr}(E_1)$, ˆ savoir l'unique sous-quotient irrŽductible sphŽrique de la sŽrie principale non ramifiŽe 
$\Sigma_y^*= i_{A_{0,H_1}}^{H_1}(\chi_{E_1,y})$ de $H_1$, o $A_{0,H_1}= (E_1^\times)^{mr}$. D'aprs le cas $E=E_1$, l'application $y\mapsto \delta(y)$ 
se quo\-tiente en une application injective 
$(\mathbb{C}^\times)^{mr}/\mathfrak{S}_{mr} \rightarrow (\mathbb{C}^\times)^n/\mathfrak{S}_n,\, y \mapsto \delta(y)\vg$ qui induit une 
application bijective $$\mathrm{Irr}^{K_{H_1}}(H_1) \rightarrow \mathrm{Irr}^K_\kappa(G)\,,\; \Sigma_y \mapsto \pi_{\delta(y)} 
\quad (y\in (\mathbb{C}^\times)^{mr}/\mathfrak{S}_{mr})\ptf$$ 
Observons que $\Sigma_y$ est aussi l'unique sous-quotient irrŽductible sphŽrique de l'induite parabolique $i_{H^\natural}^{H_1}(\Pi^\natural_y)$. 
\hfill $\blacksquare$
}
\end{remark}

\subsection{$\kappa$-relvement faible.}\label{kappa-rel faible} Pour $y\in (\mathbb{C}^\times)^{mr}$, la reprŽsentation irrŽductible sphŽrique $\kappa$-stable $\pi_{\delta(y)}$ de $G$ 
est un candidat naturel pour tre un $\kappa$-relvement de $\Pi_y$. Nous verrons ˆ la fin de la dŽmonstration qu'il en est bien ainsi, pourvu que la reprŽsentation $\Pi_y$ soit unitaire. Dans un premier temps nous allons nous contenter d'une notion plus faible: celle de \textit{$\kappa$-relvement faible}. 

On note $\mathcal{H}=\mathcal{H}(G,K)\;(\subset C^\infty_{\mathrm{c}}(G))$ l'algbre de Hecke sphŽrique formŽe des fonctions $K$-biinvariantes ˆ support compact sur $G$. On dŽfinit de la mme manire l'algbre de Hecke sphŽrique $\mathcal{H}_H=\mathcal{H}(H,K_H)$.

Si $E/F$ est une extension de corps (\ie si $E=E_1$), en \ref{lecaslocalnr} est dŽfini un homomorphisme d'algbres $b: \mathcal{H} \rightarrow \mathcal{H}_H$, qui co\"{\i}ncide avec celui dŽfini en \cite[II.3]{W1}. 
En gŽnŽral, on commence par considrer l'homomorphisme d'algbres 
$b_1: \mathcal{H} \rightarrow \mathcal{H}_{H_1}=\mathcal{H}(H_1,K_{H_1})$ dŽfini en \ref{lecaslocalnr} (relativement ˆ l'extension $E_1/F$). 
Rappelons qu'on a notŽ  
$P_1= H^\natural\ltimes U_1$ le sous-groupe parabolique standard de $H_1$ de composante de Levi $H^\natural= \GLm(E_1)^r$. L'application \guill{terme constant} 
$$\mathcal{H}_{H_1} \rightarrow \mathcal{H}_{H^\natural}= \mathcal{H}(H^\natural,K_{H^\natural})\,,\; \phi_1 \mapsto (\phi_1)^{P_1}$$ 
est elle aussi un homomorphisme d'algbres; o (rappel) 
$$(\phi_1)^{P_1}(h^\natural)= \int_{U_1}\phi_1(h^\natural u_1)\d u_1\quad (h^\natural\in H^\natural) \vgq \mathrm{vol}(U_1\cap K_{H_1},\d u_1)=1\ptf$$  
Soit $b: \mathcal{H} \rightarrow \mathcal{H}_H$ l'homomorphisme d'algbres dŽfini par 
$$b f = \iota^{-1}((b_1 f)^{P_1})\qhq{pour toute fonction} f \in  \mathcal{H}\ptf$$

\begin{definition}\label{def kapparelfaible}
\textup{Soit $\tau$ une reprŽsentation irrŽductible sphŽrique de $H$ et soit $\pi$ une reprŽsentation irrŽductible sphŽrique $\kappa$-stable de $G$. On dit que 
$\pi$ est un \textit{$\kappa$-relvement faible} de $\tau$ si pour toute fonction $f\in \mathcal{H}$, on a l'ŽgalitŽ
$$\mathrm{tr}(\pi(f)) = \mathrm{tr}(\tau(b f))\ptf$$}
\end{definition}

La notion de $\kappa$-relvement faible ne dŽpend que des classes d'isomorphisme de $\tau$ et $\pi$. 
Un $\kappa$-relvement faible, s'il existe, est unique ˆ isomorphisme prs. 

\begin{lemme}\label{description du K-rel faible dans le cas gŽnŽral}
Pour tout $y\in (\mathbb{C}^\times)^{mr}$, $\pi_{\delta(y)}$ est un $\kappa$-relvement faible de $\Pi_y$. 
\end{lemme}

\begin{demo}
Si $E=E_1$, c'est l'ŽgalitŽ  \ref{lecaslocalnr}\,(3). En particulier avec les notations de \ref{kappa-rel faible et Sigma}, $\pi_{\delta(y)}$ est un $\kappa$-relvement faible de $\Sigma_y$, au sens o
$$\mathrm{tr}(\pi_{\delta(y)}(f)) = \mathrm{tr}(\Sigma_y(b_1 f))\qhq{pour toute fonction} f\in \mathcal{H}\ptf$$ 
Puisque $\Sigma_y$ est l'unique sous-quotient irrŽductible sphŽrique de $i_{H^\natural}^{H_1}(\Pi^\natural_y)$, pour toute fonction 
$\phi_1\in C^\infty_{\mathrm{c}}(H_1)$, on a 
$$\mathrm{tr}(\Sigma_y(\phi_1))= \mathrm{tr}\left(i_{H^\natural}^{H_1}(\Pi^\natural_y)(\phi_1)\right)= \mathrm{tr}\left(\Pi^\natural_y((\phi_1)^{P_1})\right)= 
\mathrm{tr}\left( \Pi_y(\iota^{-1}((\phi_1)^{P_1}))\right)\ptf$$
Cela prouve le lemme. \hfill $\square$
\end{demo}

\begin{exemple}\label{kappa-rel triv}
\textup{Supposons que $r=1$ (\ie $E=E_1$) et soit $\mu$ un caractre non ramifiŽ de $E^\times$. Notons $\mu_m$ le caractre non ramifiŽ $\mu\circ \det_E$ de $\GLm(E)$, \ie $\mu_m = \mu \bs{1}_{E,m}$ o $\bs{1}_{E,m}$ est la reprŽsentation triviale de $\GLm(E)$. D'aprs \ref{caractres nr}, $\bs{1}_{E,m}$ correspond ˆ l'ŽlŽment $y_{E,m}= (q_E^{\smash{-\frac{m-1}{2}}}\!, \ldots , q_E^{\smash{\frac{m-1}{2}}})$ de $(\mathbb{C}^\times)^m$, au sens o $\bs{1}_{E,m}= \Pi_{y_{E,m}}$. PrŽcisŽment, $\bs{1}_{E,m}$ est le quotient de Langlands $u(\bs{1}_{E},m)$ de l'induite parabolique $\nu_E^{\smash{\frac{m-1}{2}}}\times \cdots \times \nu_E^{\smash{-\frac{m-1}{2}}}$; et $\mu_m= u(\mu,m)$. Choisissons un caractre non ramifiŽ $\xi$ de $F^\times$ tel que $\xi(\varpi)^m= \mu(\varpi)$, \ie tel que $\mu = \xi\circ \mathrm{N}_{E/F}$, et notons $\xi_m$ le caractre non ramifiŽ $\xi \circ \det_F$ de $\GLm(F)$. Puisque $q_F^d = q_E$, d'aprs \ref{description du K-rel faible dans le cas gŽnŽral}, la reprŽsentation (irrŽductible unitaire) 
$$\bs{1}_{F,m,\kappa}\bydef \bs{1}_{F,m} \times \kappa \bs{1}_{F,m} \times \cdots \times \kappa^{d-1}\bs{1}_{F,m}$$ est un $\kappa$-relvement faible de $\bs{1}_{E,m}$; et la reprŽsentation $$\xi_{m,\kappa}\bydef \xi \bs{1}_{F,m,\kappa}= \xi_m \times \kappa \xi_m \times \cdots \times \kappa^{d-1}\xi_m$$ est un 
$\kappa$-relvement faible de $\mu_m$. Observons que d'aprs \cite[6.1]{BH}, $\xi_m$ est un relvement ˆ la Shintani (\ie pour le changement de base) de $\mu_m$.\hfill $\blacksquare$
}
\end{exemple}

\begin{remark}\label{divers, LF et op norm}
\textup{
\begin{enumerate}
\item[(i)] Le lemme fondamental pour l'induction automorphe a ŽtŽ prouvŽ par Waldspurger \cite{W1} dans le cas o $p$ (la caractŽristique rŽsiduelle de $F$) est strictement supŽrieur ˆ $n$. La restriction sur $p$ a ensuite ŽtŽ otŽe dans \cite{HL2}. Ainsi d'aprs ce lemme fondamental, pour toute fonction $f\in \mathcal{H}(G,K)$, la fonction $\phi_1 = b_1 (f)$ concorde avec $f$ (relativement ˆ $E_1/F$); et d'aprs \ref{existence transfert bis}, la fonction $\iota^{-1}((\phi_1)^{P_1})= b f$ concorde avec $f$ (relativement ˆ $E/F$). 
\item[(ii)]Si $\pi$ est une reprŽsentation irrŽductible sphŽrique $\kappa$-stable de $G$, puisque le caractre $\kappa$ de $G$ est trivial sur $K$, tout opŽrateur $A\in \mathrm{Isom}_G(\kappa \pi, \pi)$ stabilise la droite $V_\pi^{K}$. Si $\pi$ est gŽnŽrique, puisque $\psi$ est de niveau $0$, 
on sait d'aprs \cite[3.1]{HL2} 
que l'opŽrateur d'entrelacement normalisŽ $ \Agen_{\pi,\psi}$ est l'identitŽ sur la droite $V_\pi^{K}$. D'aprs la propriŽtŽ d'induction parabolique et de multiplicitŽ $1$, 
ce rŽsultat reste vrai sans l'hypothse de gŽnŽricitŽ. En particulier on a l'ŽgalitŽ 
$$\mathrm{tr}(\pi(f))= \mathrm{tr}(\pi(f) \circ A_{\pi,\psi})\qhq{pour toute fonction} f\in \mathcal{H}(G,K)\ptf$$  
\item[(iii)] Soit $\pi$ une reprŽsentation irrŽductible sphŽrique $\kappa$-stable de $G$ qui soit un $\kappa$-relvement d'une reprŽsentation 
irrŽductible sphŽrique $\tau$ de $H$. Alors la constante $c(\tau,\pi, A_{\pi,\psi})$ vaut $1$. Il suffit pour cela d'appliquer l'ŽgalitŽ de $\kappa$-relvement 
aux fonctions caractŽristiques $f_0$ de $K$ et $\phi_0= b f_0$ de $H$. En effet on a 
$$\mathrm{tr}(\pi(f_0)\circ A_{\pi,\psi})= \mathrm{tr}(\pi(f_0))=1=
\mathrm{tr}(\tau(\phi_0))\ptf$$
En particulier $\pi$ est \textit{a fortiori} un $\kappa$-relvement faible de $\tau$. \hfill $\blacksquare$
\end{enumerate}
}
\end{remark}

\subsection{Le cas des reprŽsentations sphŽriques unitaires.}\label{kappa-rel faible unit} 
Soit $\tau$ une reprŽsentation irrŽductible sphŽrique unitaire de $\GLm(E)$.  On suppose tout d'abord que $r=1$ (\ie $E=E_1$). 
D'aprs \ref{classif (rappels)}, $\tau$ est isomorphe ˆ une induite parabolique (irrŽductible) de la forme 
$$\alpha_1\times \cdots \times \alpha_r  \times (\nu_E^{a_1}\beta_1\times \nu_E^{-a_1}\beta_1)\times \cdots \times (\nu_E^{a_s}\beta_s\times \nu_E^{-a_s}\beta_s)$$  
pour des caractres non ramifiŽs unitaires $\alpha_i= (\mu_i)_{m_i}\;(=\mu_i \circ \det_E)$ de $\mathrm{GL}_{m_i}(E)$, resp. $\beta_j= (\mu'_j)_{m'_j}$ de $\mathrm{GL}_{m'_j}(E)$, et des nombres rŽels $a_j\in \;]0,\frac{1}{2}[$; avec l'ŽgalitŽ 
$\sum_{i=1}^r m_i + 2 \sum_{j=1}^s m'_j = m$. Rappelons que l'on a $$\alpha_i = u(\mu_i,m_i)\qhq{et} \beta_j= u(\mu'_j,m'_j)\ptf$$ 
Pour $i= 1,\ldots , r$, choisissons un caractre non ramifiŽ (forcŽment unitaire) $\xi_i$ de $F^\times$ tel que $\mu_i= \xi_i \circ \mathrm{N}_{E/F}$; de mme pour $j=1,\ldots ,s$, choisissons un caractre non ramifiŽ $\xi'_j$ de $F^\times$ tel que $\mu'_j = \xi_j \circ \mathrm{N}_{E/F}$. D'aprs l'exemple \ref{kappa-rel triv}, pour $i=1,\ldots ,r$, 
la reprŽsentation $$\pi_i = (\xi_i)_{m_i,\kappa}\;\left(=(\xi_i)_{m_i}\times \kappa (\xi_i)_{m_i}\times \cdots \times \kappa^{d-1} (\xi_i)_{m_i}\right)$$ est un $\kappa$-relvement faible de $\alpha_i$; et pour $j=1,\ldots , s$, la reprŽsentation $$\pi'_j  =(\xi'_j)_{m'_j,\kappa}$$ est un $\kappa$-relvement faible de $\beta_j$. On en dŽduit que la reprŽsentation (irrŽductible sphŽrique unitaire) $$\pi= \pi_1\times \cdots \times \pi_r \times (\nu_F^{a_1}\pi'_1\times \nu_F^{-a_1}\pi'_1)\times \cdots \times 
(\nu_F^{a_s}\pi'_s\times \nu_F^{-a_s}\pi'_s)$$ est un $\kappa$-relvement faible de $\tau$. 

Passons au cas gŽnŽral (\ie $r\geq 1$). \'Ecrivons $\tau= \tau_1\otimes \cdots \otimes \tau_r$ o $\tau_r$ est une reprŽsentation irrŽductible sphŽrique unitaire de $\GLm(E_i)$. On commence par former l'induite parabolique $\rho= \tau_1 \times (\tau_2)^\sigma\times \cdots \times (\tau_r)^{\sigma^{r-1}}$, qui est une reprŽsentation irrŽductible sphŽrique unitaire de $\mathrm{GL}_{mr}(E_1)$, puis on considre un $\kappa$-relvement faible $\pi$ de $\rho$ (relativement ˆ $E_1/F$). Alors $\pi$ est un $\kappa$-relvement faible de $\tau$ (relativement ˆ $E/F$). Pour $i=1,\ldots , r$, la reprŽsentation $(\tau_i)^{\sigma^{i-1}}$ de $\GLm(E_1)$ est irrŽductible sphŽrique et unitaire; elle s'Žcrit donc comme un produit de caractres du type $\alpha=\mu_k\; (= \mu\circ \det_{E_1}): \mathrm{GL}_k(E_1) \rightarrow \mathbb{U}$ pour un caractre non ramifiŽ unitaire $\mu$ de $E_1^\times$ et un entier $k\geq 1$, et de reprŽsentations du type $\nu_{E_1}^a (\mu')_{l} \times \nu_{E_1}^{-a} (\mu')_l$, $\mu'_l = \mu'\circ \det_{E_1}: \mathrm{GL}_l(E_1)\rightarrow \mathbb{U}$, pour un caractre non ramifiŽ unitaire $\mu'$ de $E_1^\times$, un entier $l\geq 1$ et un nombre rŽel $a\in \;]0,\frac{1}{2}[$. La reprŽsentation $\rho$ est elle aussi un produit de reprŽsentations de ces deux types et la description de $\pi$ se ramne ˆ celle du cas $r=1$. 

\subsection{Un lemme local-global.}\label{un lemme local-global} Soit $\E/\F$ une extension finie cyclique de corps de nombres, de degrŽ $d$. On fixe un caractre $\K$ de 
$\AF^\times$ de noyau $\F^\times \mathrm{N}_{\E/\F}(\AE^\times)$. Le lemme suivant permet de sŽparer les reprŽsentations automorphes induites de discrtes (voir \ref{induite de rep disc}) de $\GLm(\AE)$ qui, en presque toute place finie $v$ de $\F$, ont le mme $\K_v$-relvement faible ˆ $\GLn(\F_v)$, $n=md$. 

Pour une place finie $v$ de $\F$, on pose ${\bf K}_v= \GLn(\mathfrak{o}_v)$ et on note $\mathcal{H}_v$ l'algbre de Hecke sphŽrique $\mathcal{H}(\GLn(\F_v),{\bf K}_v)$. On pose ${\bf K}'_v= \GLm(\mathfrak{o}_{\E_v})= \prod_{w\vert v} \GLm(\mathfrak{o}_{\E_w})$ et on note $\mathcal{H}'_v$ l'algbre de Hecke sphŽrique 
$\mathcal{H}(\GLm(\E_v),{\bf K}'_v)$. Enfin on note $b_v: \mathcal{H}_v\rightarrow \mathcal{H}'_v$ l'homomorphisme d'algbres dŽfini en \ref{kappa-rel faible} (voir aussi \ref{lecaslocalnr}).   

\begin{lemme}\label{lemme local-global}
Soient $\Pi$ et $\Pi'$ deux reprŽsentations automorphes de $\GLm(\AE)$ de la forme\footnote{On note $\Delta^{\!\times l}= \Delta\times\cdots \times\Delta$ ($l$ fois) l'induite parabolique (au sens des reprŽsentations automorphes, cf. \cite{La2}) de la reprŽsentation $\Delta\otimes\cdots \otimes \Delta$ de 
$M(\AE)= \mathrm{GL}_{a}(\AE)^l$, $a=\frac{m}{l}$, suivant le sous-groupe parabolique standard de $\GLm(\AE)$ de composante de Levi $M(\AE)$.} $\Pi= \Delta^{\!\times l}$ et $\Pi'=(\Delta')^{\!\times l'}$ pour des entiers 
$l,\, l'\geq 1$ divisant $m$ et des reprŽsentations automorphes discrtes $\Delta$ de $\mathrm{GL}_{\frac{m}{l}}(\AE)$ et $\Delta'$ de $\mathrm{GL}_{\frac{m}{l'}}(\AE)$. On suppose qu'en presque toute place finie $v$ de $\F$, on a l'ŽgalitŽ 
$$\mathrm{tr}(\Pi_v(b_vf_v))= \mathrm{tr}(\Pi'_v(b_vf_v))\qhq{pour toute fonction} f\in \mathcal{H}_v\ptf$$ 
Alors $l=l'$ et il existe un ŽlŽment $\gamma\in \Gamma(\E/\F)$ tel que ${^\gamma\!\Delta}\simeq \Delta'$.
\end{lemme}

\begin{demo}
On fixe un sous-ensemble fini de places $S\subset \V$ contenant $\Vinf$, toutes les places $v\in \Vfin$ o la $\F_v$-algbre cyclique $\E_v$ est ramifiŽe et toutes celles o l'une ou l'autre des reprŽsentations $\Pi_v$ et $\Pi'_v$ n'est pas sphŽrique. Quitte ˆ remplacer $S$ par un ensemble plus gros, on peut aussi supposer que 
pour toute place $v\in \V^S = \V \smallsetminus S$, on a  
$$\mathrm{tr}(\Pi_v(b_vf_v))= \mathrm{tr}(\Pi'_v(b_vf_v))\qhq{pour toute fonction} f\in \mathcal{H}_v\ptf\leqno{(1)}$$ 
En d'autres termes, on suppose qu'en toute place $v\in \V^S$, $\Pi_v$ et $\Pi'_v$ ont mme $\K_v$-relvement faible (d'aprs \ref{description du K-rel faible dans le cas gŽnŽral}, un tel relvement existe). 

On veut prouver que $l'=l$ et $\Delta'= {^\gamma\!\Delta}$ pour un $\gamma \in \Gamma(\E/\F)$. On commence par se ramener au cas o $\Pi$ et $\Pi'$ sont des induites de cuspidales unitaires, \ie le cas o $\Delta$ et $\Delta'$ sont des reprŽsentations automorphes cuspidales unitaires. \'Ecrivons $\Delta=u(\rho,k)$ et $\Delta'= u(\rho',k')$ pour un entier 
$k\geq 1$, resp. $k'\geq 1$, divisant $\frac{m}{l}$, resp. $\frac{m}{l'}$, et une reprŽsentation automorphe cuspidale unitaire $\rho$ de $\mathrm{GL}_{\frac{m}{kl}}(\AE)$, resp. $\rho'$ de $\mathrm{GL}_{\frac{m}{k'l'}}(\AE)$. Pour chaque place $v\in \V^S$, la composante locale $\rho_v$ est sphŽrique (gŽnŽrique unitaire) et d'aprs \cite[9.2]{J} (voir aussi plus loin la dŽmonstration de \ref{A=I}), la composante locale $\Delta_v$ est Žgale ˆ $u(\rho_v,k)$, \ie l'unique quotient irrŽductible de l'induite parabolique 
$$R(\rho_v,k)= \nu_{\E_v}^{\frac{k-1}{2}}\rho_v \times \nu_{\E_v}^{\frac{k-1}{2}-1}\rho_v \times \cdots \times \nu_{\E_v}^{-\frac{k-1}{2}}\rho_v\ptf$$ 
Si la place $v$ est inerte dans $\E$, alors d'aprs \ref{classif (rappels)}, $\rho_v$ est isomorphe ˆ un produit de caractres non ramifiŽs de $\E_v^\times$ de la forme 
$$\mu_1\times \cdots \times \mu_r  \times (\nu_{\E_v}^{a_1}\eta_1\times \nu_{\E_v}^{-a_1}\eta_1)\times \cdots \times (\nu_{\E_v}^{a_s}\eta_s\times \nu_{\E_v}^{-a_s}\eta_s)$$ 
pour des caractres non ramifiŽs unitaires $\mu_1, \ldots , \mu_r, \eta_1,\ldots , \eta_s$ de $\E_v^\times$ et des nombres rŽels $a_j\in\; ]0,\frac{1}{2}[$. 
Par consŽquent la reprŽsentation $\rho_v^{\times l} = \rho_v \times \cdots \times \rho_v $ est elle aussi un produit de caractres non ramifiŽs de $\E_v^\times$. D'autre part $\Delta_v$ est une reprŽsentation irrŽductible sphŽrique unitaire de $\mathrm{GL}_{\frac{m}{l}}(\E_v)$ qui est isomorphe ˆ un produit de caractres non ramifiŽs de $\mathrm{GL}_k(\E_v)$; prŽcisŽment on a 
$$\Delta_v \simeq \alpha_1 \times \cdots\times  \alpha_r \times  (\nu_{\E_v}^{a_1}\beta_1\times \nu_{\E_v}^{-a_1}\beta_1)\times \cdots \times (\nu_{\E_v}^{a_s}\beta_s\times \nu_{\E_v}^{-a_s}\beta_s)$$ avec 
$$\alpha_i= u(\mu_i,k)\qhq{et} \beta_j = u(\eta_j,k)\ptf$$ Par consŽquent la reprŽsentation $\Pi_v= \Delta_v^{\!\times l}$ de $\GLm(\E_v)$ est elle aussi isomorphe ˆ un produit de caractres non ramifiŽs de $\mathrm{GL}_k(\E_v)$. Si maintenant $v$ n'est pas inerte dans $\E$, alors la reprŽsentation $\Pi_v$ s'Žcrit $\Pi_v= \otimes_{w\vert v}\Pi_w$ et $\Pi_w$ est isomorphe ˆ un produit de caractres non ramifiŽs de $\mathrm{GL}_k(\E_w)$. De la mme manire, la reprŽsentation $\Pi'_v$ s'Žcrit $\Pi'_v = \otimes_{w\vert v}\Pi'_w$ et $\Pi'_w$ est isomorphe ˆ un produit de caractres non ramifiŽs de $\mathrm{GL}_k(\E_w)$. 
Par hypothse $\Pi_v$ et $\Pi'_v$ ont mme $\K_v$-relvement faible. D'aprs \ref{kappa-rel faible unit}, cela n'est possible que si $k=k'$ et si $\rho_v^{\times l}$ et $\rho'^{\times l'}_v$ ont mme $\K_v$-relvement faible. On peut donc supposer que $\Pi= \rho^{\times l}$ et $\Pi'= \rho'^{\times l'}$.

Posons $\Gamma = \Gamma(\E/\F)$. Il s'agit de prouver l'ŽgalitŽ de fonctions $L$: 
$$\prod_{\gamma \in \Gamma} L^S(s,\rho' \times {^\gamma\check{\rho}})^{l'}= \prod_{\gamma \in \Gamma} L^S(s,\rho \times {^\gamma\check{\rho}})^l\vg\leqno{(2)}$$ 
o $\check{\rho}$ est la contragrŽdiente de $\rho$. 
En effet supposons que l'ŽgalitŽ (2) soit vŽrifiŽe. D'aprs Jacquet-Shalika \cite{JS1}, son terme de droite a un p™le d'ordre au moins $1$ en $s=1$. 
Par consŽquent son terme de gauche a lui aussi un p™le en $s=1$, ce qui implique 
que $\rho'\simeq {^\gamma\rho}$ pour un $\gamma\in \Gamma$. Comme l'ordre du p™le en $s=1$ du terme de droite est le mme que celui du terme de gauche, 
on a forcŽment $l=l'$. On obtient aussi que $\Pi' \simeq ({^\gamma\rho})^{\times l} \simeq {^\gamma\Pi}$.

Prouvons l'ŽgalitŽ (2). Pour $v\in \V^S$, on a $\E_v = \prod_{w\vert v}\E_w$ o $w$ parcourt les places de $\E$ au-dessus de $v$. Notons $f_v$ le degrŽ de l'extension cyclique non ramifiŽe 
$\E_w/\F_v$ pour une (\ie pour toute) place $w\vert v$, et posons $e_v = \frac{d}{f_v}$. Choisissons 
un reprŽsentant $y_v= \prod_{w\vert v} y_w \in ((\mathbb{C}^\times)^{\frac{m}{l}})^{e_v}$ du paramtre de Satake de $\rho_v = \otimes_{w\vert v}\rho_w$; alors 
$\check{y}_v = \prod_{w\vert v}\check{y}_w$ avec $\check{y}_w= y_w^{-1}$ est un reprŽsentant du paramtre de Satake de $\check{\rho}_v= \otimes_{w\vert v} \check{\rho}_w$. Choisissons aussi un reprŽsentant $y'_v = \prod_{w\vert v}y'_w\in ((\mathbb{C}^\times)^{\frac{m}{l'}})^{e_v}$ du paramtre de Satake de 
$\rho'_v = \otimes_{w\vert v} \rho'_w$. On voit $y_v$ et $\check{y}_v$ comme des endomorphismes diagonaux de 
$(\mathbb{C}^\times)^{m_v}$ avec $m_v= \frac{m}{l}e_v$, et on voit $y'_v$ comme un endomorphisme diagonal de 
$(\mathbb{C}^\times)^{m'_v}$ avec $m'_v = \frac{m}{l'}e_v$. Le groupe $\Gamma$ opre sur 
$((\mathbb{C}^\times)^{\frac{m}{l}})^{e_v}$ par permutation des facteurs $(\mathbb{C}^\times)^{\frac{m}{l}}$: 
pour $\gamma \in \Gamma(\E/\F)$ et $z= \prod_{w\vert v}$, on pose $\gamma(z)= \prod_{w\vert v} \gamma(z)_w$ avec $\gamma(z)_w = z_{\gamma(w)}$. 
Fixons un gŽnŽrateur $\sigma$ de $\Gamma$; on a donc $\sigma^{e_v}(z)=z$ pour tout $z\in (\mathbb{C}^{\frac{m}{l}})^{e_v}$. 
Le terme de gauche de l'ŽgalitŽ (2) s'Žcrit alors 
$$\prod_{\gamma \in \Gamma} L^S(s,\rho' \times {^\gamma\check{\rho}})^{l'}=\prod_{v \in \V^S} \prod_{i=0}^{d-1} \det(1- y'_v\otimes \sigma^i(\check{y}_v)q_v^{-f_vs})^{-l'}$$ 
o $q_v= q_{\F_v}$ est le cardinal du corps rŽsiduel de $\F_v$ (ainsi pour toute place $w\vert v$, on a $q_v^{f_v}= q_{\E_{w}}$) et o l'on voit le produit tensoriel $y'_v\otimes \sigma^i(\check{y}_v)q_v^{-f_vs}$ comme un endomorphisme de $(\mathbb{C}^\times)^{m'_vm_v}$. 
Posons $\zeta_v = \K_v(\varpi_v)$; c'est une racine primitive $f_v$-ime de l'unité.
Choisissons un ŽlŽment $t_v\in ((\mathbb{C}^\times)^{\frac{m}{l}})^{e_v}$ tel que $(t_v)^{f_v}= y_v$ et posons 
$$\delta_v= (t_v, \zeta_vt_v,\ldots , \zeta_v^{f_v-1}t_v)\in ((\mathbb{C}^\times)^{\frac{m}{l}})^{d}= (\mathbb{C}^\times)^{\frac{n}{l}}\ptf$$ 
On pose aussi $$\check{\delta}_v = (\check{t}_v,\zeta_v^{-1}\check{t}_v,\ldots ,\zeta_v^{-(f_v-1)}\check{t}_v)\ptf$$
De la même manière, choisissons un ŽlŽment $t'_v\in  ((\mathbb{C}^\times)^{\frac{m}{l'}})^{e_v}$ tel que 
$(t'_v)^{f_v}= y'_v$ et posons 
$$\delta'_v= (t'_v, \zeta_vt'_v,\ldots , \zeta_v^{f_v-1}t'_v)\in ((\mathbb{C}^\times)^{\frac{m}{l'}})^{d}= (\mathbb{C}^\times)^{\frac{n}{l'}}\ptf$$  On obtient 
$$\prod_{\gamma \in \Gamma} L^S(s,\rho' \times {^\gamma\check{\rho}})^{l'}=\prod_{v\in \V^S} \det(1-\delta'_v\otimes \check{\delta}_vq^{-s})^{-d l'}\ptf$$ 
De la mme manire on obtient (terme de droite de l'ŽgalitŽ (2))
$$\prod_{\gamma \in \Gamma} L^S(s,\rho \times {^\gamma\check{\rho}})^l=\prod_{v\in \V^S} \det(1-\delta_v\otimes \check{\delta}_vq^{-s})^{-d l}\ptf$$ 
Par ailleurs pour $v\in \V^S$, d'aprs l'ŽgalitŽ (1), les représentations (irrŽductibles sphŽriques) $\Pi_v= \rho_v^{\times l}$ et $\Pi'_v = \rho'^{\times l'}_v$ ont le même $\frak{K}_v$-relèvement faible, 
ce qui Žquivaut ˆ dire que les ŽlŽments $(\delta_v, \ldots , \delta_v)$ ($l$ copies) et $(\delta'_v,\ldots , \delta'_v)$ ($l'$ copies) de $(\mathbb{C}^\times)^n$ 
sont $\mathfrak{S}_n$-conjuguŽs. D'o l'ŽgalitŽ (2).\hfill $\blacksquare$
\end{demo}

\section{La thŽorie dans le cas archimŽdien}\label{le cas a}
Dans cette section, $F$ est un corps local archimŽdien (\ie $F\simeq \mathbb{R}$ ou $\mathbb{C}$) et $E$ est une $F$-algbre cyclique de degrŽ fini $d$. On dŽcrit la version archimŽdienne des sections \ref{le cas n-a} et \ref{le cas d'une algbre cyclique}.  
Toutes les reprŽsentations sont supposŽes admissibles, ˆ valeurs dans le groupe des automorphismes d'un espace 
de Hilbert. 

\subsection{Rappels sur les reprŽsentations de $\GLn(F)$.}\label{rappels (g,K)-modules} On pose $\Gamma= \Gamma(E/F)$. 
Soit $G=\GLn(F)$ pour un entier $n\geq 1$. On reprend les notations de \ref{IP}. 

On s'intŽresse principalement aux reprŽsentations unitaires irrŽductibles (topologiquement) de $G$ et ˆ leurs $(\mathfrak{g},K)$-modules sous-jacents; o $\mathfrak{g}=\mathrm{Lie}(G)$ et $K$ est un sous-groupe compact maximal de $G$. On sait que ces $(\mathfrak{g},K)$-modules sont exactement les $(\mathfrak{g},K)$-modules irrŽductibles \textit{unitarisables} et que deux reprŽsentations unitaires irrŽductibles de $G$ sont isomorphes si et seulement si 
leurs $(\mathfrak{g},K)$-modules sous-jacents le sont.  

Les $(\mathfrak{g},K)$-modules irrŽductibles ont ŽtŽ classifiŽs par Langlands \cite{La1}; cette classification est exactement la mme que dans le cas non archimŽdien (Borel-Wallach \cite[IV, theorem~4.11]{BW}). Quant au dual unitaire de $G$, il a ŽtŽ classifiŽ par Vogan \cite{V2}; on utilisera ici la mme prŽsentation que celle de Tadi\'c pour les corps locaux non archimŽdiens (cf. \cite{T1}, voir aussi la fin de cette sous-section). D'ailleurs gr‰ce ˆ la preuve de Baruch \cite{Bar} de la conjecture de Kirillov, on a une dŽmonstration analogue dans les cas archimŽdiens et non archimŽdiens \cite{T2}. En effet la conjecture de Kirillov donne une preuve simple du rŽsultat suivant \cite{V2} (voir aussi \cite{Bar}): l'induite parabolique (normalisŽe) d'une reprŽsentation unitaire irrŽductible d'un facteur de Levi de $G$ est irrŽductible. 

Il y a cependant quelques diffŽrences entre les cas archimŽdiens et non archimŽdiens. La premire est que les reprŽsentations irrŽductibles de carrŽ intŽgrable (modulo le centre) n'existent qu'en petite dimension: $n=1$ si $F\simeq \mathbb{C}$; $n\in\{1,2\}$ si $F \simeq \mathbb{R}$. La seconde concerne les reprŽsentations gŽnŽriques. 

\vskip1mm
{\bf $\bullet$ ReprŽsentations elliptiques.} Tout comme les reprŽsentations irrŽductibles de carrŽ intŽgrable, les reprŽsentations elliptiques n'existent qu'en petite dimension (sinon il n'y a pas d'ŽlŽment semi-simple rŽgulier elliptique): $n=1$ si $F\simeq \mathbb{C}$; $n\in\{1,2\}$ si $F \simeq \mathbb{R}$. 

Si $n=1$, les reprŽsentations de $G= F^\times$ sont les caractres et elles sont toutes elliptiques. Un caractre de $\mathbb{C}^\times$ est de la forme $\xi_{s,k}(z)= \vert z\vert_{\mathbb{C}}^s z^k$ pour un $s\in \mathbb{C}$ et un $k\in \mathbb{Z}$. Ainsi $\xi_{s,k}$ est unitaire si et seulement si $\Re(s)+k=0$; et 
$\xi_{s,k}$ est invariant par la conjugaison complexe si et seulement si $k=0$ ($\xi_{k,0}= \nu_{\mathbb{C}}^s$). Soit 
$\kappa: \mathbb{R}^\times \rightarrow \{\pm 1\}$ le caractre signe. Un caractre de $\mathbb{R}^\times$ est de la forme $\chi_s=\nu_{\mathbb{R}}^s$ 
ou $\chi'_s=\kappa \nu_{\mathbb{R}}^s$ pour un $s \in \mathbb{C}$. Observons que $\chi_s \circ \mathrm{N}_{\mathbb{C}/\mathbb{R}}=\chi'_s\circ \mathrm{N}_{\mathbb{C}/\mathbb{R}}= \xi_{2s,0}$. 

Traitons maintenant le cas $n=2$ et $F= \mathbb{R}$, \ie $G= \mathrm{GL}_2(\mathbb{R})$. Reprenons brivement la description de \cite[6.4]{BH}. Pour $s,\,t\in \mathbb{C}$, on sait que:
\begin{itemize}
\item la reprŽsentation $\chi_s \times \chi_t$, resp. 
$\chi'_s \times \chi'_t,$ est rŽductible si et seulement si $s- t\in 2\mathbb{Z}+1$; 
\item la reprŽsentation $\chi_s \times \chi'_t$, resp. 
$\chi'_s \times \chi_t$, est rŽductible si et seulement si $s - t \in 2\mathbb{Z}$.
\end{itemize}
Soit $a\in \mathbb{N}^*$. On pose $$R(a)= \left\{\begin{array}{ll} \chi_{\frac{a}{2}}\times \chi_{-\frac{a}{2}} & \hbox{si $a$ est impair}\\
\chi_{\frac{a}{2}}\times \chi'_{-\frac{a}{2}} & \hbox{si $a$ est pair}
\end{array}\right..$$ La reprŽsentation $R(a)$ a exactement deux sous-quotients irrŽductibles, notŽs $\delta(a)$ et $L(a)$: 
$\delta(a)$ est l'unique sous-reprŽsentation irrŽductible de $R(a)$ et $L(a)$ est le quotient 
de Langlands de $R(a)$. La reprŽsentation $\delta(a)$ est essentiellement de carrŽ intŽgrable et $\kappa$-stable; 
la reprŽsentation $L(a)$ est de dimension finie $a$ et n'est pas $\kappa$-stable. On sait aussi que $L(1)$ est le caractre trivial de $G$. 

\begin{remark}
\textup{Soit $\mu(a)$ le caractre de $\mathbb{R}^\times$ dŽfini par 
$$\mu(a)= \left\{\begin{array}{ll} \chi_{\frac{a}{2}}(\chi_{-\frac{a}{2}})^{-1}& \hbox{si $a$ est impair}\\
\chi_{\frac{a}{2}}(\chi'_{-\frac{a}{2}})^{-1}& \hbox{si $a$ est pair}\end{array}\right..
$$ On a $\mu(a)(x)= x^a \kappa(x)$. \hfill $\blacksquare$}
\end{remark} 

Ainsi aprs semi-simplification, on a l'ŽgalitŽ dans le groupe 
de Grothendieck des reprŽsentations de longueur finie de $G$: 
$$\kappa R(a)= \delta(a)+ \kappa L(a) = \left\{ \begin{array}{ll}
\chi'_{\frac{a}{2}}\times \chi'_{-\frac{a}{2}}&\hbox{si $a$ est impair}\\
\chi'_{\frac{a}{2}}\times \chi_{-\frac{a}{2}} & \hbox{si $a$ est pair} \end{array}\right..$$ Un caractre de $G$ est une reprŽsentation de la forme $\chi_s\circ \det $ ou $\chi'_s\circ\det$ pour un $s\in \mathbb{C}$. Les reprŽsentations $\delta(a)$ et $L(a)$ pour $a\in \mathbb{N}^*$ sont toutes elliptiques; et toute reprŽsentation elliptique de $G$ est isomorphe ˆ 
$\chi_s\delta(a)= (\chi_s\circ\det)\otimes \delta(a)$, $\chi_s L(a)$ ou $\chi'_sL(a)$ pour un $s\in \mathbb{C}$ et un $a\in \mathbb{N}^*$. 

\vskip1mm
{\bf $\bullet$ ReprŽsentations  unitaires gŽnŽriques.}
Fixons un caractre additif non trivial $\psi$ de $F$. Si $(\pi,V)$ est une reprŽsentation unitaire de $G$, on note $V^\infty\subset V$ le sous-espace des vecteurs \textit{lisses}; il est naturellement muni d'une structure d'espace topologique de FrŽchet. On note $\pi^\infty$ la reprŽsentation (continue) de $G$ d'espace 
$V^\infty$. On dŽfinit comme dans le cas non archimŽdien l'espace $\ES{D}(\pi^\infty)=\ES{D}(\pi^\infty\!,\psi)$ des fonctionnelles de Whittaker 
\textit{continues} sur $V^\infty$, et on dit que $\pi$ est \textit{gŽnŽrique} si $\ES{D}(\pi^\infty)\neq \{0\}$. On sait que si $\pi$ est irrŽductible, l'espace $\ES{D}(\pi^\infty)$ est de dimension $0$ ou $1$ \cite{Sha}. De plus (toujours si $\pi$ est unitaire irrŽductible), la gŽnŽricitŽ de $\pi$ se lit sur le $(\mathfrak{g},K)$-module sous-jacent $V^{\infty,\mathrm{fin}}\subset V^\infty$ (rappelons que $V^{\infty,\mathrm{fin}}$ est dense dans $V$): d'aprs Vogan \cite[p.~98]{V1}, les rŽsultats de Kostant \cite{Ko} entra"nent qu'une reprŽsentation unitaire irrŽductible de $G$ est gŽnŽrique si et seulement si le $(\mathfrak{g},K)$-module sous-jacent est \guill{gros} (voir aussi Henniart \cite{H1}); 
auquel cas on peut dŽfinir des fonctionnelles de Whittaker continues sur 
le $(\mathfrak{g},K)$-module sous-jacent, mais en gŽnŽral elles forment un espace vectoriel de dimension $>1$. En tous cas, la classification des reprŽsentations unitaires irrŽductibles gŽnŽriques de $G$ est la mme que dans le cas non archimŽdien. 

\vskip1mm 
En conclusion, toute reprŽsentation unitaire irrŽductible $\pi$ de $G$ est un produit de reprŽsentations du type $u(\delta,k)$ et $u(\delta,k;\alpha)$ o $\delta$ est une reprŽsentation irrŽductible de carrŽ intŽgrable de $\mathrm{GL}_{n'}(F)$, $k\in \mathbb{N}^*$ et 
$\alpha \in ]0,\frac{1}{2}[$; avec $n'\in \{1,2\}$ si $F\simeq \mathbb{R}$ et $n'=1$ si $F\simeq \mathbb{C}$. La classe d'isomorphisme de $\pi$ dŽtermine les facteurs du produit ˆ permutation prs. 
Pour qu'un tel produit soit gŽnŽrique il faut et il suffit que tous les facteurs soient du type $u(\delta,k=1)$ ou $u(\delta,k=1;\alpha)$. 

\subsection{Normalisation $A_\pi$ ($E/F\simeq \mathbb{C}/\mathbb{R}$).}  
On suppose que $E/F$ est une extension isomorphe ˆ $\mathbb{C}/\mathbb{R}$. Soit $\pi$ une reprŽsentation unitaire irrŽductible $\kappa$-stable de $G$. Si 
$\pi$ est gŽnŽrique, on dŽfinit comme dans le cas non archimŽdien un opŽrateur d'entrelacement normalisŽ $I_{\pi}^{\mathrm{g\acute{e}n}}\in \mathrm{Isom}_G(\kappa\pi^\infty\!,\pi^\infty)$ sur l'espace $V_{\pi^\infty}=(V_\pi)^\infty$ de $\pi^\infty$. D'aprs \ref{op unit}\,(ii), cet opŽrateur est unitaire. Il s'Žtend donc de manire unique en un opŽrateur unitaire sur l'espace $V_\pi$ de $\pi$; ce dernier donne par restriction un opŽrateur sur le $(\mathfrak{g},K)$-module sous-jacent. On note encore $I_\pi^{\mathrm{g\acute{e}n}}\in \mathrm{Isom}_G(\kappa\pi,\pi)$ cet opŽrateur normalisŽ sur $V_\pi$, ainsi que sa restriction au $(\mathfrak{g},K)$-module sous-jacent. Si $\wt{\pi}= \nu^r\pi \;(=(\nu^r\circ \det)\otimes \pi)$ pour un $r\in \mathbb{R}$, on pose $A_{\wt{\pi}}^{\mathrm{g\acute{e}n}}= A_\pi^{\mathrm{g\acute{e}n}}$. Puisqu'une reprŽsentation irrŽductible tempŽrŽe est gŽnŽrique unitaire, gr‰ce ˆ la classification de Langlands, on dŽfinit comme dans le cas non archimŽdien un opŽrateur d'entrelacement normalisŽ $A_\pi\in \mathrm{Isom}_G(\kappa\pi,\pi)$ pour toute reprŽsentation irrŽductible $\kappa$-stable $\pi$ de $G$. 
Si $\pi$ est gŽnŽrique unitaire, on a l'ŽgalitŽ $A_\pi= A_\pi^{\mathrm{g\acute{e}n}}$. Pour obtenir cette ŽgalitŽ, on remplace les rŽsultats de Jacquet-Shalika \cite{JS2} utilisŽs dans le cas non archimŽdien, par ceux de Wallach \cite[15.6.7]{Wa}. 

\subsection{Fonctions concordantes ($E/F\simeq \mathbb{C}/\mathbb{R}$).}\label{fonctions concor archim}
On suppose toujours que $E/F$ est une extension isomorphe ˆ $\mathbb{C}/\mathbb{R}$. 
On fixe un entier $m\geq 1$ et on pose $n=2m$. 
Soient $G=\GLn(F)$ et $H= \GLm(E)$. On note $\sigma$ le gŽnŽrateur de $\Gamma$ et $\kappa$ le caractre non trivial de $F^\times/\mathrm{N}_{E/F}(E^\times)$. 

On dŽfinit comme dans le cas non archimŽdien, via les choix d'une $F$-base de $E^m$ (qui dŽfinit un plongement de $H$ dans $G$) et d'un ŽlŽment 
$e\in E^\times$ tel que $\sigma e = (-1)^{2(m-1)}$, un facteur de transfert $\Delta(\gamma)=\kappa(e\wt{\Delta}(\gamma))\in \{\pm 1\}$ pour tout ŽlŽment $\gamma\in \HGreg$. 
On note $C^\infty_{\mathrm{c}}(G)$ l'espace des fonctions lisses ˆ support compact sur $G$.  On dŽfinit la notion de fonctions $f\in C^\infty_{\mathrm{c}}(G)$ et $\phi\in C^\infty_{\mathrm{c}}(H)$ \textit{concordantes} comme dans le cas non archimŽdien (cf. \ref{fonctionsconcor}), via les  
choix de mesures de Haar $\d g$ sur $G$ et $\d h$ sur $H$. 
D'aprs Shelstad \cite[cor.~2.2]{She}, on a l'analogue de \ref{existence transfert}\,(i): pour toute fonction $f\in \CG$, il existe une fonction $\phi\in \CH$ qui concorde avec $f$. 

On a aussi la version $K$-finie de ce rŽsultat. Soit $C^{\infty,K-\mathrm{fin}}_{\mathrm{c}}(G)\subset C^\infty_{\mathrm{c}}(G)$ le sous-espace des fonctions qui sont $K$-finies ˆ droite et ˆ gauche. Soit aussi $C^{\infty,K_H-\mathrm{fin}}_{\mathrm{c}}(H) \subset \CH$ le sous-espace des fonctions qui sont $K_H$-finies ˆ droite et ˆ gauche, pour un sous-groupe compact maximal $K_H$ de $H$. D'aprs M\oe glin-Waldspurger \cite[ch.~IV, 3.4]{MW2}, pour toute fonction $f\in C^{\infty,K-\mathrm{fin}}_{\mathrm{c}}(G)$, il existe une fonction $\phi\in C^{\infty,K_H-\mathrm{fin}}_{\mathrm{c}}(H)$ qui concorde avec $f$. 

\subsection{$\kappa$-relvements ($E/F\simeq \mathbb{C}/\mathbb{R}$).}\label{kappa-rel archim} Continuons avec les hypothses et notations de \ref{fonctions concor archim}. Commenons par rappeler quelques rŽsultats bien connus (essentiellement džs ˆ Harish-Chandra). 
Si $\pi$ est une reprŽsentation unitaire irrŽductible de $G$, elle est admissible, par suite pour toute fonction $f\in C^\infty_{\mathrm{c}}(G)$, l'opŽrateur $\pi(f)= \int_G f(g)\pi(g) \d g$ est ˆ trace; et la distribution $f \mapsto \mathrm{tr}(\pi(f)))$ sur $G$ est reprŽsentŽe par une fonction analytique $\Theta_\pi$ sur $\Greg$, localement intŽgrable sur $G$ (cf. \cite[\S 10]{Kn})\footnote{Les reprŽsentations unitaires irrŽductibles de $H$ vŽrifient bien sžr les mmes propriŽtŽs. D'ailleurs par restriction des scalaires, on peut toujours voir $H$ comme un groupe de Lie rŽel.}. Si de plus $\pi$ est $\kappa$-stable 
et si $A\in \mathrm{Isom}_G(\kappa\pi,\pi)$, pour toute fonction $f\in \CG$, l'opŽrateur $\pi(f)\circ A$ est lui aussi ˆ trace; et la distribution $f\mapsto \mathrm{tr}(\pi(f)\circ A)$ sur $G$ est reprŽsentŽe par une fonction analytique $\Theta_\pi^{A}$ sur $\Greg$, localement intŽgrable sur $G$ (cf. \cite[3.6]{H1}).

On dŽfinit la notion de \textit{$\kappa$-relvement} comme dans le cas non archimŽdien (cf. \ref{defkapparel}), en se limitant reprŽsentations unitaires: une reprŽsentation unitaire irrŽductible $\kappa$-stable $\pi$ de $G$ est un $\kappa$-relvement d'une reprŽsentation unitaire irrŽductible $\tau$ de $H$ si pour tout 
$A\in \mathrm{Isom}_G(\kappa\pi,\pi)$, il existe une constante $c=c(\tau,\pi,A)>0$ telle que $\mathrm{tr}(\pi(f)\circ A)= c\, \mathrm{tr}(\tau(\phi))$ pour toute paire 
de fonctions concordantes $(f,\phi)\in \CGreg \times \CHGreg$. On a aussi la formulation Žquivalente en termes de fonctions-caractres (cf. \ref{def kappa-rel}\,(2)). 
Un tel $\kappa$-relvement, s'il existe, est unique ˆ isomorphisme prs. 

D'aprs la formulation en termes de fonctions-caractres (jointe ˆ la formule d'intŽgration de Weyl) si une reprŽsentation unitaire irrŽductible $\tau$ de $H$ admet un $\kappa$-relvement $\pi$ ˆ $G$ et si 
$A\in \mathrm{Isom}_G(\kappa\pi,\pi)$, alors l'ŽgalitŽ $$\mathrm{tr}(\pi(f)\circ A)= c(\tau,\pi,A)\mathrm{tr}(\tau(\phi))$$ est vraie pour toute paire de fonctions concordantes 
$(f,\phi)\in C^\infty_{\mathrm{c}}(G)\times C^\infty_{\mathrm{c}}(H)$. On en dŽduit aussi que pour que 
$\pi$ soit un $\kappa$-relvement de $\tau$, il faut et il suffit que l'ŽgalitŽ ci-dessus soit vŽrifiŽe pour toute 
paire de fonctions concordantes 
$(f,\phi)\in C^{\infty,K-\mathrm{fin}}_{\mathrm{c}}(G)\times C^{\infty,K_H-\mathrm{fin}}_{\mathrm{c}}(H)$. 

Henniart \cite{H1} prouve que toute reprŽsentation unitaire irrŽductible gŽnŽrique de $H$ admet un $\kappa$-relvement ˆ $G$. Il prouve aussi (\textit{loc.\,cit.}) que la constante $c(\tau,\pi, A_\pi^{\mathrm{g\acute{e}n}})$ est indŽpendante des reprŽsentations (unitaires irrŽductibles gŽnŽriques) telles que $\pi$ soit un $\kappa$-relvement de $\tau$; d'ailleurs comme dans le cas non archimŽdien, si l'on normalise le facteur de transfert comme dans Kottwitz-Shelstad \cite{KS}, cette constante vaut $1$ (Hiraga-Ichino \cite{HI}). La preuve d'Henniart 
est purement locale et consiste ˆ se ramener, gr‰ce ˆ la compatibilitŽ entre les opŽrations de $\kappa$-relvement et d'induction parabolique 
(cf. \cite[4]{H1}), au cas o $m=1$, \cad ˆ l'induction automorphe de $\mathbb{C}^\times$ ˆ $\mathrm{GL}_2(\mathbb{R})$, qui dŽcoule essentiellement des rŽsultats de Labesse-Langlands  \cite{LL}. 

\subsection{Le cas d'une $F$-algbre cyclique $E/F$ ($F\simeq \mathbb{R}$ ou $\mathbb{C}$).} Les rŽsultats de \cite{H1} sont en fait dŽmontrŽs dans le cadre plus gŽnŽral d'une $F$-algbre cyclique de degrŽ fini $d$. On Žcrit $E=E_1\times \cdots \times E_r$ o $E_i/F$ est une extension cyclique de degrŽ $s= \frac{d}{r}$. Comme en \ref{notations ext cyclique}, on suppose fixŽ un gŽnŽrateur $\sigma$ de $\Gamma$ tel que $\sigma E_i=E_{i+1}$ ($i=1,\ldots , r-1$) et $\sigma E_r= E_1$. On a deux cas possibles: $s=1$, \textit{i.e} $E=F^d$; ou bien $s=2$, \ie $E_i/F\simeq \mathbb{C}/\mathbb{R}$. On dŽfinit comme dans le cas non archimŽdien 
le facteur de transfert $\Delta(\gamma)\in \{\pm 1\}$ pour tout $\gamma\in \HGreg$, via les choix d'une $F$-base de $E^m$ (qui dŽfinit un plongement de $H$ dans $G$) et d'un ŽlŽment $e\in E^\times$ tel que $\sigma e = (-1)^{d(m-1)}e$. Observons que $\Delta \equiv 1$ si $s=1$. 
D'o (comme dans le cas $s=d=2$) la notion de fonctions $f\in C^\infty_{\mathrm{c}}(G)$ et $\phi\in C^\infty_{\mathrm{c}}(G)$ concordantes, et de $\kappa$-relvement (relativement ˆ $E/F$) pour les reprŽsentations irrŽductibles unitaires de $H$, qui s'exprime aussi comme une ŽgalitŽ de fonctions-caractres. 

Henniart \cite{H1} prouve que toute reprŽsentation unitaire irrŽductible gŽnŽrique de $H$ admet un $\kappa$-relvement ˆ $G$ (relativement ˆ $E/F$). 
PrŽcisŽment, quitte ˆ changer la $F$-base de $E^m$, on peut supposer que le plongement $\varphi$ de $H$ dans $G$ se dŽcompose en $\varphi_1\times \cdots \times \varphi_r$ o $\varphi_i$ est un plongement de $\GLm(E_i)$ dans $\mathrm{GL}_{sm}(F)$. Soit $\tau$ une reprŽsentation unitaire irrŽductible gŽnŽrique de $H$. Elle s'Žcrit comme un produit tensoriel complŽtŽ 
$\tau_1\widehat{\otimes} \cdots \widehat{\otimes}\tau_r$ o $\tau_i$ est une reprŽsentation unitaire irrŽductible gŽnŽrique de $\GLm(E_i)$. 
Si $s=1$ (\ie $E_i=F$), d'aprs \textit{loc.\,cit.}, l'induite parabolique (normalisŽe) de $\tau$ relativement au sous-groupe parabolique standard de $G$ de composante de Levi $\mathrm{GL}_m(F)^d$, est un $\kappa$-relvement de $\tau$. Si $s=2$ (\ie $E_i/F\simeq \mathbb{C}/ \mathbb{R}$), pour $i=1,\ldots ,r$, soit $\pi_i$ un $\kappa$-relvement de $\tau_i$ (relativement ˆ $E_i/F$); d'aprs \textit{loc.\,cit.}, l'induite parabolique $\pi$ de $\pi_1\widehat{\otimes} \cdots 
\widehat{\otimes} \pi_r$ relativement au sous-groupe parabolique standard de composante de Levi $\mathrm{GL}_{2m}(F)^r$ est irrŽductible et c'est un $\kappa$-relvement de $\tau$.

\begin{remark}
\textup{
La construction ci-dessus est la version archimŽdienne de la variante \ref{variante inverse} du cas non archimŽdien. Comme toutes les induites paraboliques considŽrŽes ici sont irrŽductibles, elle est Žquivalente ˆ celle proposŽe en \ref{kappa-rel bis}. En effet si $s=1$, c'est exactement la mme construction, \ie une 
induction parabolique de $\GLm(F)^d$ ˆ $G$. Si $s=2$, posons 
$H^\natural=(\GLm(E_1))^r$ et notons $\iota: H \rightarrow H^\natural$ l'isomorphisme 
$(x_1,\ldots ,x_r)\mapsto (x_1,\sigma^{-1} x_2, \ldots , \sigma^{1-r} x_r)$. Posons $H_1= \mathrm{GL}_{mr}(E_1)$ et notons $\varphi_1$ le plongement de $H_1$ dans $G$ dŽduit de $\varphi$ comme en \ref{facteurs de transfert bis}. Soit $\rho$ l'induite parabolique de $\tau^\natural= \tau\circ \iota^{-1}\;(=\tau_1^\natural\wh{\otimes}\cdots \wh{\otimes} \tau^\natural_r)$ ˆ $H_1$ suivant le sous-groupe parabolique standard de composante de Levi $H^\natural$. Puisque $\rho$ est irrŽductible, par compatibilitŽ entre les opŽrations de $\kappa$-relvement et d'induction parabolique (\ie l'analogue de \ref{compatibilitŽ IP-IA}), on a  
que $\pi$ est un $\kappa$-relvement de $\rho$ (relativement ˆ $E_1/F$). \hfill $\blacksquare$
}\end{remark}

Dans \cite{H1}, Henniart travaille exclusivement avec les fonctions-carac\-tres et n'utilise pas de fonctions concordantes. On peut cependant Žtendre les rŽsultats de \ref{fonctions concor archim} au cas d'une $F$-algbre cyclique. Soit $K_1$ un sous-groupe compact maximal de $H_1=\mathrm{GL}_{mr}(E_1)$ en 
bonne position par rapport ˆ $A_{0,H_1}= (E_1^\times)^{rm}$; on a la dŽcomposition d'Iwasawa 
$H_1=K_1P_1$ o $P_1$ est le sous-groupe parabolique standard de $H_1$ de composante de Levi $H^\natural$, et $K_{H^\natural} = K_1 \cap H^\natural$ est un sous-groupe compact maximal de $K_{H^\natural}$. On peut prendre pour $K_H$ le groupe $\iota^{-1}(K_{H^\natural})$. Pour $\phi \in C^\infty_{\mathrm{c}}(H_1)$, on dŽfinit comme en \ref{facteurs de transfert bis} le terme constant $(\phi_1)^{P_1}\in C^\infty_{\mathrm{c}}(H^\natural)$ de $\phi_1$ suivant 
$(K_1,P_1)$. Observons que si $\phi_1$ est $K_1$-finie ˆ droite et ˆ gauche, alors $(\phi_1)^{P_1}$ est $K_{H^\natural}$-finie ˆ droite et ˆ gauche. 
Le lemme \ref{existence transfert bis} est encore vrai ici. En particulier pour toute fonction $f\in C^\infty_{\mathrm{c}}(G)$, resp. $f\in C^{\infty,K-\mathrm{fin}}_{\mathrm{c}}(G)$, il existe une fonction $\phi\in C^\infty_{\mathrm{c}}(H)$, resp. $\phi\in C^{\infty,K_H-\mathrm{fin}}_{\mathrm{c}}(H)$, qui concorde avec $f$.

\section{\'EnoncŽ des rŽsultats}\label{ŽnoncŽ des rŽsultats}

Dans cette section, on Žnonce les rŽsultats locaux et globaux dŽmontrŽs dans cet article. 
Tous ces rŽsultats seront dŽmontrŽs dans la section \ref{demo}, ˆ l'exception de la proposition \ref{prop ell} qui sera dŽmontrŽe dans la section 
\ref{kappa-rel elliptiques}. En \ref{principe 1}, on donne les principaux arguments de la preuve des thŽormes \ref{theo local} (local) et \ref{theo global} (global), qui est assez tortueuse. 

\subsection{RŽsultats locaux.} Soit $E/F$ une extension finie cyclique de corps locaux (archimŽdiens ou non archimŽdiens), de degr\'e $d$. 
On note $X(E/F)$ le groupe des caractres de 
$F^\times/ \mathrm{N}_{E/F}(E^\times)$ et on fixe un gŽnŽrateur $\kappa$ de $X(E/F)$. 

Si $\pi$ est une reprŽsentation irrŽductible de $\GLn(F)$ pour un entier $n\geq 1$, on note $X(\pi)$ l'ensemble des classes 
d'isomorphisme de reprŽsentations $\kappa^i\pi= (\kappa^i\circ \det)\otimes \pi$ pour $i\in \mathbb{Z}$, et $x(\pi)$ le cardinal de $X(\pi)$. 
Observons que si $\pi$ est une reprŽsentation de Speh de la forme $u(\delta,q)$ pour un entier $q\geq 1$ divisant $n$ et une reprŽsentation irrŽductible de carrŽ intŽgrable $\delta$ de $\mathrm{GL}_{\smash{\frac{n}{q}}}(F)$, d'aprs l'unicitŽ du quotient de Langlands, on a l'ŽgalitŽ $x(\pi)=x(\delta)$. 

Si $\tau$ est une repr\'esentation irr\'eductible de $\GLm(E)$ pour un entier $m\geq 1$, on note $\Gamma(\tau)$ l'ensemble des classes d'isomorphisme de repr\'esentations ${^\gamma\tau}= \tau\circ \gamma^{-1}$ pour $\gamma\in \Gamma(E/F)$. On note $g(\tau)$ le cardinal de $\Gamma(\tau)$ et $r(\tau)$ celui du 
stabilisateur $\{\gamma\in \Gamma(E/F)\,\vert\, {^\gamma\tau} \simeq \tau\}$. On a donc l'Žgalit\'e $g(\tau)r(\tau)=d$. Observons que si 
$\tau$ est une reprŽsentation de Speh de la forme $u(\delta_E,q)$ pour un entier $q\geq 1$ divisant $m$ et une reprŽsentation irrŽductible de carrŽ intŽgrable $\delta_E$ de $\mathrm{GL}_{\smash{\frac{m}{q}}}(E)$, on a l'ŽgalitŽ $r(\tau)=r(\delta_E)$. 

Si $\delta_E$ est une reprŽsentation irrŽductible de carrŽ intŽgrable de $\GLm(E)$ pour un entier $m\geq 1$, on sait d'aprs \cite{HH} qu'il existe un $\kappa$-relvement $\pi$ de $\delta_E$. Pr\'ecisŽment, l'entier $r= r(\delta_E)$ divise $md$ et il existe une reprŽsentation irrŽductible de carrŽ intŽgrable $\kappa^r$-stable $\delta$ de $\mathrm{GL}_{n'}(F)$, $n'= \frac{md}{r}$, telle que $\pi\simeq  \delta \times \kappa\delta \times \cdots \times \kappa^{r-1}\delta$. On a 
$x(\delta)=r$, l'ensemble $X(\delta)$ est dŽterminŽ de manire unique par la classe d'isomorphisme de $\delta_E$, et $\delta$ est cuspidale si et seulement si $\delta_E$ est cuspidale. 
RŽciproquement (\textit{loc.\,cit.}), toute reprŽsentation irrŽductible de $\GLmd(F)$ de la forme $\pi=\delta\times \kappa\delta\times \cdots \times \kappa^{r-1}\delta$ pour un entier $r\geq 1$ divisant $d$ et une reprŽsentation irrŽductible de carrŽ intŽgrable $\kappa^r$-stable $\delta$ de $\mathrm{GL}_{n'}(F)$, $n'=\frac{md}{r}$, telle que $x(\delta)=r$, est un $\kappa$-relvement d'une reprŽsentation irrŽductible de carrŽ intŽgrable $\delta_E$ de $\GLm(E)$. De plus les classes d'isomorphisme de reprŽsentations irrŽductibles \textit{gŽnŽriques}\footnote{On verra plus loin qu'on peut supprimer le \guill{gŽnŽriques} ici; voir \ref{remarque sur les ell}\,(ii).} de $\GLm(E)$ ayant $\pi$ pour $\kappa$-relvement sont les ŽlŽments de $\Gamma(\delta_E)$. 

Ces rŽsultats s'Žtendent naturellement aux reprŽsentations essentiellement de carrŽ intŽgrable. Soit $\delta_E$ une reprŽsentation irrŽductible de carrŽ intŽgrable de $\GLm(E)$ ayant pour $\kappa$-relvement $\pi= \delta\times \kappa\delta\times \cdots \times \kappa^{r-1}\delta$, o $r=r(\delta_E)$ et $\delta$ est une reprŽsentation irrŽductible de carrŽ intŽgrable $\kappa^r$-stable de $\mathrm{GL}_{n'}(F)$, $n'=\frac{n}{r}$. Pour $i\in \mathbb{Z}$, la reprŽsentation (irrŽductible essentiellement de carrŽ intŽgrable)  $\tilde{\delta}_E = \nu_E^i\delta$ de $\GLm(E)$ a pour $\kappa$-relvement $\nu^i \pi \simeq \tilde{\delta}\times \kappa\tilde{\delta}\times \cdots \times 
\kappa^{r-1}\tilde{\delta}$ avec  $\tilde{\delta}=\nu^i \delta$. La description de l'image et des fibres de l'application de $\kappa$-relvement pour les reprŽsentations essentiellement de carrŽ intŽgrables est la mme que pour les reprŽsentations de carrŽ intŽgrable.

La proposition et le thŽorme ci-dessous Žtendent ces rŽsultats aux reprŽsentations \textit{elliptiques}, resp. \textit{de Speh}, de $\GLm(E)$. 

\begin{proposition}\label{prop ell}
\begin{enumerate} 
\item[(i)]  Soit $\tilde{\delta}_E$ une reprŽsentation irrŽductible elliptique de $\GLm(F)$ pour un entier $m\geq 1$. Posons $r=r(\tilde{\delta}_E)$. Il existe une reprŽsentation irrŽductible elliptique $\kappa^r$-stable $\tilde{\delta}$ de $\mathrm{GL}_{n'}(F)$, $n'=\frac{md}{r}$, telle que $\pi= \tilde{\delta} \times \kappa \tilde{\delta}\times \cdots \times \kappa^{r-1}\tilde{\delta}$ soit un $\kappa$-relvement de $\tilde{\delta}_E$. De plus on a $x(\tilde{\delta})=r$ et $\tilde{\delta}$ est essentiellement de carrŽ intŽgrable (\ie gŽnŽrique) si et seulement si $\tilde{\delta}_E$ est essentiellement de carrŽ intŽgrable.
\item[(ii)] Soit un entier $r\geq 1$ divisant $d$ et soit $\tilde{\delta}$ une reprŽsentation irrŽductible elliptique de $\mathrm{GL}_{n'}(F)$ pour un entier $n'\geq 1$, telle que 
$x(\delta)=r$. Alors $d$ divise $rn'$ et $\pi= \tilde{\delta} \times \kappa\tilde{\delta} \times\cdots \times \kappa^{r-1}\tilde{\delta}$ est le $\kappa$-relvement d'une reprŽsentation irrŽductible elliptique 
$\tilde{\delta}_E$ de $\GLm(E)$, $m= \frac{rn'}{d}$. De plus on a $r(\tilde{\delta}_E)=r$ et les classes d'isomorphisme de reprŽsentations irrŽductibles de $\GLm(E)$ ayant $\pi$ pour $\kappa$-relvement sont les ŽlŽments de $\Gamma(\tilde{\delta}_E)$. 
\end{enumerate}
\end{proposition}

\begin{remark}\label{remarque sur les ell}
\textup{\begin{enumerate}
\item[(i)]La proposition \ref{prop ell} est essentiellement dŽmontrŽe dans \cite{Fa} (cas non archimŽdien) et dans \cite{H1} (cas archimŽdien). Dans les deux cas, la preuve est locale; mais dans le cas non archimŽdien, elle utilise le $\kappa$-relvement des reprŽsentations irrŽductibles de carrŽ intŽgrable qui a ŽtŽ obtenu par voie globale \cite{HH}. 
\item[(ii)]Le point (ii) de \ref{prop ell} amŽliore lŽgrement la description des fibres de l'application de $\kappa$-relvement des reprŽsentations essentiellement de carrŽ intŽgrable dans le cas non archimŽdien \cite[5.2, cor.]{HH}.   
\end{enumerate}}
\end{remark}

\begin{theorem}\label{theo local}
\begin{enumerate}
\item[(i)]Soit un entier $m\geq 1$. Soit $u_E=u(\delta_E,q)$ une reprŽsentation de Speh de $\mathrm{GL}_m(E)$, o $q\geq 1$ est un entier divisant $m$ et $\delta_E$ est une 
reprŽsentation irrŽductible de carrŽ intŽgrable de $\mathrm{GL}_{\frac{m}{q}}(E)$. Posons $r=r(\delta_E)$. Soit 
$\delta$ une reprŽsentation irrŽductible de carrŽ intŽgrable $\kappa^r$-stable de $\mathrm{GL}_{\frac{md}{rq}}(F)$ telle que $\delta \times \kappa \delta \times \cdots \times \kappa^{r-1}\delta$ soit un $\kappa$-relvement de $\delta_E$. Posons $u= u(\delta,q)$. On a $x(u)=x(\delta)=r$ et 
$\pi = u \times \kappa u \times \cdots \times \kappa^ru$ est un $\kappa$-relvement de $u_E$.  
\item[(ii)] Soit un entier $r\geq 1$ divisant $d$ et soit $u$ une reprŽsentation de Speh de $\mathrm{GL}_{n'}(F)$ pour un entier $n'\geq 1$, telle que $x(u)=r$. 
Alors $d$ divise $rn'$ et il existe une reprŽsentation de Speh $u_E$ de $\mathrm{GL}_{m}(E)$, $m= \frac{rn'}{d}$, telle que $\pi= u\times \kappa u \times\cdots \times \kappa^{r-1}u$ soit un $\kappa$-relvement de $u_E$. De plus on a $r(u_E)=r$ et les classes d'isomorphisme de reprŽsentations irrŽductibles de $\GLm(E)$ ayant $\pi$ pour $\kappa$-relvement sont les ŽlŽments de $\Gamma(u_E)$.   
\end{enumerate}
\end{theorem}

La description des reprŽsentations irrŽductibles unitaires de $\GLm(E)$ et la compatibilitŽ entre les opŽrations de $\kappa$-relvement et 
d'induction parabolique entra"nent la 

\begin{proposition}\label{prop unit}
Soit un entier $m\geq 1$.
\begin{enumerate}
\item[(i)] Pour toute reprŽsentation irrŽductible unitaire $\tau$ de $\GLm(E)$, il existe une reprŽsentation irrŽductible unitaire $\pi$ de $\GLmd(F)$ qui soit un $\kappa$-relvement de $\tau$. De plus $\tau$ est gŽnŽrique si et seulement si $\pi$ est gŽnŽrique.
\item[(ii)] Toute reprŽsentation irrŽductible unitaire $\kappa$-stable $\pi$ de $\GLmd(F)$ est le $\kappa$-relvement d'une reprŽsentation irrŽductible unitaire $\tau$ de $\GLm(E)$. 
\item[(iii)] Soit $\tau= \tau_1\times \cdots \times \tau_s$ une reprŽsentation irrŽductible unitaire de $\GLm(E)$, 
o $\tau_i$ est une reprŽsentation du type $u(\delta_{E,i},q_i)$ ou $u(\delta_{E,i},q_i;\alpha_i)$, $\delta_{E,i}$ est  
une reprŽsentation irrŽductible de carrŽ intŽgrable $\delta_{E,i}$ de $\mathrm{GL}_{a_i}(E)$ pour un entier  $a_i\geq 1$, 
$q_i\in \mathbb{N}^*$ et $\alpha_i\in\;]0,\frac{1}{2}[$. Les reprŽsentations irrŽductibles \textup{unitaires} de $\GLm(E)$ ayant mme 
$\kappa$-relvement que $\tau$ sont, ˆ isomorphisme et permutation des facteurs prs, les reprŽsentations 
$ {^{\gamma_1}\tau_1}\times \cdots \times {^{\gamma_s}\tau_s}$ pour des ŽlŽments $\gamma_i \in \Gamma(E/F)$. 
\end{enumerate}
\end{proposition}

\subsection{RŽsultats globaux.} Soit $\bs{E}/\bs{F}$ une extension finie cyclique de corps de nombres, de degrŽ $d$. On note 
$X(\E/\F)$ le groupe des caractres de $\A^\times$ qui sont triviaux sur 
$ \F^\times \mathrm{N}_{\E/\F}(\AE^\times)$ et on fixe un gŽnŽrateur $\K$ de $X(\E/\F)$.  

On s'intŽresse aux repr\'esentations automorphes discrtes\footnote{Cf. \ref{rep aut disc}.} de $\GLm(\AE)$ pour un entier $m\geq 1$, et ˆ certaines reprŽsentations automorphes de $\GLn(\A)$ pour un entier $n\geq 1$, appelŽes \guill{induites de discrtes}\footnote{Nous reprenons ici la notion de reprŽsentation automorphe \guill{induite de cuspidale unitaire} 
de \cite[chap.~3, def.~4.1]{AC} et \cite{H2}.}: 

\begin{definition}\label{induite de rep disc} \textup{Soit $n\geq 1$ un entier. Une repr\'esentation automorphe $\pi$ de $\GLn(\A)$ est dite \textit{induite de discrte} si elle est isomorphe ˆ $\delta_1\times \cdots \times \delta_r$ pour des entiers $n_i\geq 1$ tels que $\sum_{i=1}^rn_i = n$ et des reprŽsentations automorphes discrtes $\delta_i$ de $\mathrm{GL}_{n_i}(\A)$; o $\delta_1 \times \cdots \times \delta_r$ dŽsigne l'induite parabolique (au sens des reprŽsentations automorphes, cf. \cite{La2}) de la reprŽsentation automorphe discrte $\delta_1\otimes \cdots \otimes \delta_r$ de $M(\A)= \mathrm{GL}_{n_1}(\A)\times\cdots\times \mathrm{GL}_{n_r}(\A)$ suivant 
le sous-groupe parabolique standard $P(\A)$ de $\GLn(\A)$ de composante de Levi $M(\A)$.
}\end{definition}

Observons qu'une reprŽsentation automorphe induite de discrte est par dŽfinition 
irrŽductible\footnote{Car en toute place $v\in \V$, la reprŽsentation $\delta_{1,v}\otimes \cdots \otimes \delta_{r,v}$ de $M(\F_v)$ Žtant unitaire, son induite parabolique 
$\delta_{1,v}\times \cdots \times \delta_{r,v}$ est irrŽductible (et unitaire).} et unitaire. Les reprŽsentations automorphes induites de discrtes 
vŽrifient le \guill{phŽnomne de rigiditŽ} (cf. \cite[2.2, 2.3]{H2}): 

\begin{lemme}\label{rigiditŽ}
Soit un entier $n\geq 1$. Soient $\pi$ et $\pi'$ deux reprŽsentations automorphes de $\GLn(\A)$ induites de discrtes, telles qu'en presque toute place finie $v$ de $\F$, on ait $\pi_v\simeq \pi'_v$. Alors les reprŽsentations $\pi$ et $\pi'$ sont isomorphes.
\end{lemme}

Nous aurons ˆ considŽrer principalement 
les reprŽsentations automorphes de $\GLmd(\A)$ pour un entier $m\geq 1$ de la forme 
$$\pi_{(r,\delta,\K)}=\delta \times \K \delta\times \cdots \times \K^{r-1}\delta$$ o $r\geq 1$ est un entier divisant $d$ et $\delta$ est une reprŽsentation automorphe discrte 
$\K^r$-stable de $\mathrm{GL}_{n'}(\A)$, $n'=\frac{md}{r}$. Une telle repr\'esentation $\pi_{(r,\delta,\K)}$ est induite de discrte et $\K$-stable. 

\begin{definition}\label{def rel faible aut}
\textup{
Soit un entier $m\geq 1$. Soit $\Pi$ une repr\'esentation automorphe de $\mathrm{GL}_{m}(\AE)$ induite de discrte et soit $\pi$ une reprŽsentation automorphe $\K$-stable de $\GLmd(\A)$ induite de discrte. 
\begin{enumerate}
\item[(i)]On dit que $\pi$ est un $\K$-relvement faible de $\Pi$ si en \textit{presque toute} place finie $v\in \V$, $\pi_v$ est un $\K_v$-relvement faible de $\Pi_v$.
\item[(ii)]On dit que $\pi$ est un $\K$-relvement (fort) de $\Pi$ si en \textit{toute} place $v\in \V$, $\pi_v$ est un $\K_v$-relvement de $\Pi_v$. 
\end{enumerate}
}
\end{definition}

\begin{remark}\label{unicitŽ relvement global}
\textup{
\begin{enumerate}
\item[(i)] Un $\K$-relvement est \textit{a fortiori} un $\K$-relvement faible, car en les places $v\in \Vfin$ o\`u les reprŽsentations sont sphŽriques, 
un $\K_v$-relvement est \textit{a fortiori} un $\K_v$-relvement faible. 
\item[(ii)] Un $\K$-relvement faible --- et donc \textit{a fortiori} un $\K$-relvement (fort) ---, s'il existe, est unique ˆ isomorphisme prs (\ref{rigiditŽ}). 
\end{enumerate}
}
\end{remark}

Si $\delta$ est une reprŽsentation automorphe discrte de $\mathrm{GL}_{n}(\A)$ pour un entier $n\geq 1$, on dŽfinit comme dans le cas local l'ensemble $X(\delta)$ et son cardinal $x(\delta)$. Si\footnote{Cf. \ref{rappels sur le spectre res}.} $\delta = u(\rho,q)$ pour un entier $q\geq 1$ divisant $n$ et une reprŽsentation automorphe cuspidale unitaire $\rho$ de $\mathrm{GL}_{\smash{\frac{n}{q}}}(\A)$, on a $x(\delta)=x(\rho)$. 

Si $\Pi$ est une reprŽsentation automorphe discrte de $\mathrm{GL}_{m}(\AE)$ pour un entier $m\geq 1$, on dŽfinit comme dans le cas 
local l'ensemble $\Gamma(\Pi)$ et son cardinal $g(\Pi)$; et on pose $r(\Pi)= \frac{d}{g(\Pi)}$. Si $\Pi=u(\Lambda,q)$ pour un entier $q\geq 1$ divisant $m$ et une reprŽsentation automorphe cuspidale unitaire $\Lambda$ de $\mathrm{GL}_{\smash{\frac{m}{q}}}(\AE)$, on a $r(\Pi)=r(\Lambda)$. 

\begin{theorem}\label{theo global}
\begin{enumerate}
\item[(i)] Soit un entier $m'\geq 1$ et soit $\Pi$ une repr\'esentation automorphe discrte de $\mathrm{GL}_{m'}(\AE)$. Posons $r=r(\Pi)$. 
Il existe une reprŽsentation automorphe discrte $\K^r$-stable $\delta$ de $\mathrm{GL}_{n'}(\A)$, $n'= \frac{m'\!d}{r}$, telle que $\pi_{(r,\delta,\K)}$ soit un $\K$-relvement (fort) de $\Pi$. On a $x(\delta)=r$ et l'ensemble $X(\delta)$ est dŽterminŽ de manire unique par la classe d'isomorphisme de $\Pi$. De plus $\Pi$ est cuspidale si et seulement si $\delta$ est cuspidale.
\item[(ii)] Soit un entier $r\geq 1$ divisant $d$ et soit $\delta$ une reprŽsentation automorphe discrte de $\mathrm{GL}_{n'}(\A)$ pour un entier $n'\geq 1$, telle que $x(\delta)=r$. Posons $\pi= \pi_{(r,\delta,\K)}$. Alors $d$ divise $rn'$ et il existe une reprŽsentation automorphe discrte $\Pi$ de $\mathrm{GL}_{m'}(\AE)$, $m'= \frac{rn'}{d}$, telle que $\pi$ soit un $\K$-relvement de $\Pi$. On a $r(\Pi)=r$ et les classes d'isomorphisme de reprŽsentations automorphes discrtes de $\mathrm{GL}_{m'}(\AE)$ ayant $\pi$ pour $\K$-relvement sont les ŽlŽments de $\Gamma(\Pi)$. 
\end{enumerate}
\end{theorem}

\begin{remark}
\textup{Soit $\Pi$ une reprŽsentation automorphe discrte de $\GLm(\AE)$ pour un entier $m\geq 1$. 
D'aprs \ref{theo global} et \ref{unicitŽ relvement global}, on a :
\begin{enumerate} 
\item[(i)]Ë isomorphisme prs il existe une unique reprŽsentation automorphe 
$\pi$ de $\GLmd(\A)$ induite de discrte qui soit un $\K$-relvement de $\Pi$. 
Observons que $\pi$ est discrte si et seulement si $r(\Pi)=1$,  
\ie si et seulement si $\Pi$ est $\Gamma(\E/\F)$-rŽgulire. 
\item[(ii)]S'il existe une reprŽsentation $\pi$ de $\GLmd(\A)$ induite de discrte qui soit un $\K$-relvement faible de $\Pi$, alors c'est 
un $\K$-relvement (fort).
\end{enumerate}}
\end{remark}

\begin{corollary}\label{cor global}
Soit un entier $m\geq 1$. 
\begin{enumerate}
\item[(i)] Pour toute reprŽsentation automorphe de $\mathrm{GL}_{m}(\AE)$ induite de discrte, il existe une reprŽsentation automorphe 
$\pi$ de $\GLmd(\A)$ induite de discrte qui soit un $\K$-relvement de $\pi$. PrŽcisŽment, si $\tau=\tau_1\times \cdots \times \tau_s$ o $\tau_i$ est une reprŽsentation automorphe discrte de $\mathrm{GL}_{m_i}(\AE)$, 
$\sum_{i=1}^s m_i=m$, et si pour $i=1,\ldots , s$, $\pi_i$ un $\K$-relvement de $\tau_i$ ˆ $\mathrm{GL}_{m_id}(\AE)$, alors $\pi= \pi_1\times \cdots \times \pi_r$ est un $\K$-relvement de $\tau$. 
\item[(ii)] Toute reprŽsentation automorphe de $\GLmd(\A)$ induite de discrte et $\K$-stable 
est un $\K$-relvement d'un reprŽsentation automorphe de $\GLm(\AE)$ induite de discrte. 
\item[(iii)] Si $\tau= \tau_1\times \cdots \times \tau_r$ pour des reprŽsentations automorphes discrtes $\tau_i$ de $\mathrm{GL}_{m_i}(\AE)$, $\sum_{i=1}^r m_i = m$, 
et si $\tau'$ est une reprŽsentation automorphe $\GLm(\AE)$ induite de discrte ayant mme $\K$-relvement que $\tau$, alors 
$\tau'$ est isomorphe ˆ ${^{\gamma_1}\tau_1}\times\cdots \times {^{\gamma_r}\tau_r}$ pour des ŽlŽments $\gamma_i \in \Gamma(\E/\F)$. 
\end{enumerate}
\end{corollary}

\subsection{Sur la preuve des thŽormes \ref{theo local} et \ref{theo global}.}\label{principe 1} La dŽmonstration des thŽormes \ref{theo local} (local) et \ref{theo global} (global) repose sur le principe local-global; en particulier la mŽthode conduit ˆ dŽmontrer simultanŽment les deux thŽormes. On procde par rŽcurrence sur l'entier $m'\geq 1$ de \ref{theo global}: on fixe un entier $m\geq 1$ et on suppose que le thŽorme \ref{theo global} est vrai pour tout $m'< m$. On veut dŽmontrer qu'il est vrai pour $m$, tout en dŽmontrant le thŽorme \ref{theo local}. La rŽcurrence dŽmarre bien car pour $m=1$, toute reprŽsentation automorphe discrte de $\mathrm{GL}_1(\AE)= \AE^\times$ est cuspidale (le cas cuspidal est connu d'aprs Henniart-Herb \cite{HH}). 

\vskip1mm
{\bf $\bullet$  Formule de comparaison.} 
Posons $G= \mathrm{GL}_{n/\F}\;(=\GLn\times_{\mathbb{Z}} \F)$ et $G'= \mathrm{Res}_{\E/\F}(\mathrm{GL}_{m/\E})$.   
L'ingrŽdient principal est la formule de comparaison \ref{dec endos}\,(4) 
$$I_{\mathrm{disc},\xi}^G(f,\K\circ \det)=\frac{1}{d} I_{\mathrm{disc},\mu}^{G'}(\phi)\leqno{(1)}$$ pour toute paire de fonctions $(f,\phi)$ qui sont associŽes au sens de \ref{fonctions associŽes}; ici $\mu$ est un caractre automorphe unitaire de $\AF^\times$ (identifiŽ ˆ $A_{G'}(\F)= A_G(\F)$) et $\xi= \K^{md(d-1)/2}\mu$. L'identitŽ (1) n'est autre que la dŽcomposition endoscopique de la partie discrte de la formule des traces pour $(\mathrm{GL}_{n/\F},\K\circ\det)$: 
d'aprs \ref{unicitŽ de la donnŽe elliptique}, $G'$ est ˆ isomorphisme prs l'unique groupe endoscopique elliptique de $(G,\K\circ\det)$. Partant de la dŽcomposition fine des expressions ˆ gauche et ˆ droite de l'ŽgalitŽ (1), il est assez facile de dŽcrire les reprŽsentations automorphes irrŽductibles pouvant y donner une contribution non triviale. On obtient 
l'ŽgalitŽ : 
$$\sum_{k\vert n} \frac{1}{k^2} \sum_{\delta} c_{(k,\delta,\K)} \mathrm{tr}(\pi_{(k,\delta,r)}(f)\circ I_{\K})=\frac{1}{d}\sum_{l\vert m}\frac{1}{l^2}\sum_{\Delta} c_{(l,\Delta)} \mathrm{tr}(\Delta^{\!\times l}(\phi))\leqno{(2)}$$ 
o:
\begin{itemize}
\item $\delta$ parcourt les classes d'isomorphisme de reprŽsentations automorphes discrtes $\K^k$-stables de $\mathrm{GL}_{\frac{n}{k}}(\A)$ de caractre central $\omega_\delta:\A^\times \rightarrow \mathbb{C}^\times$ vŽrifiant $(\omega_\delta)^k \K^{n(k-1)/2}=\xi$; 
\item $\Delta$ parcourt les classes d'isomorphisme de reprŽsentations automorphes discrtes de $\mathrm{GL}_{\frac{m}{l}}(\AE)$, de caractre central $\omega_\Delta: \AE^\times \rightarrow \mathbb{C}^\times$ vŽrifiant $(\omega_\Delta)^l\vert_{\A^\times}= \mu$;
\item $I_{\K}$ est l'opŽrateur \guill{physique} donnŽ par $I_\K(\varphi) (g)= \K(g)\varphi(g)$ pour toute fonction $\varphi\in L^2(\GLn(\F)\backslash \GLn(\A), \xi)$ et tout 
$g\in \GLn(\A)$; 
\item $c_{(k,\delta,\K)}$ et $c_{(l,\Delta)}$ sont des constantes complexes non nulles.
\end{itemize}  

\vskip1mm
{\bf $\bullet$ La \guill{mŽthode standard}.}
FixŽs un entier $l\geq 1$ divisant $m$ et une reprŽsentation automorphe discrte $\Delta$ de $\mathrm{GL}_{\frac{m}{l}}(\AE)$, la 
\textit{mŽthode standard} (cf. Flath \cite{Fl2}) et le lemme local-global (\ref{lemme local-global}) 
permettent de sŽparer ˆ droite de l'ŽgalitŽ (2) les reprŽsentations automorphes de $\GLm(\AE)$ ayant mme $\K$-relvement faible que $\Pi=\Delta^{\!\times l}$ en dehors d'un sous-ensemble fini $S\subset \V$ (contenant $\Vinf$ toutes les places $v\in \Vfin$ o la reprŽsentation $\Delta_v$ n'est pas sphŽrique). Pour $l=1$ et $\Pi=\Delta$, on obtient ainsi l'existence d'un $\K$-relvement faible $\pi$ de $\Pi$ de la forme $\pi=\pi_{(k,\delta,\K)}$ pour un entier $k\geq 1$ divisant $n$ et une reprŽsentation automorphe discrte $\K^k$-stable $\delta$ de $\mathrm{GL}_{\frac{n}{k}}(\A)$. Plus prŽcisŽment, on prouve (gr‰ce ˆ la description explicite de $\Pi$ donnŽe par \cite{MW1}) que $r(\Pi)=k$. De plus en toute place $v\in S$, $\pi_v$ est un $\K_v$-relvement de $\Pi_v$; 
et (par construction) en toute place $v\in \V\smallsetminus S$, $\pi_v$ est un $\K_v$-relvement faible de $\Pi_v$. \`A ce stade de la dŽmonstration, 
on utilise l'hypothse de rŽcurrence pour obtenir l'ŽgalitŽ $x(\delta)=k$. 

\vskip1mm
{\bf $\bullet$ $\kappa$-relvement des reprŽsentations de Speh.} Nous ne dŽcrivons ici que le cas local non archimŽdien, le cas local archimŽdien Žtant dŽjˆ connu d'aprs \cite{LL} (cf. \ref{preuve theo local i}). L'Žtape suivante consiste ˆ rŽaliser une reprŽsentation de Speh (locale) $\Pi_0= u(\Lambda_0,q)$ de $\GLm(E)$, pour un entier $q\geq 1$ divisant $m$ et une reprŽsentation irrŽductible de carrŽ intŽgrable $\Lambda_0$ de $\mathrm{GL}_{\frac{m}{q}}(E)$, comme une composante locale $\Pi_{v_0}$ d'une reprŽsentation automorphe discrte $\Pi$ de $\GLm(\AE)$. Pour cela on choisit l'extension $\E/\F$ (telle que $\E_{v_0}/\F_{v_0}\simeq E/F$) et la reprŽsentation $\Pi$ (telle que $\Pi_{v_0}\simeq \Pi_0$) de manire la plus favorable possible; en particulier on s'arrange pour qu'en une place finie $v_1$ de $\F$ distincte de $v_0$, inerte et non ramifiŽe dans $\E$, la composante locale $\Pi_{v_1}$ soit $\Gamma(\E/\F)$-rŽgulire. La reprŽsentation $\Pi$ est alors \textit{a fortiori} $\Gamma(\E/\F)$-rŽgulire, \ie $r(\Pi)=1$. Par la mŽthode standard dŽcrite dans la paragraphe prŽcŽdent, on obtient l'existence d'un $\kappa$-relvement pour la reprŽsentation de Speh $\Pi_0= u(\Lambda_0,q)$), \ie \ref{theo local}\,(i). Le point (ii) de \ref{theo local} est une consŽquence de la description (dŽjˆ connue) de l'image et des fibres de l'application de $\kappa$-relvement pour les reprŽsentations de carrŽ intŽgrable. 

\vskip1mm
{\bf $\bullet$ Suite et fin de la dŽmonstration.} Gr‰ce ˆ la description du dual unitaire de $\GLm(E)$, on en dŽduit, par compatibilitŽ entre les applications de $\kappa$-relvement et d'induction parabolique locales, l'existence d'un $\kappa$-relvement pour toutes les reprŽsentations irrŽductibles unitaires de $\GLm(E)$; \ie \ref{prop unit}\,(i). Les mme arguments permettent de prouver \ref{prop unit}\,(ii) (surjectivitŽ) et \ref{prop unit}\,(iii) (description des fibres). 
Si l'extension de corps locaux non archimŽdiens $E/F$ est non ramifiŽe, cela prouve que le $\kappa$-relvement faible d'une reprŽsentation irrŽductible unitaire sphŽrique de $\GLm(E)$ est en fait un $\kappa$-relvement. En particulier le $\K$-relvement faible $\pi$ de la reprŽsentation automorphe discrte $\Pi$ de $\GLm(\AE)$ construit plus haut, est un $\K$-relvement (fort); on a donc prouvŽ \ref{theo global}\,(i). Quant ˆ \ref{theo global}\,(ii), la \textit{mŽthode standard} permet aussi de dŽcrire l'image de l'application de $\K$-relvement pour les reprŽsentations automorphes discrtes de $\GLm(\AE)$; et la description des fibres rŽsulte du lemme local-global \ref{lemme local-global}. L'hypothse de rŽcurrence (pour $m'<m$) est donc vraie pour $m$, ce qui achve la dŽmonstration. 

\section{$\kappa$-relvement des reprŽsentations elliptiques}\label{kappa-rel elliptiques}
Dans cette section, $E/F$ est une extension finie cyclique de corps locaux, de degrŽ $d$. 
On fixe un entier $m\geq 1$ et on pose $n= md$. On note $G$ le groupe $\GLn(F)$ et $H$ le groupe $\GLm(E)$. 
On dŽmontre la proposition \ref{prop ell}. PrŽcisŽment on complte les arguments de \cite{Fa} dans le cas non archimŽdien   
et on rappelle brivement les rŽsultats de \cite{LL,KS} (cf. \cite{H1}) dans le cas archimŽdien.

\subsection{Le cas non archimŽdien.} Rappelons la description des reprŽsentations irrŽductibles elliptiques de $H$ (cf. \ref{classif (rappels)}). Pour un entier $k\geq 1$ divisant $m$, on note $P'_k$ le sous-groupe parabolique standard de $H$ de composante de Levi standard $L'_k=\mathrm{GL}_{\frac{m}{k}}(E)\times \cdots \times \mathrm{GL}_{\frac{m}{k}}(E)$ ($k$ fois), et $\ES{L}_{\mathrm{st}}(L'_k)$ l'ensemble des facteurs de Levi standard de $H$ contenant $L'_k$. Soit $\rho_E$ une reprŽsentation irrŽductible cuspidale de $\mathrm{GL}_{\frac{m}{k}}(E)$. Pour $M'\in \ES{L}_{\mathrm{st}}(L'_k)$, on forme l'induite parabolique 
$$\wt{R}^{M'\!}(\rho_E,k)= i_{P'_k \cap M'}^{M'}(\nu_E^{k-1}\rho_E \otimes \nu_E^{k-2}\rho_E \otimes \cdots \otimes \rho_E)$$ suivant le sous-groupe 
parabolique standard $P_k\cap M'$ de $M'$. Cette reprŽsentation $\wt{R}^{M'\!}(\rho_E,k)$ de $M'$ a une unique sous-reprŽsentation irrŽductible, notŽe $\tilde{\delta}^{M'\!}(\rho_E,k)$, qui est essentiellement de carrŽ intŽgrable. On forme ensuite l'induite parabolique $$\wt{R}(\rho_E,k;M')= i_{P'}^H(\tilde{\delta}^{M'\!}(\rho_E,k))$$ suivant le sous-groupe parabolique standard $P'$ de $H$ de composante de Levi $M'$, et on note 
$\wt{u}(\rho_E,k;M')$ l'unique quotient irrŽductible de $\wt{R}(\rho_E,k;M')$. L'application $$M' \mapsto \wt{u}(\rho_E,k;M')$$ est une bijection de $\ES{L}_{\mathrm{st}}(L'_k)$ sur l'ensemble des sous-quotients irrŽductibles de $\wt{R}(\rho_E,k)= \wt{R}^H(\rho_E,k)$, chacun apparaissant avec multiplicitŽ $1$. Pour $M'=H$, la reprŽsentation 
$$\wt{u}(\rho_E,k;H)= \wt{R}(\rho_E,k;H)= \tilde{\delta}^H(\rho_E,k)=\tilde{\delta}(\rho_E,k)$$ est l'unique sous-quotient irrŽductible de $\wt{R}(\rho_E,k)$ qui soit essentiellement de carrŽ intŽgrable. D'aprs la propriŽtŽ de multiplicitŽ $1$, pour tout sous-quotient irrŽductible 
$\tau$ de $\wt{R}(\rho_E,k)$, on a $r(\tau)=r(\rho_E)$. 

L'entier $r= r(\rho_E)$ divise $\frac{n}{k}$ et on sait d'aprs \cite{HH} que la reprŽsentation $\rho_E$ de $\mathrm{GL}_{\frac{m}{k}}(E)$ admet un $\kappa$-relvement ˆ $\mathrm{GL}_{\frac{n}{k}}(F)$ de la forme $\rho \times \kappa\rho \times\cdots \times \kappa^{r-1}\rho$ pour une reprŽsentation irrŽductible cuspidale $\kappa^r$-stable $\rho$ de $\mathrm{GL}_{\frac{n}{kr}}(F)$; de plus on a $x(\rho)=r$. Alors $\wt{\delta}(\rho,k)$ est une reprŽsentation irrŽductible essentiellement de carrŽ intŽgrable $\kappa^r$-stable de $\mathrm{GL}_{\frac{n}{r}}(F)$ telle que $x(\tilde{\delta}(\rho,k))=r$, et 
$$\tilde{\delta}(\rho,k)\times \kappa\tilde{\delta}(\rho,k) \times\cdots \times \kappa^{r-1}\tilde{\delta}(\rho,k)$$ est un 
$\kappa$-relvement de $\tilde{\delta}(\rho_E,k)$. 

Posons $g= \frac{d}{r}$. Pour $M' \in \ES{L}_{\mathrm{st}}(L_k)$, on dŽfinit comme suit un facteur de Levi standard $\iota_g(M')$ de $\mathrm{GL}_{gm}(F)$, $gm=\frac{n}{r}$. On Žcrit $M'=\mathrm{GL}_{m_1}(E)\times \cdots \times \mathrm{GL}_{m_s}(E)$ avec $k \vert m_i$ et $\sum_{i=1}^s m_i = m$, et on pose 
$\iota_g(M')=\mathrm{GL}_{gm_1}(F)\times \cdots \times \mathrm{GL}_{gm_s}(F)$.  
Observons que $\iota_g(M')$ contient $\iota_g(L'_k)= \mathrm{GL}_{\frac{gm}{k}}(F)\times \cdots \times \mathrm{GL}_{\frac{gm}{k}}(F)$ ($k$ fois). L'application 
$M'\mapsto \iota_g(M')$ est une bijection de $\ES{L}_{\mathrm{st}}(L'_k)$ sur l'ensemble $\ES{L}_{\mathrm{st}}(\iota_g(L'_k))$ des facteurs de Levi standard 
de $\mathrm{GL}_{gm}(F)$ qui contiennent $\iota_g(L'_k)$. On dŽfinit comme en \ref{classif (rappels)} la reprŽsentation irrŽductible essentiellement de carrŽ intŽgrable $\tilde{\delta}^{\iota_g(M')}(\rho,k)$ de $\iota_g(M')$, puis la reprŽsentation irrŽductible elliptique $\wt{u}(\rho, k; \iota_g(M'))$ de $\mathrm{GL}_{gm}(F)$. L'application $M' \mapsto \wt{u}(\rho,k;\iota_g(M'))$ est une bijection de $\ES{L}_{\mathrm{st}}(L'_k)$ sur l'ensemble des sous-quotients irrŽductibles de la reprŽsentation $\wt{R}(\rho,k)= \nu^{k-1}\rho\times \nu^{k-2}\rho\times\cdots \times \rho$ de $\mathrm{GL}_{gm}(F)$; ou, ce qui revient au mme, sur l'ensemble des classes d'isomorphisme de reprŽsentations irrŽductibles elliptiques de $\mathrm{GL}_{gm}(F)$ de suppport cuspidal le segment $[\rho,\nu^{k-1}\rho]$. Observons que (d'aprs la propriŽtŽ de multiplicitŽ $1$) tout sous-quotient irrŽductible $\pi$ de $\wt{R}(\rho,k)$ est $\kappa^r$-stable et vŽrifie $x(\pi)=r$. Pour $M'\in \ES{L}_{\mathrm{st}}(L'_k)$ et 
$M=\iota_g(M')$, posons 
$$\tilde{\delta}_{E,M'}= \wt{u}(\rho_E,k;M')\qhq{et} \tilde{\delta}_{M}= \wt{u}(\rho,k;M)\ptf$$
D'aprs \cite{Fa}, la reprŽsentation $$\pi_{M}= \tilde{\delta}_{M}\times \kappa \tilde{\delta}_{M}\times\cdots \times \kappa^{r-1}\tilde{\delta}_{M}$$  de $G=\mathrm{GL}_n(F)$ est un $\kappa$-relvement de $\tilde{\delta}_{E,M'}$. 
La reprŽsentation (irrŽductible elliptique) $\tilde{\delta}_{M}$ de $\mathrm{GL}_{\frac{n}{r}}(F)$ est $\kappa^r$-stable et vŽrifie $x(\tilde{\delta}_{M})=r=r(\tilde{\delta}_{E,M'})$. 
On voit que les conditions suivant sont Žquivalentes: 
\begin {itemize}
\item $\tilde{\delta}_{E,M'}$ est essentiellement de carrŽ intŽgrable, \ie $\tilde{\delta}_{E,M'} = \tilde{\delta}(\rho_E,k)$;
\item $M'=H$, \ie $M= \mathrm{GL}_{gm}(F)$; 
\item $\tilde{\delta}_{M}$ est essentiellement de carrŽ intŽgrable, \ie $\tilde{\delta}_{M}= \tilde{\delta}(\rho,k)$.
\end{itemize} 
Cela prouve \ref{prop ell}\,(i).

RŽciproquement, soit $r\geq 1$ un entier divisant $n$ et $\tilde{\delta}$ une reprŽsentation irrŽductible elliptique $\kappa^r$-stable de $\mathrm{GL}_{n'}(F)$, $n'=\frac{n}{r}$,  telle que 
$x(\tilde{\delta})=r$. Posons $\pi=\tilde{\delta}\times \kappa \tilde{\delta} \times \cdots \times \kappa^{r-1}\tilde{\delta}$; c'est une reprŽsentation (irrŽductible $\kappa$-stable)  de $G=\GLn(F)$. On veut prouver que $\pi$ est un $\kappa$-relvement d'une reprŽsentation irrŽductible elliptique $\tilde{\delta}_E$ de $\GLm(E)$. 
D'aprs la description  des reprŽsentations irrŽductibles elliptiques de $G$ (cf. \ref{classif (rappels)}), il existe un entier $k\geq 1$ divisant $n'$, une reprŽsentation irrŽductible cuspidale $\rho$ de $\mathrm{GL}_{n''}(F)$, $n''=\frac{n'}{k}=\frac{n}{kr}$, et un facteur de Levi standard $M$ de $\mathrm{GL}_{n'}(F)$ contenant $\mathrm{GL}_{n''}(F)\times\cdots \times \mathrm{GL}_{n''}(F)$ ($k$ fois), tels que $\tilde{\delta}\simeq \wt{u}(\rho,k;M)$. On a $x(\rho)=r$, ce qui assure que $g=\frac{d}{r}$ divise $n''$ (observons que $k \frac{n''}{g}=\frac{n}{d}=m$). On en dŽduit que: 
\begin{itemize}
\item la reprŽsentation $\rho\times \kappa \rho \times \cdots \times \kappa^{r-1}\rho$ de 
$\mathrm{GL}_{n''r}(F)$ est un $\kappa$-relvement d'une reprŽsentation irrŽductible cuspidale $\rho_E$ de $\mathrm{GL}_{\frac{n''}{g}}(E)$; 
\item $M= \iota_g(M')$ pour un (unique) facteur de Levi standard $M'\in \ES{L}_{\mathrm{st}}(L'_k)$; 
\item $\pi$ est un $\kappa$-relvement de $\tilde{\delta}_{E,M'}= \wt{u}(\rho_E,k;M')$ de $\GLm(E)$.
\end{itemize}
Observons que la classe d'isomorphisme de $\pi$ dŽtermine les entiers $k$ et $r$, l'ensemble $\Gamma(\rho_E)$ --- ou, ce qui revient au mme, l'ensemble $X(\rho)$ --- et le facteur de Levi $M'\in \ES{L}_{\mathrm{st}}(L'_k)$. Quant ˆ la description des fibres, elle rŽsulte de cette dernire observation pourvu que toute reprŽsentation irrŽductible $\tau$ de $H$ ayant $\pi$ pour $\kappa$-relvement, soit elliptique. En effet, une telle reprŽsentation $\tau$ est alors isomorphe ˆ $\wt{u}({^\gamma(\rho_E)},k,M')\simeq {^\gamma(\wt{u}(\rho_E,k,M'))}$ pour un ŽlŽment $\gamma\in \Gamma(E/F)$. 
Reste ˆ prouver que $\tau$ est elliptique. Il suffit pour cela, d'aprs la formule de caractres dŽfinissant la notion de $\kappa$-relvement (\ref{def kappa-rel}\,(2) ou \ref{rem fonc-car kappa-rel}\,(ii) pour les ŽlŽments elliptiques), que $\pi$ soit \textit{$\kappa$-elliptique} au sens suivant. 

\begin{definition}
\textup{
Une reprŽsentation irrŽductible $\kappa$-stable $\pi$ de $G$ est dite \textit{$\kappa$-elliptique} si pour un (\ie pour tout) 
$A\in \mathrm{Isom}_G(\kappa\pi,\pi)$, la fonction-caractre $\Theta_\pi^A$ n'est pas identiquement nulle sur l'ouvert de $G$ formŽ des ŽlŽments semi-simples rŽguliers elliptiques. }
\end{definition}

\begin{lemme}
Soient $r\geq 1$ un entier divisant $n$ et  $\tilde{\delta}$ une reprŽsentation irrŽductible elliptique $\kappa^r$-stable 
de $\mathrm{GL}_{n'}(F)$, $n'= \frac{n}{r}$, telle que $x(\tilde{\delta})=r$. La reprŽsentation 
$\pi= \tilde{\delta}\times \kappa\tilde{\delta}\times \cdots \times \kappa^{r-1}\tilde{\delta}$ de $G$ est $\kappa$-elliptique. 
\end{lemme}

\begin{demo}
Pour $\tilde{\delta}$ essentiellement de carrŽ intŽgrable, le rŽsultat est connu: en effet dans ce cas $\pi$ est essentiellement $\kappa$-discrte au sens de \cite[5.3]{HH}, donc $\kappa$-elliptique \cite[6.2, cor.~1]{HH}. En gŽnŽral, 
Žcrivons $\tilde{\delta}= \tilde{\delta}(\rho,k; L)$ pour un entier $k\geq 1$ divisant $n'$, une reprŽsentation irrŽductible cuspidale $\rho$ de $G_{n''}=\mathrm{GL}_{n''}(F)$, $n''= \frac{n'}{k}$, et un facteur de Levi standard $L$ de $G_{n'}=\mathrm{GL}_{n'}(F)$ contenant $L_k^{n'}=G_{n''}\times\cdots\times G_{n''}$ ($k$ fois). 
Soit $\ES{G}(G_{n'})$ le groupe de Grothendieck des reprŽsentations de longeur finie de $G_{n'}$. Les semi-simplifiŽes 
des reprŽsentations standard de $G_{n'}$ forment une $\mathbb{Z}$-base de $\ES{G}(G_{n'})$. On peut dŽcomposer $\tilde{\delta}$ suivant cette 
base; seules les reprŽsentations standard de $G_{n'}$ ayant mme support cuspidal (multi-segment au sens de \cite{Z}) que $\tilde{\delta}$ 
interviennent dans cette decomposition. 
Ces reprŽsentations sont prŽcisŽment les $\wt{R}(\rho,k; M_1)$ o $M_1$ parcourt l'ensemble $\ES{L}_{\mathrm{st}}(L_k^{n'})$ 
des facteurs de Levi standard de $G_{n'}$ contenant $L_k^{n'}$: $$\tilde{\delta}= \sum_{M_1} a_{M_1}(\tilde{\delta}) \wt{R}(\rho,k;M_1)\qhq{avec} a_{M_1}(\tilde{\delta})\in \mathbb{Z}\pvg$$ o (par abus de notation) on note de la mme manire une reprŽsentation de longeur finie et sa semi-simplifiŽe. 
Rappelons que $\wt{R}(\rho,k;M_1)= \iota_{P_1}^{G_{n'}}(\tilde{\delta}^{M_1}(\rho,k))$ o $P_1$ est le sous-groupe parabolique standard de $G_{n'}$ de composante de Levi $M_1$. Puisque la reprŽsentation $ \wt{R}(\rho,k;G_{n'})= \tilde{\delta}^{G_{n'}}(\rho,k)$ est elliptique et que pour $M_1\neq G_{n'}$, $\wt{R}(\rho,k;M_1)$ est une induite parabolique stricte, on a $a_{G_{n'}}(\tilde{\delta})\neq 0$. 
Pour $M_1\in \ES{L}_{\mathrm{st}}(L_k^{n'})$, posons 
$$X_{M_1}=\wt{R}(\rho,k;M_1)\times \kappa \wt{R}(\rho,k;M_1)\times \cdots \times \kappa^{r-1}\wt{R}(\rho,k;M_1)\ptf$$ C'est une reprŽsentation $\kappa$-stable de 
$G=G_n$.  \'Ecrivons $M_1= G_{n_1}\times \cdots \times G_{n_s}$ avec $k\vert n_i$, $\sum_{i=1}^sn_i=n'$, et notons $M=\iota_r(M_1)$ le facteur de Levi standard 
de $G$ dŽfini par $M= G_{rn_1} \times \cdots \times G_{rn_s}$. Pour $i=1,\ldots ,s$, notons $\pi_i$ la reprŽsentation (irrŽductible, gŽnŽrique) $\kappa$-stable de $G_{rn_i}$ dŽfinie 
par $$\pi_i= \tilde{\delta}^{M_1}(\rho,k)\times \kappa \tilde{\delta}^{M_1}(\rho,k)\times\cdots \times \kappa^{r-1}\tilde{\delta}^{M_1}(\rho,k)\ptf$$ 
La reprŽsentation $\pi^M= \pi_1\otimes \cdots \otimes \pi_s$ de $M$ est  irrŽductible, gŽnŽrique et $\kappa$-stable. Notons $Y_{M_1}$ 
l'induite parabolique $ \pi_1\times \cdots \times \pi_s= \iota_{P_M}^G(\pi^M)$, $P=P_M$; elle est isomorphe ˆ $X_{M_1}$. On a la dŽcomposition 
dans le groupe de Grothendieck $\ES{G}(G)$ des reprŽsentations de longueur finie de $G$: 
$$\pi= \sum_{M_1} a_{M_1}(\tilde{\delta}) X_{M_1} =\sum_{M_1} a_{M_1}(\tilde{\delta}) Y_{M_1} \ptf$$ 

Pour $M_1\in \ES{L}_{\mathrm{st}}(L_k^{n'})$, $M= \iota_r(M_1)$ et $P=P_M$, 
on dispose de l'opŽrateur d'entrelacement normalisŽ $B_{M_1} = A_{\pi^M}^{\mathrm{g\acute{e}n}}\in \mathrm{Isom}_M(\kappa\pi^M\!,\pi^M)$. 
On pose $$A_{M_1}=i_P^G(B_{M_1})\in \mathrm{Isom}_G(\kappa Y_{M_1}, Y_{M_1})\ptf $$ 
Cela dŽfinit une reprŽsentation $Y_{M_1}^\kappa$ de $(G,\kappa\circ\det)$ 
au sens de \cite[5.3]{HL4}: pour $g\in G$, on pose $$Y_{M_1}^\kappa(g)= Y_{M_1}(g)\circ A_{M_1}\ptf$$ La distribution $f \mapsto \mathrm{tr}(Y_{M_1}^\kappa(f))= \mathrm{tr}(Y_{M_1}(f)\circ A_{M_1})$ sur $G$ est donnŽe sur $\Greg$ par une fonction localement constante $\Theta_{Y_{M_1}^\kappa}= \Theta_{Y_{M_1}}^{\smash{A_{M_1}}}$ qui, d'aprs l'analogue tordu de la formule de descente de Van Djik \cite[7.3.9]{HL4}, est nulle sur les ŽlŽments elliptiques si $M_1\neq G_{n'}$ (\ie si $M\neq G$). De plus cette distribution 
$\mathrm{tr}(Y_{M_1}^\kappa)$ ne \guill{voit} que les sous-quotients irrŽductibles $\kappa$-stables de $Y_{M_1}$ (cf. \cite[5.9]{HL4}). Notons $\pi^\kappa$ l'image 
de la somme virtuelle $\sum_{M_1} a_{M_1}(\tilde{\delta})Y_{M_1}^\kappa$ dans le $\mathbb{C}$-espace vectoriel $\ES{G}_{\mathbb{C}}(G,\kappa\circ\det)$ dŽfini en \cite[5.7]{HL4}. C'est une combinaison linŽaire de (classes d'isomorphisme de) reprŽsentations $G$-irrŽductibles\footnote{Au sens ou la reprŽsentation sous-jacente ($\kappa$-stable) de $G$ est irrŽductible.} de $(G,\kappa\circ\det)$, de caractre-distribution $\mathrm{tr}(\pi^\kappa)= \sum_{M_1} a_{M_1}(\tilde{\delta}) \mathrm{tr}(Y_{M_1}^\kappa)$. Cette distribution $\mathrm{tr}(\pi^\kappa)$ est donnŽe sur $\Greg$ par la fonction 
$\sum_{M_1}a_{M_1}(\tilde{\delta}) \Theta_{Y_{M_1}^\kappa}$. Seul le terme d'indice $M_1=G_{n'}$, \cad la reprŽsentation irrŽductible $\kappa$-stable essentiellement de carrŽ intŽgrable 
$$Y_{G_{n'}}= \tilde{\delta}(\rho,k)\times \kappa \tilde{\delta}(\rho,k)\times\cdots\times\kappa^{r-1}\tilde{\delta}(\rho,k)\vg$$ 
contribue non trivialement ˆ la restriction de cette fonction aux ŽlŽments (semi-simples rŽguliers) elliptiques de $G$. Comme $a_{G_{n'}}(\tilde{\delta})\neq 0$, on a $\pi^\kappa\neq 0$. Par construction, $\pi^\kappa$ est la classe d'isomorphisme d'une reprŽsentation $G$-irrŽductible de $(G,\kappa\circ\det)$ qui prolonge $\pi$. 
Le lemme est dŽmontrŽ.  
\hfill$\square$ 
\end{demo}

\vskip2mm
Cela achve la preuve de \ref{prop ell}\,(ii). La proposition \ref{prop ell} est compltement dŽmontrŽe dans le cas non archimŽdien.
 
\subsection{Le cas archimŽdien.}\label{rel ell cas archimŽdien} On peut supposer que $d=2$ et $E/F= \mathbb{C}/\mathbb{R}$. On peut aussi supposer que $m=1$ sinon la proposition \ref{prop ell} est vide. Toutes les reprŽsentations de $H=\mathbb{C}^\times$, \ie les caractres $\xi=\xi_{s,k}$ pour $s\in \mathbb{C}$ et $k\in \mathbb{Z}$ (cf. \ref{rappels (g,K)-modules}), sont elliptiques (et gŽnŽriques). Le rŽsultat 
rappelŽ ci-dessous est essentiellement contenu dans \cite{LL} et explicitement dŽcrit dans \cite[5.3]{KS} (voir aussi \cite[4.2]{H1}). Reprenons les notations de \ref{rappels (g,K)-modules}. 

Si $\xi$ est stable par la conjugaison complexe, \ie si $\xi= \xi_{s,0}$ pour un $s\in \mathbb{C}$, alors 
$\chi_{\frac{s}{2}}\times \chi'_{\frac{s}{2}}= (\chi_{\frac{s}{2}}\circ \det)\otimes (1\times \kappa)$ est un $\kappa$-relvement de $\xi$ ˆ $G=\mathrm{GL}_2(\mathbb{R})$. 

Supposons maintenant que $\xi = \xi_{s,k}$ pour un  $s\in \mathbb{C}^\times$ et un $k\in \mathbb{Z}\smallsetminus \{0\}$. Quitte ˆ remplacer $\xi$ par $\bar{\xi}= \xi_{s,-k}$, on peut supposer $k>0$. \'Ecrivons $\xi= \xi_{s+k,0}\theta^k$ o $\theta$ est le caractre (unitaire) $z \mapsto \frac{z}{\vert z \vert_{\mathbb{C}}}$ de $\mathbb{C}^\times$; on a donc $\theta^k= \xi_{-k,k}$. Au caractre $\theta^k$ est associŽ dans \cite{JL} une reprŽsentation irrŽductible $\pi(\theta^k)$ 
de $G$. On a $\pi(\theta^k)=\delta(k)$, o (rappel) $\delta(k)$ est l'unique sous-reprŽsentation irrŽductible de la sŽrie principale $\chi_{\frac{k}{2}} \times \chi_{-\frac{k}{2}}$ si $k$ est impair, resp. $\chi_{\frac{k}{2}} \times \chi'_{-\frac{k}{2}}$ si $k$ est pair, de $G$. Le noyau de $\kappa\circ \det$ est le sous-groupe $G_+$ de $G$ formŽ des matrices 
de dŽterminant strictement positif. Ce groupe $G_+$ est le produit de $\mathrm{SL}_2(\mathbb{R})$ par le sous-groupe de $G$ formŽ des matrices centrales de dŽterminant strictement positif. Fixons un caractre additif non trivial $\psi$ de $\mathbb{R}$. La restriction de $\delta(k)$ ˆ $G_+$ se casse en deux reprŽsentations irrŽductibles non isomorphes, notŽes $\delta(k)^+$ et $\delta(k)^-$, o $\delta(k)^+= \delta(k,\psi)^+$ est le composant gŽnŽrique (relativement ˆ $\psi)$. 
L'opŽrateur d'entrelacement normalisŽ $A= A_{\delta(k)}^{\mathrm{g\acute{e}n}}$ est l'identitŽ sur l'espace de $\delta(k)^+$ et opre par $v\mapsto -v$ sur l'espace de $\delta(k)^-$. En notant $\Theta_{\delta(k)^+}$ et $\Theta_{\delta(k)^-}$ les fonctions-caractres de $\delta(k)^+$ et $\delta(k)^-$, qui sont des fonctions analytiques sur $G_+\cap \Greg$, localement intŽgrables sur $G$, on a l'ŽgalitŽ $\Theta_{\delta(k)}^A = \Theta_{\delta(k)^+} - \Theta_{\delta(k)^-}$. D'aprs 
\cite[5.3]{KS}, la reprŽsentation $\delta(k)$ est un $\kappa$-relvement de $\theta^k$. Revenons au caractre $\xi=\xi_{s,k}$ de $\mathbb{C}^\times$. Puisque $\xi_{s+k,0}= \chi_{\frac{1}{2}(s+k)}\circ \mathrm{N}_{\mathbb{C}/\mathbb{R}}$, on en dŽduit que la reprŽsentation $\chi_{\frac{1}{2}(s+k)}\delta(k)= (\chi_{\frac{1}{2}(s+k)}\circ\det)\otimes \delta(k)$ est un $\kappa$-relvement de $\xi$; c'est bien sžr aussi un $\kappa$-relvement de $\bar{\xi}$. 

En conclusion, tout caractre $\xi$ de $\mathbb{C}^\times$ a un $\kappa$-relvement ˆ $\mathrm{GL}_2(F)$; et si $\xi_{s,k}$ et $\xi_{s'\!,k'}$ ont mme $\kappa$-relvement ˆ $\mathrm{GL}_2(\mathbb{R})$ alors $\vert k\vert = \vert k' \vert $ et $s=s'$. Cela prouve \ref{prop ell} dans le cas archimŽdien. 

\section{Comparaison des parties discrtes des formules des traces}\label{comparaison part disc}

Dans cette section, $\F$ est un corps de nombres. En \ref{rep aut disc}--\ref{variante avec car aut}, $G$ est un groupe rŽductif connexe dŽfini sur $\F$ et 
$\omega$ est un caractre unitaire de $G(\A)$ trivial sur $G(\F)$. En \ref{l'ŽgalitŽ}--\ref{dec endos}, $\E/\F$ est une extension finie cyclique de corps de nombres, de 
degrŽ $d$, $n=md$, $G=\mathrm{GL}_{n/\F}$ et $G'= \mathrm{Res}_{\E/\F}(\mathrm{GL}_{m/\E})$; et 
$\omega= \K\circ\det$ pour un caractre $\K$ de $\A^\times$ trivial sur $\F^\times \mathrm{N}_{\E/\F}(\AE^\times)$. 
En \ref{lecaslocalnr}, les donnŽes sont locales: $E/F$ est une extension finie non ramifiŽe de corps locaux non archimŽdiens, de degrŽ $d$, $n=md$, 
$G=\mathrm{GL}_{n/F}$ et $G'=\mathrm{Res}_{E/F}(\mathrm{GL}_{m/E})$. 

\subsection{ReprŽsentations automorphes discrtes.}\label{rep aut disc} 
Soit $G$ un groupe rŽductif connexe dŽfini sur $\F$. On note $Z_G$ le centre de $G$ et $A_G$ le tore $\F$-dŽployŽ maximal de $Z_G$. 
On fixe un  caractre unitaire $\omega$ 
de $G(\A)$ trivial sur $G(\F)$. On s'intŽresse aux reprŽsentations automorphes $\pi$ de $G(\A)$ telles 
que $\omega\otimes \pi \simeq \pi$ et plus particulirement ˆ la partie discrte de la formule des traces pour $(G,\omega)$. 
On suppose donc que $\omega$ est trivial sur $Z_G(\A)$ sinon la thŽorie est vide. 

Fixons un sous-groupe compact maximal ${\bf K}= \prod_{v}{\bf K}_v$ de $G(\A)$ qui soit \guill{bon} relativement ˆ un tore $\F$-dŽployŽ maximal $A_0$ de $G$, 
cf. \cite[3.1]{LW}.   

On ne peut pas parler de 
reprŽsentation automorphe discrte sans fixer un caractre automorphe unitaire 
de $Z_G(\A)$ ou d'un sous-groupe co-compact $\mathfrak{A}_G$ de $Z_G(\F)\backslash Z_G(\A)$. Comme il est d'usage (cf. \cite[1.1]{LW}), on prend pour $\mathfrak{A}_G$ 
la composante neutre $A_{G_{\mathbb{Q}}}(\mathbb{R})^\circ$ du groupe des points sur $\mathbb{R}$ du tore central $\mathbb{Q}$-dŽployŽ maximal 
$A_{G_{\mathbb{Q}}}$ du $\mathbb{Q}$-groupe alg\'ebrique $G_{\mathbb{Q}}= \mathrm{Res}_{\F/\mathbb{Q}}(G)$. Alors $\mathfrak{A}_G$ 
est un sous-groupe de Lie connexe de $Z_G(\F_{\!\infty})\;(\subset Z_G(\F)\backslash Z_G(\A))$, 
co-compact dans $Z_G(\F)\backslash Z_G(\A)$. On note $X_{\F}(G)$ le groupe des caractres algŽbriques de $G$ qui sont dŽfinis sur $\F$. 
L'homomorphisme 
de Harish-Chandra (cf. \cite[1.1]{LW}) $${\bf H}_G: G(\A) \rightarrow \mathfrak{a}_G = \mathrm{Hom}(X_\F(G), \mathbb{R})$$ se restreint en un isomorphisme 
$\mathfrak{A}_G \buildrel\simeq\over{\longrightarrow} \mathfrak{a}_G$. En posant $G(\A)^1=\ker({\bf H}_G)$, on obtient que l'application naturelle 
$G(\A)^1 \rightarrow \mathfrak{A}_G\backslash G(\A)$ est un isomorphisme. Posons 
$$\bs{X}_{\!G} = \mathfrak{A}_G G(\F)\backslash G(\A)\simeq G(\F)\backslash G(\A)^1\ptf$$

\vskip1mm
Si $\alpha$ est un caractre unitaire de $\mathfrak{A}_G$, on note:
\begin{itemize}
\item $L^2(G(\F)\backslash G(\A),\alpha)$ l'espace des fonctions sur $G(\F)\backslash G(\A)$ qui se transforment suivant 
$\alpha$ par translations sous $\mathfrak{A}_G$ et sont $L^2$ sur $\bs{X}_{\!G}$; 
\item $\rho_{G,\alpha}$ la reprŽsentation rŽgulire droite de $G(\A)$ sur $L^2(G(\F)\backslash G(\A),\alpha)$;
\item $\rho_{G,\mathrm{disc},\alpha}$ la somme directe hilbertienne des reprŽsentations irrŽductibles qui apparaissent
discr\`etement dans la dŽcomposition spectrale de $\rho_{G,\alpha}$; 
\item $L^2_{\mathrm{disc}}(G(\F)\backslash G(\A),\alpha)\subset L^2(G(\F)\backslash G(\A),\alpha)$ l'espace de $\rho_{G,\mathrm{disc},\alpha}$;
\item $C^\infty_{\mathrm{c}}(G(\A),\alpha)=C^\infty_{\mathrm{c}}(G(\F_{\!\infty}),\alpha)\otimes C^\infty_{\mathrm{c}}(G(\Afin))$ l'espace des fonctions lisses sur $G(\A)$ 
qui se transforment suivant $\alpha^{-1}$ par translations sous $\mathfrak{A}_G$ et dont le support est d'image compacte dans $\mathfrak{A}_G \backslash G(\A)$. 
\end{itemize} 

On appelle \textit{sŽrie discrte} de $G(\A)$ une sous-reprŽsentation irrŽductible $\pi$ de $\rho_{G,\mathrm{disc},\alpha}$ (donnŽe par un sous-espace de Hilbert $G(\A)$-stable irrŽductible $V_\pi$ de $L^2(G(\F)\backslash G(\A),\alpha)$) pour un caractre unitaire $\alpha$ de $\mathfrak{A}_G$. Une sŽrie discrte $\pi$ de $G(\A)$ se dŽcompose en --- au sens o elle est est isomorphe ˆ --- un produit tensoriel restreint complŽtŽ $\widehat{\otimes}_{v} \pi_v$ o $\pi_v$ est une reprŽsentation unitaire irrŽductible de $G(\F_v)$; 
et o le \guill{restreint} est par rapport au choix pour presque toute place $v\in \Vfin$ d'un vecteur non nul ${\bf K}_v$-invariant $x_v\in V_{\pi_v}$ (cf. \cite{Fl1}). 
Pour chaque place $v\in \V$, on note $\pi_v^\infty$ la reprŽsentation de $G(\F_v)$ sur le sous-espace de $V_{\pi_v}^\infty \subset V_{\pi_v}$ formŽ des vecteurs \textit{lisses}. Si $v\in \Vfin$, $\pi_v^\infty$ est donc une reprŽsentation (lisse) au sens de la section \ref{le cas n-a}. Si $v\in \Vinf$, on peut aussi considŽrer le sous-espace $V_{\pi_v}^{\infty,\mathrm{fin}}\subset V_{\pi_v}^\infty$ formŽ des vecteurs ${\bf K}_v$-finis, \ie l'espace du $(\mathfrak{g}(\F_v),{\bf K}_v)$-module $\pi_v^{\infty,\mathrm{fin}}$ sous-jacent ˆ $\pi_v$  (cf. \ref{rappels (g,K)-modules}); o l'on a posŽ $\mathfrak{g}=\mathrm{Lie}(G)$. Pour toute place $v\in \V$ et toute fonction 
$f_v\in C^\infty_{\mathrm{c}}(G(\F_v))$, on a l'ŽgalitŽ $\mathrm{tr}(\pi_v(f_v))= \mathrm{tr}(\pi_v^\infty(f_v))$; o la trace est dŽfinie via le choix d'une mesure de Haar $\d g_v$ sur $G(\F_v)$. Si de plus $v\in \Vinf$, 
alors pour toute fonction $f_v\in C^{\infty, {\bf K}_v-\mathrm{fin}}_{\mathrm{c}}(G(\F_v))$, on a l'ŽgalitŽ $\mathrm{tr}(\pi_v(f_v))= \mathrm{tr}(\pi_v^{\infty, \mathrm{fin}}(f_v))$.

Si $\pi$ est une sŽrie discrte de $G(\A)$, la reprŽsentation automorphe $\pi^{\sharp}$ sous-jacente (cf. \cite{BJ}) est appelŽe \textit{reprŽsentation automorphe discrte} de $G(\A)$. Elle est irrŽductible et se dŽcompose en un produit tensoriel restreint $\otimes_v \pi_v^\sharp$ 
o $\pi_v^{\sharp}= \pi_v^\infty$ si $v\in \Vfin$ et $\pi_v^{\sharp}=\pi_v^{\infty,\mathrm{fin}}$ si $v\in \Vinf$. Observons que $\pi^{\sharp}$ --- et plus gŽnŽralement toute reprŽsentation automorphe irrŽductible de $G(\A)$ --- n'est donc pas ˆ proprement parler un $G(\A)$-module; en revanche c'est toujours un 
$(\mathfrak{g}(\F_\infty) \times {\bf K}_\infty) \times G(\Afin)$-module, o l'on a posŽ ${\bf K}_\infty = \prod_{v\in \Vinf} {\bf K}_v$.

Deux sŽries discrtes $\pi_1$ et $\pi_2$ de $G(\A)$ sont isomorphes si et seulement si les reprŽsentations automorphes discrtes $\pi_1^\sharp$ et $\pi_2^\sharp$ le sont. MalgrŽ qu'il est d'usage d'appeler de la mme manire une sŽrie discrte et sa reprŽsentation automorphe sous-jacente (cf. \cite[4.6]{BJ}), nous Žviterons de le faire ici. 

\subsection{Partie discrte de la formule des traces pour $(G,\omega)$.}\label{la FdT pour G} 
Continuons avec les hypothses et notations de \ref{rep aut disc}. 
Par commoditŽ nous utiliserons comme rŽfŽrence le livre \cite{LW} mme si les rŽsultats rappelŽs ici se dŽduisent facilement des travaux d'Arthur 
(la torsion par $\omega$ n'est pas une \guill{vraie} torsion).

Fixons une mesure de Haar $\d g$ sur $G(\A)$. Si $(\pi,V)$ est une reprŽsentation unitaire de $G(\A)$, pour toute fonction $f\in C^\infty_{\mathrm{c}}(G(\A))$, on note $\pi(f)= \pi(f\d g)$ l'opŽrateur sur l'espace $V$ dŽfini par l'intŽgrale (au sens de la topologie forte sur l'espace des opŽrateurs)
$$\pi(f)= \int_{G(\A)}f(g)\pi(g)\d g\ptf$$

Fixons un caractre unitaire $\alpha$ de $\mathfrak{A}_G$. Fixons aussi une mesure de Haar $\d a$ sur $\mathfrak{A}_G$ et notons 
$\d \dot{g}$ la mesure quotient $\frac{\d g}{\d a}$ sur $\mathfrak{A}_G \backslash G(\A)$. On dispose d'une application 
surjective $C^\infty_{\mathrm{c}}(G(\A))\rightarrow C^\infty_{\mathrm{c}}(G(\A),\alpha),\, f\mapsto f^\alpha$ donnŽe par 
$$f^\alpha(g)= \int_{\mathfrak{A}_G}f(ag)\alpha(a)\d a\ptf$$ 
Si $(\pi,V)$ est une reprŽsentation unitaire de $G(\A)$ se transformant suivant $\alpha$ par translations sous $\mathfrak{A}_G$, pour 
toute fonction $f\in C^\infty_{\mathrm{c}}(G(\A),\alpha)$, on dŽfinit l'opŽrateur $\pi(f) = \pi(f \d \dot{g})$ sur $V$ par une intŽgrale sur $\mathfrak{A}_G \backslash G(\A)$ comme ci-dessus. Pour $f\in C^\infty_{\mathrm{c}}(G(\A))$, on a donc $\pi(f)= \pi(f^\alpha)$.  

Pour $f\in C^\infty_{\mathrm{c}}(G(\A))$, l'opŽrateur $\rho_{G,\alpha}(f)= \rho_{G,\alpha}(f^\alpha)$ est reprŽsentŽ par le noyau 
$$K_{G,\alpha}(f;x,y)= \sum_{\gamma \in G(\F)} f^\alpha(x^{-1}\gamma y)\qhq{pour tous} x,\, y \in G(\A)\ptf$$
On dispose aussi sur l'espace $L^2(G(\F)\backslash G(\A),\alpha)$ de l'opŽrateur d'entrelacement \guill{physique} $I_\omega$ donnŽ par la multiplication par le caractre 
$\omega$ de $G(\A)$: $$(I_w \Phi)(x)= \omega(x)\Phi(x)\ptf$$ 
Pour $f\in C^\infty_{\mathrm{c}}(G(\A))$, l'opŽrateur $\rho_{G,\alpha}(f,\omega)=\rho_{G,\alpha}(f)\circ I_\omega$ est reprŽsentŽ par le noyau 
$$K_{G,\alpha}(f,\omega;x,y)= \omega(y)K_{G,\alpha}(f;x,y)\qhq{pour tous}x,\,y\in G(\A)\ptf$$
L'opŽrateur $I_\omega$ stabilise le sous-espace $L^2_{\mathrm{disc}}(G(F)\backslash G(\A),\alpha)$: il envoie une sous-repr\'esentation irr\'eductible $\pi$ de $\rho_{G,\mathrm{disc},\alpha}$ sur la reprŽsentation 
$\omega^{-1}\pi= \omega^{-1}\otimes\pi$; ou, ce qui revient au même, il envoie $\omega \pi$ sur $\pi$. 

Si $\pi$ est une sous-reprŽsentation irrŽductible $\omega$-stable de $\rho_{G,\mathrm{disc},\alpha}$, on peut choisir un prolongement $\wt{\pi}$ de $\pi$ ˆ $(G(\A),\omega)$, 
\ie un opŽrateur d'entrelacement $A \in \mathrm{Isom}_G(\omega\pi,\pi)$, cf. \cite[2.3]{LW}: on a $\wt{\pi}(g,\omega)= \pi(g)\circ A$. On dŽfinit alors comme en \cite[2.4]{LW} la \guill{multiplicitŽ tordue} 
$m(\pi,\wt{\pi})$ de $\pi$ dans la reprŽsentation $(g,\omega)\mapsto  \rho_{G,\mathrm{disc},\alpha}(g)\circ I_\omega$ de $(G(\A),\omega)$. C'est un nombre complexe qui dŽpend du choix de $\wt{\pi}$, mais l'expression 
$m(\pi,\wt{\pi}) \mathrm{trace}(\wt{\pi}(f, \omega))$ n'en dŽpend pas. 

D'aprs M\"{u}ller \cite{Mu}, pour $f\in C^\infty_{\mathrm{c}}(G(\A))$, l'opŽrateur $\rho_{G,\mathrm{disc},\alpha}(\rho,f)$ est \textit{ˆ trace} (cf. \cite[8.1.3]{LW}). 
La contribution de $L^2_{\mathrm{disc}}(G(\F)\backslash G(\A),\alpha)$ ˆ la formule des traces 
pour $(G,\omega)$ est donnŽe par l'expression, : 
$$\mathrm{tr}(\rho_{G,\mathrm{disc},\alpha}(f)\circ I_\omega)= \sum_\pi m(\pi,\wt{\pi})\mathrm{tr}(\wt{\pi}(f,\omega))$$ 
o $\pi$ parcourt un systme de reprŽsentants des classes d'isomorphisme de sous-reprŽsentations irrŽductibles $\omega$-stables 
de $\rho_{G,\mathrm{disc},\alpha}$ et 
o, pour chaque $\pi$, on a choisi un prolongement $\wt{\pi}$ de $\pi$ ˆ $(G(\A),\omega)$.

La partie discrte du dŽveloppement spectral de la formules des traces pour $(G,\omega)$ est une somme de termes, parmi lesquels on a $\mathrm{tr}(\rho_{G,\mathrm{disc},\alpha}(f)\circ I_\omega)$. 
Mais il y a aussi d'autres termes discrets, ne faisant apparaître aucune intŽgrale dans leur expression bien que provenant du spectre continu;  
nous les dŽcrivons ci-dessous.

Rappelons que l'on a choisi un tore $\F$-dŽployŽ maximal $A_0$ de $G$. On note $M_0= Z^G(A_0)$, resp. $N^G(A_0)$, le centrali\-sateur, resp. normalisateur, de $A_0$ dans $G$, $W^G_0$ le groupe de Weyl $N^G(A_0)/M_0$ et   
$\ES{L}(A_0)$ l'ensemble des $\F$-facteurs de Levi de $G$ qui contiennent $A_0$. Pour $L\in \ES{L}(A_0)$, on note $\EuScript{P}(L)$ l'ensemble 
des $\F$-sous-groupes paraboliques de $G$ de composante de Levi $L$. 

Soit $L\in \EuScript{L}(A_0)$. On note $W_0^L\subset W_0^G$ le groupe de Weyl $N^L(A_0)/M_0$ et 
on pose $$W(L)=N^G(A_L)/L=N^{W^G_0}(A_L)/W_0^L$$ Le groupe $W(L)$ opre naturellement sur le noyau $\mathfrak{a}_L^G$ du morphisme  
$\mathfrak{a}_L \rightarrow \mathfrak{a}_G$ dŽduit de l'injection naturelle (restriction des caractres) $X_{\F}(G)\subset X_{\F}(L)$. On pose 
$$W(L)_{\mathrm{r\acute{e}g}}= \{s\in W(L)\,\vert\, \det(s-1\vert \mathfrak{a}_L^G)\neq 0\}\ptf$$ 
On dŽfinit comme plus haut le sous-groupe 
$$\mathfrak{A}_L= A_{L_{\mathbb{Q}}}(\mathbb{R})^\circ\subset A_L(\F_\infty)\ptf$$ Puisque $\mathfrak{A}_G\simeq \mathfrak{a}_G$ et $\mathfrak{A}_L\simeq \mathfrak{a}_L$, l'inclusion 
$\mathfrak{A}_G\subset \mathfrak{A}_L$ fournit une section $\mathfrak{a}_G \rightarrow \mathfrak{a}_L$ du morphisme surjectif 
$\mathfrak{a}_L \rightarrow \mathfrak{a}_G$ et donc une dŽcomposition $\mathfrak{a}_L= \mathfrak{a}_L^G \oplus \mathfrak{a}_G$. On note $\mathfrak{A}_L^G$ le sous-groupe de $\mathfrak{A}_L$ 
prŽ-image de $\mathfrak{a}_L^G$ par l'isomorphisme $\mathfrak{A}_L\simeq \mathfrak{a}_L$. Le caractre unitaire $\alpha$ de $\mathfrak{A}_G$ se prolonge de manire unique en caractre (unitaire) de $\mathfrak{A}_L$ trivial sur $\mathfrak{A}_L^G$, notŽ $\alpha_L$. Observons que le caractre $\omega$ de $G(\A)$ est trivial sur $\mathfrak{A}_L^G$. 

Pour $L\in \ES{L}(A_0)$ et $P\in \EuScript{P}(L)$, on peut former la 
reprŽsentation $\pi_{P,\alpha}$ de $G(\A)$ induite parabolique suivant $P(\A)$ de la reprŽsentation $\rho_{L,\mathrm{disc},\alpha_L}$, \ie la somme directe (hilbertienne) des induites paraboliques ˆ $G(\A)$ suivant $P(\A)$ des sŽries discrtes de $L(\A)$ de caractre central prolongeant $\alpha$. Notons $\ES{V}_{P,\alpha}$ l'espace de $\pi_{P,\alpha}$. 
Si $Q\in \EuScript{P}(L)$, l'opŽrateur d'entrelacement\footnote{DonnŽ par une formule intŽgrale pour $\lambda\in \mathfrak{a}_L^*\otimes \mathbb{C}$ assez rŽgulier, cf. \cite[5.2]{LW}.} $M_{P\vert Q}(1,\lambda)$ dŽfinit par prolongement mŽromorphe un opŽrateur d'entrelacement 
$$M_{P\vert Q}(1,0): \ES{V}_{Q,\alpha}\rightarrow \ES{V}_{P,\alpha}\ptf$$ Pour $Q={^s\!P}\;(=sPs^{-1})$ avec $s\in W(L)$ et $y\in G(\A)$, on dŽfinit un opŽrateur $$\pi_{P,\alpha}(s,y,\omega): \ES{V}_{P,\alpha} \rightarrow \ES{V}_{Q,\alpha}$$ en posant 
$$(\pi_{P,\alpha}(s,y,\omega)\Phi)(x)= \omega(n_s^{-1}xy)\Phi(n_s^{-1}xy),$$ o $n_s$ est un relvement de $s$ dans $N^G(L)(\F)$. En notant $\bs{s}_P: \ES{V}_{P,\alpha}\rightarrow \ES{V}_{Q,\alpha}$ 
l'opŽ\-rateur d'entrelacement dŽfini par $(\bs{s}_P\Phi)(x)= \Phi(n_s^{-1}x)$, on a 
$$\pi_{P,\alpha}(s,y ,\omega)= \bs{s}_P \circ \pi_{P,\alpha}(y) \circ I_\omega\ptf$$ 
 
Pour $s\in W(L)_{\mathrm{r\acute{e}g}}$ et $y\in G(\A)$, notons $\rho_{P,\mathrm{disc},s,\alpha}(y,\omega)$ l'opŽrateur sur $\ES{V}_{P,\alpha}$ dŽfini par 
$$\rho_{P,\mathrm{disc},s,\alpha}(y,\omega)= M_{P,{^s\!P}}(1,0) \circ  \pi_{P,\alpha}(s,y,\omega)\ptf$$ 
Il dŽfinit une reprŽsentation de $(G(\A),\omega)$ d'espace $\ES{V}_{P,\alpha}$ (cf. \cite[2.3]{LW}), \cad que $$y\mapsto \rho_{P,\mathrm{disc},s,\alpha}(y)= M_{P,{^s\!P}}(1,0) \circ \bs{s}_{P}\circ  \pi_{P,\alpha}(y)$$ est une reprŽsentation 
$\omega$-stable de $G(\A)$. 

\vskip1mm
La partie discrte du 
dŽveloppement spectral de la formule des traces pour $(G,\omega)$ est donnŽe par l'expression \cite[14.3.2]{LW}, pour $f\in C^\infty_{\mathrm{c}}(G(\A))$:  
$$I_{\mathrm{disc},\alpha}^G(f,\omega)
=\!\! \sum_{L\in \ES{L}(A_0)} \!\frac{\vert W^L_0 \vert}{\vert W_0^G\vert}\!\sum_{s\in W(L)_{\mathrm{r\acute{e}g}}} \!\!
\vert \det(s-1\vert \mathfrak{a}_L^G )\vert^{-1}\mathrm{tr}\left(\rho_{P,\mathrm{disc},s,\alpha}(f,\omega) \right)\leqno{(1)}$$
o, pour chaque $L$, on a choisi un $P\in \EuScript{P}(L)$; le choix de $P$ n'affecte pas la trace. Comme plus haut, on a notŽ  
$\rho_{P,\mathrm{disc},s,\alpha}(f,\omega)=\rho_{P,\mathrm{disc},s,\alpha}(f^\alpha\!,\omega)$ l'opŽrateur 
$$\int_{G(\A)} f(g) \rho_{G,\mathrm{disc},s,\alpha}(g,\omega)\d g=\int_{\mathfrak{A}_G\backslash G(\A)}f^\alpha(g)\rho_{P,\mathrm{disc},s,\alpha}(g,\omega)\d \dot{g}\ptf$$

\begin{remark}\label{travaux de FLM} 
\textup{\begin{enumerate}
\item[(i)]On sait d'aprs Finis, Lapid et M\"{u}ller \cite{FL,FLM} que le dŽveloppement spectral est absolument convergent (cf. \cite[14.3.1]{LW}). 
En particulier pour tous $L\in \ES{L}(A_0)$ et $s\in W(L)_{\mathrm{r\acute{e}g}}$, et pour toute fonction $f\in C^\infty_{\mathrm{c}}(G(\A))$, l'opŽrateur $\rho_{P,\mathrm{disc},s,\alpha}(f,\omega) $ est ˆ trace. 
\item[(ii)] Soient $L\in \ES{L}(A_0)$, $s\in W(L)_{\mathrm{r\acute{e}g}}$ et $P$ un $\F$-sous-groupe parabolique de $G$ de composante de Levi 
$L$. Pour $\tau$ une composante irrŽductible de $\rho_{L,\mathrm{disc},\alpha_L}$, notons $\pi_\tau$ l'induite 
parabolique de $\tau$ ˆ $G$ suivant $P$. Pour que $\pi_\tau$ donne une contribution non triviale ˆ la trace $\mathrm{tr}(\rho_{P,\mathrm{disc},s,\alpha}(f,\omega))$, il faut que $\omega\vert_{L(\A)}\otimes s(\tau) \simeq \tau$. Observons que cela implique que la restriction ˆ $\mathfrak{A}_L$ du caractre central $\omega_\tau$ de $\tau$ vŽrifie $\omega_\tau(a^{-1}s(a))= \omega(a)$ pour tout $a\in \mathfrak{A}_L$, \ie que $\omega_\tau$ est trivial sur $\mathfrak{A}_L^G$.
\end{enumerate}
}
\end{remark}

L'expression (1) n'est pas exactement celle de \cite[14.3.2]{LW} o $\alpha$ est le caractre trivial de $\mathfrak{A}_G$, mais elle s'en dŽduit aisŽment. 
En effet le caractre unitaire $\alpha$ de $\mathfrak{A}_G$ s'Žtend en un caractre unitaire de $G(\A)$ trivial sur $G(\A)^1$, encore notŽ $\alpha$. 
Pour $L\in \EuScript{L}(A_0)$, $P\in \EuScript{P}(L)$ et $s\in W(L)_{\mathrm{r\acute{e}g}}$, on a 
$$\rho_{P,\mathrm{disc},s,\alpha}\simeq \rho_{P,\mathrm{disc},s,1}\otimes \alpha$$ o 
\guill{$1$} est le caractre trivial de $\mathfrak{A}_G$; et pour $f\in C^\infty_{\mathrm{c}}(G(\A))$, on a l'ŽgalitŽ $$\mathrm{tr}\left(\rho_{P,\mathrm{disc},s,\alpha}(f^\alpha\!,\omega) \right)=\mathrm{tr}\left(\rho_{P,\mathrm{disc},s,1}(f^1,\omega) \right)\ptf$$ D'o l'ŽgalitŽ, pour toute fonction $f\in C^\infty_{\mathrm{c}}(G(\A))$: 
$$I^G_{\mathrm{disc},\alpha}(f,\omega)= I^G_{\mathrm{disc},1}(f,\omega)\ptf$$

\begin{remark}\label{point central}
\textup{
\begin{enumerate}
\item[(i)]Si le \guill{point central} $T_0\in \mathfrak{a}_{L_0}^G$ (cf. \cite[3.3.3]{LW}) est $T_0=0$ --- \textit{e.g.} si $G= \mathrm{GL}_{n/\F}$ ---, alors d'aprs \cite[5.2.1, 5.2.2\,(1)]{LW}, on a 
$$M_{P,{^s\!P}}(1,0)\circ \bs{s}_P = M_{P,{^s\!P}}(1,0) \circ M_{{^s\!P},P}(s,0)= M_{P\vert P}(s,0)\ptf$$ 
Dans ce cas la reprŽsentation $\rho_{P,\mathrm{disc},s,\alpha}$ de $(G(\A),\omega)$ s'Žcrit aussi $$\rho_{P,\mathrm{disc},s,\alpha}(y,\omega)=M_{P\vert P}(s,0)\circ \pi_{P,\alpha}(y)\circ I_\omega\ptf$$
\item[(ii)]\`A une constante (strictement positive) prs, l'expression (1) dŽpend du choix de la mesure de Haar $\d g$ sur $G(\A)$. 
En revanche elle ne dŽpend pas du choix de la mesure de Haar $\d a$ sur $\mathfrak{A}_G$ (servant ˆ dŽfinir l'application $f\mapsto f^\alpha$ et la mesure quotient $\d\dot{g}$ sur $\mathfrak{A}_G\backslash G(\A)$).  
\end{enumerate}}
\end{remark}

\subsection{Variante avec un caractre automorphe de $A_G(\A)$.}\label{variante avec car aut} On peut aussi dans l'expression \ref{la FdT pour G}\,(1) remplacer le caractre unitaire $\alpha$ de $\mathfrak{A}_G$ par un caractre automorphe unitaire de $A_G(\A)$. PrŽcisons cette variante. 

Soit $\Xi(G,\alpha)$ l'ensemble des caractres de $A_G(\A)$ qui sont triviaux sur $A_G(\F)$ et prolongent $\alpha$. C'est un espace principal homogne sous le dual de Pontryagin $\Xi(G,1)$ du groupe abŽlien compact $$\bs{X}_{\!A_G}=\mathfrak{A}_GA_G(\F)\backslash A_G(\A)\simeq A_G(\F)\backslash A_G(\A)^1\ptf$$ 
Observons que tout caractre $\xi\in \Xi(G,\alpha)$ est unitaire. Pour $f\in C^\infty_{\mathrm{c}}(G(\A))$ et $a\in A_G(\A)$, notons $f_a$ la fonction $g\mapsto f(ag)$. Pour $\xi\in \Xi(G,\alpha)$, posons 
$$I^G_{\mathrm{disc},\xi}(f,\omega)= \int_{\bs{X}_{\!A_G}}\xi(a) I^G_{\mathrm{disc},\alpha}(f_{a} ,\omega )\d \dot{a}$$
o $\d \dot{a}$ est la mesure normalisŽe sur $\bs{X}_{\!A_G}$ --- \ie $\mathrm{vol}(\bs{X}_{\!A_G},\d \dot{a})= 1$ --- et o, pour chaque $\dot{a}\in \bs{X}_{\!A_G}$, on a choisi un reprŽsentant $a$ de $\dot{a}$ dans $A_G(\A)^1$ pour dŽfinir $f_a$.  La fonction $f_a$ dŽpend bien sžr du choix de $a$ mais l'expression 
$I^G_{\mathrm{disc},\alpha}(f_a,\omega)$ n'en dŽpend pas. Par inversion de Fourier, on a la dŽcomposition 
$$I^G_{\mathrm{disc},\alpha}(f,\omega)=\sum_{\xi \in \Xi(G,\alpha)} I^G_{\mathrm{disc},\xi}(f,\omega)\ptf \leqno{(1)}$$ 

Soit $\xi\in \Xi(G,\alpha)$. On note $\rho_{G,\mathrm{disc},\xi}$ la somme directe hilbertienne des sous-reprŽsentations irrŽductibles de $\rho_{G,\mathrm{disc},\alpha}$ de caractre central prolongeant $\xi$. 
Pour $L\in \ES{L}(A_0)$, puisque $\mathfrak{A}_L\cap A_G(\A)= \mathfrak{A}_G$ (et donc $\mathfrak{A}_L^G \cap A_G(\A)=\{1\}$), le caractre $\xi$ de $A_G(\A)$ se prolonge de manire unique en un caractre (unitaire) $\xi_L$ de $\mathfrak{A}_L^GA_G(\A)$ trivial sur $\mathfrak{A}_L^G$. 
On note $\rho_{L,\mathrm{disc},\xi_L}$ la somme directe hilbertienne des sous-reprŽsentations irrŽductibles de $\rho_{L,\mathrm{disc},\alpha_L}$ de caractre central prolongeant $\xi_L$, \ie prolongeant $\xi$. Si $P$ est un $\F$-sous-groupe parabolique de $G$ de composante de Levi $L$, on note $\pi_{P,\xi}$ l'induite parabolique de $\rho_{L,\mathrm{disc},\xi_L}$ ˆ $G(\A)$ suivant $P(\A)$. Alors on a 
$$I_{\mathrm{disc},\xi}^G(f,\omega)
=\!\! \sum_{L\in \ES{L}(A_0)} \!\frac{\vert W^L_0 \vert}{\vert W_0^G\vert}\!\sum_{s\in W(L)_{\mathrm{r\acute{e}g}}} \!\!
\vert \det(s-1\vert \mathfrak{a}_L^G )\vert^{-1}\mathrm{tr}\left(\rho_{P,\mathrm{disc},s,\xi}(f,\omega) \right)\leqno{(2)}$$ 
o, pour $y\in G(\A)$, l'opŽrateur $\rho_{P,\mathrm{disc},s,\xi}(y,\omega)$ sur l'espace $\ES{V}_{P,\xi}$ de $\pi_{P,\xi}$ est dŽfini comme en \ref{la FdT pour G}.

\subsection{L'ŽgalitŽ ˆ Žtablir $(G=\mathrm{GL}_{n/\F},\omega=\K\circ \det)$.}\label{l'ŽgalitŽ} 
Soit $\E/\F$ une extension finie cyclique de corps de nombres, 
de degrŽ $d$. On pose $\Gamma=\Gamma(\E/\F)$. On fixe un caractre 
$\K$ de $\A^\times$ de noyau $\F^\times \mathrm{N}_{\E/\F}(\AE^\times)$ et un gŽnŽrateur $\sigma$ de $\Gamma$.

Soit $G= \mathrm{GL}_{n/\F}$ pour un entier $n\geq 1$. Posons  
$$\omega=\K\circ \det : G(\A) \rightarrow \mathbb{C}^\times\ptf$$ Observons que $\omega$ est trivial sur $Z_G(\A)$. D'autre part 
le centre $Z_G$ de $G$ co\"{\i}ncide avec $A_G$ et s'identifie naturellement au groupe multiplicatif $\mathbb{G}_{\mathrm{m}/\F}\;(=\mathbb{G}_{\mathrm{m}}\times_{\mathbb{Z}}\F)$. Ainsi on a l'identification 
$$\mathfrak{A}_G= (\mathbb{R}^\times)^\circ= \mathbb{R}_{>0} \subset (\F\otimes_{\mathbb{Q}}\mathbb{R})^\times \subset \A^\times \;( = Z_G(\A))\ptf$$ 
On prend pour $A_0$ le tore diagonal $(\mathbb{G}_{\mathrm{m}/\F})^n$ de $G$. En toute place 
$v\in \Vfin$, on prend ${\bf K}_v= \GLn(\mathfrak{o}_F)$. 

Supposons que $d=[\E:\F]$ divise $n$. Posons $m= \frac{n}{d}$ et $G'= \mathrm{Res}_{\E/\F}(\mathrm{GL}_{m/\E})$. On a 
encore des identifications naturelles $$A_{G'}= \mathbb{G}_{\mathrm{m}/\F}\subset \mathrm{Res}_{\E/\F}(\mathbb{G}_{\mathrm{m}/\E})=Z_{G'}\ptf$$ 
D'o les identifications $A_{G'}=A_G$ et $\mathfrak{A}_{G'}= \mathfrak{A}_G$.  
On note $A'_{0,\E}$ le tore diagonal $(\mathbb{G}_{\mathrm{m}/\E})^m$ de $\mathrm{GL}_{m/\E}$ et $A'_0$ 
le tore $\F$-dŽployŽ maximal de $G'$ dŽfini par $$A'_0 = (\mathbb{G}_{\mathrm{m}/\F})^m\subset M'_0= Z^{G'\!}(A'_0)\;(=\mathrm{Res}_{\E/\F}(A'_{0,\E})). $$ 
En toute place $v\in \Vfin$, on prend ${\bf K}'_v=\GLm(\mathfrak{o}_{\E_v})=\prod_{w\vert v}\GLm(\mathfrak{o}_{\E_w})$. 
On note $\ES{L}(A'_0)$ l'ensemble des $\F$-facteurs de Levi de $G'$ qui contiennent $A'_0$. Observons que $\ES{L}(A'_0)$ est naturellement en bijection avec l'ensemble des $\E$-facteurs de Levi de $\mathrm{GL}_{m/\E}$ qui contiennent $A'_{0,\E}$. Pour $L'\in \EuScript{L}(A'_0)$, on note 
$\EuScript{P}(L')$ l'ensemble des $\F$-sous-groupes paraboliques de $G'$ de composante de Levi $L'$. 

Fixons un caractre unitaire $\beta$ de $\mathfrak{A}_{G'}=\mathfrak{A}_G$. Comme en \ref{la FdT pour G}, on considre la reprŽsentation $\rho_{G'\!, \beta}$ d'espace $L^2(G'(\F)\backslash G'(\mathbb{A}_F),\beta)$. La partie discrte du dŽveloppement spectral de la formule des traces pour $G'$ est donnŽe par l'expression \ref{la FdT pour G}\,(1) (avec $\omega=1$, que l'on omet ici), 
pour $f'\in C^\infty_{\mathrm{c}}(G'(\A))$: 
$$I_{\mathrm{disc},\beta}^{G'}(f')=\sum_{L'\in \ES{L}(A'_0)} \frac{\vert W^{L'}_0 \vert}{\vert W_0^{G'}\vert}
\sum_{s\in W(L')_{\mathrm{r\acute{e}g}}}\!\!\! \vert \det(s-1\vert \mathfrak{a}_{L'}^{G'} )\vert^{-1}\mathrm{tr}\left(\rho_{P'\!,\mathrm{disc},s,\beta}(f')
\right)\leqno{(1)}$$ 
o, pour chaque $L'$, on a choisi un $P'\in \ES{L}(L')$.  
Observons que puisque le \guill{point central} pour $G'$ est Žgal ˆ $0$, l'opŽrateur 
$\rho_{P'\!,\mathrm{disc},s,\beta}(f')= \rho_{P'\!,\mathrm{disc},s,\beta}(f'^\beta)$ est donnŽ par la formule (remarque \ref{point central}\,(i))
$$\rho_{P'\!,\mathrm{disc},s,\beta}(f')=M_{P'\vert P'}(s,0) \circ \pi_{P'\!,\beta}(f')\ptf$$ 
Pour dŽfinir ces opŽrateurs, on a choisi comme en \ref{la FdT pour G} une mesure de Haar $\d h$ sur $G'(\A)$.  

Pour allŽger l'Žcriture, posons $$\G= G(\F)= \GLn(\F) \qhq{et} \G'= G'(\F)=\GLm(\E)\ptf$$ 
De mme pour $v\in \V$, posons $\G_v=\GLn(\F_v)$ et $\G'_v=\GLm(\E_v)$.  
On fixe une $\F$-base de $\E^m$, d'o un plongement $\varphi$ de $\G'$ dans $\G$. 
Pour $v\in \V$, cette base donne une $\F_v$-base de $(\E_v)^m$, d'o un plongement $\varphi_v$ de $\G'_v$ dans $\G_v$. Pour presque toute place $v\in \Vfin$, on obtient en fait une $\mathfrak{o}_v$-base de $(\mathfrak{o}_{\E_v})^m$, et $\varphi_v$ se restreint en un plongement de 
${\bf K}'_v = \GLm(\mathfrak{o}_{\E_v})$ dans ${\bf K}_v= \GLn(\mathfrak{o}_v)$. 

On dŽfinit comme dans le cas local un facteur de transfert $\wt{\bs{\Delta}}(\gamma)\in \E^\times$ pour tout $\gamma \in \G' \cap \G_{\mathrm{r\acute{e}g}}$; on a aussi, pour $v\in \V$, un facteur de transfert $\wt{\bs{\Delta}}_v(\gamma)\in \E_v^\times$ pour tout $\gamma\in \G'_v\cap \G_{v,\mathrm{r\acute{e}g}}$, et $\wt{\bs{\Delta}}$ est la restriction de $\wt{\bs{\Delta}}$ aux ŽlŽments $\F$-rationnels de 
$\G'_v\cap \G_{v,\mathrm{r\acute{e}g}}$. On choisit un ŽlŽment $\bs{e}\in \E^\times$ tel que $\sigma \bs{e} = (-1)^{m(d-1)}\bs{e}$. Pour presque toute place 
$v\in\Vfin$, $\bs{e}$ appartient ˆ $\mathfrak{o}_{\E_v}^\times$. Pour $v\in \V$, on note $\bs{\Delta}_v$ la fonction sur $\G'_v\cap \G_{v,\mathrm{r\acute{e}g}}$ dŽfinie par 
$$\bs{\Delta}_v(\gamma)= \K_v(\bs{e}\wt{\bs{\Delta}}(\gamma))\ptf$$ 
Si $\gamma\in \G' \cap \G_{\mathrm{r\acute{e}g}}$, puisque $\K$ est trivial sur 
$\F^\times$, on a $$\prod_v \bs{\Delta}_v(\gamma) = \K(\bs{e}\wt{\bs{\Delta}}(\gamma))=1\ptf$$ 
On a aussi, par la formule du produit, 
$$\prod_v \vert D_{\G_v}(\gamma)\vert_{\F_v}=1\qhq{et} \prod_v \vert D_{\G'_v}(\gamma)\vert_{\E_v}=1 \ptf $$ 

La mesure de Haar $\d g$ sur $G(\A)$ se dŽcompose en un produit de mesures de Haar locales $\d g_v$ sur $\G_v$, avec la propriŽtŽ que 
pour presque toute place $v\in \Vfin$, on a $\mathrm{vol}({\bf K}_v,\d g_v)=1$. De mme, la mesure de Haar $\d h$ sur $G'(\A)$ se dŽcompose en un produit de mesures de Haar 
$\d h_v$ sur $\G'_v$, avec $\mathrm{vol}({\bf K}'_v, \d h_v)=1$ pour presque toute place $v\in \Vfin$. 

Pour $v\in \V$, on a donc une notion de concordance pour les fonctions $f_v\in C^\infty_{\mathrm{c}}(\G_v)$ et 
$\phi_v\in C^\infty_{\mathrm{c}}(\G'_v)$. 

Reprenons les notations de \ref{un lemme local-global}. En particulier pour $v\in \Vfin$, on a un 
homomorphisme d'algbres 
$$b_v: \mathcal{H}_v= \mathcal{H}(\G_v,{\bf K}_v) \rightarrow \mathcal{H}'_v= \mathcal{H}(\G'_v,{\bf K}'_v)\ptf$$ 

Soit $C^{\infty,\bs{K}-\mathrm{fin}}_{\mathrm{c}}(G(\A))\subset C^\infty_{\mathrm{c}}(G(\A))$ 
le sous-espace formŽ des fonctions qui sont $\bs{K}$-finies ˆ droite et ˆ gauche. On a 
$$C^{\infty,\bs{K}-\mathrm{fin}}_{\mathrm{c}}(G(\A))= \left(\otimes_{v\in \Vinf} C^{\infty, \bs{K}_v-\mathrm{fin}}_{\mathrm{c}}(\bs{G}_v)\right)\otimes C^\infty_{\mathrm{c}}(G(\Afin))\ptf$$

\begin{definition}\label{fonctions associŽes}
\textup{(Une version globale de la notion de concordance.) Deux fonctions $f\in C^{\infty,\bs{K}-\mathrm{fin}}_{\mathrm{c}}(G(\A))$ et $\phi\in C^{\infty,\bs{K}'-\mathrm{fin}}_{\mathrm{c}}(G'(\A))$ sont dites \textit{associŽes} si elles sont dŽcomposŽes suivant les places de $\F$, \ie 
si $f= \prod_v f_v$ et $\phi= \prod_v \phi_v$, et s'il existe un sous-ensemble fini $S\subset \V$ contenant $\Vinf$ 
et tel que: 
\begin{itemize} 
\item pour toute place $v\in S$, les fonctions $f_v\in C^\infty_{\mathrm{c}}(\G_v)$ et $\phi_v\in C^\infty_{\mathrm{c}}(\G'_v)$ sont concordantes;
\item pour toute place $v\in \V^S=\V\smallsetminus S$, on a $f_v\in \mathcal{H}_v$ et $\phi_v= b_vf_v$. 
\item pour presque toute place $v\in \V^S$, on a $f_v=f_{v,0}$ (l'ŽlŽment unitŽ de $\mathcal{H}_v$).
\end{itemize} 
On note $$\bs{\mathfrak{F}}=\bs{\mathfrak{F}}^{\bs{K},\bs{K'}}_{G(\A),G'(\A)}$$ le sous-ensemble de $ C^{\infty,\bs{K}-\mathrm{fin}}_{\mathrm{c}}(G(\A))\times  C^{\infty,\bs{K}'-\mathrm{fin}}_{\mathrm{c}}(G'(\A))$ formŽ des paires de fonctions $(f,\phi)$ qui sont associŽes. }
\end{definition}

Soit $\alpha$ le caractre unitaire de $\mathfrak{A}_G$ dŽfini par 
$$\alpha = \K^{md(d-1)/2}\beta\ptf$$
L'ŽgalitŽ ˆ Žtablir est la suivante, pour toute paire de fonctions $(f,\phi)\in \bs{\mathfrak{F}}$: 
$$I^G_{\mathrm{disc},\alpha}(f,\omega)=\frac{1}{d}I^{G'}_{\mathrm{disc},\beta}(\phi)\pvg\leqno{(2)}$$ soit encore
$$\sum_{L\in \ES{L}(A_0)} \frac{\vert W^L_0 \vert}{\vert W_0^G\vert}\sum_{s\in W(L)_{\mathrm{r\acute{e}g}}} \vert \det(s-1\vert \mathfrak{a}_L^G )\vert^{-1}\mathrm{tr}\left(M_{P\vert P}(s,0)\circ 
 \pi_{P,\alpha}(f)\circ I_\omega\right)\leqno{(3)}$$
$$= \frac{1}{d}\sum_{L'\in \ES{L}(A'_0)} \frac{\vert W^{L'}_0 \vert}{\vert W_0^{G'}\vert}\sum_{s\in W(L')_{\mathrm{r\acute{e}g}}} \vert \det(s'-1\vert \mathfrak{a}_{L'}^{G'} )\vert^{-1}\mathrm{tr}\left(M_{P'\vert P'}(s,0)
\circ \pi_{P'\!,\beta}(\phi)\right)\ptf$$ 
Bien sžr pour que l'ŽgalitŽ (3) ait un sens, il faut prŽciser les mesures de Haar $\d g$ sur $G(\A)$ et $\d h$ sur $G'(\A)$, ce que nous ferons plus loin (voir \ref{dec endos}).

Nous verrons que l'ŽgalitŽ (3) n'est autre que la dŽcomposition endoscopique de la partie discrte du dŽveloppement spectral de la formule des traces pour $(G,\omega)$: en effet $G'$ est, ˆ Žquivalence prs, le seul groupe endoscopique elliptique de $(G,\omega)$. 

\begin{remark}
\textup{
Pour $L\in \ES{L}(A_0)$, on a notŽ $\pi_{P,\alpha}$ l'induite parabolique de $\rho_{L,\mathrm{disc},\alpha}$ ˆ $G(\A)$ suivant $P(\A)$. 
Comme la fonction test $f$ est supposŽe $\bs{K}$-finie ˆ droite et ˆ gauche, pour $s\in W(L)_{\mathrm{r\acute{e}g}}$, la trace 
$\mathrm{tr}\left(M_{P\vert P}(s,0)\circ \pi_{P,\alpha}(f)\circ I_\omega\right)$ ne \guill{voit} que les reprŽsentations automorphes sous-jacentes 
aux composantes irrŽductibles de $\pi_{P,\alpha}$, \ie les induites paraboliques ˆ $G(\A)$ suivant $P(\A)$ des 
reprŽsentations automorphes sous-jacentes aux composantes irrŽductibles de $\rho_{L,\mathrm{disc},\alpha_L}$. La mme remarque s'applique bien sžr ˆ 
la trace $\mathrm{tr}\left(M_{P'\vert P'}(s,0)\circ \pi_{P'\!,\beta}(\phi)\right)$ pour $L'\in \ES{L}(A'_0)$. }
\end{remark}

\subsection{Endoscopie pour $(\mathrm{GL}_{n/\F},\K\circ \det)$.}\label{endoscopie} On fixe une clôture algŽbrique $\overline{\F}$ de $\F$ contenant $\E$. On note $\Gamma_{\F}$ le groupe de Galois de $\overline{\F}/\F$ et $W_{\F}$ son groupe de Weil. 
On a donc $W_{\F}/W_{\E}\simeq \Gamma_\F/ \Gamma_\E = \Gamma\;(=\Gamma(\E/\F))$. Soit $\bs{a}$ le caractre de $W_{\F}$ obtenu en composant $\K$ 
avec le morphisme $W_{\F}\rightarrow W_{\F}^{\mathrm{ab}}\simeq \A^\times/\F^\times$ fourni par 
la thŽorie du corps de classes. Il se factorise ˆ travers $W_{\F}/W_{\E}$ et donc dŽfinit un caractre de $\Gamma$, encore notŽ $\bs{a}$, que l'on peut voir comme un caractre de $\Gamma_\F$ trivial sur $\Gamma_\E$. 

On suppose toujours que $(G,\omega)=(\mathrm{GL}_{n/\F},\K\circ \det)$. 
On utilise ici la notion de donnŽe endoscopique 
dŽfinie dans \cite[I.1.5, VI.3.1]{MW2}. On note $\wh{G}= \mathrm{GL}_n(\mathbb{C})$ le groupe dual de $G$ et 
${^LG}= \wh{G}\times W_{\F}$ son $L$-groupe. On appelle \textit{donnŽe endoscopique de $(G,\bs{a})$} un 
triplet $\underline{G}'= (G',\ES{G}',s)$ o:
\begin{itemize}
\item $G'$ est un groupe rŽductif connexe dŽfini et quasi-dŽployŽ sur $\F$;
\item $\ES{G}'$ est un sous-groupe fermŽ de ${^LG}$; 
\item $s$ est un ŽlŽment semi-simple de $\wh{G}$ tel que $\ES{G}' \cap \wh{G}= \wh{G}_s$ o $\wh{G}_s$ est le centra\-lisateur de $s$ dans $\wh{G}$ (on sait qu'il est connexe);
\item pour tout $(g,\sigma)\in \ES{G}'$, on a $s(g,\sigma)s^{-1}= (\bs{a}(\sigma)g,\sigma)$.
\end{itemize}
On a donc une suite $$1 \rightarrow \wh{G}_s \rightarrow \ES{G}' \rightarrow W_{\F} \rightarrow 1\ptf$$ On demande que cette suite soit exacte et scindŽe, \cad qu'il existe une section $W_{\F}\rightarrow \ES{G}'$ qui 
soit un morphisme continu de groupes. Fixons une paire de Borel ŽpinglŽe $\wh{\ES{E}}'=\left(\wh{B}'\!,\wh{T}'\!, (\wh{E}'_{\alpha'})_{\alpha'\in \Delta'}\right)$ de $\wh{G}_s$. Pour $\sigma \in W_{\F}$, on choisit un 
ŽlŽment $g_\sigma=(g(\sigma),\sigma)\in \ES{G'}$ tel que l'automorphisme $\mathrm{Int}_{g_\sigma}$ de $\wh{G}_s$ conserve $\wh{\ES{E}}'$. 
L'application $\sigma\mapsto \mathrm{Int}_{g_\sigma}$ s'Žtend en une action galoisienne sur $\wh{G}_s$. On suppose que $\wh{G}_s$ muni de cette action est un groupe dual de $G'$. On peut donc poser 
$$\wh{G}'= \wh{G}_s \qhq{et} {^LG'}= \ES{G}'\ptf$$ 

Deux donnŽes endoscopiques $\underline{G}'_1=(G'_1,\ES{G}'_1,s_1)$ et $\underline{G}'_2=(G'_2,\ES{G}'_2,s_2)$ de $(G,\omega)$ sont dites \textit{Žquivalentes} s'il existe un ŽlŽment $x\in \wh{G}$ tel que 
$$x\ES{G}'_1x^{-1}= \ES{G}'_2\quad \hbox{et}\quad xs_1x^{-1} \in Z_{\wh{G}}s_2\ptf$$ De l'isomorphisme 
$\mathrm{Int}_{x^{-1}}: \wh{G}'_2\rightarrow \wh{G}'_1$ se dŽduit par dualitŽ un $\F$-isomorphisme de groupes algŽbriques $\iota_x: G_1 \rightarrow G_2$, bien dŽfini 
modulo l'action de $\mathrm{Ad}(G'_1)(\F)$ ou, ce qui revient au même, de $\mathrm{Ad}(G'_2)(\F)$. En particulier pour une seule donnŽe endoscopique 
$\underline{G}'=(G'\!,\ES{G}'\!,s)$ de $(G,\bs{a})$, le groupe $\mathrm{Aut}(\underline{G}')$ des Žquivalences entre cette donnŽe et elle-même contient $\wh{G}'$. 
On pose $$\mathrm{Out}(\underline{G}')= \mathrm{Aut}(\underline{G}')/\wh{G}'\ptf$$ Soit $\mathrm{Out}_\F(G')= \mathrm{Aut}_\F(G')/\mathrm{Ad}(G')(\F)$ le groupe des 
$\F$-automor\-phismes extŽrieurs de $G'$. L'application $\mathrm{Aut}(\underline{G}')\rightarrow \mathrm{Out}_\F(G'),\, x \mapsto \iota_x$ se quotiente en un isomorphisme 
$$\mathrm{Out}(\underline{G}')\buildrel\simeq\over{\longrightarrow} \mathrm{Out}_\F(G')\ptf$$

\vskip1mm 
Une donnŽe endoscopique $\underline{G}'=(G'\!,\ES{G}'\!,s)$ de $(G,\bs{a})$ est dite \textit{elliptique} si elle vŽrifie l'inclusion 
$$(Z_{\wh{G}'})^{\Gamma_\F,\circ} \subset Z_{\wh{G}}$$ o $(Z_{\wh{G}'})^{\Gamma_\F,\circ}$ est la composante neutre du sous-groupe du centre $Z_{\wh{G}'}$ 
de $\wh{G}'$ formŽ des points fixes sous $\Gamma_\F$. 
Nous allons dans un premier temps exhiber une donnŽe endoscopique elliptique de $(G,\bs{a})$. Nous verrons ensuite que c'est la seule ˆ Žquivalence prs. 

Fixons un relvement $\phi\in W_{\F}$ de $\sigma$ --- rappelons que $\sigma$ est un gŽnŽrateur de $\Gamma$ --- et notons $\zeta\in \mathbb{C}^\times$ la racine primitive $d$-ime de l'unitŽ dŽfinie par $\bs{a}(\phi)=\zeta$. 
Soit $\wh{T}=(\mathbb{C}^\times)^n$ le tore diagonal de $\wh{G}$. Notons $\eta: (\mathbb{C}^\times)^d \rightarrow \wh{T}$ le plongement dŽfini par 
$$\eta(z_1,\ldots , z_d) = (z_1,\ldots , z_1, z_2,\ldots ,z_2,\ldots , z_d, \ldots , z_d)$$ o chaque $z_i$ apparaît $m$ fois. Soit $s$ l'ŽlŽment de $\wh{T}$ dŽfini par 
$$s= \eta(1,\zeta, \ldots , \zeta^{d-1})\ptf$$
Le centralisateur $\wh{G}_s$ dans $\wh{G}$ est donnŽ par 
$$\wh{G}_s = \mathrm{GL}_m(\mathbb{C})\times \cdots \times \mathrm{GL}_m(\mathbb{C})\ptf$$ 
Soit $h\in \widehat{G}=\mathrm{GL}_{md}(\mathbb{C})$ l'ŽlŽment dŽfini par $$h= \left( \begin{array}{ccc}0_{m,m(d-1)} &  & 1_{m,m}\\
1_{m(d-1),m(d-1)} &  & 0_{m(d-1),m} 
\end{array}\right)\ptf$$ 
Pour $g'=(g_1,\ldots , g_d)\in \wh{G}_s$, on a $hg'h^{-1}=(g_d,g_1,\ldots , g_{d-1})$. En particulier 
$$hsh^{-1}= \eta(\zeta^{d-1},1, \zeta,\ldots , \zeta^{d-2})=\zeta^{-1}s \ptf$$
Notons $h_{\phi}$ l'ŽlŽment $(h,\phi)$ de ${^LG}= \wh{G}\times W_\F$. Soit $\ES{G}'$ le sous-groupe fermŽ de ${^LG}$ engendrŽ par 
$\wh{G}_s \times W_{\E}$ et $h_\phi$. 
On a bien $\ES{G}'\cap \wh{G}= \wh{G}_s$ et $s(g,\sigma)s^{-1}=(\bs{a}(\sigma)g,\sigma)$ pour tout $(g,\sigma)\in \ES{G}'$. Le triplet 
$\underline{G}'=(G'=\mathrm{Res}_{\E/\F}(\mathrm{GL}_{m/\E}),\ES{G}'\!,s)$ est une donnŽe endoscopique de $(G,\bs{a})$. Elle est elliptique 
car le groupe $(Z_{\wh{G}'})^{\Gamma_\F}$ co\"{\i}ncide avec le commutant de $h_\phi$ dans $Z_{\wh{G}_s}$, qui est $Z_{\wh{G}}$. 

\begin{remark}\label{remark endos}
\textup{
\begin{enumerate}
\item[(i)]Puisque $h_\phi^d= (1,\phi^d)$ avec $\phi^d\in W_{\E}$, l'application $$\ES{G}' \rightarrow {^LG'},\, (g,w)h_{\phi}^k \mapsto (g, w \phi^k)\,, \quad (g,w)\in \widehat{G}_s\times W_{\E}\,,\; k\in \mathbb{Z}\,,$$ 
est un isomorphisme. 
\item[(ii)]Pour $\gamma \in \Gamma$, Žcrivons $\gamma = \sigma^k$ avec $k\in \{0,\ldots ,k-1\}$ et posons $h(\gamma)= h^k$. L'application $\gamma \mapsto h(\gamma)$ se quotiente 
en un isomorphisme $\Gamma \buildrel\simeq \over{\longrightarrow} \mathrm{Out}(\underline{G}')$.  
En effet si $xsx^{-1} = zs$ pour un $z\in Z_{\wh{G}}= \mathbb{C}^\times$, puisque $s$ et $zs$ ont les mêmes valeurs propres, on a nŽcessairement 
$z= \zeta^k$ pour un $k\in \{0,\ldots ,d-1\}$; par consŽquent $xsx^{-1}= h^{k}sh^{-k}$ et $x \in h^{k} \wh{G}'$. 
\end{enumerate}}
\end{remark}

\begin{lemme}\label{unicitŽ de la donnŽe elliptique}
\textit{\`A Žquivalence prs, le triplet $\underline{G}'=(G'=\mathrm{Res}_{\E/\F}(\mathrm{GL}_{m/\E}),\ES{G'},s)$ construit ci-dessus est l'unique donnŽe endoscopique elliptique de $(G,\bs{a})$.}
\end{lemme}

\begin{demo} Soit $\underline{G}''=(G''\!,\ES{G}''\!, s_*)$ une donnŽe endoscopique 
de $(G,\bs{a})$. \`A Žquivalence prs on peut supposer que $s_*$ appartient ˆ $\wh{T}$. Soit $(g,\phi)\in \ES{G}''$. On a 
$s_*(g,\phi)s_*^{-1}= (\bs{a}(\phi)g,\phi)$ et donc $g s_* g^{-1}= \zeta^{-1}s_*$. Puisque $s_*$ et $\zeta s_*\;(=g^{-1}s_*g)$ ont les mêmes valeurs propres, les valeurs propres distinctes 
de $s_*$ sont de la forme $$\mu_1, \zeta \mu_1,\ldots , \zeta^{d-1}\mu_1, \mu_2,\zeta\mu_2,\ldots , \zeta^{d-1}\mu_2, \ldots , \mu_k,\zeta\mu_k, \ldots , \zeta^{d-1}\mu_k$$ pour des $\mu_i\in \mathbb{C}^\times$ tels que 
$\mu_i \not \in \langle \zeta \rangle \mu_j$ si $i\neq j$. En notant $m_i$ la multiplicitŽ de $\mu_i$ et $\tilde{\mu}_{i}$ l'ŽlŽment $\mu_i 1_{m_i}$ de $Z_{\mathrm{GL}_{m_i}(\mathbb{C})}$, 
on peut, ˆ Žquivalence prs, supposer que $s_*$ est de la forme 
$$s_*= (\tilde{\mu}_1, \zeta\tilde{\mu}_1,\ldots , \zeta^{d-1}\tilde{\mu}_1,\ldots , \tilde{\mu}_k ,\zeta \tilde{\mu}_k,\ldots \zeta^{d-1}\tilde{\mu}_k )\ptf$$ 
Observons que $\sum_{i=1}^k m_i = m \;(= \frac{n}{d})$ et $$\wh{G}_{s_*} = \mathrm{GL}_{m_1}(\mathbb{C})\times \cdots \times \mathrm{GL}_{m_1}(\mathbb{C})\times \cdots \times  \mathrm{GL}_{m_k}(\mathbb{C})\times \cdots \times \mathrm{GL}_{m_k}(\mathbb{C})\ptf$$ 
Pour $i=1,\ldots , k$, notons $h_{m_i}\in \mathrm{GL}_{m_id}(\mathbb{C})$ la matrice dŽfinie en remplaçant 
$m$ par $m_i$ dans la dŽfinition de $h\in \mathrm{GL}_{md}(\mathbb{C})$. L'ŽlŽment $h_* = \mathrm{diag}(h_{m_1},\ldots ,h_{m_k})\in \mathrm{GL}_n(\mathbb{C})$ vŽrifiant 
$h_* s_* h_*^{-1}= \zeta^{-1} s_*$, on a $g^{-1}h_*\in \wh{G}_{s_*}$ et $(h_*,\phi)\in \ES{G}''$. Le groupe 
$(Z_{\wh{G}_{s_*}})^{\Gamma_{\F}}$ est connexe et il co\"{\i}ncide avec le commutant de $h_*$ dans $Z_{\wh{G}_{s_*}}$, qui est le tore $(\mathbb{C}^\times)^k $ naturellement identifiŽ au centre de 
$\mathrm{GL}_{m_1d}(\mathbb{C})\times \cdots \times \mathrm{GL}_{m_k d}(\mathbb{C})$. En particulier la donnŽe $\underline{G}''$ est elliptique si et seulement si $k=1$, 
auquel cas elle est Žquivalente ˆ $\underline{G}'$. \hfill $ \square$
\end{demo}

\begin{remark}
\textup{
Il rŽsulte de la preuve de \ref{unicitŽ de la donnŽe elliptique} qu'une donnŽe endoscopique de $(G,\bs{a})$ est Žquivalente ˆ un produit $\underline{G}'_1\times \cdots \times \underline{G}'_k$ o 
$\underline{G}'_i$ est une donnŽe endoscopique elliptique de $(\mathrm{GL}_{m_id/\F},\bs{a})$ avec $\sum_{i=1}^km_i = m$. \hfill $\blacksquare$
}
\end{remark}

\subsection{DŽcomposition endoscopique de la partie discrte de la formule des traces pour $(\mathrm{GL}_{n/\F}, \K\circ \det$).}\label{dec endos} Continuons avec les notations de \ref{l'ŽgalitŽ} et \ref{endoscopie}. En particulier on a  l'\guill{unique} donnŽe endoscopique elliptique $\underline{G}'=(G'\!,\ES{G}'\!,s)$ de $(G,\bs{a})$ construite en \ref{endoscopie}. 
En \cite[VI.5.1]{MW2} est dŽfini un invariant\footnote{La dŽfinition de $i(G,\underline{G}')$ contient de nombreux autres termes, qui valent tous $1$ ici.} 
$$i(G,\underline{G}')= \vert \mathrm{Out}(\underline{G}')\vert^{-1}\ptf$$ 
Comme d'aprs la remarque \ref{remark endos}\,(ii), on a $\mathrm{Out}(\underline{G}')\simeq \Gamma$, on obtient que 
$$i(G,G')= \vert \Gamma \vert^{-1}= d^{-1}.\leqno{(1)}$$

On fixe un caractre unitaire $\beta$ de $\mathfrak{A}_G$ et on pose $\alpha = \K^{md(d-1)/2}\beta$. 
Comme en \cite[VII.4.1]{MW2}, on munit les groupes $G(\A)= \mathrm{GL}_n(\A)$ et $G'(\A)= \mathrm{GL}_m(\AE)$ des mesures de Tamagawa $\d g$ et $\d h$. 
On fixe une mesure de Haar arbitraire $\d a$ sur le groupe $\mathfrak{A}_G\simeq (\mathbb{R},+)$ et on note $\d \dot{g}$ et $\d \dot{h}$ les mesures quotients sur les groupes 
$\mathfrak{A}_G \backslash G(\A)$ et $\mathfrak{A}_G \backslash G'(\A)$. Cela permet de dŽfinir les opŽrateurs $\rho_{G,\mathrm{disc},\alpha,s}(f^\alpha,\omega)$ et 
$\rho_{G'\!,\mathrm{disc},\beta,s}(f'^\beta)$ apparaissant dans les expressions \ref{la FdT pour G}\,(1) et \ref{l'ŽgalitŽ}\,(1). Rappelons que ces expressions ne dŽpendent pas du choix de $\d a$ (\ref{point central}\,(ii)). 

Soit $(f,\phi)\in \bs{\mathfrak{F}}$, \ie (rappel) $f \in C^{\infty,\bs{K}-\mathrm{fin}}_{\mathrm{c}}(G(\A))$ et 
$\phi\in C^{\infty,\bs{K}'-\mathrm{fin}}_{\mathrm{c}}(G'(\A))$ sont deux fonctions dŽcomposŽes suivant les places de $\F$ et associŽes au sens de \ref{fonctions associŽes}. D'aprs le lemme fondamental (pour toutes les fonctions de l'algbre de Hecke sphŽrique $\mathcal{H}_v$, en toute place $v\in \Vfin$), la fonction $f^{\underline{G}'}=\phi$ est un transfert (global) de $f$. Puisque $\underline{G}'$ est ˆ Žquivalence prs  la seule donnŽe endoscopique elliptique de $(G,\omega)$, l'expression $I_{\mathrm{disc},\alpha}^G(f,\omega)$ est dŽjˆ \guill{stable} et d'aprs \cite[X.8.1]{MW2}, on a l'ŽgalitŽ 
$$I_{\mathrm{disc},\alpha}^G(f,\omega)= i(G,\underline{G}')SI^{\underline{G}'}_{\mathrm{disc},\beta}(f^{\underline{G}'}).\leqno{(2)}$$ 
L'expression $SI^{\underline{G}'}_{\mathrm{disc},\beta}(f^{\underline{G}'})$ est la stabilisation de l'expression $I^{G'}_{\mathrm{disc},\beta}(f^{\underline{G}'})$, mais 
comme ˆ Žquivalence prs la seule donnŽe endoscopique elliptique de $G'$ est la donnŽe triviale $(G'\!, {^LG'},1)$, l'expression $I^{G'}_{\mathrm{disc},\beta}(f^{\underline{G}'})$ est elle 
aussi dŽjˆ stable et donc elle co\"{\i}ncide avec $SI^{\underline{G}'}_{\mathrm{disc},\beta}(f^{\underline{G}'})$. On obtient l'ŽgalitŽ 
cherchŽe $$I_{\mathrm{disc},\alpha}^G(f,\omega)=\frac{1}{d} I_{\mathrm{disc},\beta}^{G'}(\phi).\leqno{(3)}$$ 
Observons que tout comme pour l'expression \ref{la FdT pour G}\,(1), l'ŽgalitŽ (2) n'est pas exactement celle de \cite[X.8.1]{MW2} mais elle s'en dŽduit facilement (cf. \ref{la FdT pour G}). 

On a aussi la variante de (3) avec un caractre automorphe unitaire de $A_G(\A)$. Soit $\mu$ un caractre de $A_G(\A)$ trivial sur $A_G(\F)$ et prolongeant $\beta$; \ie $\mu \in \Xi(G,\beta)$ avec les notations de \ref{variante avec car aut}. Le caractre $\xi = \K^{md(d-1)/2}\mu$ appartient ˆ $\Xi(G,\alpha)$. Pour $\phi \in C^\infty_{\mathrm{c}}(G'(\A))$, on dŽfinit l'expression $I^{G'}_{\mathrm{disc},\mu}(\phi)$ comme en \ref{variante avec car aut}. Alors pour toute paire de fonctions 
$(f,\phi)\in \bs{\mathfrak{F}}$, on a l'ŽgalitŽ 
$$I_{\mathrm{disc},\xi}^G(f,\omega)=\frac{1}{d} I_{\mathrm{disc},\mu}^{G'}(\phi).\leqno{(4)}$$ 
En effet on a $$I_{\mathrm{disc},\xi}^G(f,\omega)= \int_{\bs{X}_{\!A_G}}\xi(a) I_{\mathrm{disc},\alpha}^G(f_a,\omega)\d a$$ 
avec (rappel) $f_a(g)=f(ag)$ ($g\in G(\A)$). Pour $a=(a_v)\in \A^\times$ et $v\in \V$, puisque le facteur de transfert local $\Delta_v$ 
vŽrifie (\ref{facteurs de transfert, propriŽtŽ et choix}\,(i)) $$\Delta_v(a_v\gamma_v)= \K_v^{md(d-1)/2}(a_v) \Delta_v(\gamma_v)\qhq{pour tout}\gamma_v \in \G'_v \cap \G_{v,\mathrm{r\acute{e}g}}\vg$$ 
la paire de fonctions $(f_a, \K^{-md(d-1)/2}(a)\phi_a)$ appartient ˆ $\bs{\mathfrak{F}}$. Par consŽquent 
$$\xi(a)I_{\mathrm{disc},\alpha}^G(f_a,\omega) = \xi(a)\frac{1}{d} I^{G'}_{\mathrm{disc},\beta}(\K^{-md(d-1)/2}(a)\phi_a)= \frac{1}{d}\mu(a)I^{G'}_{\mathrm{disc},\beta}(\phi_a)$$ et 
$$ \int_{\bs{X}_{\!A_G}}\xi(a) I_{\mathrm{disc},\alpha}^G(f_a,\omega)\d a = 
\frac{1}{d} \int_{\bs{X}_{\!A_G}}\mu(a) I_{\mathrm{disc},\beta}^{G'}(\phi_a)\d a = \frac{1}{d} I^{G'}_{\mathrm{disc},\mu}(\phi)\ptf$$

\subsection{Le cas local non ramifiŽ: dŽfinition de l'homomorphisme $b$.}\label{lecaslocalnr} Soient $F$ une extension finie de $\mathbb{Q}_p$, $E/F$ une extension finie non ramifiŽe, de degrŽ $d$, et $n=md$ pour un entier $m\geq 1$. On pose $\Gamma=\Gamma(E/F)$, $G=\mathrm{GL}_{n/F}$ et $G'= \mathrm{Res}_{E/F}(\mathrm{GL}_{m/E})$. On pose aussi $K= \mathrm{GL}_n(\mathfrak{o}_F)$ et $K'= \mathrm{GL}_m(\mathfrak{o}_E)$.   
On note $\mathrm{d}g$ et $\mathrm{d}g'$ les mesures de Haar sur $G(F)$ et $G'(F)$ normalisŽes par 
$\mathrm{vol}(K,\mathrm{d}g)= 1= \mathrm{vol}(K'\!, \mathrm{d}g')$. 

On a une notion de \textit{$\kappa$-relvement faible} pour les reprŽsentations irrŽductibles sphŽriques de $G'(F)= \mathrm{GL}_m(E)$, cf. \ref{def kapparelfaible}. 
On dŽcrit ici explicitement cette opŽration en termes des paramtres de Satake. 

On fixe un caractre $\kappa$ de $F^\times$ de noyau $\mathrm{N}_{E/F}(E^\times)$. Il lui correspond par la thŽorie du corps de classes local un caractre 
$\bs{a}$ du groupe de Weil $W_F$ de $F$, de noyau $W_E$, que l'on voit aussi comme un caractre de $\Gamma$. 
On construit comme en \ref{endoscopie} \guill{la} donnŽe endoscopique elliptique $\underline{G}'=(G'\!,\ES{G}'\!, s)$ de $(G,\bs{a})$. On peut supposer que $\phi\in W_F$ est un ŽlŽment de Frobenius, de sorte que 
$\bs{a}(\phi)=\zeta = \kappa(\varpi)$ o $\varpi$ est une uniformisante de $F$. Rappelons qu'on a notŽ $h_\phi$ l'ŽlŽment $(h,\phi)$ de ${^LG}$, o $h\in \wh{G}=\mathrm{GL}_n(\mathbb{C})$ est dŽfini comme 
en \ref{endoscopie}. Alors $$\ES{G}'= (\wh{G}_s \times I_F)\rtimes \langle h_\phi \rangle$$ o $I_F$ est le sous-groupe d'inertie de $W_F$ 
(toutes les actions se factorisent ici par $W_F \rightarrow W_F/I_F \rightarrow \Gamma$).

Notons $\widehat{\mathcal{H}}_\phi$ l'algbre des fonctions polynomiales sur $\widehat{G}\times \phi$ 
qui sont invariantes par conjugaison par $\widehat{G}$ et 
$\widehat{\mathcal{H}}'_\phi$ l'algbre des fonctions polynomiales sur $\ES{G}'\cap (\widehat{G}\times \phi)$ qui sont 
invariantes par conjugaison par $\widehat{G}'= \widehat{G}_s$. 
On note $$\widehat{b}: \widehat{\mathcal{H}}_\phi \rightarrow \widehat{\mathcal{H}}'_\phi$$ l'homomorphisme de restriction. Notons 
$\mathcal{H}$ et $\mathcal{H}'$ les algbres de Hecke sphŽriques $\mathcal{H}(G(F),K)$ et $\mathcal{H}(G'(F),K')$. Soit 
$b:\mathcal{H} \rightarrow \mathcal{H}'$ l'homomorphisme dŽduit de $\widehat{b}$ via les isomorphismes de Satake  
$\widehat{S}: \mathcal{H}\buildrel\simeq\over{\longrightarrow} \widehat{\mathcal{H}}_\phi$ et $ \widehat{S}':\mathcal{H}'\buildrel\simeq\over{\longrightarrow} \widehat{\mathcal{H}}'_\phi$.

On rappelle la construction de $\wh{S}$ et de $\wh{S}'$. Les classes d'isomorphisme de reprŽsentations irr\'eductibles 
sphŽriques de $G(F)$ sont paramŽtrŽes par les clas\-ses de conjugaison semi-simples dans $\widehat{G}$ ou, ce qui revient au même, par les classes de $\mathfrak{S}_n$-conjugaison dans $(\mathbb{C}^\times)^n$, cf. \ref{param satake}. Pour $y\in (\mathbb{C}^\times)^n$, on dŽfinit $\pi_y^*$ et $\pi_y$ comme en \textit{loc.~cit.}; la reprŽsentation $\pi_y$ de $G(F)$ est irrŽductible et sphŽrique, et sa classe d'isomorphisme ne dŽpend que de la classe de $\mathfrak{S}_n$-conjugaison de $y$. On peut donc dŽfinir $\pi_y$ ˆ isomorphisme prs pour tout ŽlŽment $y\in \widehat{G}$ semi-simple et poser, pour toute fonction 
$f\in \mathcal{H}$: $$\widehat{S}f(y,\phi)  = {\rm trace}(\pi_y(f))\leqno{(1)}$$ 
o l'opŽrateur $\pi_y(f)$ est dŽfini ˆ l'aide de la mesure $\d g$. 
Observons que le $\phi$ ne sert ˆ rien ici, on pourrait le supprimer; il est commode pour dŽfinir $\widehat{b}$ par restriction. 
Deux ŽlŽments de $\wh{G}$ sont dits \textit{ss-conjuguŽs} si leurs parties semi-simples (pour la dŽcomposition de Jordan) sont conjuguŽes dans $\wh{G}$. 
Pour $f\in \mathcal{H}$, l'application $y \mapsto \widehat{S}f(y,\phi)$ --- d\'efinie pour $y\in\widehat{G}$ semi-simple --- s'Žtend de manire unique en une fonction sur 
$\wh{G}$ constante sur les classes de ss-conjugaison, qui est 
un ŽlŽment de $\wh{\mathcal{H}}_\phi$. L'application $\widehat{S}: \mathcal{H} \rightarrow \widehat{\mathcal{H}}_\phi$ ainsi dŽfinie est un isomor\-phisme d'algbres. 

Pour $G'$, la construction de $\wh{S}'$ est l\'eg\`erement plus compliquŽe 
car le groupe $G'$ n'est pas dŽployŽ sur $F$. Pour $y\in \mathrm{GL}_m(\mathbb{C})$ semi-simple, on note $\Pi_y$ la classe d'isomorphisme de reprŽsentations irrŽductibles sphŽriques de $\mathrm{GL}_m(E)$ associ\'ee \`a 
la classe de conjugaison de $y$ dans $\mathrm{GL}_m(\mathbb{C})$. Pour 
un ŽlŽment $\bs{x}=(x,1)h_\phi\in \ES{G}'\cap (\widehat{G}\times \phi)$ avec $x=(x_1,\dots , x_d)\in \widehat{G}'=\widehat{G}_s$, on pose 
$$N(\bs{x})=x_1x_2\cdots x_d\in \GLm(\mathbb{C})\ptf$$
Observons que pour $g=(g_1,\ldots ,g_d)\in \widehat{G}'$, on a $N(g\bs{x}g^{-1})= g_1N(\bs{x})g_1^{-1}$. 
Pour $f'\in \mathcal{H}'$ et $\bs{x}\in \ES{G}'\cap (\widehat{G}\times \phi)$ tel que $N(\bs{x})$ soit semi-simple, on pose: 
$$\widehat{S}' f' (\bs{x})= {\rm trace}(\Pi_{N(\bs{x})}(f'))\leqno{(2)}$$ o l'opŽrateur $\Pi_{N(\bs{x})}(f')$ est dŽfini ˆ l'aide de la mesure $\d g'$. 
Deux \'el\'ement $\bs{x}$ et $\bs{x}'$ de $\ES{G}'\cap (\widehat{G}\times \phi)$ sont dit 
\textit{$N$-ss-conjugu\'es} si les \'el\'ements $N(\bs{x})$ et $N(\bs{x}')$ sont ss-conjugu\'es dans $\mathrm{GL}_m(\mathbb{C})$. 
L'application $\bs{x}\mapsto \widehat{S}'f'(\bs{x})$ --- d\'efinie pour $\bs{x}\in \ES{G}'\cap (\widehat{G}\times \phi)$ tel que $N(\bs{x})$ soit semi-simple ---   
s'Žtend de manire unique 
en une fonction sur $\ES{G}'\cap (\widehat{G}\times \phi)$ constante sur les classes de $N$-ss-conjugaison, qui par construction est un \'el\'ement de $\widehat{\mathcal{H}}'_\phi$. L'application $\widehat{S}': \mathcal{H}' \rightarrow \widehat{\mathcal{H}}'_\phi$ ainsi dŽfinie est  
un isomorphisme d'algbres.

\begin{remark}
\textup{
L'homomorphisme d'algbres $b= \widehat{S}'^{-1}\circ \widehat{b} \circ \widehat{S}:\mathcal{H} \rightarrow \mathcal{H}'$ co\"{\i}ncide avec celui dŽfini en 
\cite[II.3]{W1}. \hfill $\blacksquare$}
\end{remark}

\vskip1mm
Soit $y=(y_1,\ldots , y_m) \in (\mathbb{C}^\times)^m$. On choisit un ŽlŽment $t=(t_1,\ldots , t_m)\in (\mathbb{C}^\times)^m$ tel que 
$t^d=y$, \ie tel que $t_i^d = y_i$ pour $i=1,\ldots ,m$. Soit $x\in (\mathbb{C}^\times)^{n}$ l'ŽlŽment $(t,\ldots ,t)$, o $t$ est rŽpŽtŽ $d$ fois. 
Posons $$\bs{x}= (x,1)h_\phi=(xh,\phi)\in \ES{G}'\cap (\widehat{G}\times \phi)\ptf$$ 
On a $N(\bs{x})= t^d=y$ et par dŽfinition, pour toute fonction $f\in \mathcal{H}$, 
$$\widehat{S}f (\bs{x})=  \widehat{b}(\widehat{S}f)(\bs{x})= \widehat{S}' (bf) (\bs{x}) = {\rm trace}(\Pi_y(bf))\,\ptf$$ 
Or $xh$ est conjuguŽ dans $\widehat{G}$ ˆ la matrice diagonale 
$$\delta (y)= (t, \zeta t, \ldots , \zeta^{d-1} t)\in (\mathbb{C}^\times)^n\,\ptf$$ 
On a donc aussi $$\widehat{S}f(\bs{x})= \widehat{S}f(\delta(y),\phi)= {\rm trace}(\pi_{\delta(y)}(f))\,\ptf$$ 
D'o l'ŽgalitŽ, pour toute fonction $f\in \mathcal{H}$:
$${\rm trace}(\pi_{\delta(y)}(f))= {\rm trace}(\Pi_{y}(bf ))\ptf\leqno{(3)}$$ 

\section{Sur l'opŽrateur d'entrelacement $I_\K$}\label{sur l'op d'entrel IK} 

Dans cette section, $\E/\F$ est une extension finie cyclique de corps de nombres, de degrŽ $d$. On fixe un caractre 
$\K$ de $\A^\times$ de noyau $\F^\times \mathrm{N}_{\E/\F}(\AE^\times)$. On fixe aussi un entier $m\geq 1$. On pose $G=\mathrm{GL}_{n/\F}$ avec $n=md$.

\subsection{OpŽrateur global issu du local.}\label{op issu du local} On reprend les notations de \ref{l'ŽgalitŽ} ˆ la seule diffŽrence suivante prs: 
le caractre $\omega=\K\circ \det$ de $G(\A)$ est ici notŽ $\K$. 
On fixe un caractre additif non trivial $\psi$ de $\A$ trivial sur $\F$. Il dŽfinit pour chaque place $v$ de $\F$ un caractre additif non trivial $\psi_v$ de $\F_v$. 

Soit $\pi$ une sŽrie discrte $\K$-stable de $G(\A)$. Pour chaque place $v$ de $\F$, sa composante locale $\pi_v$ est $\K_v$-stable et on dispose d'un opŽrateur d'entrelacement normalisŽ $A_{\pi_v}=A_{\pi_v,\psi_v}\in \mathrm{Isom}_{{\bf G}_v}(\K_v\pi_v,\pi_v)$. Rappellons que pour $v\in \EuScript{V}_\infty$, cet opŽrateur est obtenu par extension de l'opŽrateur d'entrelacement normalisŽ $A_{\pi_v^\infty}=A_{\pi_v^\infty,\psi_v}\in \mathrm{Isom}_{{\bf G}_v}(\K_v\pi_v^\infty,\pi_v^\infty)$; o 
$\pi_v^\infty$ est la reprŽsentation de $\G_v$ sur le sous-espace $V_{\pi_v}^\infty \subset V_{\pi_v}$ formŽ des vecteurs lisses. Fixons un isomorphisme $\iota$ du produit tensoriel restreint complŽtŽ $\widehat{\otimes}_v \pi_v$ sur $\pi$; o le \guill{restreint} est par rapport au choix pour presque toute place $v\in \Vfin$ d'un vecteur non nul ${\bf K}_v$-invariant $x_v\in V_{\pi_v}$. Pour chaque place $v$ de $\F$, l'opŽrateur $A_{\pi_v}$ est unitaire (\ref{op unit}\,(ii)); et pour presque toute place $v\in \Vfin$, il fixe le vecteur $x_v$ (\ref{divers, LF et op norm}\,(ii)). L'opŽrateur $\otimes_v A_{\pi_v^\infty}: \otimes_v \K_v\pi_v^\infty \rightarrow  \otimes_v \pi_v^\infty$ s'Žtend de manire unique en un opŽrateur unitaire
$\widehat{\otimes}_v A_{\pi_v} : \widehat{\otimes}_v \K_v\pi_v \rightarrow \widehat{\otimes}_v \pi_v$. On note $A_\pi = A_{\pi,\psi}\in \mathrm{Isom}_{G(\A)}(\K \pi, \pi)$ l'opŽrateur global unitaire dŽfini par$$A_{\pi}= \iota \circ \left( \widehat{\otimes}_v A_{\pi_v}\right)\circ \iota^{-1} \ptf$$ 

On note $A_{\pi^\sharp}= A_{\pi^\sharp,\psi}\in \mathrm{Isom}_{G(\A)}(\K \pi^\sharp, \pi^\sharp)$ l'opŽrateur global dŽduit de $A_\pi$ par restriction; o (rappel) $\pi^\sharp$ est la reprŽsentation automorphe de $G(\A)$ sous-jacente ˆ $\pi$, cf. \ref{rep aut disc}. 

On Žtend naturellement cette dŽfinition aux reprŽsentations automorphes de $G(\A)$ qui sont induites de discrtes. Si $\pi$ est une reprŽsentation automorphe de $G(\A)$ de la forme $\pi_1\times \cdots \times \pi_r$ pour des reprŽsentations automorphes discrtes $\pi_i$ de $\mathrm{GL}_{n_i}(\A)$ avec $\sum_{i=1}^rn_i=n$, on note $A_\pi$ l'opŽrateur global dŽfini comme dans le cas local par $A_\pi= A_{\pi_1}\times \cdots \times A_{\pi_r}$.

\subsection{Rappels sur le spectre rŽsiduel.}\label{rappels sur le spectre res} La description du spectre rŽsiduel de $G(\A)$ est due ˆ M\oe glin-Waldspurger \cite{MW1}. 
On note $\nu = \nu_{\A}$ le caractre de $\A^\times$ trivial sur $\F^\times$ dŽfini par $\nu(x)= \vert x \vert_{\A}=\prod_v \vert x_v\vert_{\F_v}$ pour $x=\prod_v x_v \in \A^\times$. 
Soient $k\geq 1$ un entier divisant $n$ et $\rho$ une reprŽsentation automorphe cuspidale unitaire de $\mathrm{GL}_{\frac{n}{k}}(\A)$. On commence par former la reprŽsentation automorphe cuspidale 
$$\nu^{\frac{k-1}{2}}\rho \otimes \nu^{\frac{k-3}{2}}\rho\otimes\cdots \otimes 
\nu^{\frac{1-k}{2}}\rho$$ du facteur de Levi standard $L_k(\A)= (\mathrm{GL}_{\frac{n}{k}}(\A))^k$ de $G(\A)$. On note 
$$R(\rho,k)=\nu^{\frac{k-1}{2}}\rho \times \nu^{\frac{k-3}{2}}\rho\times\cdots \times 
\nu^{\frac{1-k}{2}}\rho$$ son induite parabolique suivant le sous-groupe parabolique standard $P_k(\A)$ de $G(\A)$ de composante de Levi $L_k(\A)$. D'aprs \cite{MW1}, $R(\rho,k)$ a un unique quotient irrŽductible\footnote{On renvoie ˆ Langlands \cite{La2} pour les notions d'induite parabolique et de quotient d'une reprŽsentation automorphe.}, notŽ $u(\rho,k)$, qui est une reprŽsentation automorphe discrte (rŽsiduelle si $k>1$) de $G(\A)$. De plus (\textit{loc.\,cit.}) ˆ isomorphisme prs toute reprŽsentation automorphe discrte de $G(\A)$ est obtenue de cette manire, et deux telles reprŽsentations automorphes discrtes 
$u(\rho,k)$ et $u(\rho'\!,k')$ de $G(\A)$ sont isomorphes si et seulement si $k=k'$ et $\rho\simeq \rho'$. En particulier $u(\rho,k)$ est $\K$-stable si et seulement si $\rho$ est $\K$-stable. 

On dŽduit de cette classification le thŽorme de multiplicitŽ $1$ fort pour les reprŽsentations automorphes discrtes de $G(\A)$, hŽritŽ du mme thŽorme pour les reprŽsentations automorphes cuspidales unitaires (cf. \cite[4]{BH}): 

\begin{proposition}
\begin{enumerate}
\item[(i)] Soient $\pi=u(\rho,k)$ et $\pi'=u(\rho'\!,k')$ deux reprŽsentations automorphes discrtes de $G(\A)$ telles que $\pi_v\simeq \pi'_v$ pour presque toute place $v\in \V$. Alors $\pi \simeq \pi'$, \ie $k=k'$ et $\rho\simeq \rho'$. 
\item[(ii)] Soit $\pi=u(\rho,k)$ une reprŽsentation automorphe discrte de $G(\A)$. Alors $\delta$ appara"t avec multiplicitŽ $1$ dans $L^2_{\mathrm{disc}}(G(\F)\backslash G(\A),\alpha)$ o $\alpha$ est la restriction ˆ $\mathfrak{A}_G$ du caractre central de $\pi$. 
\end{enumerate}
\end{proposition} 

Soient $k\vert n$ et $\rho$ une reprŽsentation automorphe cuspidale unitaire de $L_k(\A)$. L'espace de $u(\rho,k)$ est obtenu comme rŽsidu de sŽries d'Eisenstein. Rappelons brivement la construction de Jacquet \cite{J}. 

On rappelle que l'homomorphisme ${\bf H}_{L_k}: L_k(\A) \rightarrow \mathfrak{a}_{L_k}$ se restreint en un isomomorphisme de $\mathfrak{A}_{L_k}\;(\simeq (\mathbb{R}_{>0})^k)$ sur $\mathfrak{a}_{L_k}$. Supposons dans un premier temps 
que le caractre central $\omega_{\rho_{L_k}}$ de $\rho_{L_k}= \rho\otimes \cdots \otimes \rho$ ($k$ fois) soit trivial sur $\mathfrak{A}_{L_k}$, ce qui Žquivaut ˆ supposer que le caractre central $\omega_{\rho}$ est trivial sur $\mathfrak{A}_{\mathrm{GL}_{a/ \F}}$, $a=\frac{n}{k}$. 
Notons $V= V_\rho\otimes\cdots \otimes V_\rho$ l'espace de $\rho_{L_k}$. 
Pour $s=(s_1,\ldots ,s_k)\in \mathbb{C}^k$, on note $\rho[s]$ la reprŽsentation automorphe cuspidale 
$\nu^{s_1}\rho\otimes \cdots \otimes \nu^{s_k}\rho$ de $L_k(\A)$ et on pose 
$\pi_{\rho[s]}=\nu^{s_1}\rho\times \cdots \times \nu^{s_k}\rho$ (induite parabolique de $\rho[s]$ ˆ $G(\A)$ suivant $P_k(\A)$). Ainsi 
$R(\rho,k)=\pi_{\rho[e]}$ avec $e= (\frac{k-1}{2}, \frac{k-3}{2},\ldots , \frac{1-k}{2})$. 
Pour dŽcrire l'unique quotient irrŽductible $u(\rho,k)$ de $R(\rho,k)$, on rŽalise toutes les induites paraboliques $\pi_{\rho[s]}$ ˆ l'aide d'un mme espace de fonctions, ˆ savoir l'espace $\ES{F}$ formŽ des fonctions $\Phi: G(\A) \rightarrow V$ telles que: 
\begin{itemize}
\item $\Phi(u\gamma a g) = \Phi(g)$ pour tous $u\in U_{P_k}(\A)$, $\gamma \in P_k(\F)$, $a\in \mathfrak{A}_{L_k}$ et $g\in G(\A)$; 
\item $\Phi$ est $\bs{K}$-finie et $\Phi(lk) = \rho_{L_k}(l)(\Phi(k))$ pour tous $l\in L_k(\A)$ et $k\in \bs{K}$. 
\end{itemize}
Pour $s\in \mathbb{C}^k$, notons $\ES{F}_s$ l'espace des fonctions sur $G(\A)$ de la forme 
$$g\mapsto e^{\langle s + \bs{\rho}_{P_k},{\bf H}_{P_k}(g) \rangle }\Phi(g)\vgq \Phi\in \ES{F}\vg$$ 
o $\bs{\rho}_{P_k}$ (ˆ ne pas confondre avec la reprŽsentation $\rho_{L_k}$!) est la demi-somme des racines de $A_{L_k}$ dans $U_{P_k}$. Ici $s$ et $\bs{\rho}_{P_k}$ sont naturellement identifiŽs ˆ des ŽlŽments du dual $\mathfrak{a}_{L_k,\mathbb{C}}^*$ de $\mathfrak{a}_{L_k,\mathbb{C}}= \mathrm{Hom}(X(A_{L_k}),\mathbb{C}^\times)$.
La reprŽsentation rŽgulire droite de $G(\A)$ sur $\ES{F}_s$ est par dŽfinition l'induite parabolique $\pi_{\rho[s]}$. 
Pour $\Phi\in \ES{F}$ et $s\in \mathbb{C}^k$, on forme la sŽrie d'Eisenstein 
$$E(g,s,\Phi) = \sum_{\gamma\in P_k(\F)\backslash G(\F)} e^{\langle s + \bs{\rho}_{P_k}, {\bf H}_{P_k}(\gamma g) \rangle} \Phi(\gamma g)\vgq g\in G(\A)\pvg$$  
La sŽrie $E(g,s,\Phi)$ converge absolument pour $s$ dans un c™ne ouvert de $\mathfrak{a}_{L_k,\mathbb{C}}^*$ et elle admet un prolongement mŽromorphe ˆ l'espace $\mathfrak{a}_{L_k,\mathbb{C}}^*$ tout entier. Le point est qu'elle admet un p™le en $s=e$. Pour attraper le rŽsidu en $s=e$, Jacquet multiplie $E(s,g,\Phi)$ par une fonction mŽromorphe $s\mapsto u(s)$  indŽpendante de $\Phi$ (dont la valeur exacte n'a pas d'importance ici) telle que la fonction $s\mapsto u(s)E(s,g,\Phi)$ soit holomorphe en $s=e$ et que les fonctions $g\mapsto u(s)E(g,s,\Phi)$ pour $\Phi\in \ES{F}$, ŽvaluŽes en $s=e$, engendrent l'espace de $u(\rho,k)$. La construction fournit un opŽrateur de $G(\A)$-entrelacement non nul $$\mathrm{Res}= \mathrm{Res}_{\rho,k}: \ES{F}_e=V_{R(\rho,k)} \rightarrow V_{u(\rho,k)}\ptf$$ 

Si le caractre central $\omega_\rho$ de $\rho$ n'est pas trivial sur $\mathfrak{A}_{\mathrm{GL}_{a/\F}}$, il existe un $t\in \mathbb{C}$ (en fait $t\in \mathbb{U}$ puisque $\rho$ est unitaire) tel que le caractre central de $\rho' = \nu^{-t}\rho$ soit trivial sur $\mathfrak{A}_{\mathrm{GL}_{a/\F}}$. On a alors $R(k,\rho)= \nu^t R(k,\rho')$ et $u(k,\rho)= \nu^t u(k,\rho')$, et la construction de $u(k,\rho)$ se dŽduit aisŽment de la construction de $u(k,\rho')$. En particulier, 
la reprŽsentation $u(k,\rho)$ se rŽalise dans l'espace des fonctions $$g \mapsto (\nu^t\Phi')(g)=\vert \det(g)\vert_{\mathbb{A}}^t\Phi'(g)\vgq \Phi'\in V_{\delta(\rho'\!,k)}\vg$$ et l'opŽrateur de $G(\A)$-entrelacement $\mathrm{Res}'=\mathrm{Res}_{\rho'\!,k}:  V_{R(\rho'\!,k)}\rightarrow V_{u(\rho'\!,k)}$ fourni par la construction de Jacquet dŽfinit un opŽrateur de $G(\A)$-entrelacement 
$$\mathrm{Res}= \mathrm{Res}_{\rho,k}: V_{R(\rho,k)}\rightarrow V_{u(\rho,k)}\vgq  
\Phi \mapsto \nu^t\mathrm{Res}'(\nu^{-t}\Phi)\ptf$$

\subsection{L'ŽgalitŽ $A_\pi= I_{\K,\pi}$.} Si $\alpha$ est un caractre unitaire de $\mathfrak{A}_G$, on dŽfinit comme en \ref{la FdT pour G} l'opŽrateur d'entrelacement \guill{physique} $I_\K=I_{\K\circ \det}$ sur l'espace $L^2(G(F)\backslash G(\A),\alpha)$. Si $\pi$ est 
une sŽrie discrte de $G(\A)$ de caractre central prolongeant $\alpha$, puisque $\pi$ appara"t avec multiplicitŽ $1$ dans 
le spectre de $L^2_{\mathrm{disc}}(G(\F)\backslash G(\A),\alpha)$, l'opŽrateur physique $I_\K$ laisse invariant 
l'espace $V_\pi$ de $\pi$. On note $I_{\K,\pi}$ la restriction de $I_{\K}$ ˆ $V_\pi$.

\begin{lemme}\label{A=I}
Soit $\pi$ une sŽrie discrte $\K$-stable de $G(\A)$. Alors $A_\pi= I_{\K,\pi}$. 
\end{lemme}

\begin{demo} 
D'aprs \ref{rappels sur le spectre res}, la reprŽsentation automorphe $\pi^\sharp$ sous-jacente ˆ $\pi$ s'Žcrit 
$\pi^\sharp = u(\rho^\sharp\!,k)$ o $k\geq 1$ est un entier divisant $n$ et $\rho$ est une sŽrie discrte cuspidale de $G(\A)$. 
Fixons un isomorphisme $\iota$ du produit tensoriel restreint complŽtŽ $\widehat{\otimes}_v \rho_v$ sur $\rho$. 

Supposons tout d'abord que $\pi$ soit cuspidale, \ie $k=1$ et $\pi=\rho$. On a sur l'espace $V^\infty \;(\subset V= V_\pi)$ de la reprŽsentation lisse $\pi^\infty = \iota(\otimes_v\pi^\infty_v)$ sous-jacente ˆ $\pi$ une fonctionnelle de Whittaker globale $\Lambda$ donnŽe par 
$$\Lambda (\Phi )= \int_{U_0(\F)\backslash U_0(\A)} \Phi(u)\overline{\theta_\psi(u)} \d u\vgq \Phi\in V^\infty \pvg$$ o $\theta_\psi$ est le caractre de $U_0(\A)$ dŽfini comme en local (cf. \ref{norm op entrelac}) et $\d u$ est une mesure de Haar sur $U_0(\A)$. Puisque $\pi$ est cuspidale, ses composantes locales $\pi_v^\infty$ sont gŽnŽriques. On peut donc Žcrire $\Lambda\circ \iota$ comme un produit tensoriel $\otimes_v \Lambda_v$ o, pour chaque place $v$ de $\F$, $\Lambda_v$ est une fonctionnelle de Whittaker pour $\pi_v$ relativement ˆ $\psi_v$ et l'opŽrateur d'entrelacement normalisŽ $A_{\pi_v^\infty}= A_{\pi_v^\infty , \psi_v}$ est caractŽrisŽ par l'ŽgalitŽ $\Lambda_v \circ A_{\pi_v^\infty} = \Lambda_v$. On a donc $\Lambda\circ \iota (\otimes_v A_{\pi_v^\infty})= \Lambda \circ \iota$, \ie 
$$\Lambda \circ A_{\pi^\infty} = \Lambda \qhq{avec} A_{\pi^\infty}= \iota \circ (\otimes_v A_{\pi_v^\infty})\circ \iota^{-1}\ptf$$ Le caractre $\K=\K \circ \det$ de $G(\A)$ Žtant trivial sur $U_0(\A)$, on a aussi $$\Lambda \circ I_\K= \Lambda\ptf$$ Par consŽquent $A_{\pi^\infty}=I_\K\vert_{V_{\pi^\infty}}$ ce qui par densitŽ entra"ne l'ŽgalitŽ 
cherchŽe: $A_\pi = I_{\K,\pi}$. 

Supposons maintenant que $\pi$ soit rŽsiduelle, \ie $k>1$. Supposons aussi dans un premier temps que le caractre central $\omega_{\rho^\sharp}$ de $\rho^\sharp$ soit trivial sur $\mathfrak{A}_{\mathrm{GL}_{a/\F}}$, $a= \frac{n}{k}$. Reprenons la construction de Jacquet \cite{J} (cf. \ref{rappels sur le spectre res}). 
On note $V=V_{\rho^\sharp} \otimes \cdots \otimes V_{\rho^\sharp}$ l'espace de la reprŽsentation automorphe cuspidale $(\rho^\sharp)_{L_k}= \rho^\sharp\otimes \cdots \otimes \rho^\sharp$ de $L_k(\A)$. Soit $R=R(\rho^\sharp,k)$ la reprŽsentation automorphe de $G(\A)$ dŽfinie par $R= \nu^{\frac{k-1}{2}}\rho^\sharp \times \nu^{\frac{k-3}{2}}\rho^\sharp\times \cdots \times \nu^{\frac{1-k}{2}}\rho^\sharp$. D'aprs Jacquet (\textit{loc.\,cit.}), pour $v\in \ES{V}$, la composante locale $\pi^\sharp_v$ de $\pi^\sharp=u(\rho^\sharp,k)$ en $v$ est le quotient de Langlands $\overline{R}_v$ de la composante locale $R_v=\nu_v^{\frac{k-1}{2}}\rho^\sharp_v \times \cdots \times \nu_v^{\frac{1-k}{2}}\rho_v^\sharp$ de $R$ en $v$, \ie la reprŽsentation notŽe $u(\rho^\sharp_v,k)$ en \ref{classif (rappels)}. PrŽcisons cette affirmation, qui nous permettra de dŽcrire l'action de l'opŽrateur physique $I_{\K}$ sur l'espace de $\pi$. L'espace $V_{R_v}$ de $R_v$ est formŽ des fonctions $f_v: G(\F_v)\rightarrow V_{\smash{\rho^\sharp_v}} \otimes \cdots \otimes V_{\smash{\rho^\sharp_v}}$ ($k$ fois) qui sont $\bs{K}_{\!v}$-finies (ˆ droites) et 
telles que 
pour tous $l\in L_k(\F_v)$, $u\in U_{P_k}(\F_v)$ et $g\in G(\F_v)$, on ait $$f_v(lug)= \bs{\delta}_{P_k}(l)^{\frac{1}{2}}(\rho_v^\sharp[e](l))f(g)\ptf$$ 
Ë l'aide de l'isomorphisme $\iota:\widehat{\otimes}_v \rho_v\rightarrow \rho$, on dŽfinit un isomorphisme $\mathfrak{I}$ du 
produit tensoriel restreint $\otimes_vR_v$ sur $R$: ˆ une fonction $f= \otimes_v f_v$ dans l'espace $\otimes_v V_{R_v}$, on associe 
une fonction $\Phi_f: G(\A) \rightarrow V$ en posant $$\Phi_f(g) = 
(\iota \otimes \cdots \otimes \iota)(\otimes_v f_v(g_v))\vgq g= (g_v)\in G(\A)\ptf$$ 
Par construction $\Phi_f$ appartient ˆ l'espace notŽ $\ES{F}_e$ en \ref{rappels sur le spectre res}, \ie l'espace de $R$, et l'application $\mathfrak{I}:f \mapsto \Phi_f$ dŽfinit un isomorphisme de $\otimes_v R_v $ sur $R$. En composant $\mathfrak{I}$ avec le morphisme (surjectif) $\mathrm{Res}:R\rightarrow \pi^\sharp$, on obtient un morphisme surjectif de $\otimes_v R_v$ sur $\pi^\sharp$. Cela entra"ne qu'en toute place $v$ de $\F$, la composante locale $\pi^\sharp_v$ de $\pi^\sharp$ ne peut tre que le quotient de Langlands $\overline{R}_v$ de $R_v$ et que le morphisme surjectif $\mathrm{Res}\circ \mathfrak{I}:f \mapsto \mathrm{Res}(\Phi_f)$ se factorise en un isomorphisme $\otimes_v \overline{R}_v\rightarrow \pi^\sharp$. Puisque l'espace $V_{\rho^\sharp}$ de $\rho^\sharp$ est formŽ de fonctions sur $\mathrm{GL}_a(\A)$, l'espace $V$ de $(\rho^\sharp)_{L_k}$ est formŽ de fonctions sur $L_k(\A)$ que l'on peut Žvaluer en $1$ (l'ŽlŽment unitŽ de $L_k(\A)$). Pour $f = \otimes_vf_v\in \otimes_v V_{R_v}$, on note $\varphi_f: G(\A)\rightarrow \mathbb{C}$ la fonction dŽfinie par 
$\varphi_f(g) = \Phi_f(g)(1)$. Observons que si $f$ se dŽcompose en $f_1\otimes \cdots \otimes f_k$ avec $f_i=\otimes_vf_{i,v}: G(\A)\rightarrow \otimes_v V_{\smash{\rho^\sharp_v}}$, on a $\varphi_f(g)= \prod_{i=1}^k \iota(f_i(g))(1)$.

Puisque $\pi$ est $\kappa$-stable, $\rho$ l'est aussi (d'aprs \cite{MW1}). L'opŽrateur physique $I_\K$ stabilise l'espace $\ES{F}_e$ de $R$ (et bien sžr aussi l'espace de $\pi^\sharp$); 
et pour $\Phi\in \ES{F}_e$, on a $$\mathrm{Res}(I_\K\Phi (g))=I_\K \mathrm{Res}(\Phi)\ptf$$ 

Pour $v\in \ES{V}$, notons $B_v\in \mathrm{Isom}_{G(\F_v)}(\K_v\rho_v^\sharp,\rho_v^\sharp)$ l'opŽrateur d'entrelacement normalisŽ $A_{\smash{\rho^\sharp_v}}=A_{\smash{\rho^\sharp_v,\psi_v}}$. Soit 
$A_{R_v}= B_v \times\cdots \times B_v  \in \mathrm{Isom}_{G(\F_v)}(\K_vR_v,R_v)$ l'opŽrateur obtenu par induction parabolique ˆ partir de l'opŽrateur
$$B_v \otimes \cdots \otimes B_v\in \mathrm{Isom}_{L_k(\F_v)} (\K_v \rho^\sharp_v[e], \rho^\sharp_v[e])\ptf$$ L'opŽrateur d'entrelacement normalisŽ $A_{\smash{\pi^\sharp_v}}= A_{\smash{\pi^\sharp_v,\psi_v}}\in \mathrm{Isom}_{G(\F_v)}(\K_v\pi^\sharp_v,\pi^\sharp_v)$ est celui dŽduit de $A_{R_v}$ par passage au quotient. Puisque l'opŽrateur physique $I_\K$ commute ˆ l'opŽrateur $\mathrm{Res}: R \rightarrow \pi^\sharp$, pour montrer que l'opŽrateur global $A_{\pi^\sharp}= \otimes_v A_{\smash{\pi^\sharp_v}}$ co\"{\i}ncide avec $I_{\K,\pi^\sharp}= I_{\K}\vert_{V_{\pi^\sharp}}$, il suffit de montrer que l'opŽrateur global $\otimes_v A_{R_v}$ co\"{\i}ncide avec $I_{\K,R}=I_\K\vert_{V_R}$, \cad que pour toute fonction $f= \otimes_v f_v$ dans l'espace $\otimes_v V_{R_v}$, on a l'ŽgalitŽ 
$$\Phi_{\otimes_v A_{R_v}f_v}(g)= \K(g)\Phi_f(g)\qhq{pour tout} g\in G(\A)\pvg$$ ou, ce qui revient au mme, on a l'ŽgalitŽ (Žvaluation en $1$)
$$\varphi_{\otimes_v A_{R_v}f_v}(g)= \K(g)\varphi_f(g)\qhq{pour tout} g\in G(\A)\ptf$$
Pour $f_v \in V_{R_v}$ et $g_v\in G(\F_v)$, on a par dŽfinition 
$$(A_{R_v}f_v)(g)= \K_v(g) (B_v\otimes\cdots \otimes B_v)(f_v(g))\ptf$$  
Par suite pour $f = \otimes_vf_v$ dans l'espace $\otimes_v R_v$ et $g=(g_v)\in G(\A)$, on a  
$$\Phi_{\otimes_v A_{R_v}f_v}(g)=\K(g)(\iota\otimes \cdots\otimes  \iota )((\otimes_vB_v)\otimes \cdots \otimes (\otimes_vB_v))(\otimes_v f_v(g_v))\pvg$$ et 
puisque $\iota \circ (\otimes_v B_v)= I_{\K,\rho^\sharp}\circ \iota$ (d'aprs le rŽsultat montrŽ prŽcŽdemment pour les sŽries discrtes cuspidales), on obtient 
$$\Phi_{\otimes_v A_{R_v}f_v}(g)= \K(g) (I_{\K,\rho^\sharp}\otimes \cdots \otimes I_{\K,\rho^\sharp})(\Phi_{f}(g))\ptf$$ 
Il ne reste plus qu'ˆ Žvaluer en $1$: 
$$\varphi_{\otimes_vA_{R_v}f_v}(g)= \K(g) (I_{\K,\rho^\sharp}\otimes \cdots \otimes I_{\K,\rho^\sharp})(\Phi_{f}(g))(1)\ptf$$ Comme 
$(I_{\K,\rho^\sharp} f^\sharp)(1)=f^\sharp(1)$ pour toute fonction $f^\sharp\in V_{\rho^\sharp}$, on obtient l'ŽgalitŽ cherchŽe: 
$\varphi_{\otimes_vA_{R_v}\phi_v}(g)= \K(v) \varphi_{f}(g)$ pour tout $g\in G(\A)$. On a donc prouvŽ l'ŽgalitŽ $A_{\pi^\sharp}= I_{\K,\pi^\sharp}$. Par densitŽ on obtient l'ŽgalitŽ $A_\pi= I_{\K,\pi}$. 

Reste ˆ supprimer l'hypothse que le caractre central de $\rho^\sharp$ est trivial sur $\mathfrak{A}_{\mathrm{GL}_{a/\F}}$. On procde pour cela comme dans le dernier paragraphe de \ref{rappels sur le spectre res}: il existe un $t\in \mathbb{U}$ tel que le caractre central de $\nu^{-t}\rho^\sharp$ soit trivial sur $\mathfrak{A}_{\mathrm{GL}_{a/\F}}$. Alors $\pi^\sharp = \nu^t u(\nu^{-t}\rho^\sharp,k)$ et d'aprs le rŽsultat prouvŽ ci-dessus pour la sŽrie discrte $\pi'$ dont $u(\nu^{-1}\rho^\sharp,k)$ est la reprŽsentation automorphe sous-jacente, on obtient le rŽsultat pour $\pi$. \hfill $\square$
\end{demo}

\section{DŽmonstration des rŽsultats de la section \ref{ŽnoncŽ des rŽsultats}}\label{demo}

Continuons avec les hypothses et les notations de la sous-section \ref{l'ŽgalitŽ}; avec $G= \mathrm{GL}_{n/\F}$, $G'= \mathrm{Res}_{\E/\F}(\mathrm{GL}_{m/\E})$, $\bs{G}_v = G(\F_v)$, $\bs{G}'_v= G'(\F_v)$, etc. Comme dans la section \ref{sur l'op d'entrel IK}, on note $\K$ et non $\omega$ le caractre $\K\circ\det$ de $G(\A)$.  
On fixe aussi des mesures de Haar $\d g$ sur $G(\A)$ et $\d g'$ sur $G'(\A)$. 

\subsection{Formule de comparaison; premires simplifications.}\label{prem simplif} 
 Identifions $A_{G'}(\A)=A_G(\A)$ ˆ $\A^\times$. Fixons un caractre unitaire $\mu$ de $\A^\times$ trivial sur $\F^\times$ et posons 
$\xi = \K^{md(d-1)/2}\mu$. Rappelons que $$\bs{\mathfrak{F}}=\bs{\mathfrak{F}}^{\bs{K},\bs{K'}}_{G(\A),G'(\A)}$$ est le sous-ensemble de $C^{\infty,\bs{K}-\mathrm{fin}}_{\mathrm{c}}(G(\A))\times C^{\infty,\bs{K}'-\mathrm{fin}}_{\mathrm{c}}(G'(\A))$ formŽ des paires de fonctions $(f,\phi)$ qui sont associŽes au sens de \ref{fonctions associŽes}. On part de l'ŽgalitŽ \ref{dec endos}\,(4), pour toute paire de fonctions $(f,\phi)\in \bs{\mathfrak{F}}$: 
$$I^G_{\mathrm{disc},\xi}(f,\K)= \frac{1}{d} I^{G'}_{\mathrm{disc},\mu}(\phi)\ptf\leqno{(1)}$$

Pour $k\in \mathbb{N}^*$ divisant $n$, on note $L_k\in \ES{L}(A_0)$ le $\F$-facteur de Levi standard de $G$ dŽfini par $L_k = \mathrm{GL}_{\frac{n}{k}/\F}\times \cdots \times \mathrm{GL}_{\frac{n}{k}/\F}$ ($k$ fois) et on pose $P_k= P_{L_k}$. Seuls les $\F$-facteurs de Levi $L\in \ES{L}(A_0)$ qui sont $G(\F)$-conjuguŽs ˆ un tel $L_k$ peuvent donner une contribution non triviale ˆ l'expression 
$$I_{\mathrm{disc},\xi}^G(f,\K)
=\!\! \sum_{L\in \ES{L}(A_0)} \!\frac{\vert W^L_0 \vert}{\vert W_0^G\vert}\!\sum_{s\in W(L)_{\mathrm{r\acute{e}g}}} \!\!
\vert \det(s-1\vert \mathfrak{a}_L^G )\vert^{-1}\mathrm{tr}\left(\rho_{P,\mathrm{disc},s,\xi}(f,\K) \right)\ptf$$ 
En effet pour les autres, l'ensemble $W(L)_{\mathrm{r\acute{e}g}}$ est vide. D'autre part deux ŽlŽments de $\ES{L}(A_0)$ qui sont conjuguŽs dans $W^G_0$ donnent la mme contribution. On peut donc remplacer dans l'expression ci-dessus la somme $\sum_{L\in \ES{L}(A_0)}  \!\frac{\vert W^L_0 \vert}{\vert W_0^G\vert}$ par une somme $\sum_{k\vert n}\frac{1}{\vert W(L_k)\vert}$ avec (rappel) $W(L_k)= N^{W^G_0}(A_{L_k})/ W^{L_k}_0$. Comme 
$W(L_k)\simeq \mathfrak{S}_k$, on a 
$\vert W(L_k)\vert = \vert \mathfrak{S}_k \vert =k!$. D'o l'ŽgalitŽ 
$$ I_{\mathrm{disc},\xi}^G(f,\K)
=\!\! \sum_{k\vert n} \!\frac{1}{k!}\!\sum_{s\in W(L_k)_{\mathrm{r\acute{e}g}}} \!\!
\vert \det(s-1\vert \mathfrak{a}_{L_k}^G )\vert^{-1}\mathrm{tr}\left(\rho_{P_k,\mathrm{disc},s,\xi}(f,\K) \right)\ptf\leqno{(2)}$$ 
Les ŽlŽments de $W(L_k)_{\mathrm{r\acute{e}g}}$ sont les $k$-cycles (ou cycles de longueur $k$) dans le groupe de permutation $W(L_k) \simeq \mathfrak{S}_k$. 
Ils forment une unique orbite sous $W(L_k)$ et tous les ŽlŽments de cette orbite sont conjuguŽs au $k$-cycle 
$$s_k= (k, k-1, \cdots , 1)\ptf$$ L'ensemble $W(L_k)_{\mathrm{r\acute{e}g}}$ contient $(k-1)!$ ŽlŽments et pour tout $s\in W(L_k)_{\mathrm{r\acute{e}g}}$, on a l'ŽgalitŽ 
$$\vert \det(s-1\vert \mathfrak{a}_{L_k}^G )\vert\;(=\vert \det(s_k-1\vert \mathfrak{a}_{L_k}^G )\vert)=k\ptf$$ 
L'ŽgalitŽ (2) se rŽcrit
$$I_{\mathrm{disc},\xi}^G(f,\K)
=\!\! \sum_{k\vert n} \!\frac{1}{k! k }\!\sum_{s\in W(L_k)_{\mathrm{r\acute{e}g}}} \!\! 
\mathrm{tr}\left(\rho_{P_k,\mathrm{disc},s,\xi}(f,\K) \right)\ptf\leqno{(3)}$$ 
avec (rappel) $$\rho_{P_k,\mathrm{disc},s,\xi}(f,\K)=M_{P_k\vert P_k}(s,0)\circ \pi_{P_k,\xi}(f)\circ I_\K\ptf$$ 

Soit $k\geq 1$ un entier divisant $n$. Pour une reprŽsentation automorphe discrte $\tau= \tau_1\otimes \cdots \otimes \tau_k$ de $L_k(\A)$, on note $\pi_\tau = \tau_1\times \cdots \times \tau_k$ l'induite parabolique de $\tau$ ˆ $G(\A)$ suivant $P_k(\A)$. Rappelons que $\pi_\tau$ est irrŽductible (et unitaire). 
On a $$\mathrm{tr}\left(M_{P_k\vert P_k}(s,0)\circ \pi_{P_k,\xi}(f)\circ I_\K\right)= \sum_\tau \mathrm{tr}\left(M_{P_k\vert P_k}(s,0)\circ \pi_\tau(f)\circ I_\K\right)$$
o $\tau$ parcourt l'ensemble $\bs{\Pi}_{\xi_{L_k}}(L_k(\A))$ des reprŽsentations automorphes discrtes de caractre central prolongeant $\xi_{L_k}$, \ie les reprŽsentations automorphes sous-jacentes ˆ une composante irrŽductible de 
$\rho_{L_k,\mathrm{disc},\xi_{L_k}}$; rappelons que $\xi_{L_k}$ est l'unique caractre de $\mathfrak{A}_{L_k}^GA_G(\A)$ trivial  sur $\mathfrak{A}_{L_k}^G$ et prolongeant $\xi$. 

\begin{remark}
\textup{Observons que l'opŽrateur $I_\K$ entrelace $\pi_\tau$ et $\K^{-1}\pi_{\K\tau}$ avec $\K\tau = \K\tau_1\otimes \cdots \otimes \K \tau_r$. 
Comme l'opŽrateur $M_{P_k\vert P_k}(s,0)$ entrelace $\pi_\tau$ et $\pi_{s(\tau)}$ avec $s(\tau)= \tau_{s^{-1}(1)}\otimes \cdots \otimes \tau_{s^{-1}(r)}$, l'opŽrateur 
$I_\K\circ M_{P_k\vert P_k}(s,0)$ entrelace $\pi_\tau$ et $\K^{-1}\pi_{s(\K \tau)}$ avec $s(\K\tau)=\K s(\tau)$. 
En particulier seuls les $\tau$ tels que $\K s(\tau)= \tau$ contribuent ˆ la somme ci-dessus (voir la preuve de \ref{simplif op entrel}).
}\end{remark}

\begin{lemme}\label{simplif op entrel}
Soient $k\vert n$, $s\in W(L_k)_{\mathrm{r\acute{e}g}}$ et $w\in W(L_k)$. On a l'ŽgalitŽ
$$\mathrm{tr}\left(\rho_{P_k,\mathrm{disc},wsw^{-1},\xi}(f,\K)\right) =\mathrm{tr}\left(\rho_{P_k,\mathrm{disc},s,\xi}(f,\K)\right) \ptf$$
\end{lemme}

\begin{demo}
Pour allŽger l'Žcriture, supprimons l'indice $P_k\vert P_k$ dans l'opŽrateur $M_{P_k\vert P_k}(s,\lambda)$. Posons $s'= wsw^{-1}$. Il s'agit de prouver l'ŽgalitŽ 
$$\sum_\tau \mathrm{tr}\left(M(s',0)\circ \pi_\tau(f)\circ I_\K\right)= \sum_\tau \mathrm{tr}\left(M(s,0)\circ \pi_\tau(f)\circ I_\K\right)$$
o (dans les deux sommes) $\tau$ parcourt l'ensemble $\bs{\Pi}_{\xi_{L_k}}(L_k(\A))$. 
L'application $\tau\mapsto \tau'= w(\tau)\;(= \tau \circ w^{-1})$ est une bijection de l'ensemble 
$\bs{\Pi}_{\xi_{L_k}}(L_k(\A))$ sur lui-mme. 
Il suffit donc de prouver que pour chaque $\tau\in \bs{\Pi}_{\xi_{L_k}}(L_k(\A))$, on a l'ŽgalitŽ 
$$\mathrm{tr}\left(M(s'\!,0)\circ \pi_{\tau'}(f)\circ I_\K\right)= \mathrm{tr}\left(M(s,0)\circ \pi_\tau(f)\circ I_\K\right)$$
ou, ce qui revient au mme, l'ŽgalitŽ 
$$\mathrm{tr}\left(I_\K\circ M(s'\!,0)\circ \pi_{\tau'}(f)\right) = \mathrm{tr}\left(I_\K\circ M(s,0)\circ \pi_{\tau}(f)\right)\ptf$$ 

Fixons $\tau=\tau_1\otimes\cdots \otimes \tau_k$. Pour $\lambda =(\lambda_1,\ldots ,\lambda_k)\in \mathbb{C}^k=\mathfrak{a}_{L_k,\mathbb{C}}^*$, 
on note $\tau_\lambda$ 
la reprŽsentation $\nu^{\lambda_1}\tau_1\otimes \cdots \otimes \nu^{\lambda_k}\tau_k$ de $L_k(\A)$ --- \ie $\tau_\lambda (x)= e^{\langle \lambda, {\bf H}_{L_k}(x)\rangle}\tau(x)$ --- et $\pi_{\tau_\lambda}$ la reprŽsentation $\nu^{\lambda_1}\tau_1\times \cdots \times \nu^{\lambda_k}\tau_k$ (induite parabolique de $\tau_\lambda$ ˆ $G(\A)$ suivant $P_k(\A)$)\footnote{Observons que pour que le caractre central de $\tau_\lambda$ prolonge $\xi$, il faut et il suffit que $\lambda_1+ \cdots + \lambda_k=0$, \ie 
que $\lambda$ appartienne au dual $ \mathfrak{a}_{L_k,\mathbb{C}}^{G,*}$ de $\mathfrak{a}_{L_k,\mathbb{C}}^G= \mathfrak{a}_{L_k}^G\otimes_{\mathbb{R}}\mathbb{C}$.}. La reprŽsentation $\pi_{\tau_\lambda}$ n'est en gŽnŽral ni irrŽductible ni unitaire. Pour $\lambda=0$, on a donc $\tau_0=\tau$ et $\pi_{\tau_0} = \pi_{\tau}$. 

Pour  $\lambda \in \mathfrak{a}_{L_k,\mathbb{C}}^*$, on a l'Žquation fonctionnelle (cf. \cite[5.2.2\,(1)]{LW})
$$M(ws,\lambda)= M(w,s(\lambda)\circ M(s,\lambda)\ptf\leqno{(4)}$$ On en dŽduit que 
$$M(wsw^{-1},\lambda)= M(w,sw^{-1}(\lambda))M(s,w^{-1}(\lambda))M(w^{-1},\lambda)\vg$$ d'o en particulier 
$$M(wsw^{-1},0)= M(w,0)M(s,0)M(w^{-1},0)\ptf\leqno{(5)}$$ 
Notons $M(s,\lambda,\tau)$ la restriction de l'opŽrateur $M(s,\lambda)$ ˆ l'espace de $\pi_{\tau_\lambda}$; 
il entrelace $\pi_{\tau_\lambda}$ et $\pi_{s(\tau_\lambda)} $, avec 
$s(\tau_\lambda ) =\tau'_{\smash{s(\lambda)}}$. 
D'aprs (4), on a l'ŽgalitŽ 
$$M(w^{-1},\lambda,\tau)^{-1}= M(w,w^{-1}(\lambda),w^{-1}(\tau))\vg$$ d'o en particulier $$M(w^{-1},0,\tau)^{-1}=M(w^{-1},0,w^{-1}(\tau))\ptf$$
Pour $\lambda$ assez rŽgulier, d'aprs la dŽfinition de 
$M(s,\lambda)$ (cf. \cite[5.2]{LW}), on a  
$$I_\K \circ M(s,\lambda,\tau)= M(s,\lambda, \K\tau )\circ I_\K\ptf$$ Par prolongement mŽromorphe on en dŽduit que 
$$I_\K\circ M(s,0,\tau)= M(s,0,\K\tau)\circ I_\K\ptf\leqno{(6)}$$ 
L'opŽrateur (6) entrelace $\pi_{\tau}$ et $\K^{-1}\pi_{\K s(\tau)}= \K^{-1}(\K \tau_{s^{-1}(1)}\times \cdots \times \K \tau_{s^{-1}(k)})$. Or, ou bien les espaces 
$V_{\pi_\tau}$ et $V_{\K^{-1}\pi_{s(\K\tau)}}= V_{\pi_{s(\K\tau)}}$ sont disjoints, auquel cas la trace $\mathrm{tr}(I_\K\circ M(s,0,\tau)\circ \pi_{\tau}(f^\xi))$ 
est nulle; ou bien ils sont Žgaux, \ie $\K s(\tau)= \tau$. On suppose donc de plus que 
$$\K s(\tau)= \tau\ptf\leqno{(7)}$$ 
Rappelons que l'on a posŽ $s'= wsw^{-1}$ et  $\tau'= w(\tau)$. Observons que la condition $\K s'(\tau')=\tau'$ est Žquivalente ˆ (7). 
D'aprs (5) on a l'ŽgalitŽ 
$$ M(s',0,\tau')= M(w,0,sw^{-1}(\tau'))\circ M(s,0,\tau)\circ M(w^{-1},0,\tau')\ptf $$ 
Comme $$I_\K \circ M(w,0,sw^{-1}(\tau'))=M(w,0,\K s(\tau))\circ I_\K=M(w,0,\tau)\circ I_\K $$ et 
$$M(w^{-1},0,\tau')= M(w,0,w^{-1}(\tau'))^{-1}= M(w,0,\tau)^{-1}\vg$$ on en dŽduit que 
$$ I_\K \circ M(s',0,\tau') = M(w,0,\tau)\circ \left(I_\K \circ M(s,0,\tau) \right) M(w,0,\tau)^{-1}\ptf$$ 
Comme d'autre part $$M(w,0,\tau)^{-1}\circ \pi_{\tau'}(f)= \pi_{\tau}(f)\circ M(w,0,\tau)^{-1}\vg$$ on obtient l'ŽgalitŽ 
$$\mathrm{tr}(I_\K \circ M(s',0,\tau')\circ \pi_{\tau'}(f))= \mathrm{tr}(I_\K\circ M(s,0,\tau)\circ \pi_\tau(f))\ptf$$ 
Cela achve la dŽmonstration du lemme. \hfill $\square$ 
\end{demo}

\vskip2mm
D'aprs \ref{simplif op entrel}, les $(k-1)!$ ŽlŽments de $W(L_k)_{\mathrm{r\acute{e}g}}$ donnent la mme contribution 
ˆ l'expression (3). Puisque $\frac{(k-1)!}{k! k } =\frac{1}{k^2}$, l'ŽgalitŽ (3) se rŽcrit 
$$I_{\mathrm{disc},\alpha}^G(f,\K)
=\!\! \sum_{k\vert n} \!\frac{1}{k^2} 
\mathrm{tr}\left(\rho_{P_k,\mathrm{disc},s_k,\xi}(f,\K) \right)\vg\leqno{(8)}$$ 
soit encore 
$$I_{\mathrm{disc},\xi}^G(f,\K)
=\!\! \sum_{k\vert n} \!\frac{1}{k^2}  \sum_{\tau} \mathrm{tr}(M_{P_k\vert P_k}(s_k,0) \circ \pi_\tau(f)\circ I_\K)$$ 
o $\tau$ parcourt l'ensemble $\bs{\Pi}_{\xi_{L_k}}(L_k(\A))$. D'aprs la dŽmonstration de \ref{simplif op entrel}, seuls les $\tau=\tau_1 \otimes \cdots \otimes \tau_k$ vŽrifiant la condition 
$$\K s_k(\tau)=\tau$$ peuvent donner une contribution non triviale ˆ la somme ci-dessus. Cette condition est vŽrifiŽe si et seulement si 
$\tau_{i+1} = \K\tau_i$ ($i=1,\ldots , k-1$) et $\K \tau_k= \tau_1$, \ie si et seulement si 
$\tau= \tau_1 \otimes\K\tau_1\otimes \cdots \otimes \K^{k-1}\tau_1$ et $\K^k\tau_1= \tau_1$. On a aussi une condition sur le caractre central $\omega_{\tau_1}:\mathbb{A}^\times \rightarrow \mathbb{U}$ 
de la reprŽsentation automorphe discrte $\tau_1$ de $\mathrm{GL}_a(\A)$, $a=\frac{n}{k}$, ˆ savoir  
$(\omega_{\tau_1})^k\K^{a\frac{k(k-1)}{2}}=\xi$. Observons que le caractre $\K^{a\frac{k(k-1)}{2}}= \K^{\frac{n(k-1)}{2}}$ de $\A^\times$ est trivial ou d'ordre $2$. On note:
\begin{itemize}
\item $\Xi_{k,\K}(\xi)$ l'ensemble des caractres automorphes (forcŽment unitaires) $\xi_1$ de $\A^\times$ tels que 
$(\xi_1)^k =\K^{a\frac{k(k-1)}{2}}\xi$;
\item $\bs{\Pi}_{\Xi_{k,\K}(\xi)}(\mathrm{GL}_a(\A))= \bigcup_{\xi_1\in \Xi_{k,\K}(\xi)} \bs{\Pi}_{\xi_1}(\mathrm{GL}_a(\A))$ l'ensemble des reprŽsentations automorphes discrtes de $\mathrm{GL}_a(\A)$ de caractre central dans $\Xi_{k,\K}(\xi)$;
\item $\bs{\Pi}_{\Xi_{k,\K}(\xi)}(\mathrm{GL}_a(\A),\K^k)\subset \bs{\Pi}_{\Xi_{k,\K}(\xi)}(\mathrm{GL}_a(\A))$ le sous-ensemble formŽ des reprŽsentations 
qui sont $\K^k$-stables. 
\end{itemize} 
Pour un entier $k\geq 1$ divisant $n$ et une reprŽsentation automorphe discrte $\K^k$-stable $\delta$ de $\mathrm{GL}_{\frac{n}{k}}(\A)$, on note: 
\begin{itemize}
\item $\tau_{(k,\delta,\K)}$ la reprŽsentation automorphe discrte $\delta \otimes \K\delta \otimes \cdots \otimes \K^{k-1}\delta$ de $L_k(\A)$; 
\item  
$\pi_{(k,\delta,\K)}= \pi_{\tau_{(k,\delta,\K)}}=\delta \times \K \delta \times \cdots \times \K^{k-1}\delta$ l'induite parabolique de $\tau_{(k,\delta,\K)}$ 
ˆ $G(\A)$ suivant $P_k(\A)$.
\end{itemize} L'ŽgalitŽ (8) se rŽcrit 
$$ I_{\mathrm{disc},\xi}^G(f,\K)
=\!\! \sum_{k\vert n} \!\frac{1}{k^2} \sum_\delta 
\mathrm{tr}\left(M_{P_k\vert P_k}(s_k,0) \circ \pi_{(k,\delta,\K)}(f)\circ I_\K\right)\leqno{(9)}$$ 
o $\delta$ parcourt l'ensemble $\bs{\Pi}_{\Xi_{k,\K}(\xi)}(\mathrm{GL}_{\frac{n}{k}}(\A),\K^k)$. 

On obtient de la mme manire une dŽcomposition de l'expression $I_{\mathrm{disc},\mu}^{G'}(\phi)$ analogue ˆ (9). 
Pour $l\in \mathbb{N}^*$ divisant $m$, on note $L'_l\in \ES{L}(A_0)$ le $\F$-facteur de Levi standard de $G'$ dŽfini par 
$L'_l=\mathrm{Res}_{\E/\F}(\mathrm{GL}_{\frac{m}{l}/\E})\times \cdots \times \mathrm{Res}_{\E/\F}(\mathrm{GL}_{\frac{m}{l}/\E})$; $P'_l= P_{L'_l}$ le $\F$-sous-groupe parabolique standard de $G'$ de composante de Levi $L'_l$;
$s'_l\in \mathfrak{S}_l$ le $l$-cycle $(l,l-1,\ldots , 1)$; et $\Xi_l(\mu)$ l'ensemble des caractres automorphes $\mu_1$ de $\mathbb{A}^\times$ tels que $(\mu_1)^l=\mu$. 
Pour une reprŽsentation automorphe discrte $\Delta$ de $\mathrm{GL}_{\frac{m}{l}}(\A_\E)$, on note $\Delta^{\otimes l}$ la reprŽsentation automorphe discrte 
$\Delta \otimes \cdots \otimes \Delta$ de $L'_l(\A)$ et $\Delta^{\!\times l}=\Delta \times \cdots \times \Delta$ 
l'induite parabolique de $\Delta^{\otimes l}$ ˆ $G'(\A)$ suivant $P'_{l}(\A)$. On obtient 
$$I^{G'}_{\mathrm{disc},\mu}(\phi)= \sum_{l\vert m} \frac{1}{l^2} \sum_{\Delta} \mathrm{tr}(M_{P'_l\vert P'_l}(s'_l,0)\circ \Delta^{\!\times l}(\phi))\leqno{(10)}$$ 
o $\Delta$ parcourt l'ensemble $\bs{\Pi}_{\Xi_l(\mu)} (\mathrm{GL}_{\frac{m}{l}}(\AE)=\bigcup_{\mu_1\in \Xi_l(\mu)} \bs{\Pi}_{\mu_1}(\mathrm{GL}_{\frac{m}{l}}(\A_E))$. 

\subsection{Un lemme technique.}\label{un lemme technique} 
Pour $k\in \mathbb{N}^*$ divisant $n$ et $\delta$ une reprŽsentation automorphe discrte $\K^k$-stable de 
$\mathrm{GL}_{\frac{n}{k}}(\A)$, d'aprs la dŽmonstration de \ref{simplif op entrel} (ŽgalitŽ (6)), l'opŽrateur 
$I_\K \circ M_{P_k\vert P_k}(s_k,0)$ entrelace $\K \pi_{(k,\delta,\K)}$ et $\pi_{(k,\delta,\K)}$. La reprŽsentation (irrŽductible) 
$\pi_{(k,\delta,\K)}$ est donc $\K$-stable, et d'aprs le lemme de Schur, il existe une constante $c_{(k,\delta,\K)}= c_{\pi_{(k,\delta,\K)}}\in \mathbb{C}^\times$ telle que 
$$ I_\K \circ M_{P_k\vert P_k}(s_k,0)\vert_{V_{\pi_{(k,\delta,\K)}}}= c_{(k,\delta,\K)}A_{\pi_{(k,\delta,\K)}}\pvg$$ o (rappel) 
$A_{\pi_{(k,\delta,\K)}}$ est l'opŽrateur global obtenu par induction parabolique ˆ partir de l'opŽrateur global $A_{\delta}\otimes A_{\K\delta} \otimes\cdots \otimes A_{\K^{k-1}\delta}$ sur l'espace de $\tau_{(k,\delta,\K)}$ dŽduit des opŽrateurs locaux normalisŽs $A_{\K^i_v \delta_v}= A_{\delta_v}$ comme en \ref{op issu du local}. 
D'aprs \ref{A=I}, pour $k=1$, \ie pour $\pi=\delta$ discrte, on a $c_\pi=1$. Pour toute fonction $f\in C^{\infty,\bs{K}-\mathrm{fin}}_{\mathrm{c}}(G(\A))$, on a donc
\begin{align*}
\mathrm{tr}(M_{P_k\vert P_k}(s_k,0)\circ \pi_{(k,\delta,\K)}(f)\circ I_\K)&= \mathrm{tr}(\pi_{(k,\delta,\K)}\circ I_\K\circ M_{P_k\vert P_k}(s_k,0))\\
&= c_{(k,\delta,\K)}\, \mathrm{tr}\left(\pi_{(k,\delta,\K)}(f) \circ A_{\pi_{(k,\delta,\K)}}\right)
\end{align*}
et l'ŽgalitŽ \ref{prem simplif}\,(9) se rŽcrit $$I_{\mathrm{disc},\xi}^G(f,\K)= \sum_{k\vert n} \!\frac{1}{k^2} \sum_\delta 
c_{(k,\delta,\K)}\mathrm{tr}\left(\pi_{(k,\delta,\K)}(f)\circ A_{\pi_{(k,\delta,\K)}}\right)\leqno{(1)}$$ o $\delta$ parcourt l'ensemble $\bs{\Pi}_{\Xi_{k,\K}(\xi)}(\mathrm{GL}_{\frac{n}{k}}(\A),\K^k)$. 

\begin{lemme}\label{ŽgalitŽ des c}
Soient $k\in \mathbb{N}^*$ divisant $n$ et $\delta$ une reprŽsentation automorphe discrte de $\mathrm{GL}_{\frac{n}{k}}(\A)$. 
On a l'ŽgalitŽ $c_{(k,\K\delta, \K)}= c_{(k,\delta, \K)}$. 
\end{lemme}

\begin{demo}
Posons $\tau= \tau_{(k,\delta,\K)}$, $\tau'= \K\tau\;(= \tau_{(k,\K\delta,\K)})$, $\pi= \pi_\tau$ et $\pi'= \pi_{\tau'}$. Reprenons les notations de la dŽmonstration de 
\ref{simplif op entrel}. 
Rappelons que $M(s_k,0,\tau)$ dŽsigne la restriction de l'opŽrateur $M(s_k,0)$ ˆ l'espace de $\pi_\tau$, qui est aussi l'espace de $\K\pi_\tau$. On a l'ŽgalitŽ (6) de \textit{loc.\,cit.}: 
$$I_\K \circ M(s_k,0,\tau)= M(s_k,0,\tau')\circ I_\K\vg$$ soit encore 
$$ I_\K \circ M(s_k,0,\tau) = I_\K^{-1} \circ (I_\K \circ M(s_k,0,\tau'))\circ I_\K\ptf $$
Or (par dŽfinition) $I_\K \circ M(s_k,0,\tau)= c_{\pi} A_{\pi}$ et $I_\K \circ M(s_k,0,\tau')= c_{\pi'}A_{\pi'}$. 
En notant $u$ l'opŽrateur d'entrelacement entre $\K \pi$ et $\pi'$, ou ce qui revient au mme entre $\pi$ et $\K^{-1}\pi'$, donnŽ par la restriction de l'opŽrateur \guill{physique} $I_\K$, on obtient l'ŽgalitŽ 
$$c_\pi A_{\pi} = c_{\pi'}\, u^{-1}\circ  A_{\pi'} \circ u \ptf $$
Comme l'opŽrateur $A_{\pi'}$ est construit ˆ partir des opŽrateurs locaux normalisŽs (cf. \ref{op issu du local}) et que ces derniers sont insensibles aux isomorphisme 
(cf. \ref{norm op entrelac}), 
a l'ŽgalitŽ $u^{-1} \circ A_{\pi'}\circ u = A_{\pi}$. D'o l'ŽgalitŽ cherchŽe: $c_\pi= c_{\pi'}$. \hfill $\square$ 
\end{demo}

\vskip2mm
De la mme manire pour $l\geq 1$ un entier divisant $m$ et $\Delta$ une reprŽsentation automorphe discrte de $\mathrm{GL}_{\frac{m}{l}}(\AE)$, l'opŽrateur 
$M_{P'_l\vert P'_l}(s'_l,0)$ entrelace la reprŽsentations $\Delta^{\!\times l}$ de $G'(\A)= \GLm(\AE)$ avec elle-mme; et d'aprs le lemme de Schur, il existe une constante 
$c_{(l,\Delta)}=c_{\Delta^{\!\times l}}\in \mathbb{C}^\times$ telle que pour toute fonction $\phi\in C^{\infty,\bs{K}'-\mathrm{fin}}_{\mathrm{c}}(G'(\A))$, on ait 
$$\mathrm{tr}(M_{P'_l\vert P'_l}(s'_l,0)\circ \Delta^l(\phi))= c_{(l,\Delta)} \mathrm{tr}(\Delta^l(\phi))\ptf $$ 
L'ŽgalitŽ \ref{prem simplif}\,(10) se rŽcrit: 
$$I^{G'}_{\mathrm{disc},\mu}(\phi) =
\sum_{l\vert m}\frac{1}{l^2}\sum_\Delta c_{(l,\Delta)}\mathrm{tr}\left(\Delta^{\!\times l}(\phi) \right)\leqno{(2)}$$ o $\Delta$ parcourt l'ensemble 
$\bs{\Pi}_{\Xi_{l}(\mu)}(\mathrm{GL}_{\frac{m}{l}}(\AE))$. 

\subsection{SŽparation des reprŽsentations par la \guill{mŽthode standard}.}\label{mŽthode standard}
On a vu que seules les reprŽsentations automorphes de $G'(\A)$ de la forme $\Delta^{\!\times l}$ pour un entier $l\geq 1$ divisant $m$ et une reprŽsentation 
$\Delta\in \bs{\Pi}_{\Xi_l(\mu)}(\mathrm{GL}_{\frac{m}{l}}(\AE))$ peuvent donner une contribution non triviale ˆ l'expression 
$I^{G'}_{\mathrm{disc},\mu}(\phi)$. Fixons une telle reprŽsentation $\Pi= \Delta^{\!\times l}$. Posons $b= \frac{m}{l}$. On fixe un caractre additif non trivial 
$\psi$ de $\A$ trivial sur $\F$ et un sous-ensemble fini $S\subset \V$ contenant $\Vinf$ et toutes les places $v\in \Vfin$ vŽrifiant au moins l'une des trois conditions 
suivantes: 
\begin{itemize}
\item la $\F_v$-algbre cyclique $\E_v$ est ramifiŽe; 
\item la composante locale $\Delta_v$ de $\Delta$ n'est pas sphŽrique; 
\item le caractre additif $\psi_v$ de $\F_v$ n'est pas de niveau $0$. 
\end{itemize}
On pose $\V^S = \V\smallsetminus S$. Rappelons que pour toute paire de fonctions $(f,\phi)\in \bs{\mathfrak{F}}$, on a l'ŽgalitŽ $$I^G_{\mathrm{disc},\xi}(f)= \frac{1}{d} I^{G'}_{\mathrm{disc},\mu}(\phi)\ptf$$ 

Afin d'isoler la reprŽsentation $\Pi$ dans le terme de droite de l'ŽgalitŽ ci-dessus, notons $\bs{\mathfrak{R}}_\E^S(\Pi)$ l'ensemble des reprŽsentations automorphes de $G'(\A)$ de la forme $\Pi_1= \Delta_1^{\!\times l_1}$ pour un entier $l_1\geq 1$ divisant $m$ et une reprŽsentation automorphe discrte $\Delta_1$ de $\mathrm{GL}_{\frac{m}{l_1}}(\AE)$ telles qu'en toute place $v\in \V^S$, on ait l'ŽgalitŽ 
$$\mathrm{tr}(\Pi_{1,v}(b_v f_v))= \mathrm{tr}(\Pi_{v}(b_vf_v))\qhq{pour toute fonction} f\in \mathcal{H}_v= \mathcal{H}(\G_v,\K_v)\ptf$$ 
Autrement dit on demande qu'en toute place $v\in \V^S$, les reprŽsentations (irrŽductibles sphŽriques unitaires) $\Pi_v$ et $\Pi_{1,v}$ de $\G'_v$ ait le mme $\K_v$-relvement faible. D'aprs \ref{lemme local-global}, l'ensemble $\bs{\mathfrak{R}}^S_\E(\Pi)$ ne contient que des reprŽsentations de la forme $\Delta_1^{\!\times l}$ avec 
$\Delta_1\in \Gamma(\Delta) = \{{^\gamma\Delta}\,\vert \, \gamma \in \Gamma\}$. Observons aussi que  si $l=1$, \ie si $\Pi= \Delta$, alors $\bs{\mathfrak{R}}_\E^S(\Pi)= \Gamma(\Pi)$. 

Notons $\bs{\mathfrak{r}}^S_\F(\Pi)$ l'ensemble des reprŽsentations automorphes de $G(\A)$ de la forme $\pi= \pi_{(k,\delta,\K)}$ pour un entier 
$k\geq 1$ divisant $n$ et une reprŽsentation automorphe discrte de $\mathrm{GL}_{\frac{n}{k}}(\A)$ telles qu'en toute place $v\in \V^S$, on ait l'ŽgalitŽ 
$$\mathrm{tr}(\pi_v(f_v))= \mathrm{tr}(\Pi_v(b_vf_v)) \qhq{pour toute fonction} f\in \mathcal{H}_v\ptf$$ 
Autrement dit on demande qu'en toute place $v\in \V^S$, $\pi_v$ soit un $\K_v$-relvement faible de $\Pi_v$. On Žcrit $\bs{\mathfrak{r}}^S_\F(\Pi)= 
\bigcup_{k\vert n} \bs{\mathfrak{r}}_{k,\F}^S(\Pi)$ avec une dŽfinition Žvidente pour $\bs{\mathfrak{r}}_{k,\F}^S(\Pi)$. 

Notons  $\bs{\mathfrak{F}}_S\subset 
C^\infty_{\mathrm{c}}(\F_S)\times C^\infty_{\mathrm{c}}(\F_S)$ l'image de $\bs{\mathfrak{F}}$ par l'application naturelle $(f,\phi)\mapsto (f_S= \otimes_{v\in S}f_v,\phi_S= \otimes_{v\in S}\phi_v)$, \ie  l'ensemble des paires de fonctions $(f_S= \otimes_vf_v,\phi_S=\otimes_v\phi_v )\in C^\infty_{\mathrm{c}}(\F_S)\times C^\infty_{\mathrm{c}}(\F_S)$ telles que $f_\infty = \otimes_{v\in \Vinf}f_v$ soit $\bs{K}_{\!\infty}$-finie ˆ droite et ˆ gauche, $\phi_\infty = \otimes_{v\in \Vinf}\phi_v$ soit $\bs{K}'_{\!\infty}$-finie ˆ droite et ˆ gauche, et les fonctions $f_v$ et $\phi_v$ soient concordantes en toute place $v\in S$.  
D'autre part soit $\mathcal{H}^S\subset C^\infty_{\mathrm{c}}(\Afin)$ le produit tensoriel restreint $\otimes_{v\in \V^S} \mathcal{H}_v$; et de la mme manire soit $\mathcal{H}'^S= \otimes_v \mathcal{H}'_v\subset C^\infty_{\mathrm{c}}(G'(\Afin))$. 

Le lemme suivant est une simple application de la \textit{mŽthode standard}.

\begin{lemme}\label{sŽp mŽthode standard}Pour toute paire de fonctions $(f_S=\otimes_{v\in S}f_v,\phi_S= \otimes_{v\in S}\phi_v)\in \bs{\mathfrak{F}}_S$, on a l'ŽgalitŽ
$$ \sum_{k\vert n} \frac{1}{k^2} \sum_{\pi \in \bs{\mathfrak{r}}_{k,\F}^S(\Pi)} c_\pi \prod_{v\in S}\mathrm{tr}(\pi_v(f_v) \circ A_{\pi_v})= \frac{1}{d} \sum_{\Pi' \in \bs{\mathfrak{R}}_\E^S(\Pi)} c_{\Pi'} \prod_{v\in S}\mathrm{tr}(\Pi'_v((\phi_v))\ptf$$ 
\end{lemme}

\begin{demo} 
Notons $d_\Pi(f_S,\phi_S)$ la diffŽrence entre le terme ˆ gauche et le terme ˆ droite de l'ŽgalitŽ de l'ŽnoncŽ. On veut montrer que $d_\Pi(f_S,\phi_S)=0$. Pour cela on applique la \textit{mŽthode standard} due ˆ la Langlands. 
L'ŽgalitŽ $$I^G_{\mathrm{disc},\alpha}(f,\K)=\frac{1}{d} I^{G'}_{\mathrm{disc}}(\phi)$$ appliquŽe aux paires de fonctions $(f,\phi)\in \bs{\mathfrak{F}}$ telles que $(f^S\!,\phi^S)\in \mathcal{H}^S \times \mathcal{H}'^S$ ne \guill{voit} que les reprŽsentations automorphes qui sont sphŽriques en dehors de $S$. FixŽe une paire $(f_S,\phi_S)\in \bs{\mathfrak{F}}_S$, l'application qui ˆ $f^S= \otimes_v f_v \in \mathcal{H}^S$ associe l'expression 
$$I^G_{\mathrm{disc},\xi}(f_S\otimes f^S,\K)- \frac{1}{d}I^{G'}_{\mathrm{disc},\mu}(\phi_S \otimes (\otimes_{v\in \V^S}b_v f_v))\leqno{(1)}$$ 
dŽfinit une forme linŽaire sur $\mathcal{H}^S$. Dans l'expression $\frac{1}{d}I^{G'}_{\mathrm{disc},\mu}(\phi_S \otimes (\otimes_{v\in \V^S}b_v f_v))$ ne contribuent que des reprŽsentations automorphes de $G'(\A)$ de la forme $\Pi_1=\Delta_1^{\!\times l_1}$ avec $l_1\vert m$ et $\Delta_1\in \bs{\Pi}_{\Xi_{l_1}(\mu)}(\mathrm{GL}_{\frac{m}{l_1}}(\AE))$ telles que $\Delta_1$ soit sphŽrique en dehors de $S$; et dans l'expression $I^G_{\mathrm{disc},\xi}(f_S\otimes f^S,\K)$ ne contribuent que des reprŽsentations automorphes de la forme $\pi=\pi_{(k_1,\delta_1,\K)}$ pour $k_1\vert n$ et $\delta_1\in \bs{\Pi}_{\Xi_{k_1,\K}(\xi)}(\mathrm{GL}_{\frac{n}{k_1}}(\A),\K^{k_1})$ telles que $\delta_1$ soit sphŽrique en dehors de $S$. Pour toute reprŽsentation automorphe $\Pi_1$ de $G(\A)$ de la forme ci-dessus et toute reprŽsentation $\pi\in\bs{\mathfrak{r}}_\F^S(\Pi_1)$, on a 
l'ŽgalitŽ $$\prod_{v\in \V^S} \mathrm{tr}(\pi_v(f_v))\;\left(=\prod_{v\in \V^S}\mathrm{tr}(\pi_v(f_v)\circ A_{\pi_v})\right)=\prod_{v\in \V^S}\mathrm{tr}(\Pi_{1,v}(b_vf_v))\leqno{(2)}$$ 
Par consŽquent l'ŽgalitŽ (1) se rŽcrit 
$$\sum_{\Pi_1}d_{\Pi_1}(f_S,\phi_S) \prod_{v\in \V^S}\mathrm{tr}(\pi_{v}(f_v))\leqno{(3)}$$ o $\Pi_1$ 
parcourt un systme de reprŽsentants des reprŽsentations automorphes de $G'(\A)$ de la forme ci-dessus modulo la relation d'Žquivalence $\approx$ 
donnŽe par $$\Pi_1\approx \Pi'_1\qhq{si et seulement si} \bs{\mathfrak{R}}_\E^S(\Pi_1)= \bs{\mathfrak{R}}_\E^S(\Pi'_1)$$ et o, pour chaque $\Pi_1$, on a choisi une reprŽsentation $\pi\in \bs{\mathfrak{r}}_{\F}^S(\Pi_1)$. 
Observons que si $\Pi_1 \not\approx \Pi'_1$, alors pour $\pi\in \bs{\mathfrak{r}}_{\F}^S(\Pi_1)$ et $\pi'\in \bs{\mathfrak{r}}_{\F}^S(\Pi'_1)$, les $\mathcal{H}^S$-modules $\otimes_{v\in \V^S}V_{\pi_v}^{\bs{K}_v}$ et 
$\otimes_{v\in \V^S}V_{\pi'_v}^{\bs{K}_v}$ sont par dŽfinition non isomorphes. Puisque l'ŽgalitŽ (2) est vraie pour toute fonction $f\in \mathcal{H}^S$, 
d'aprs l'indŽpendance linŽaire des caractres de $\mathcal{H}^S$ (cf. \cite[theorem~4.2]{Fl2}), toutes les constantes $d_{\Pi_1}(f_S,\phi_S)$ sont nulles. 
Cela prouve le lemme.  \hfill$\square$
\end{demo}

\subsection{Construction d'un $\K$-relvement faible.}\label{construction; hypo} 
Dans cette sous-section, on prouve que toute reprŽsentation automorphe discrte $\Pi$ 
de $G'(\A)$ admet un $\K$-relvement faible de la forme $\pi_{(k,\delta,\K)}$, avec $k=r(\Pi)$. Sous une certaine hypothse de rŽcurrence sur le $\F$-rang de $G'$ (\ref{hypo rec}), 
on prouve aussi que $k= x(\delta)$. 
 
Appliquons le lemme \ref{sŽp mŽthode standard} dans le cas o particulier o la reprŽsentation $\Pi$ est discrte, \ie $l=1$ et $\Pi= \Delta$. Alors 
$\bs{\mathfrak{R}}_\E^S(\Pi)= \Gamma(\Pi)$ et pour toute 
reprŽsentation $\Pi'\in \Gamma(\Pi)$, la constante $c_{\Pi'}$ vaut $1$ (il n'y a aucun opŽrateur d'entrelacement du c™tŽ $G'$ dans ce cas). L'ŽgalitŽ de \textit{loc.\,cit.} s'Žcrit, 
pour toute paire de fonctions $(f_S,\phi_S)\in \bs{\mathfrak{F}}_S$:
$$\sum_{k\vert n} \frac{1}{k^2}\!\!\!\sum_{\pi \in \bs{\mathfrak{r}}_{k,\F}^S(\Pi)}\!\! c_\pi \prod_{v\in S}\mathrm{tr}(\pi_v(f_v)\circ A_{\pi_v})= 
\frac{1}{d} \sum_{\Pi'\in \Gamma(\Pi)} \prod_{v\in S}\mathrm{tr}(\Pi'_v(\phi_v))\ptf\leqno{(1)}$$ 
Observons que pour toute reprŽsentation $\Pi'\in \Gamma(\Pi)$ et toute place $v\in S$, on a l'ŽgalitŽ 
$$\mathrm{tr}(\Pi'_v(\phi_v))= \mathrm{tr}(\Pi_v(\phi_v))\ptf$$ En effet puisque $\phi_v$ et $f_v$ concordent, les intŽgrales semi-simples rŽgulires de $\phi_v$ sont invariantes sous l'action du groupe $\Gamma=\Gamma(\E_v/\F_v)$ --- c'est la remarque \ref{propriŽtŽs fonct concor}\,(i) dans le cas o la place $v$ est inerte dans $\E$, et le cas gŽnŽral s'en dŽduit gr‰ce au lemme \ref{existence transfert bis} ---; et pour tout $\gamma\in \Gamma$, 
cette invariance entra"ne l'ŽgalitŽ $\mathrm{tr}(({^\gamma\Pi)_v}(\phi_v))= \mathrm{tr}(\Pi_v(\phi_v))$ d'aprs la formule d'intŽgration de Weyl. 
Le terme de droite de l'ŽgalitŽ (1) se rŽcrit donc 
$$\frac{g(\Pi)}{d} \prod_{v\in S}\mathrm{tr}(\Pi_v(\phi_v))=\frac{1}{r(\Pi)} \prod_{v\in S}\mathrm{tr}(\Pi_v(\phi_v))\ptf$$ Puisqu'on peut toujours choisir la fonction $\phi_S = \prod_{v\in S}\phi_v$ (avec $(f_S,\phi_S)\in \bs{\mathfrak{F}}_S$ pour une fonction $f_S= \prod_v f_v\in C^\infty_{\mathrm{c}}(\F_S)$) de telle manire que l'expression ci-dessus soit non nulle, l'expression a gauche de l'ŽgalitŽ (1) n'est pas identiquement nulle et donc il existe un $k\vert n$ tel que l'ensemble $\bs{\mathfrak{r}}_{k,\F}^S(\Pi)$ ne soit pas vide. Soit $\pi= \pi_{(k,\delta,\K)}\in \bs{\mathfrak{r}}_{k,\F}^S(\Pi)$ pour un tel $k$. Par dŽfinition 
$\pi$ est un $\K$-relvement faible de $\Pi$. 
D'aprs \cite{JS1}, on a le thŽorme de multiplicitŽ $1$ fort dans le spectre automorphe, qui dit que deux reprŽsentations automorphes irrŽductibles de $G(\A)$ qui sont isomorphes en presque toute  place de $\F$ ont le mme support cuspidal. Par consquent si $\pi'= \pi_{(k'\!,\delta'\!,\K)}\in \bs{\mathfrak{r}}_{k'\!,\F}^S(\Pi)$, puisque $\pi$ et $\pi'$ ont le mme support cuspidal, on a forcŽment $k'=k$ et $\delta' = \kappa^i \delta$ pour un $i\in \{0,\ldots , x(\delta)-1\}$. On a donc $$\bs{\mathfrak{r}}_\F^S(\Pi)= \bs{\mathfrak{r}}_{k,\F}^S(\Pi)= \{\pi_{(k,\K^i\delta,\K)}\,\vert \, i = 0,\ldots x(\delta)-1\}\ptf$$ En particulier 
$\bs{\mathfrak{r}}_\F^S(\Pi)$ est de cardinal $x(\delta)$. Observons que pour $\pi'= \pi_{(k,\K\delta,\K)}$, en notant (comme dans la preuve de \ref{ŽgalitŽ des c}) $u\in \mathrm{Isom}_{G(\A)}(\K\pi, \pi')$ l'isomorphisme donnŽ par la restriction de l'opŽrateur 
physique $I_\K$ ˆ l'espace de $\pi$, on a les ŽgalitŽs 
$A_{\pi}= u^{-1}\circ A_{\pi'}\circ u$ et $\mathrm{tr}(\pi(f)\circ A_\pi)= \mathrm{tr}(\pi'(f)\circ A_{\pi'})$. Compte-tenu de l'ŽgalitŽ $c_\pi=c_{\pi'}$ (\ref{ŽgalitŽ des c}), le terme de gauche de l'ŽgalitŽ (1) se rŽcrit 
$$\frac{x(\delta)c_\pi}{k^2}\prod_{v\in S} \mathrm{tr}(\pi_v(f_v)\circ A_{\pi_v})\ptf$$ 
Pour toute paire de fonctions $(f_S,\phi_S)\in \bs{\mathfrak{F}}_S$, on a donc simplifiŽ l'ŽgalitŽ (1) en 
$$\frac{x(\delta)c_\pi}{k^2}\prod_{v\in S} \mathrm{tr}(\pi_v(f_v))\circ A_{\pi_v})= \frac{1}{r(\Pi)} \prod_{v\in S}\mathrm{tr}(\Pi_v(\phi_v))
\ptf\leqno{(2)}$$

\begin{lemme}\label{rel en v dans S}
Pour chaque place $v\in S$, $\pi_v$ est un $\K_v$-relvement de $\Pi_v$. 
\end{lemme}

\begin{demo} Par construction, l'expression (2) n'est pas identiquement nulle. Fixons une paire de fonctions $(f_S,\phi_S)\in \bs{\mathfrak{F}}_S$ telle que 
$\prod_{v\in S}\mathrm{tr}(\Pi_v(\phi_v))\neq 0$; on a donc aussi $\prod_{v\in S} \mathrm{tr}(\pi_v(f_v))\circ A_{\pi_v})\neq 0$. Fixons une place $v_0\in S$ et faisons varier la paire de fonctions concordantes $(f_{v_0},\phi_{v_0}$): pour toute paire de fonctions concordantes $(f'_{v_0},\phi'_{v_0})\in C^\infty_{\mathrm{c}}(\G_{v_0})\times C^\infty_{\mathrm{c}}(\G'_{v_0})$ telle que si $v_0\in \Vinf$, $f_{v_0}$ soit $\bs{K}_{\!v_0}$-finie ˆ droite et ˆ gauche et $\phi_{v_0}$ soit $\bs{K}'_{\! v_0}$-finie ˆ droite et ˆ gauche, on a l'ŽgalitŽ 
$$\mathrm{tr}(\pi_{v_0}(f'_{v_0})\circ A_{\pi_{v_0}})= c_0 \mathrm{tr}(\Pi_{v_0}(\phi'_{v_0}))$$ o la constante $c_0\in \mathbb{C}^\times$ est dŽfinie par 
$$c_0= \left(\frac{r(\Pi)x(\delta)c_\pi}{k^2}\right)^{-1}\!\!\left(\prod_{v\in S\smallsetminus \{v_0\}}\!\!\! \mathrm{tr}(\pi_v(f_v)\circ A_{\pi_v})\right)^{-1}\!\!\! \prod_{v\in S\smallsetminus \{v_0\}} \mathrm{tr}(\Pi_{v}(\phi_{v})\ptf$$ 
En d'autres termes, $\pi_{v_0}$ est un $\K_{v_0}$-relvement de $\Pi_{v_0}$. \hfill $\square$
\end{demo}

\begin{lemme}\label{r=k}
On a l'ŽgalitŽ $r(\Pi)=k$.
\end{lemme}

\begin{demo}
On Žcrit $\Pi=u(\Lambda,q)$ pour un entier $q\geq 1$ divisant $m$ et une reprŽsentation automorphe cuspidale unitaire $\Lambda$ de 
$\mathrm{GL}_{\frac{m}{q}}(\A)$, cf. \ref{rappels sur le spectre res}. Rappelons que $\Pi$ est l'unique quotient irrŽductible de la reprŽsentation automorphe 
$$R(\Lambda,q)= \nu_{\AE}^{\frac{q-1}{2}}\Lambda \times \nu_{\AE}^{\frac{q-1}{2}-1}\Lambda \times \cdots \times \nu_{\AE}^{-\frac{q-1}{2}}\Lambda$$ de $\GLm(\AE)$. En toute place $v$ de $\F$, la composante locale $\Pi_v $ est le quotient de Langlands $u_v(\Lambda_v,q)$ de $\nu_v^{\frac{q-1}{2}}\Lambda_v \times \cdots \times \ni_v^{-\frac{q-1}{2}}\Lambda_v$ (voir la preuve de \ref{A=I}). Par consŽquent dans la dŽfinition du sous-ensemble fini $S\subset \V$, la condition \guill{$\Pi_v$ est sphŽrique en toute place $v\in \V^S$} est Žquivalente ˆ \guill{$\Lambda_v$ est sphŽrique en toute place $v\in \V^S$}. Posons $s=r(\Lambda)$. D'aprs le thŽorme \ref{theo global}\,(i) connu pour les reprŽsentations automorphes cuspidales unitaires \cite{H2}\footnote{PrŽcisŽment, Henniart prouve dans \cite{H2} l'existence d'un tel $\lambda$ qui soit un $\K$-relvement \textit{faible} de $\Lambda$. Comme en toute place $v$ de $\F$, la composante locale $\Lambda_v$ est gŽnŽrique unitaire et que l'existence d'un $\K_v$-relvement est Žtabli pour toutes les reprŽsentations irrŽductibles gŽnŽriques unitaires \cite{HH,H1}, $\lambda$ est un $\K$-relvement fort.}, $s$ divise $d\frac{m}{q}=\frac{n}{q}$ et il existe un $\K$-relvement $\lambda$ de $\Lambda$ ˆ $\mathrm{GL}_{\frac{n}{q}}(\A)$ de la forme $\lambda = \beta \times \K\beta\times \cdots \times \K^{s-1}\beta$ pour une reprŽsentation automorphe cuspidale unitaire $\K^s$-stable $\beta$ de $\mathrm{GL}_{\frac{n}{qs}}(\A)$ telle que $x(\beta)=s$. Notons $\alpha$ la reprŽsentation automorphe discrte de $\mathrm{GL}_{\frac{n}{s}}(\A)$ dŽfinie par $\alpha= u(\beta,q)$ et posons $\theta= \alpha \times \K \alpha \times \cdots \times \K^{s-1}\alpha$. 
Pour chaque place $v\in \V^S$, d'aprs la description explicite de l'application de $\K_v$-relvement faible (\ref{description du K-rel faible dans le cas gŽnŽral}), 
on a que $\theta_v$ est un $\K_v$-relvement faible de $\Pi_v$. Par consŽquent $\theta$ est isomorphe ˆ un ŽlŽment de $\bs{\mathfrak{r}}_{\F}^S(\Pi)= 
\{\pi_{(k,\K^i\delta,\K)}\,\vert \, i = 0,\ldots x(\delta)-1\}$. Comme $\theta$ et $\pi= \pi_{(k,\delta,\K)}$ ont le mme support cuspidal, on a que $\delta$ appartient ˆ $X(\alpha)$. Le calcul des exposants centraux montre ensuite que $s=k$. Puisque $r(\Pi)=r(\Lambda)$, le lemme est dŽmontrŽ. 
 \hfill$\square$ 
\end{demo}

\vskip2mm 
Puisque la reprŽsentation automorphe discrte $\delta$ de $\mathrm{GL}_{\frac{n}{k}}(\A)$ est $\K^k$-stable, le cardinal $x(\delta)$ de $X(\delta)$  
divise $k$. On veut prouver l'ŽgalitŽ $x(\delta)=k$. Pour cela introduisons l'hypothse de rŽcurrence suivante: 

\begin{hypothese}\label{hypo rec}Supposons par r\'ecurrence que le thŽorme \ref{theo global} soit vrai pour tout entier $m'\geq 1$ tel que $m'<m$ (cela signifie en particulier que \ref{theo global}\,(iii) est vrai pour tout entier $r\geq 1$ divisant $d$ et tout entier $n'\geq 1$ tel que $rn'= m'd < md =n$).
\end{hypothese}

Le cas $m'=1$ est dŽmontrŽ dans \cite{H2} puisque toute reprŽsentation automorphe discrte de $\mathrm{GL}_1(\AE)$ est cuspidale unitaire 
(c'est un caractre automorphe unitaire de $\AE^\times$).

\begin{lemme}[sous l'hypothse \ref{hypo rec}]\label{x=k}
On a l'ŽgalitŽ $x(\delta)=k$.
\end{lemme}

\begin{demo}
Posons $b= \frac{k}{x(\delta)}$. 
En notant $\pi_1$ la reprŽsentation automorphe $\delta \times \K\delta\times\cdots  \times \K^{x(\delta)-1}$ de $\mathrm{GL}_{x(\delta)\frac{n}{k}}(\A)$, on a 
$\pi \simeq \pi_1^{\times b}$. Si $b>1$, puisque $r=x(\delta)$ divise $d$ et $n'=\frac{n}{k}<n$, on peut appliquer le thŽorme \ref{theo global}\,(iii): il existe un entier $m'\geq 1$ divisant $rn'= \frac{n}{b}$ avec $m'd= rn'$ et une reprŽsentation automorphe discrte $\Delta_1$ de $\mathrm{GL}_{m'}(\AE)$ tels que $\pi_1$ soit un $\K$-relvement de $\Delta_1$; de plus on $r(\Delta_1)= r$. Observons que $m'= \frac{m}{b}$. Par compatibilitŽ, en toute place $v\in \V$, entre les applications de $\K_v$-relvement et d'induction parabolique locales, 
on obtient que $\pi$ est un $\K$-relvement (fort) de la reprŽsentation automorphe $\Delta_1^{\!\times b}$ de 
$\GLm(\AE)$. En particulier en toute place $v\in \V^S$, $\pi_v$ est un $\K_v$-relvement --- et donc \textit{a fortiori} un $\K_v$-relvement faible --- 
de $\Delta_{1,v}$. D'aprs \ref{lemme local-global}, cela entra"ne que $b=1$ et $\Delta_1\in \Gamma(\Delta)$; contradiction. Donc $b=1$.\hfill$\square$
\end{demo}

\vskip2mm
Sous l'hypothse \ref{hypo rec}, d'aprs \ref{x=k} et \ref{r=k}, on a $x(\delta)=k=r(\Pi)$ et l'ŽgalitŽ (2) se rŽcrit 
$$ c_\pi \prod_{v\in S} \mathrm{tr}(\pi_v(f_v)\circ A_{\pi_v})=  \prod_{v\in S}\mathrm{tr}(\Pi_v(\phi_v))
\ptf\leqno{(3)}$$

\vskip2mm
RŽcapitulons. Rappelons que $\Pi=u(\Lambda,q)$. Posons $r=r(\Lambda)=r(\Pi)$. 
On a prouvŽ (sous l'hypothse \ref{hypo rec}) que si $\lambda$ est un $\kappa$-relvement (fort) de $\Lambda$, 
en Žcrivant $\lambda = \delta \times \kappa \delta\times \cdots \times \kappa^{r-1} \delta$ pour une reprŽsentation automorphe cuspidale unitaire $\kappa^r$-stable $\delta$ de $\mathrm{GL}_{\smash{\frac{n}{r}}}(\A)$, la reprŽsentation automorphe $$\pi= u(\lambda,q) \times \kappa u(\lambda,q)\times \cdots \times \kappa^{r-1} u(\lambda,q)$$ de $G(\A)$ est 
un $\K$-relvement faible de $\Pi$. De plus on a $r=r(\lambda)=r(u(\lambda,q))$ et il existe un sous-ensemble fini $S\subset \V$ contenant 
$\Vinf$ tel que
\begin{itemize}
\item en toute place $v\in S$, $\pi_v$ est un $\K_v$-relvement de $\Pi_v$;
\item en toute place $v\in \V^S$, les reprŽsentations $\Pi_v$ et $\pi_v$ sont sphŽriques et $\pi_v$ est un $\K_v$-relvement faible de $\Pi_v$. 
\end{itemize}

\subsection{ThŽorme \ref{theo local} (local): preuve du point (i).}\label{preuve theo local i} Soit $E/F$ une extension finie cyclique de corps locaux de caractŽristique nulle, de degrŽ $d$. Soit $\kappa$ un caractre de $F^\times$ de noyau $\mathrm{N}_{E/F}(E^\times)$. 
On veut prouver que toute reprŽsentation de Speh $\Pi_0=u(\Lambda_0,q)$ de $\GLm(E)$ admet un $\kappa$-relvement ˆ $\GLn(F)$; o $q\geq 1$ est un entier divisant $m$ et $\Lambda_0$ est une reprŽsentation de carrŽ intŽgrable de $\GLm(E)$. On peut supposer $d> 1$ sinon il n'y a rien ˆ dŽmontrer. 
On peut aussi supposer $q>1$, le rŽsultat Žtant dŽjˆ connu pour $q=1$ \cite{HH}. Observons que si $F$ est archimŽdien, alors $d=2$ et $q=m$ (car si $q<m$, il n'y pas de reprŽsentation de carrŽ intŽgrable de $\mathrm{GL}_{\frac{m}{q}}(E)\simeq \mathrm{GL}_{\frac{m}{q}}(\mathbb{C})$). 

\vskip1mm
On suppose tout d'abord que $F$ est non archimŽdien. Soit $\Pi_0=u(\Lambda_0,q)$ comme ci-dessus. Choisissons une extension finie cyclique de corps de nombres $\E/\F$ redonnant l'extension de corps locaux $E/F$ en une place finie de $\F$. Plus prŽcisŽment, on suppose qu'il existe une place finie $v_0$ de $\F$ inerte dans $\E$ et un isomorphisme de corps topologiques $\iota:\E_{v_0}\buildrel\simeq \over{\longrightarrow} E$ qui induit un isomorphisme 
de $\F_{v_0}$ sur $F$. On note $\K$ l'unique caractre 
de $\A^\times=\A^\times_\F$ de noyau $\F^\times \mathrm{N}_{\E/\F}(\AE^\times)$ tel que $\K_{v_0} = \kappa\circ \iota\vert_{\F_{v_0}^\times}$. 
Choisissons une place finie $v_1\neq v_0$ de $\F$ inerte et non ramifiŽe dans $\E$, de caractŽristique rŽsiduelle $>n$. 
Choisissons aussi une place finie $v_2$ de $\F$ dŽployŽe dans $\E$; on a donc $v_2\neq v_0$ et $v_2\neq v_1$. 
On peut aussi comme dans \cite{HL2} supposer que toutes les places $v\in \Vinf$ sont dŽployŽes dans $\E$ (\eg si $\F$ est totalement imaginaire). 
D'aprs \cite{HL2}, il existe une reprŽsentation automorphe cuspidale unitaire $\Lambda$ de $\mathrm{GL}_{\frac{m}{q}}(\AE)$ telle que: 
\begin{itemize}
\item $\Lambda_{v_0}\simeq \Lambda_0$; 
\item $\Lambda_{v_1}$ soit cuspidale, $\Gamma(\E/\F)$-rŽgulire et de niveau $0$; 
\item $\Lambda_{v}$ soit sphŽrique en toute place finie $v\notin \{v_0,v_1,v_2\}$.
\end{itemize}
Puisque la composante locale $\Lambda_{v_1}$ est $\Gamma(\E/\F)$-rŽgulire, la reprŽsentation $\Lambda$ est \textit{a fortiori} $\Gamma(\E/\F)$-rŽgulire, \ie $r(\Lambda)=1$. Soit $S$ un ensemble fini de places de $\F$ contenant $\Vinf \cup\{v_0,v_1,v_2\}$, toutes les places $v\in \Vfin$ o la $\F_v$-algbre cyclique 
$\E_v$ est ramifiŽe et toutes les places $v\in \Vfin$ o le caractre additif $\psi_v$ de $\F_v$ n'est pas de niveau $0$. Posons $m'= \frac{m}{q}$ et $n'= m'd= \frac{n}{q}$. Puisque $r(\Lambda)=1$, d'aprs le thŽorme \ref{theo global}\,(i) connu pour les reprŽsentations automorphes cuspidales unitaires \cite{H2}, il existe une reprŽsentation automorphe cuspidale unitaire $\K$-stable $\lambda$ de $\mathrm{GL}_{n'}(\A)$ qui soit un $\K$-relvement (fort) de $\Lambda$. Posons $\Pi=u(\Lambda,q)$ et 
$\pi=u(\lambda,q)$. On a $r(\Pi)=r(\Lambda)=1=x(\lambda)$. D'aprs la preuve de \ref{A=I}, pour $v\in \V$, on a $\Pi_v= u(\Lambda_v,q)$ et $\pi_v=u(\lambda_v,q)$. D'autre part d'aprs \ref{construction; hypo} (sous l'hypothse \ref{hypo rec}) $\pi$ est un $\K$-relvement faible de $\Pi$ et pour $v\in S$, $\pi_v$ est un $\K_v$-relvement de $\Pi_v$. En particulier $\pi_{v_0}$ est un $\K_{v_0}$-relvement de $\Pi_{v_0}\simeq \Pi_0$. 

Soit $r= r(\Lambda_0)$. D'aprs le thŽorme \ref{theo local}\,(i) connu pour les reprŽsentations irrŽductibles de carrŽ intŽgrable, le $\kappa$-relvement $\lambda_0\simeq\lambda_{v_0}$ de $\Lambda_0$ s'Žcrit $\lambda_0= \delta_0 \times \kappa \delta_0 \times \cdots \times \kappa^{r-1}\delta_0$ pour un entier $r\geq 1$ divisant $n'$ et une reprŽsentation irrŽductible de carrŽ intŽgrable $\kappa^{r}$-stable $\delta_0$ de $\mathrm{GL}_{\frac{n}{r}}(F)$ telle que $x(\delta_0)= r$. Alors $\pi_{v_0} \simeq u_0\times \kappa u_0 \times \cdots \times \kappa^{r-1}u_0$ avec $u_0=u(\delta_0,q)$. 
Cela achve la preuve de \ref{theo local}\,(i) dans le cas non archimŽdien (sous l'hypothse \ref{hypo rec}). 
  
\vskip1mm
On suppose maintenant que $F$ est archimŽdien. Puisque $d=2$, on peut supposer que $E/F= \mathbb{C}/\mathbb{R}$. 
Soit $\Pi_0= u(\Lambda_0,m)$ pour un caractre unitaire $\Lambda_0$ de $\mathbb{C}^\times$. Rappelons que (avec les notations de \ref{rappels (g,K)-modules} et 
\ref{rel ell cas archimŽdien}) 
$\Lambda_0$ est de la forme 
$$\xi_{s-k,k}=\xi_{s,0}\theta^k: z \mapsto \vert z\vert_{\mathbb{C}}^{s-k} z^k$$ pour un $s \in \mathbb{U}$ et un $k\in \mathbb{Z}$, et que $\Lambda_0$ est $\Gamma(\mathbb{C}/\mathbb{R})$-rŽgulier si et seulement si $k\neq 0$. 
Choisissons une extension quadratique de corps de nombres $\E/\F$ et un caractre automorphe unitaire $\Gamma(\E/\F)$-rŽgulier 
$\Lambda$ de $\A^\times=\AF^\times$ tels qu'en une place $v_0\in \Vinf$, il existe un isomorphisme de corps topologiques $\iota: \E_{v_0} \buildrel\simeq \over{\longrightarrow} \mathbb{C}$ induisant un isomorphisme de $\F_{v_0}$ sur $\mathbb{R}$ tel que $\Lambda_{v_0}= \Lambda_0\circ \iota\vert_{\E_{v_0}^\times}$. Soit $\K$ le gŽnŽrateur de $X(\E/\F)$ et soit $\kappa$ le gŽnŽrateur de $X(\mathbb{C}/\mathbb{R})$, \ie le caractre signe de $\mathbb{R}^\times$; 
on a $\K_{v_0}=  \kappa \circ \iota \vert_{\F_{v_0}^\times}$. On sait \cite{H2} qu'il existe une reprŽsentation automorphe cuspidale unitaire $\K$-stable $\lambda$ de  $\mathrm{GL}_{2}(\A)$ qui soit un $\K$-relvement (fort) de $\Lambda$. Observons que pour $v\in \V$, la composante locale $\lambda_v$ de $\lambda$ en $v$ est: 
\begin{itemize} 
\item une reprŽsentation irrŽductible de carrŽ intŽgrable $\K_v$-stable de $\mathrm{GL}_2(\F_v)$ si $\Lambda_v$ est $\Gamma(\E/\F)$-rŽgulier; 
\item une induite parabolique irrŽductible de la forme $\chi_v \times \K_v\chi_v$ pour un caractre unitaire 
$\chi_v$ de $\F_v^\times$ sinon.
\end{itemize}
Posons $\Pi= u(\Lambda,m)$ et $\pi= u(\lambda,m)$. On procde ensuite comme dans le cas non archimŽdien. On obtient que 
$\pi_{v_0}= u(\lambda_{v_0},m)$ est un $\kappa$-relvement de $\Lambda_{v_0}\simeq \Lambda_0$. Si $\Lambda_0$ est $\mathrm{Gal}(\mathbb{C}/\mathbb{R})$-rŽgulier, alors $x(\pi_{v_0})=x(\lambda_{v_0})=1$; et si $\Lambda_0$ est invariant par la conjugaison complexe, alors $\lambda_{v_0}= \chi_{v_0} \times \kappa \chi_{v_0}$ pour un caractre unitaire de $\F_{v_0}^\times$ et $\pi_{v_0}\simeq u(\chi_{v_0},m)\times \kappa u(\chi_{v_0},m)$ avec $x(u(\chi_{v_0},m))= x(\chi_{v_0})=1$. Cela achve la preuve de \ref{theo local}\,(i) dans le cas archimŽdien.

\subsection{ThŽorme \ref{theo local} (local): preuve du point (ii).}\label{preuve theo local ii}
Soient $E/F$ et $\kappa$ comme en \ref{preuve theo local i}. Soit un entier $r\geq 1$ divisant $d$ et soit $u$ une reprŽsentation de Speh de $\mathrm{GL}_{n'}(F)$ pour un entier $n'\geq 1$, telle que $x(u)=r$. Posons $n=rn'$ et $\pi=u\times \kappa u \times \cdots \times \kappa^{r-1}u$. Puisque 
$\pi$ est une reprŽsentation irrŽductible $\kappa$-stable de $\GLn(F)$, $d$ divise $n$ (le caractre central de $\kappa\pi$ est Žgal ˆ celui de $\pi$, par suite 
$\kappa\circ \det $ est trivial sur le centre de $\GLn(F)$, \ie $\kappa^n =1$). \'Ecrivons $u=u(\delta,q)$ pour un entier $q\geq 1$ divisant $n'$ et une reprŽsentation irrŽductible de carrŽ intŽgrable $\delta$ de $\mathrm{GL}_{n''}(F)$, $n''=\frac{n'}{q}$. Puisque $x(\delta)=x(u)=r$, d'aprs le thŽorme \ref{theo local}\,(ii) connu dans ce cas, $d$ divise $rn''$ et il existe une reprŽsentation irrŽductible de carrŽ intŽgrable $\delta_E$ de $\mathrm{GL}_{m'}(E)$, $m'= \frac{rn''}{d}$, telle que $r(\delta_E)=r$ et $\pi'=\delta \times \kappa \delta \times \cdots \times \kappa^{r-1}\delta$ soit un $\kappa$ relvement de $\delta_E$ ˆ $\mathrm{GL}_{rn''}(F)$.  D'aprs le point (i) du thŽorme \ref{theo local} (prouvŽ en \ref{preuve theo local i}), $\pi$ est un $\kappa$-relvement de $u_E=u(\delta_E,q)$. De plus on a $r(u_E)=r(\delta_E)= x(\delta)=x(u)$. Quant ˆ la description des fibres, elle rŽsulte du point (iii) de \ref{prop unit}, qui sera prouvŽ en \ref{preuve prop unit}.

\subsection{Proposition \ref{prop unit} (local).}\label{preuve prop unit} 
Toute reprŽsentation irrŽductible unitaire $\tau$ de $\GLm(E)$ s'Žcrit comme un produit (irrŽductible) de reprŽsentations unitaires irrŽductibles du type $u(\delta_E,q)$ ou $u(\delta_E,q; \alpha)= \nu_E^\alpha u(\delta_E,q)\times \nu_E^{-\alpha}u(\delta_E,q)$; o $\delta_E$ est une reprŽsentation irrŽductible de carrŽ intŽgrable de $\mathrm{GL}_a(E)$ pour un entier $a\geq 1$, $q\in \mathbb{N}^*$ et $\alpha \in\;  ]0,\frac{1}{2}[$. On a vu (\ref{preuve theo local i}) que toute reprŽsentation de $\mathrm{GL}_a(E)$ 
du type $u_E=u(\delta_E,q)$ admet un $\kappa$-relvement ˆ $\mathrm{GL}_{ad}(F)$, disons $u$, du type $u(\delta,q)\times \kappa u(\delta,q)\times \cdots \times \kappa^{r-1}u(\delta,q)$. Par compatibilitŽ entre les opŽrations de $\kappa$-relvement et d'induction parabolique, cela entra"ne que pour $\alpha\in\; ]0,\frac{1}{2}[$, la reprŽsentation $\nu_F^\alpha u \times \nu_F^{-\alpha}u$ est un $\kappa$-relvement de 
$\nu_E^{\alpha}u_E \times \nu_E^{-\alpha} u_E=u(\delta_E,q;\alpha)$. On en dŽduit, ˆ nouveau par compatibilitŽ entre les opŽrations de $\kappa$-relvement et d'induction parabolique, que toute reprŽsentation irrŽductible unitaire $\tau$ de $\GLm(E)$ admet un $\kappa$-relvement ˆ $\GLmd(F)$. 
Puisque $\tau$ est gŽnŽrique si et seulement si tous les facteurs du produit sont du type $u(\delta_E,q=1)$ ou $u(\delta_E,q=1;\alpha)$, on en dŽduit aussi que $\tau$ est gŽnŽrique si et seulement si son $\kappa$-relvement est gŽnŽrique. Cela prouve \ref{prop unit}\,(i). 

Soit $\pi$ une reprŽsentation irrŽductible $\kappa$-stable de $\GLmd(F)$. On dŽcompose $\pi$ en un produit de reprŽsentations irrŽductibles unitaires du type $u(\delta,q)$ ou $u(\delta,q;\alpha)$; o $\delta$ est une reprŽsentation irrŽductible de carrŽ intŽgrable de $\mathrm{GL}_a(F)$ pour un entier $a\geq 1$, $q\in \mathbb{N}^*$ et $\alpha \in\; ]0,\frac{1}{2}[$. Puisque $\kappa u(\delta,q) \simeq u(\kappa\delta,q)$ et $\kappa u(\delta,q;\alpha)\simeq u(\kappa\delta,q; \alpha)$, par unicitŽ de la dŽcomposition ˆ isomorphisme et permutation prs des facteurs du produit, on en dŽduit que $\pi$ s'Žcrit comme un produit de reprŽsentations irrŽductibles unitaires $\kappa$-stables du type 
$\pi_{(x(u),u,\kappa)}= u \times \kappa u \times \cdots \times \kappa^{x(u)-1}u$ ou $\nu^\alpha \pi_{(x(u),u,\kappa)} \times \nu^{-\alpha}\pi_{(x(u),u,\kappa)}$ avec $u=u(\delta,q)$ et $\alpha$ comme ci-dessus; rappelons que $x(u(\delta,q)))=x(\delta)$. On a vu (\ref{preuve theo local ii}) que chacune de ces reprŽsentations 
$\pi_{(x(u),u,\kappa)}$ est un $\kappa$-relvement. Par compatibilitŽ entre les opŽrations de $\kappa$-relvement et d'induction parabolique, on en dŽduit que 
chacune des reprŽsentations $\nu^\alpha \pi_{(x(u),u,\kappa)} \times \nu^{-\alpha}\pi_{(x(u),u,\kappa)}$ est un $\kappa$-relvement, puis que $\pi$ est un $\kappa$-relvement. Cela prouve \ref{prop unit}\,(ii). 

Quant ˆ la description des fibres (point (iii)), soient $\tau= \tau_1\times \cdots \times \tau_s$ une reprŽsentation irrŽductible unitaire de $\GLm(E)$, o $\tau_i$ est du type $u(\delta_{E,i},q_i)$ (type (I)) ou $u(\delta_{E,i},q_i;\alpha_i)$ (type (II)) comme dans l'ŽnoncŽ de \ref{prop unit}\,(iii). Soit $\tau'$ une reprŽsentation irrŽductible unitaire de $\GLm(E)$ ayant mme $\kappa$-relvement que $\tau$. DŽcomposons $\tau'$ en $ \tau'_1\times \cdots \times \tau'_{t}$ o $\tau'_i$ est du type (I) ou (II) comme dans l'ŽnoncŽ de \ref{prop unit}\,(iii).  D'aprs la forme du $\kappa$-relvement des $\tau_i$ (cf. la preuve de \ref{prop unit}\,(i) ci-dessus), on a forcŽment $t=s$; et quitte ˆ permuter les $\tau'_i$, on peut supposer que $\tau_i$ et $\tau'_i$ ont mme $\kappa$-relvement ($i=1,\ldots ,s$). Alors 
pour $i=1,\ldots , s$, $\tau_i$ et $\tau'_i$ sont du mme type; plus prŽcisŽment 
$$\tau'_i = \left\{\begin{array}{ll}u(\delta'_{E,i}, q'_i=q_i) & \hbox{si $\tau_i$ est du type (I)}\\
u(\delta'_{E,i},q'_i=q_i;\alpha'_i=\alpha_i)& \hbox{si $\tau_i$ est du type (II)}\end{array}\right.$$
avec dans les deux cas $r(\delta'_{E,i})=r(\delta_{E,i})$. 
On est donc ramenŽ au cas o $s=1$, \ie $\tau=u(\delta_E,q)$ ou $\tau=u(\delta_E,q;\alpha)$. Traitons d'abord le cas (I): 
$\tau= u(\delta_E,q)$ et 
$\tau'=u(\delta'_E,q)$. Posons $r=r(\delta_E)=r(\delta'_E)$. Soit $\delta$ une reprŽsentation irrŽductible de carrŽ intŽgrable de $\mathrm{GL}_{\frac{md}{rq}}(F)$ telle que $\delta\times \kappa \delta\times \cdots \times \kappa^{r-1}\delta$ soit un $\kappa$-relvement de $\delta_E$. Posons $u=u(\delta,q)$. Alors $\pi= u \times \kappa u \times \cdots \times \kappa^{r-1}u$ est un $\kappa$-relvement de $\tau$. De la mme manire, soit $\delta'$ une reprŽsentation  irrŽductible de carrŽ intŽgrable de 
$\mathrm{GL}_{\frac{md}{rq}}(F)$ telle que $\delta' \times \kappa \delta'\times \cdots \times \kappa^{r-1}\delta'$ soit un $\kappa$-relvement de $\delta'_E$. 
Posons $u'=u(\delta',q)$. Puisque $\pi\simeq u'\times \kappa u' \times \cdots \times \kappa^{r-1}u'$, il existe un 
$i\in \{0,\ldots ,r-1\}$ tel que $u' \simeq \kappa^i u \simeq u(\kappa^i \delta,q)$. Par consŽquent $\delta' \simeq \kappa^i \delta$; et $\delta_E$ et $\delta'_E$ ont mme $\kappa$-relvement. D'aprs \ref{prop ell}\,(ii) (prouvŽ dans la section \ref{kappa-rel elliptiques}), cela n'est possible que si $\delta'_E ={^\gamma \delta_E}$ pour un $\gamma\in \Gamma(E/F)$. Comme $u({^\gamma{\delta_E}},k)\simeq {^\gamma\tau}$, on obtient que $\tau'\simeq {^\gamma\tau}$. Traitons maintenant le cas 
(II): $\tau=u(\delta_E,q;\alpha)$ et $\tau'= u(\delta'_E,q; \alpha)$. Posons $u_E= u(\delta_E,q)$ et $u'_E=u(\delta'_E,q)$. Puisque $\tau = \nu_E^\alpha u_E\times \nu_E^{-\alpha}u_E$ et $\tau'= \nu_E^\alpha u'_E \times \nu_E^{-\alpha} u'_E$ ont mme $\kappa$-relvement, $u_E$ et $u'_E$ ont forcŽment mme $\kappa$-relvement. D'aprs le cas (I) traitŽ ci-dessus, on a $u'_E \simeq {^\gamma{u_E}}$ pour un $\gamma\in \Gamma(E/F)$. Comme $u({^\gamma{u_E}},q;\alpha)\simeq {^\gamma\tau}$, on obtient que $\tau'\simeq {^\gamma\tau}$. Cela achve la preuve de \ref{prop unit}\,(iii). 

\subsection{ThŽorme \ref{theo global} (global): preuve du point (i).}\label{preuve theo global i} Soit $\Pi$ une reprŽsentation automorphe discrte de $\GLm(\AE)$. En \ref{construction; hypo} on a construit (sous l'hypothse de rŽcurrence \ref{hypo rec}) un $\K$-relvement faible de $\Pi$ de la forme $\pi= \pi_{(k,\delta,\K)}$ pour un entier $k\geq 1$ divisant $n$ et une reprŽsentation automorphe discrte $\K^k$-stable $\delta$ de $\mathrm{GL}_{\frac{n}{k}}(\A)$ telle que $x(\delta)=k = r(\Pi)$. De plus il existe un ensemble fini $S$ de places de $\F$ contenant $\Vinf$ tel que: 
\begin{itemize} 
\item pour $v\in S$, $\pi_v$ est un $\K_v$-relvement de $\Pi_v$; 
\item pour $v\in \V^S$, la $\F_v$-algbres $\E_v$ est non ramifiŽe, les reprŽsentations $\Pi_v$ de $\G'_v$ et $\pi_v$ de $\G_v$ 
sont sphŽriques, et $\pi_v$ est un $\K_v$-relvement faible de $\Pi_v$.
\end{itemize}
Or pour toute place $v$ de $\F$, et en particulier pour $v\in \V^S$, toute reprŽsentation irrŽductible unitaire de $\G'_v$ admet un $\K_v$-relvement ˆ $\G_v$ 
(d'aprs \ref{preuve theo local i}). Par suite pour 
$v\in \V^S$, la reprŽsentation $\pi_v$ est un $\K_v$-relvement de $\Pi_v$. Donc $\pi$ est un $\K$-relvement (fort) de $\Pi$. Comme on sait dŽjˆ (d'aprs \cite{H2}) que toute reprŽsentation automorphe cuspidale unitaire de $\GLm(\AE)$ admet un $\K$-relvement ˆ $\GLn(\A)$ qui est une induite de cuspidale unitaire, on a que $\Pi$ est cuspidale si et seulement si $\delta$ est cuspidale. Cela achve la preuve 
de \ref{theo global}\,(i). 

\subsection{ThŽorme \ref{theo global} (global): preuve du point (ii).}\label{preuve theo global ii} 
Soit $\pi$ une reprŽsentation automorphe de $\GLn(\A)$ de la forme 
$\pi= \delta \times \K \delta \times \cdots \times \K^{r-1}\delta$ pour un entier $k\geq 1$ divisant $n$ et une reprŽsentation automorphe discrte 
$\delta$ de $\mathrm{GL}_{n'}(\A)$, $n'=\frac{n}{r}$, telle que $x(\delta)=r$. On veut prouver que $\pi$ est un $\K$-relvement (fort) d'une reprŽsentation automorphe discrte $\Pi$ de $\GLm(\AE)$ telle que $r(\Pi)=r$. Pour cela on reprend la \textit{mŽthode standard} de \ref{mŽthode standard}, en inversant les r™les de $\Pi$ et $\pi$.

Le caractre additif non trivial $\psi$ de $\A$ trivial sur $\F$ Žtant fixŽ, on fixe un sous-ensemble fini $S\subset \V$ contenant $\Vinf$ et toutes les places $v\in \Vfin$ vŽrifiant au moins l'une des trois conditions 
suivantes: 
\begin{itemize}
\item la $\F_v$-algbre cyclique $\E_v$ est ramifiŽe; 
\item la composante locale $\delta_v$ de $\delta$ n'est pas sphŽrique; 
\item le caractre additif $\psi_v$ de $\F_v$ n'est pas de niveau $0$. 
\end{itemize}
On pose $\V^S = \V\smallsetminus S$. Rappelons que pour toute paire de fonctions $(f,\phi)\in \bs{\mathfrak{F}}$, on a l'ŽgalitŽ $$I^G_{\mathrm{disc},\xi}(f)= \frac{1}{d} I^{G'}_{\mathrm{disc},\mu}(\phi)$$ 
qui, d'aprs \ref{prem simplif}  et \ref{un lemme technique}, se simplifie en 
$$\sum_{k\vert n} \!\frac{1}{k^2} \sum_\delta 
c_{(k,\delta,\K)}\mathrm{tr}\left(\pi_{(k,\delta,\K)}(f)\circ A_{\pi_{(k,\delta,\K)}}\right)
=\frac{1}{d}\sum_{l\vert m}\frac{1}{l^2}\sum_\Delta c_{(l,\Delta)}\mathrm{tr}\left(\Delta^{\!\times l}(\phi) \right)\leqno{(1)}$$ o 
$\delta$ parcourt l'ensemble $\bs{\Pi}_{\Xi_{k,\K}(\xi)}(\mathrm{GL}_{\frac{n}{k}}(\A),\K^k)$ et $\Delta$ parcourt l'ensemble  
$\bs{\Pi}_{\Xi_{l}(\mu)}(\mathrm{GL}_{\frac{m}{l}}(\AE))$. Notons:
\begin{itemize}
\item $\bs{\mathfrak{r}}_\F^S(\pi)= \bigcup_{k\vert n} \bs{\mathfrak{r}}_{k,\F}^S(\pi)$ l'ensemble des reprŽsentations automorphes $\pi'$ de $\GLn(\A)$ de la forme 
$\pi'=\pi_{(k,\delta,\K)}$ apparaissant ˆ gauche de l'ŽgalitŽ (1) et telles qu'en toute place $v\in \V^S$, on ait $\pi'_v \simeq \pi_v$; 
\item $\bs{\mathfrak{R}}_\E^S(\pi)= \bigcup_{l\vert m} \bs{\mathfrak{R}}_{l,\E}^S(\pi)$ l'ensemble des reprŽsentations automorphes $\Pi$ de $\GLm(\AE)$ de la forme $\Pi= \Delta^{\!\times l}$ apparaissant ˆ droite de l'ŽgalitŽ (1) et telles qu'en toute place $v\in \V^S$, $\pi_v$ soit un $\K_v$-relvement faible de $\Pi_v$. 
\end{itemize}
D'aprs le lemme local-global \ref{lemme local-global}, si l'ensemble $\bs{\mathfrak{R}}_\E^S(\pi)$ est non vide (ce que l'on veut dŽmontrer), alors il est de la forme 
$$\bs{\mathfrak{R}}_\E^S(\pi)= \bs{\mathfrak{R}}_{\E,a}^S(\pi)= \{({^{\sigma^i\!\!}\Delta})^{\!\times a}\,\vert \, i=0,\ldots , r(\Delta)-1\}$$ pour un entier $a\geq 1$ divisant $m$ et une reprŽsentation automorphe discrte $\Delta$ de $\mathrm{GL}_{\frac{m}{a}}(\AE)$; o $\sigma$ est un gŽnŽrateur de $\Gamma(\E/\F)$. Supposons que $\bs{\mathfrak{R}}_\E^S(\pi) $ soit non vide et soit $\Pi=\Delta^{\!\times a}\in \bs{\mathfrak{R}}_\E^S(\pi)$. Posons $r'=r(\Delta)$. 
Si $a>1$, alors d'aprs l'hypothse de rŽcurrence \ref{hypo rec}, il existe une reprŽsentation automorphe discrte $\K^r$-stable 
$\delta'$ de $\mathrm{GL}_{\smash{\frac{md}{ra}}}(\A)$ telle que $\pi'= \delta' \times \K \delta' \times \cdots \times \K^{r'-1}\delta'$ soit une $\K$-relvement de $\Delta$; et on a $r'=x(\delta')$. Par compatibilitŽ, en toute place $v\in \V$, entre les opŽrations locales de $\K_v$-relvement et d'induction parabolique, on en dŽduit que 
$\pi'^{\times a}$ est un $\K$-relvement de $\Pi$; ce qui est impossible. Donc $a=1$ et $\Pi=\Delta$. D'autre part d'aprs le \guill{phŽnomne de rigiditŽ} \ref{rigiditŽ}, on a 
$$\bs{\mathfrak{r}}_\F^S(\pi)= \{\pi_{(r,\K^i \delta, \K)} \,\vert \, i=0,\ldots ,r-1\}\ptf$$ 

Pour toute paire de fonctions $(f_S=\otimes_{v\in S}f_v,\phi_S= \otimes_{v\in S}\phi_v)\in \bs{\mathfrak{F}}_S$, la \textit{mŽthode standard} fournit l'ŽgalitŽ
$$\frac{1}{r^2}\sum_{\pi'\in R_{\F}^S(\pi)} c_{\pi'} \prod_{v\in S}\mathrm{tr}(\pi'_v(f_v) \circ A_{\pi'_v})= \frac{1}{d} \sum_{\Pi' \in \bs{\mathfrak{R}}_\E^S(\pi)} c_{\Pi'}\prod_{v\in S}\mathrm{tr}(\Pi'_v((\phi_v))\ptf\leqno{(2)}$$ 
Comme en \ref{construction; hypo}, on obtient (compte-tenu du lemme \ref{ŽgalitŽ des c}) que l'expression ˆ gauche de l'ŽgalitŽ (2) est Žgale 
ˆ $$\frac{c_\pi}{r} \prod_{v\in S}\mathrm{tr}(\pi_v(f_v) \circ A_{\pi_v})\ptf$$
On peut choisir les fonctions $f_v$ ($v\in S$) de telle manire que l'expression ci-dessus soit non nulle.  
Par consŽquent l'expression ˆ droite de l'ŽgalitŽ (2) n'est pas identiquement nulle et l'ensemble $\bs{\mathfrak{R}}_\E^S(\pi)$ n'est pas vide. Si $\Pi\in \bs{\mathfrak{R}}_\E^S(\pi)$, 
alors $\Pi$ est discrte (voir plus haut) et comme en \ref{construction; hypo}, on obtient que 
l'expression ˆ droite de l'ŽgalitŽ (2) est Žgale ˆ $$\frac{1}{r(\Pi)}\prod_{v\in S}\mathrm{tr}(\Pi_v(\phi_v))\ptf$$ 
On en dŽduit alors comme en \ref{rel en v dans S} qu'en toute place $v\in S$, $\pi_v$ est un $\K_v$-relvement de $\Pi_v$. 
En toute place $v\in \V^S$, $\pi_v$ est par construction un $\K_v$-relvement faible de $\Pi_v$; et d'aprs \ref{prop unit}\,(i) (prouvŽ en \ref{preuve prop unit}), c'est un 
$\K_v$-relvement de $\Pi_v$. Donc $\pi$ est un $\K$-relvement (fort) de $\pi$. D'aprs \ref{theo global}\,(i) (prouvŽ en \ref{preuve theo global i}), on a 
$r(\Pi)= x(\delta)=r$. 

Quant ˆ la description des fibres, elle rŽsulte du lemme local-global \ref{lemme local-global}: si $\Pi'$ est une reprŽsentation 
automorphe discrte de $\GLm(\AE)$ ayant mme $\K$-relvement que $\Pi$, alors il existe un $\gamma\in \Gamma(\E/\F)$ tel que $\Pi' = {^\gamma\Pi}$. 
Comme d'autre part pour tout $\gamma\in \Gamma(\E/\F)$, $\pi$ est un $\K$-relvement de ${^\gamma\Pi}$ (d'aprs \ref{[HH] 4.5}\,(iv)), 
cela achve la dŽmonstration \ref{theo global}\,(ii). 

\subsection{\,Corollaire \ref{cor global} (global).}\label{preuve cor global} Soit $\tau$ une reprŽsentation automorphe de $\GLm(\AE)$ de la forme 
$\tau = \tau_1\times \cdots \times \tau_s$ o $\tau_i$ est une reprŽsentation automorphe discrte de $\mathrm{GL}_{m_i}(\AE)$, 
$\sum_{i=1}^s m_i=m$. Pour $i=1,\ldots , s$, soit $\pi_i$ un $\K$-relvement de $\tau_i$ ˆ $\mathrm{GL}_{m_id}(\AE)$. 
Par compatibilitŽ, en toute place $v\in \V$, entre les applications de $\K_v$-relvement et d'induction parabolique locales, on obtient que 
$\pi= \pi_1\times \cdots \times \pi_s$ est un $\K$-relvement de $\tau$. Cela prouve \ref{cor global}\,(i). 

Soit $\pi$ une reprŽsentation automorphe $\K$-stable de $\mathrm{GL}_{md}(\A)$ de la forme $\pi=\pi_1\times \cdots \times \pi_s$ 
o $\pi_i$ est une reprŽsentation automorphe discrte $\pi_i$ de $\mathrm{GL}_{n_i}(\A)$, $\sum_{i=1}^s n_i = md$. Puisque 
$\pi$ est $\kappa$-stable et que $\kappa \pi \simeq \kappa \pi_1 \times \cdots \times \kappa\pi_s$, il existe un sous-ensemble 
$\{\delta_1, \ldots ,\delta_t\}$ de $\{\pi_1,\ldots , \pi_s\}$ tel que 
$$\pi \simeq \pi'_1\times \cdots \times \pi'_t \qhq{avec} \pi'_i = \delta_i \times \K\delta_i \times\cdots \times \K^{x(\delta_i)-1}\delta_i\ptf$$ 
Pour $i=1,\ldots ,t$, $\delta_i$ est une reprŽsentation automorphe discrte de 
$\mathrm{GL}_{a_i}(\A)$ pour un entier $a_i\geq 1$ et $\pi'_i $ est une reprŽsentation automorphe (irrŽductible) $\K$-stable de 
$\mathrm{GL}_{n'_i}(\A)$, $n'_i= x(\delta_i)a_i$. Puisque $\K\pi'_i\simeq \pi'_i$, $\K$ est trivial sur le centre de $\mathrm{GL}_{n'_i}(\A)$, \ie $d$ divise $n'_i$. D'aprs \ref{theo global}\,(ii) (prouvŽ en \ref{preuve theo global ii}), $\pi'_i$ est un $\K$-relvement d'une reprŽsentation automorphe discrte $\Pi_i$ de $\mathrm{GL}_{m_i}(\AE)$, $m_i = \frac{n'_i}{d}$.  Par compatibilitŽ, en toute place $v\in \V$, entre les applications de $\K_v$-relvement et d'induction parabolique locales, on obtient que 
$\pi$ est un $\K$-relvement de $\Pi_1\times \cdots \times \Pi_t$. Cela prouve \ref{cor global}\,(ii). 

Quant ˆ la description des fibres, soient $\Pi$ et $\Pi'$ deux reprŽsentations automorphes de $\GLm(\AE)$ de la forme $\Pi= \Pi_1\times \cdots \times \Pi_s$ et $\Pi'=\Pi'_1\times \cdots \times \Pi'_t$ o $\Pi_i$, resp. $\Pi'_j$, est une reprŽsentation automorphe discrte de $\mathrm{GL}_{m_i}(\AE)$, resp. $\mathrm{GL}_{m'_j}(\AE)$. On suppose que $\Pi$ et $\Pi'$ ont mme $\K$-relvement ˆ $\GLmd(\A)$. Cela entra"ne que $s=t$ et que, quitte ˆ permuter les $\Pi'_j$, on peut supposer que pour $i=1,\ldots ,s$, on a: 
\begin{itemize}
\item $m_i=m'_i$; 
\item $\Pi_i$ et $\Pi'_i$ ont mme $\K$-relvement ˆ $\mathrm{GL}_{m_id}(\A)$.  
\end{itemize}
On conclut gr‰ce ˆ la description des fibres dans \ref{theo global}\,(ii) (prouvŽ en \ref{preuve theo global ii}). Cela prouve \ref{cor global}\,(iii). 

\section{CompatibilitŽ avec le changement de base.}\label{comp CB} 
Dans \cite{H2}, Henniart dŽduit l'existence et les propriŽtŽs du $\K$-relvement des reprŽsentations automorphes cuspidales de $\GLm(\AE)$ ˆ $\GLmd(\A)$ des rŽsultats analogues 
pour le $\sigma$-relvement (\ie le relvement pour le changement de base) de $\GLn(\A)$ ˆ $\GLn(\AE)$. Les applications de $\K$-relvement et de $\sigma$-relvement Žtant dŽsormais Žtablies pour les reprŽsentation automorphes induites de discrtes, on peut les relier l'une ˆ l'autre comme dans \textit{loc.~cit.} 

\subsection{Rappels sur le changement de base.} Soit un entier $n\geq 1$. Soit $E/F$ une extension finie cyclique de corps locaux de caractŽristique nulle. Fixons un gŽnŽrateur $\sigma$ de $\Gamma(E/F)$. Dans \cite{AC}, Arthur et Clozel construisent une application de changement de base, que nous appellerons ici \textit{$\sigma$-relvement}, qui ˆ une reprŽsentation irrŽductible tempŽrŽe $\pi$ de $\GLn(F)$ associe une reprŽsentation irrŽductible tempŽrŽe 
$\sigma$-stable $\pi_E$ de $\GLn(E)$. Cette application est caractŽrisŽe par une identitŽ de caractre, appelŽe \textit{identitŽ de Shintani}, entre le caractre de $\pi_E$ tordu par un isomorphisme entre $\pi_E$ et $\pi_E^\sigma= \pi_E\circ \sigma$ et le caractre de $\pi$. En \cite{BH}, Badulescu et Henniart Žtendent ce rŽsultat ˆ toutes les reprŽsentation irrŽductibles elliptiques, resp. unitaires, de $\GLn(F)$. Ils obtiennent aussi un rŽsultat analogue dans le cas archimŽdien, 
\ie pour $E/F\simeq \mathbb{C}/\mathbb{R}$. 

La dŽmonstration de \cite{BH} est, comme celle du prŽsent article qui s'en inspire largement, basŽe sur le principe local-global. 
Comme consŽquence, 
Badulescu et Henniart obtiennent, pour une extension finie cyclique $\E/\F$ de corps de nombres et un gŽnŽrateur $\sigma$ de $\Gamma(\E/\F)$, que 
toute reprŽsentation automorphe discrte $\pi$ de $\GLn(\AF)$ admet un $\sigma$-relvement (fort) ˆ $\GLn(\AE)$; au sens 
o il existe une reprŽsentation automorphe $\pi_E$ de $\GLn(\AE)$ (induite de discrte) telle qu'en toute place $v$ de $\F$, la composante 
local $\pi_{\E,v}$ de $\pi_\E$ en $v$ soit un $\sigma$-relvement de $\pi_v$. Bien sžr pour que cela ait un sens,  il faut au prŽalable avoir Žtendu la notion de $\sigma$-relvement local au cas d'une $F$-algbre cyclique de degrŽ fini $E$, ce qui est fait dans \cite{BH} et \cite{H1}. 
PrŽcisons la forme de ce $\sigma$-relvement $\pi_E$. Rappelons qu'on a notŽ $x(\pi)$ le cardinal de $X(\pi)$; ainsi $r= \frac{d}{x(\pi)}$ est le cardinal du stabilisateur de la classe d'isomorphisme de $\pi$ dans $X(\E/\F)$. Puisque $\kappa^{x(\pi)}\pi \simeq \pi$, l'entier $d$ divise $x(\pi)n$ (d'aprs l'argument habituel sur les caractres centraux), \ie $r$ divise $n$. Alors (d'aprs \cite{BH}) il existe une reprŽsentation automorphe discrte $\sigma^r$-stable $\delta_E$ de $\mathrm{GL}_{\frac{n}{r}}(\AE)$ telle que la reprŽsentation automorphe $\pi_E = \delta_E \times \delta_E^\sigma\times \cdots \times \delta_E^{\sigma^{r-1}}$ de $\GLn(\AE)$ soit un $\sigma$-relvement de $\pi$. De plus on a $r(\delta_E)= r$ et $\pi$ est cuspidale si et seulement si $\delta_E$ est cuspidale. RŽciproquement, toute reprŽsentation automorphe $\pi_E$ de $\GLn(\AE)$ de la forme $\pi_E= \delta_E \times \cdots \times \delta_E^{\sigma^{r-1}}$ pour un entier $r\geq 1$ divisant $n$ et une reprŽsentation automorphe discrte $\delta_E$ de $\mathrm{GL}_{\frac{n}{r}}(\AE)$ telle que $r(\delta_E)=r$, est le $\sigma$-relvement d'une reprŽsentation automorphe discrte de $\GLn(\AF)$. 

Le rŽsultat global ci-dessus s'Žtend naturellement ˆ toutes les reprŽsentations automorphes de $\GLn(\AF)$ qui sont induites de discrtes, par compatibilitŽ, 
en toute place $v$ de $\F$, entre les applications de $\sigma$-relvement et d'induction parabolique locales. En d'autres termes, on a l'analogue du corollaire \ref{cor global} pour le changement de base.  

\subsection{$\sigma$-relvement faible.}\label{sigma-rel faible} Dans cette sous-section, $F$ est une extension finie de $\mathbb{Q}_p$ et $E$ est 
une $F$-algbre cyclique \textit{non ramifiŽe} de degrŽ fini $d$, de groupe $\Gamma$. 
On Žcrit $E=E_1\times \cdots \times E_r$ o $E_i/F$ est une extension non ramifiŽe de degrŽ $s= \frac{d}{r}$. On reprend les hypothses et les notations de \ref{notations ext cyclique}. En particulier, on a 
$\sigma E_i = E_{i+1}$ ($i=1,\ldots , r-1$) et $\sigma E_r = E_1$ pour un gŽnŽrateur $\sigma$ de $\Gamma$. 
On fixe aussi un caractre $\kappa$ de $F^\times$ de noyau $\mathrm{N}_{E/F}(E^\times)= \mathrm{N}_{E_1/F}(E_1^\times)$ et on pose $\zeta=\kappa(\varpi)$ pour une (i.e. pour toute) uniformisante $\varpi$ de $F^\times$. Pour $m\in \mathbb{N}^*$ et $y\in (\mathbb{C}^\times)^{mr}= ((\mathbb{C}^\times)^m)^r$, on a dŽfini en \ref{param satake} une reprŽsentation irrŽductible sphŽrique $\Pi_y$ de $\GLm(E)= \prod_{i=1}^r\GLm(E_i)$ et une reprŽsentation irrŽductible sphŽrique $\kappa$-stable $\pi_{\delta(y)}$ de $\GLmd(F)$ qui (d'aprs \ref{description du K-rel faible dans le cas gŽnŽral}) est un $\kappa$-relvement faible de $\Pi_y$. Si de plus $\Pi_y$ est unitaire, on sait (d'aprs \ref{prop unit}\,(i)) que $\pi_{\delta(y)}$ est un $\kappa$-relvement de $\Pi_y$. 

L'entier $n\geq 1$ Žtant fixŽ, posons $K_E=\GLn(\mathfrak{o}_E)= \prod_{i=1}^r\GLn(\mathfrak{o}_{E_i})$ et notons $\mathrm{Irr}^{K_E}_\sigma(\GLn(E))$ le sous-ensemble 
$\mathrm{Irr}^{K_E}(\GLn(E))$ formŽ des reprŽsentations qui sont $\sigma$-stables. Observons que si $E=E_1$, une reprŽsentation irrŽductible sphŽrique de $\GLn(E)$ est toujours $\sigma$-stable.

Supposons pour commencer que $E/F$ soit une extension de corps (\ie $E=E_1$). Pour $y\in (\mathbb{C}^\times)^n$, on a dŽfini en \ref{param satake} une reprŽsentation irrŽductible sphŽrique $\pi_y$ de $\GLn(F)$ dont la classe d'isomorphisme ne dŽpend que de l'image de $y$ dans $(\mathbb{C}^\times)^n/ \mathfrak{S}_n$. On dŽfinit de la mme manire la reprŽsentation sphŽrique $\pi_{E,y}$ de $\GLn(E)$, ˆ partir du caractre non ramifiŽ $\chi_{E,y}$ de $A_{E,0} = (E^\times)^n$ donnŽ par $\chi_{E,y}(\varpi^{a_1},\ldots ,\varpi^{a_n})= \prod_{i=1}^n y_i^{a_i}$. L'application $y\mapsto \pi_{E,y}$ induit une bijection de 
$(\mathbb{C}^\times)^n/\mathfrak{S}_n$ sur $\mathrm{Irr}^{K_E}(\GLn(E))= \mathrm{Irr}^{K_E}_\sigma(\GLn(E)$. Pour $y=(y_1,\ldots ,y_n)\in (\mathbb{C^\times})^n$, posons $y^d =(y_1^d, \ldots , y_n^d)$. Toute reprŽsentation irrŽductible sphŽrique de $\GLn(E)$ isomorphe ˆ $\pi_{E,y^d}$ est appelŽe un \textit{$\sigma$-relvement faible} de $\pi_y$. L'application de $\sigma$-relvement faible correspond par dualitŽ ˆ l'homomorphisme de changement de base $$c:\mathcal{H}(\GLn(F),K) \rightarrow \mathcal{H}(\GLn(E),K_E)$$ dŽduit comme en \cite[ch.~I, \S4.2]{AC} des isomorphismes de Satake. 

\begin{remark}\label{rem comp IA/CB local}
\textup{
\begin{enumerate}
\item[(i)] L'application surjective $(\mathbb{C}^\times )^n \rightarrow (\mathbb{C}^\times)^n, \, y\mapsto y^d$ se factorise en une application surjective 
$(\mathbb{C}^\times )^n/\mathfrak{S}_n \rightarrow (\mathbb{C}^\times)^n/ \mathfrak{S}_n$, qui n'est pas injective si $d>1$. L'application 
de $\sigma$-relvement faible $\mathrm{Irr}^K(\GLn(F))\rightarrow \mathrm{Irr}^{K_E}(\GLn(E))$ est donc surjective, et deux  
reprŽsentations irrŽductibles sphŽriques $\pi_y$ et $\pi_{y'}$ de $\GLn(F)$ ont mme $\sigma$-relvement faible si et seulement s'il existe une permutation 
$i \mapsto \nu(i)$ de $\{1,\ldots ,n\}$ telle que $y'_i\in \langle \zeta \rangle y_{\nu(i)}$ ($i=1,\ldots ,n)$. 
\item[(ii)] Supposons que $d$ divise $n$ et posons $m=\frac{n}{d}$. Pour $y\in (\mathbb{C}^\times)^m$, on a dŽfini en \ref{param satake} un ŽlŽment $\delta(y)\in (\mathbb{C}^\times)^n$. Cet ŽlŽment dŽpend des choix de $\zeta$ (\ie de $\kappa$) et d'un ŽlŽment $t\in (\mathbb{C}^\times)^m$ tel que $t^d=y$, mais l'image de $\delta(y)$ dans $(\mathbb{C}^\times)^n/\mathfrak{S}_n$ est bien dŽfinie (elle ne dŽpend pas de ces choix). Observons que l'ŽlŽment $\delta(y)^d \in (\mathbb{C}^\times)^n$ est lui aussi bien dŽfini et que son image dans $(\mathbb{C}^\times)^n/\mathfrak{S}_n$ co\"{\i}ncide avec celle de $(y,\ldots ,y)$ ($d$ fois). En d'autres termes, la reprŽsentation irrŽductible sphŽrique $\Pi_y \times \cdots \times \Pi_y$ ($d$ fois) de $\GLn(E)$ est un $\sigma$-relvement faible de $\pi_{\delta(y)}$. \hfill $\blacksquare$
\end{enumerate}
}
\end{remark}

Passons au cas o $E$ est une $F$-algbre cyclique quelconque (de degrŽ fini $d$). Dans ce cas l'homomorphisme de changement de base  
$$c:\mathcal{H}(\GLn(E),K_E)= \bigotimes_{i=1}^r \mathcal{H}(\GLn(E_i),K_{E_i})\rightarrow  \mathcal{H}(\GLn(F),K)$$ est donnŽ par\footnote{Ici $*$ est le produit de convolution dans $\mathcal{H}(\GLn(E_1),K_{E_1})$ dŽfini ˆ l'aide de la mesure de Haar sur $\GLn(E_1)$ qui donne le volume $1$ ˆ $K_{E_1}$.} 
$$c(\phi)= c_{E_1}(\phi_1 * \phi_2^{\sigma} * \cdots * \phi_r^{\sigma^{r-1}})$$ pour toute fonction $\phi= \phi_1\otimes \cdots \otimes \phi_r\in \mathcal{H}(\GLn(E),K_E)$. L'application de $\sigma$-relvement faible correspondant ˆ $c$ par dualitŽ est dŽfinie comme suit: pour $y\in (\mathbb{C}^\times)^n$, toute reprŽsentation irrŽductible sphŽrique de $\GLn(E)= \prod_{i=1}^r \GLn(E_i)$ isomorphe ˆ $\pi_{E_1,y^s}\otimes \cdots \otimes \pi_{E_r,y^s}$ est appelŽe un 
\textit{$\sigma$-relvement faible} de $\pi_y$.

\begin{remark}
\textup{La remarque \ref{rem comp IA/CB local}\,(ii) est encore vraie ici. Supposons que $d$ divise $n$ et posons $m=\frac{n}{d}$. Pour $y\in (\mathbb{C}^\times)^{mr}$, la reprŽsentation $\Pi_y \times \cdots \times \Pi_y$ ($d$ fois) de 
$\GLn(E)$ est un $\sigma$-relvement faible de $\pi_{\delta(y)}$. \hfill $\blacksquare$}
\end{remark}

\subsection{CompatibilitŽ induction automorphe -- changement de base.} 
Soit $\E/\F$ une extension finie cyclique de corps de nombres. On fixe un gŽnŽrateur $\sigma$ de $\Gamma(\E/\F)$ et un caractre $\K$ de $\A^\times$ de noyau $F^\times \mathrm{N}_{\E/\F}(\AE)$. On dŽfinit comme en \ref{def rel faible aut}\,(i) la 
notion de $\sigma$-relvement faible  pour les reprŽsentations automorphes de $\GLn(\A)$ qui sont induites de discrte. D'aprs \cite{BH}, toute reprŽsentation automorphe $\pi$ de $\GLn(\A)$ induite de discrte admet un $\sigma$-relvement (fort); en particulier (d'aprs \ref{rigiditŽ}) tout $\sigma$-relvement faible de $\pi$ est un 
$\sigma$-relvement. 

Si $\pi$ est une reprŽsentation automorphe discrte de $\GLn(\A)$ pour un entier $n\geq 1$, on note 
$\wt{\pi}$ la reprŽsentation de $\mathrm{GL}_{x(\pi)n}(\A)$ dŽfinie par 
$$\wt{\pi}= \pi \times \K\pi\times\cdots \times\K^{x(\pi)-1}\pi\ptf$$ 
Si $\Pi$ est une reprŽsentation automorphe discrte de $\GLm(\AE)$ pour un entier $m\geq 1$, on note 
$\wt{\Pi}$ la reprŽsentation automorphe de $\mathrm{GL}_{g(\Pi)m}(\AE)$ induite de discrte dŽfinie par 
$$\wt{\Pi}=\Pi\times \Pi^\sigma \times \cdots \times \Pi^{\sigma^{g(\Pi)-1}}\ptf$$ 

\begin{theorem}\label{compatibilitŽ IA/BC}Soient $m\geq 1$ et $n\geq 1$ deux entiers.
\begin{enumerate}
\item[(i)] Soit $\Pi$ une reprŽsentation automorphe discrte de $\GLm(\AE)$. Posons $r=r(\Pi)$ et $g=g(\Pi)\;(=\frac{d}{r})$. La repr\'esentation 
$\wt{\Pi}$ de $\mathrm{GL}_{gm}(\AE)$ est un $\sigma$-relvement 
d'une reprŽsentation automorphe discrte $\delta$ de $\mathrm{GL}_{gm}(\A)$. On a $x(\delta)=r$ et la reprŽsentation 
$\wt{\delta}$ de $\GLmd(\A)$ est un $\K$-relvement de $\Pi$.
\item[(ii)] Soit $\pi$ une reprŽsentation automorphe discrte de $\GLn(\A)$. Posons $r=x(\pi)$. 
Alors $d$ divise $rn$ et la reprŽsentation $\wt{\pi}$ de $\mathrm{GL}_{rn}(\A)$ est 
un $\K$-relvement d'une reprŽsentation automorphe discrte $\Pi$ de $\mathrm{GL}_{\frac{rn}{d}}(\AE)$. On a $r(\Pi)=r$ et la 
reprŽsentation $\wt{\Pi}$ de $\GLn(\AE)$ est un $\sigma$-relvement de $\pi$. 
\end{enumerate}
\end{theorem}

\begin{demo}
Seules la dernire assertion du point (i) et la dernire assertion du point (ii) sont ˆ prouver. Pour (i), en presque toute place $v\in \Vfin$, $\widetilde{\delta}_v$ est un $\K_v$-relvement faible de $\Pi_v$ (d'aprs \ref{description du K-rel faible dans le cas gŽnŽral}). Donc $\wt{\delta}$ est un $\K$-relvement (fort) de $\pi$. 
Pour (ii), en presque toute place $v\in \Vfin$, $\wt{\Pi}$ est un  $\sigma$-relvement faible de $\pi$ (d'aprs \ref{sigma-rel faible}). Donc $\wt{\Pi}$ est un $\sigma$-relvement (fort) de $\pi$. \hfill $\square$
\end{demo}

\vskip2mm
\begin{corollary}Soient $m\geq 1$ et $n\geq 1$ deux entiers.
\begin{enumerate}
\item[(i)] Soit $\Pi$ une reprŽsentation automorphe de $\GLm(\AE)$ induite de discrte. La repr\'esentation 
$\Pi\times \Pi^\sigma \times \cdots \times \Pi^{\sigma^{d-1}}$ 
de $\GLmd(\AE)$ est un $\sigma$-relvement d'une reprŽsentation automorphe $\pi$ de $\GLmd(\A)$ induite de discrte 
et $\K$-stable, qui est un $\K$-relvement de $\Pi$.
\item[(ii)] Soit $\pi$ une reprŽsentation automorphe de $\GLn(\A)$ induite de discrte. La reprŽsentation 
$\pi \times \K\pi \times\cdots \times \K^{d-1}\pi$ de $\mathrm{GL}_{dn}(\A)$ est un 
$\K$-relvement d'une reprŽsentation automorphe $\Pi$ de $\GLn(\AE)$ induite de discrte et $\sigma$-stable, qui 
est un $\sigma$-relvement de $\pi$. 
\end{enumerate}
\end{corollary}

\begin{demo}
Par compatibilitŽ entre les applications (globales) de $\K$-relvement, resp. $\sigma$-relvement, et d'induction parabolique (\ref{cor global}\,(i) et \cite{BH}). 
 \hfill $\square$
\end{demo}

\section{L'Žtat des lieux pour les corps de fonctions}\label{Žtat des lieux CF}

Dans cette section, $F$ est un corps local non archimŽdien de caractŽristique $p>0$, \ie $F\simeq \mathbb{F}_{q}((\varpi))$ avec $q=p^r$ pour un entier $r\geq 1$; et $\F$ est un corps global de caractŽristique $p>0$, \ie le corps des fonctions d'une courbe projective lisse et gŽomŽtriquement connexe sur $\mathbb{F}_q$. 

\subsection{Transfert, lemme fondamental et $\kappa$-relvement.} Soit $E$ une extension cyclique de $F$, ou plus gŽnŽralement une $F$-algbre cyclique, de degrŽ fini $d$. 
Fixons un entier $m\geq 1$ et posons $n=md$. Soient $G=\GLn(F)$ et $H= \GLm(E)$. On fixe un caractre $\kappa$ de $F^\times$ de noyau $\mathrm{N}_{E/F}(E^\times)$. On dŽfinit comme en caractŽristique nulle un facteur de transfert $\Delta(\gamma)$ pour $\gamma\in \HGreg$ et la notion de fonctions concordantes (\ref{fonctionsconcor}). Le transfert (thŽorme \ref{existence transfert} en caractŽristique nulle) n'est pas dŽmontrŽ en caractŽristique $p>0$, 
mais on sait construire suffisamment de couples de fonctions concordantes $(f,\phi)$. PrŽcisŽment (\cite[3.8]{HH}; voir aussi \ref{facteurs de transfert bis}), pour toute fonction $f\in C^\infty_{\mathrm{c}}(\Greg)$, il existe une fonction $\phi\in C^\infty_{\mathrm{c}}(\HGreg)$ qui concorde avec $f$; et pour toute fonction $\phi\in C^\infty_{\mathrm{c}}(\HGreg)$ vŽrifiant la condition de \ref{existence transfert}\,(ii), il existe une fonction $f\in C^\infty_{\mathrm{c}}(\Greg)$ qui concorde avec 
$\phi$. 

Le lemme fondamental pour l'induction automorphe est prouvŽ en toute gŽnŽralitŽ dans \cite{HL1}, ˆ partir du rŽsultat de Waldspurger pour $p>n$ \cite{W1}. On devrait pouvoir en dŽduire le transfert comme le fait Waldspurger en caractŽristique nulle  \cite{W2, W3}, mais cela reste ˆ faire. 

La notion de $\kappa$-relvement est dŽfinie comme en \ref{defkapparel}. L'absence du transfert conduit ˆ privilŽgier la formulation en termes de fonctions-caractres (\ref{def kappa-rel}\,(2)). L'existence d'un $\kappa$-relvement $\pi$ pour toute reprŽsentation irrŽductible gŽnŽrique unitaire $\tau$ de $\GLm(E)$, est prouvŽe dans \cite[ch.~III]{HL3}. On prouve aussi, comme en caractŽristique nulle \cite{HL2}, que la constante de proportionalitŽ $c(\tau,\pi,\psi)= c(\tau,\pi, A_{\pi,\psi}^{\mathrm{g\acute{e}n}})$ ne dŽpend pas de $\tau$, o $\psi$ est un caractre additif non trivial de $F$ fixŽ au dŽpart (\textit{loc.\,cit.}). 

\subsection{$\K$-relvement (global).} 
Soit $\E$ une extension finie cyclique de $\F$, de degrŽ $d$, et soit $\K$ un caractre de $\AF^\times$ de noyau $\F^\times \mathrm{N}_{\E/\F}(\AE)$. 
La correspondance de Langlands (globale) Žtablie par L.~Lafforgue \cite{Laf} permet d'associer ˆ toute reprŽsentation automorphe cuspidale $\Pi$ de $\GLm(\AE)$ une reprŽsentation automorphe $\pi$ de $\GLn(\AF)$, qui est un $\K$-relvement faible de $\Pi$ au sens de \ref{def kapparelfaible}. On prouve dans \cite[ch.~IV]{HL3} que cette application de relvement globale est compatible aux applications de relvement locales: en toute place $v$ de $\F$, $\pi_v$ est un $\K_v$-relvement de $\Pi_v$. On prouve aussi (\textit{loc.\,cit.}) que pour chaque place $v$ de $\F$, l'application de $\K_v$-relvement est compatible ˆ la correspondance de Langlands locale Žtablie par laumon, Rapoport et Stuhler \cite{LRS}. 

\subsection{Stabilisation de la formule des traces pour $(G,\omega)$.} Le seul outil qui nous manque pour faire fonctionner la dŽmonstration en caractŽristique $p>0$ est la formule de comparaison \ref{dec endos}\,(4), \ie la stablisation de la formule des traces pour $(\mathrm{GL}_{n/\F},\K\circ\det)$. On dispose aujourd'hui du dŽveloppement spectral fin de la formule des traces tordue pour les corps globaux de caractŽristique $p>0$ (\cite{LLe}). Pour stabiliser ce dŽveloppement spectral, ou mme seulement la partie discrte qui nous intŽresse, la stratŽgie d'Arthur et Langlands reprise par M\oe glin et Walspurger dans le cas tordu \cite{MW2} conduit ˆ mener de front des rŽductions pour les expressions spectrales et gŽomŽtriques de cette stabilisation. C™tŽ gŽomŽtrique, les identitŽs de transfert nŽcessaires entre intŽgrales orbitales pondŽrŽes semi-locales ne sont pas Žtablies en caractŽristiques $p$. D'ailleurs on ne dispose mme pas encore d'un dŽveloppement gŽomŽtrique fin. Ce dernier fait l'objet d'un travail en cours \cite{Le1, Le2}. Les difficultŽs nouvelles sont essentiellement dues ˆ la prŽsence d'ŽlŽments primitifs insŽparables  (qui n'existent pas en caractŽristique nulle). Bien sžr si l'on suppose que la caractŽristique $p$ de $\F$ est strictement supŽrieure ˆ $n$, ou mme plus gŽnŽralement que $p$ ne divise pas $n$, alors ces difficultŽs disparaissent et il ne devrait pas tre trs difficile d'Žtendre les rŽsultats connus en caractŽristique nulle ˆ la caractŽristique $p>0$ (du moins pour le groupe $\GLn$). 

En conclusion, si l'on admet la formule de comparaison \ref{dec endos}\,(4) pour les corps globaux de caractŽristique $p>0$, 
alors les rŽsultats de cet article s'Žtendent  ˆ la caractŽristique $p>0$.

\end{document}